\documentclass[phd,tocprelim]{cornell}
\usepackage{graphicx,pstricks}
\usepackage{graphics}
\usepackage{moreverb}
\usepackage{subfigure}
\usepackage{epsfig}
\usepackage{subfigure}
\usepackage{hangcaption}
\usepackage[font=small,labelsep=none]{caption}
\usepackage{xcolor}

\renewcommand{\caption}[1]{\singlespacing\hangcaption{#1}\normalspacing}

\usepackage{amsfonts}
\usepackage{amsmath}
\usepackage{amssymb,bm}
\usepackage{amsthm,color}
\usepackage{tikz, tikz-cd}
\usepackage{physics}
\usepackage{blkarray, array}
\usepackage{stmaryrd}
\usepackage{hyperref}
\usepackage{pdflscape}
\usepackage{longtable}
\usepackage{changepage}
\usepackage{geometry}
\usepackage{mathtools}

\newtheorem{Thm}{Theorem}

\newtheorem{Cor}{Corollary}
\newtheorem{Lem}{Lemma}
\newtheorem{Prop}{Proposition}

\newtheorem{Def}{Definition}

\newtheorem{Exa}{Example}
\newtheorem*{Proof}{Proof}

\newcommand{\C}{\mathbb{C}}
\newcommand{\Z}{\mathbb{Z}}

\def\+{\includegraphics[scale=0.4]{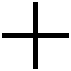}}
\def\elbow{\includegraphics[scale=0.4]{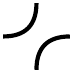}}

\def\rin{\rotatebox[origin=c]{90}{$\in$}}

\def\+{\includegraphics[scale=0.3]{crosstile.pdf}}
\def\-{\includegraphics[scale=0.3]{elbowtile.pdf}}

\newcommand{\ra}{\ensuremath{\xrightarrow}}
\newcommand{\A}{\ensuremath{\mathbb{A}}}
\newcommand{\R}{\ensuremath{\mathbb{R}}}

\newcommand{\PP}{\ensuremath{\mathbb{P}}}

\newcommand{\del}{\ensuremath{\partial}}

\newcommand{\bs}{\ensuremath{\backslash}}

\newcommand\rd[1]{{\color{red}#1}}
\newcommand\teal[1]{{\color{teal}#1}}

\newcommand{\mf}{\ensuremath{\mathfrak}}

\newcommand{\ov}{\ensuremath{\overline}}
\newcommand{\ep}{\ensuremath{\epsilon}}

\DeclareMathOperator{\codim}{codim}

\newcommand\minip[3]{\noindent
  \begin{minipage}{#1\linewidth}
    \hskip #2
    #3
  \end{minipage}
}

\title {Singular Points on Positroid Varieties and Planar N=4 Supersymmetric Yang-Mills Theory (Revised 2026)}
\author {Joseph Fluegemann}
\conferraldate {August}{2024}
\degreefield {Ph.D.}
\copyrightholder{Joseph Fluegemann}
\copyrightyear{2024}

\begin{document}

\maketitle
\makecopyright

\begin{abstract}

Positroid varieties $\Pi_f$ are geometric shapes that give us a way of cutting up the Grassmannian $Gr(k,n)$ of $k$-dimensional planes in $n(\geq k)$ dimensions. They have nice combinatorial properties, since they can be enumerated using special permutations called bounded affine permutations (denoted, say by $f$), which interestingly, also characterize juggling patterns \cite{KLS}. Furthermore, (the real positive part of) positroid varieties parameterize the space that is integrated over when calculating leading singularities in planar $N=4$ supersymmetric Yang-Mills theory (leading singularities can be used to build up the amplitude in certain quantum field theories: see Section \ref{sec:leadsing} for a quick introduction). 

The main question we answer in this thesis is whether a positroid variety  is smooth (everywhere); that is, whether a positroid variety $\Pi_f$ has any geometric singularities. 
We first show that it is sufficient to check singularity at the ${n \choose k}$ special torus $T$-fixed points (denoted $\lambda$), since if $\Pi_f$ has any singularities, then it must be singular at one of these points $\lambda$ (Theorem \ref{thm:testonnchoosek}).
Then, we show that we can check singularity at these points $\lambda$ by calculating the equivariant cohomology restricted to the point: $[\Pi_f]|_\lambda$ and then restricting to a circle $C^\times$ within the torus $T$; since this will give us the multiplicity (degree of the tangent cone) of $\Pi_f$ at $\lambda$ (Section \ref{sec:smoothness}). We calculate $[\Pi_f]|_\lambda$ in 2 ways:

(1) A diagrammatically nice and more intuitive way using affine pipe dreams (Subsection \ref{subsec:AJSBilley}) 

(2) A computational method (Chapter \ref{sec:compformula})

For applications, (2) is easier to use: the code can be found at \url{josephflueg.github.io} (along with other resources); how to use the code is described in Section \ref{sec:usingcode}: simply input $f$ and $\lambda$, and the code will spit out the multiplicity, telling you whether $\Pi_f$ is singular at $\lambda$ or not.  

Our proof of method (1) shows that the equivariant cohomology can be read off of these affine pipe dream diagrams, and the multiplicity is equal to the number of such diagrams. In fact, the Main Theorem then says it is sufficient to just draw a single diagram and check a criterion for determining whether the multiplicity is greater than 1 (so $\Pi_f$ singular) or not. Thus, method (1) gives a nice way of visualizing and doing this calculation relatively quickly by hand.

Using this code, we have included tables in Chapter \ref{ch:catalogue} of this thesis listing the multiplicities of all the points on positroid varieties up to $n=6$, including whether each positroid variety is smooth or not. 

Section \ref{ch:affinefinitepipes} describes the relationship between what we call finite pipe dreams (found in e.g. \cite{Post} and \cite{ARW}) and the affine pipe dreams we work with primarily in this thesis. It shows that the results on smoothness for affine pipe dreams apply for finite pipe dreams as well. 

It is also possible to put an ordering on the set of pairs (positroid, point)$=(\Pi_f,\lambda)$ given by deletion and contraction. A pictorial example of this ordering can be found in Figure \ref{fig:atomicorderexa}. See Sections \ref{sec:DelContr} and \ref{section:theorder} for the background on this ordering and what these operations are, Section \ref{section:delcontrpipedream} for their relationship to affine pipe dreams and Section \ref{sec:smoothnessordering} for their relationship to smoothness or singularity. Using the results in Section \ref{sec:smoothnessordering}, we can give another (but nonconstructive) characterization of a singular point $\lambda$ on a positroid variety $\Pi_f$  as being a pair $(\Pi_f,\lambda)$ that can reach an atomic configuration (a minimal singular pair) under the ordering given by deletion and contraction. 

In Part II of this thesis, we have attempted to connect the mathematical content of the first part with the physics of planar $\mathcal{N}=4$ supersymmetric Yang-Mills Theory. We start in Chapter \ref{ch:QFT} by giving a brief introduction to quantum field theory, what leading singularities are, where Grassmannians and then positroids come into play in the story, etc. Then in Chapter \ref{ch:BCFW}, we discuss the Britto-Cachazo-Feng-Witten recursion relations. We begin by introducing the recursion, before focusing on the notion of a BCFW bridge in on-shell diagrams, and the construction of a permutation using these BCFW bridges first described in \cite{ABCGPT}. We end by describing how to build an on-shell diagram on a pipe dream using BCFW bridges (Section \ref{sec:pipesplabicexample} does an example in detail going all the way from permutation to plabic diagram). We noted above that deletion and contraction operations are involved in the ordering on atomic pairs; having seen that a deletion-like operation appears in the inverse soft construction in the physics literature, in Chapter \ref{ch:softlimits} we work out some simple combinatorics related to these inverse soft factors, showing that the inverse soft factors are a particular choice of right inverse for the deletion and contraction operators, that produce ``2-step" switches when deletion or contraction are applied. Finally in Chapter \ref{ch:smoothphysics}, we explore whether singularities in positroid varieties have connections to amplitudes in planar N=4 SYM, concluding that they do not. 

A more in-depth overview and summary of the contents of this thesis is in the Introduction.

Note: This is an revised and cleaned up version of my thesis originally written in August 2024. 

\end{abstract}

\begin{acknowledgements} The intellectual content of this thesis owes the most to two individuals: first and foremost to my advisor Allen Knutson, and secondly to Nima Arkani-Hamed.  

Thank you Allen for countless conversations discussing many different topics in math, including of course, positroid varieties. Thank you for your guidance and support over many years.

Thank you Nima for introducing me to the field of Amplitudes during the A.D. White Cornell Professor joint math-physics colloquim (organized by Liam McAllister and Mike Stillman) in early 2019. Thank you for teaching me much of what I know about Amplitudes during the Harvard class, for answering my questions and bringing me to Harvard that semester for the workshop.

Thank you to my thesis committee members Allen, Liam, and Karola Meszaros.

I'm thankful to the Professors who taught me physics at Cornell: Liam McAllister, Csaba Csaki, Maxim Perelstein, and Tom Hartman, as well as my homework buddies Andrew Gomes, Ameen Ismail, Colin Bundschu, Zebang Xu, and Daniel Longenecker, and Mijo Ghosh and Aaron Hui for your tutoring during quantum field theory.

I am incredibly grateful to many members of the Amplitudes community for their help during graduate school. Thank you Jake Bourjaily for inviting me to Penn State and giving me the first opportunity to present the results in this thesis in December 2021, as well as for all your explanations. Thank you Jaroslav Trnka for introducing me to my first resources in Amplitudes in 2019, and to you and Shruti Paranjape for involving me in a project. Thank you Sebastian Mizera for several illuminating discussions and introducing me to the Amplitudes Journal Club. Thank you Nick Early for lots of helpful career advice and inspiration. Thank you Aidan Herderschee for teaching me about research and reading papers. Thank you Tomek Lukowski for discussing research with me. Thank you Lauren Williams for pointing out to me several helpful ideas. Thank you Thomas Lam and David Speyer for bringing me to University of Michigan to speak and for a couple helpful conversations. Thank you Karen Yeats for your support and encouragement. Thank you to Lecheng Ren for your hospitality while I was at Brown.

I'm also very thankful for the Amplitudes and Geometry Zoom Journal Club, which I was able to co-organize with Jonah Stalknecht for awhile, and which gave me an opportunity to interact with the Amplitudes community from all over the world. Thank you Robert Moerman for teaching me about your research. Thank you Gabriele Dian for your many helpful question and insights. Thank you Carolina Figuereido, Hadleigh Frost, Zhang Yong, Nick Early, and Ran Tessler for giving amazing lectures.

I also want to acknowledge the support of the Munich Institute for Astro-, Particle and BioPhysics (MIAPbP) for the Scattering Amplitudes Program in July 2022, which provided an excellent environment for thinking about and working on physics, as well as thank June Huh for allowing me to participate as visiting scholar in the Special Program on Algebraic and Geometric Combinatorics at the Institute for Advanced Study in the Fall of 2024.

I'm thankful for the love and support of many friends and family members.

I'm grateful to my parents and my brother for their continued love and support and always being available when I need help. Thank you to all my cousins for your support, and especially Michael Fluegemann for being such an amazing friend. Thank you Uncle Dave and Aunt Betty and Uncle Herb and Aunt Diane for opening up your homes to me when I'm visiting.

I'm especially grateful to several friends whose companionship was instrumental in getting me through the program. Thank you Mark Dalthorp for being the best friend, roommate and brother someone could ask for. Thank you Daniel Longenecker for countless hangouts discussing physics and life, and the opportunity to stimulate each other's interest in Amplitudes while you were at Cornell. Thank you Kam Black for always staying in touch and your great advice. Thank you David Kim for deep conversations and always being available. Thank you Tim Hogue for inspiring me with your excellence in both work and personal life. Thank you Olu Olorode for your wise insights and example of discipline. Thank you Wesley Cheung for our fun catchups multiple times a year.

Thank you to graduate field coordinators Melissa Totman, Jim Utz, Elly O'Brien, and Betsy Collins for all your administrative help and all your efforts behind the scenes. Thank you Scott and Joni Hathorn, David and Mavis Ng, Jianhua Chen and Xiaohong Wang, J.W. and Stephanie Betts, Nic and Christy Brenner, Ben and Hannah Hutton for your hospitality and kindness. Thank you Lucas Raley and Brennon Dunbar and all my friends in Cru. Thank you Inbeom Lee, Leo Moon, David Chung, James and Elisha, Yurong, Hongpyo Ha, Cheol Seung, Yun, Caiwei, Haokun, Jialiang for your hospitality. Thank you Carl and Vivek and everyone at Chesterton. Thank you Mavis Ng, Christine Trumble, Joy O'Connor, and Nicole Cornell, and Nicholas Serniak for your health-related support.

\end{acknowledgements}

\contentspage
\figurelistpage

\normalspacing \setcounter{page}{1} \pagenumbering{arabic}
\pagestyle{cornell} \addtolength{\parskip}{0.5\baselineskip}

\chapter{Introduction (Overview/Summary)}

In this introduction, we give an overview/summary of the thesis, pointing to specific sections and results that may interest the reader. 

\section{Math}

The space in which positroid varieties live is the Grassmannian $Gr(k,n)$, a geometric space parameterizing the space of $k$-dimensional planes in $n$-dimensional space. For example, if we are over the real $\R$, the lines in 2D space forms a circle.

This Grassmannian $Gr(k,n)$ can be subdivided into smaller pieces, with particular subdivisions like the Schubert decomposition having interesting combinatorial structure and connections to enumerative geometry. Positroid varieties give us a finer decomposition than Schubert varieties (in fact, you can think of positroid varieties what one gets from rotating Schubert varieties in the flag and intersecting them). Positroid varieties have very nice combinatorial properties; they can be defined by cyclic rank conditions (i.e. taking $k \times n$ matrices mod $GL(k)$ action on the left and imposing that consecutive columns, say columns $i$ to $i+j$, have rank less than $a$) and can be enumerated by bounded affine permutations which describe juggling patterns \cite{KLS}. In fact, \cite{Post} showed that there are numerous mathematical objects that are in bijection with these bounded affine permutations, including positroids (hence the name), Le-diagrams, Grassmann necklaces, plabic networks (up to equivalence moves), etc. 

We can ask questions about the geometry of these positroid varieties. One natural question to ask is whether they are smooth, or if they have singularities. This thesis answers this question as follows.
First, Theorem \ref{thm:testonnchoosek} shows that if a positroid variety $\Pi_f$ has any singularities, then it must be singular at one of the $n\choose k$ torus-fixed points $\{\lambda\}$; this implies that if we want to know that $\Pi_f$ has at least one singular point, all we need to do is to check on these $n\choose k$ points $\{\lambda\}$. Thus, we calculate the multiplicity (the degree of the tangent cone) of $\Pi_f$ at these points $\lambda$; if the multiplicity is greater than or equal to $2$, then $\Pi_f$ is singular at $\lambda$, if multiplicity $1$ then smooth, and if multiplicity $0$ then $\lambda \not\in \Pi_f$ - c.f. Theorem \ref{thm:smoothone}). We show in Section \ref{sec:smoothness} that we can calculate the multiplicity by calculating the $T$-equivariant cohomology of $\Pi_f$ restricted to a point $\lambda$: $[\Pi_f]|_\lambda$, and then restricting to the appropriate $\C^\times$-equivariant cohomology. We can do this in 2 ways:

(1) We can calculate $[\Pi_f]|_\lambda$ directly. This is shown in Section \ref{sec:compformula}. The formula is:
\begin{align*}
[\Pi_f]|_\lambda &= \sum_{u\in S_n, u\mapsto \lambda} [X_v^w]|_u \frac{[\lambda]|_\lambda}{[uB/B]|_u} \\
&= \sum_{u\in S_n, u\mapsto \lambda} \frac{[X_v]|_u  \,\, w_0\cdot([X_{w_0*w}]|_{w_0*u})}{(\Pi_{i<j\leq k} (y_{u(i)}-y_{u(j)}) \Pi_{k<i< j} (y_{u(i)}-y_{u(j)})) }
\end{align*}
where $w_0 *$performs: $y_i\mapsto y_{n+1-i}$, and $X^w=w_0\cdot x_{w_0w}$ (this is just left-multiplying by the matrix $w_0$ and is effectively taking Schubert varieties to opposite Schubert varieties).

In fact, this formula not only tells us which ones are singular, but what the multiplicities are, which gives a very rough measure of how singular they are at the given points. 

This formula is relatively easy to program. We have code written in Macaulay2 on \url{josephflueg.github.io}. How to use this code is described in Section \ref{sec:usingcode}/ Using this code, we can put together some tables that list which positroid varieties are smooth and which ones are singular (including their multiplicities) for any $Gr(k,n)$. The tables are provided in Chapter \ref{ch:catalogue}.

(2) We show that $[\Pi_f]|_\lambda$ actually is a sum of terms $t_1 +\cdots t_m$, where $m$ is the multiplicity, and each term $t_i$ actually comes from a diagram called an affine pipe dream consisting of cross and elbow tiles (c.f. discussion in Section \ref{subsec:AJSBilley}). An affine pipe dream looks like:

 \begin{figure}[htbp] \centering
	\includegraphics[scale=0.6,clip=true]{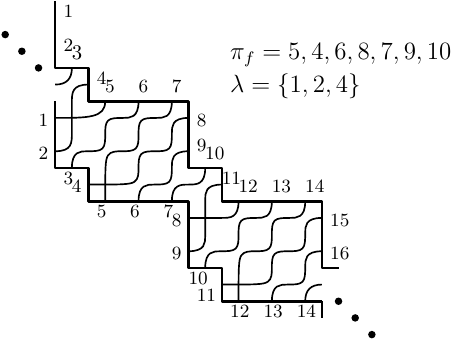}
\end{figure}

Since this is perhaps the central result of this thesis, we pause to provide a quick example, especially for the reader who may not read beyond the introduction. 

\begin{Exa}

First, we spoke about ``multiple pipe dreams" because they are nonunique (\textit{nonrigid}) due to the possibility of \textit{moves}; in other words, each pipe dream gives rise to a permutation by following the pipes from the bottom to the top boundary, but there can be nonuniqueness because sometimes elbow and cross tiles can be exchanged in such a way that the permutation does not change, e.g.:

 \begin{figure}[htbp] \centering
	\includegraphics[scale=0.8,clip=true]{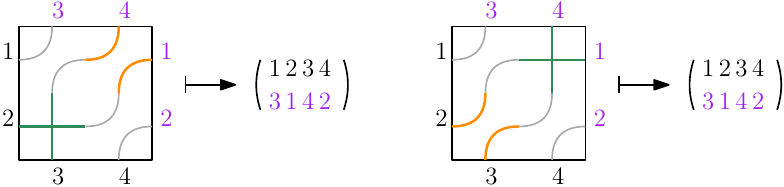}
\end{figure}

We showed that the number of affine pipe dreams counts the multiplicity by showing the the multiplicity equals $[\Pi_f]|_\lambda$ in T-equivariant cohomology $H_T^*$ restricted to $H_{\C^\times}^*$, and proving that we can obtain $[\Pi_f]|_\lambda$ from affine pipe dreams diagrammatically. We next illustrate this in an example:

 \begin{figure}[htbp] \centering
	\includegraphics[scale=0.9,clip=true]{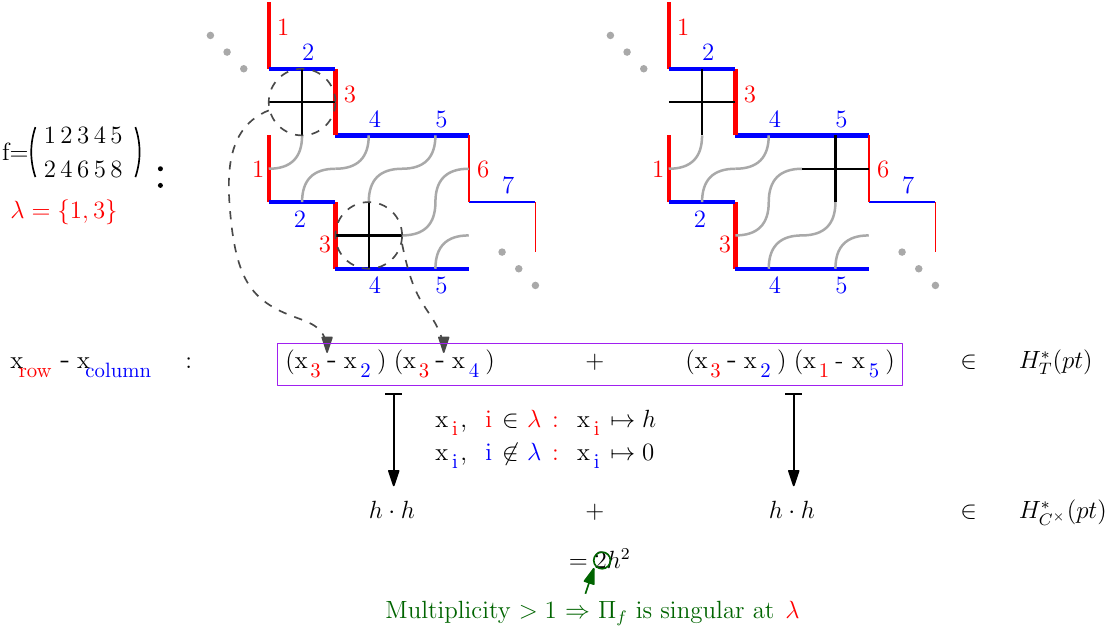}
\end{figure}

In this example, we start with the permutation $f=(24658)$ and we want to consider whether the point $\lambda=\{1,3\}$ is smooth. $\lambda$ gives us the outer boundary (the \textit{shape}) of the pipe dream by drawing the $\lambda$ number as vertical line segments. In this case, there are 2 affine pipe dreams (details of how to draw them are in Section \ref{affinepipedreams}). Each affine pipe dream produces a term $t_i$  in $H_T^*$: this term is a product of $x_{row}-x_{col}$ for each cross in the affine pipe dream. The final restriction to $H_{\C^\times}^*$ simply sends each term to $1\cdot h^{codim(\Pi_f)} \in \Z[h]=H_{\C^\times}^*$, so that the final result is $\# \text{of pipe dreams} \cdot h^{codim(\Pi_f)}$, so we can read out the multiplicity as the coefficient. We get this simply by choosing $\C^\times\subseteq T$ such that the map sends each $x_{row}\mapsto h$ and $x_{col}\mapsto 0$. 

Note: in the figure, the polynomial boxed in purple is the same as what one would obtain by using the formula $\sum_{u\in S_n, u\mapsto \lambda} \frac{[X_v]|_u  \,\, w_0\cdot([X_{w_0*w}]|_{w_0*u})}{(\Pi_{i<j\leq k} (y_{u(i)}-y_{u(j)}) \Pi_{k<i< j} (y_{u(i)}-y_{u(j)})) }$ above.

\end{Exa}

So whereas method (1) is perhaps the superior one for ease of computation, method (2) is superior for intuition and visualization. With (2), we can literally draw out a picture to check for smoothness of $\Pi_f$ at $\lambda$.

In fact, the Main Theorem \ref{thm:mainthm} says that we don't need to draw out all diagrams to check for smoothness - it is sufficient to just draw 1. Then, we check the criterion given in Theorem \ref{thm:mainthm}, and this will tell us whether the point is smooth or not. 

Note that we answer the question of whether a positroid variety $\Pi_f$ has a singularity anywhere: whether there exist some singular point $p\in \Pi_f$. We do not answer the question of what the singular locus $Sing(\Pi_f)$ looks like; that is, we do not consider the set of all singularities $Sing(\Pi_f)$ on a positroid variety $\Pi_f$. The best result that we know of regarding this was stated in my 2024 Research Statement on \url{josephflueg.github.io}, where the following proposition was noted:
 \begin{Prop}
For a given positroid variety $\Pi_f$, the singular locus $Sing(\Pi_f)$ of $\Pi_f$ is a union of smaller positroid varieties: $Sing(\Pi_f) = \cup_i^m \Pi_{f_i}, \forall \Pi_{f_i} \subseteq \Pi_f$  
\end{Prop}

Additionally, we can put an ordering on the set of pairs (positroid, point)$=(\Pi_f,\lambda)$ where we go down in the order via deletion and contraction operations. It may be helpful to visualize what this ordering looks like, as depicted in Figure \ref{fig:atomicorderexa}. The key fact that allows us to connect this ordering with smoothness is Proposition \ref{prop:delsmooth}, which says that if $\lambda$ is smooth on $\Pi_f$, then any further deletions or contractions will continue to keep the point smooth. This implies that there are minimal singular pairs $(\Pi_f,\lambda)$ in the ordering, where any additional deletion or contraction is smooth. We call these minimal singular pairs \textit{atomic positroid pairs} and characterize them using affine pipe dreams. This gives us a (nonconstructive) characterization of a singular point $\lambda$ on a positroid variety $\Pi_f$  as being a pair $(\Pi_f,\lambda)$ that can reach an atomic configuration (a minimal singular pair) under the ordering given by deletion and contraction. 

Finally in Part I, although we work primarily with affine pipe dreams in this thesis, the reader may have prior experience with what we called \textit{finite pipe dreams}, which appear for example in \cite{Post} and \cite{ARW}). In Section \ref{ch:affinefinitepipes}, we give the explicit mapping between the two, showing that finite pipe dreams are embedded in affine ones by letting everything below the distinguished path be filled with crosses. Using this fact, we then show that our results about smoothness from earlier in the thesis (e.g. the main theorem, moves, deletion and contraction on pipe dreams etc.) apply to finite pipe dreams as well. Finally, we write out the proof that finite pipe dreams are on a patch given by $\lambda=I_1$, the first Grassmann necklace element. 

\section{Physics}

One of the original motivations for this research was to make contact with the field of Amplitudes, in particular the material in the book \cite{ABCGPT}. We have attempted to do that in Part II of this thesis.

Perhaps the key quantity that is calculated in quantum field theory is the amplitude $A$. We can think of the amplitude as a quantity $A(p_1,...,p_n)$ depending on $n$ external momenta $p_i$ and their spin (or helicity) such that $\abs{A(p_1,...,p_n)}^2$ gives the probability for this scattering process to occur (the probability that if $p_1$ and $p_2$ are collided, that $p_3,...,p_n$ are detected as emerging). The \textit{theory} is encapsulated by the choice of Lagrangian, an expression which states the particles/fields and what interactions are allowed between them. 

Planar $\mathcal{N}=4$ supersymmetric Yang-Mills theory is a particularly nice toy-model theory. The supersymmetry means that that there are both bosonic and fermionic fields such that the Lagrangian is symmetric under exchange of these fields. The symmetries are captured by their Lie algebra generators; in massless $N=4$, there are 4 sets of such supersymmetry generators (the maximum allowed). There are techniques for calculating amplitudes called unitarity methods - these unitarity methods allow one to obtain the amplitudes in certain theories as a sum over what are called \textit{leading singularities} (more details in Section \ref{sec:leadsing} or the notes on \url{josephflueg.github.io)}. It was discovered that in planar $N=4$ SYM, the leading singularities could be calculated using the following \textit{link formula} for $n$ external momenta with helicity $k$, and $\Gamma_f$ the positroid cell with positroid $f$:
\[L_{n,k}(f) =\oint_{C\subset \Gamma_f} \frac{d^{k\times n} C}{vol(GL(k))}\frac{1}{(1\cdots k)\cdots(n\cdots k-1)} \delta^{k\times 4}(C\cdot \tilde{\eta})\delta^{k\times 2}(C\cdot \tilde{\lambda}) \delta^{2\times (n-k)}(\lambda \cdot C^\perp)\]
In other words, in this theory, leading singularities are calculated by putting in the information of $n$ massless external momenta as well as the matrix $C$, a parameterization of a positroid cell (a positroid variety restricted to its nonnegative real part; we comment a bit on the reason for the nonnegativity in Section \ref{sec:positivereal}). In terms of which particular combinations of positroid cells (that is, which permutations $f$) make up an amplitude, one must go through a recursive method such as the Britto Cachazo Feng Witten (BCFW) relations, or find a set of positroid cells that tile the amplituhedron. 

In Chapter \ref{ch:BCFW}, we discuss the Britto-Cachazo-Feng-Witten (BCFW) relations. We begin with some background, where we derive the BCFW relations. We then turn out attention to the BCFW bridge for on-shell diagrams. We prove a couple of statements regarding the BCFW decomposition method in \cite{ABCGPT}, wring out proofs on how the bridge decomposition algorithm works and that the formula for how the C-matrix changes upon addition of a bridge changes a permutation in the expected way. Finally in this section, we show how to put BCFW bridges on a pipe dream to build the plabic diagram for a permutation $f$; this gives us a way to associate a plabic diagram to a pipe dream. The connection comes from the fact that a pipe dream is essentially a (rotated) wiring diagram, and what the BCFW bridge decomposition method in \cite{ABCGPT} does is it creates a fully crossed wiring diagram, while picking a certain order for adding crossings and keeping track of this order. Thus, we just have to see what that ordering looks like on a pipe dream. The process is not difficult, but it is not simpler than simply running through the bridge decomposition algorithm - the main complication/rule is that when a leg $i$ becomes ``locked" because $f(i)=n$, we may have to eliminate this row, changing the pipe dream and relabeling (a sort of full contraction plus projection, in the language of Section \ref{sec:DelContr}), if there are any non-cross tiles in rows higher up in the pipe dream. Perhaps the key section in this chapter is Section \ref{sec:pipesplabicexample}. It provides a detailed example going all the way from a permutation $f$ and point $\lambda$, to creating the pipe dream, figuring out the order corresponding to the bridge decomposition, and finally adding the bridges.

In Chapter \ref{ch:softlimits}, we discuss soft factors. Soft factors describe how the amplitude changes if one lets one of the external momentum go to zero: in nice cases, the amplitude breaks up into a smaller point amplitude, times a predictable ``soft factor." The notion of sending a momentum to zero seems similar to applying a deletion operation that sends one of the columns to zero. We focus in particular on the notion of inverse soft factors, where a higher-point amplitude can be built up from these smaller soft factors, which has appeared in the physics literature in say Section 16 of \cite{ABCGPT} or \cite{HeZh}. In this chapter, we begin by going to the non-patch-dependent language of the \textit{diagonal diagrams} in the Appendix of \cite{Flu}. We do some basic combinatorics to show that the inverse soft factors are a sort of right inverse to deletion and contraction on positroid varieties: they are a particular choice of going to a larger $n$ positroid variety, such that deletion (or contraction) goes back to the original positroid variety via a 2-step flip on the diagonal diagrams. We describe how the inverse soft factors act on siteswaps. In \cite{HeZh}, the authors describe how to obtain the expression for an on-shell diagram by applying a series of deletions and contractions. Our previous results in Section \ref{sec:DelContr}, that deletion and contraction are intersections with subGrassmannians, show that this process is commutative; however, it was already shown in the cited papers that the addition of inverse soft factors is order-dependent. 

Finally, in Chapter \ref{ch:smoothphysics} we talk about why the singularity structure of positroid varieties does not matter for the calculation of amplitudes. We do a couple of simple examples, showing that an identity for amplitudes in terms of positroid cells does not hold for equivariant cohomology of the positroid varieties restricted to a point, and that the singularities of positroid varieties do not correspond to spurious poles in a simple amplituhedorn example.

The appendices contain a couple of fun small exercises I did during grad school (2023ish), on generalized color orderings diagrams and (attempts at) deletion and contraction on nonplanar on-shell varieties. These probably should not be taken too seriously as their correctness cannot be guaranteed, but are included in case someone might take some inspiration from them.

\section{Other Notes}

Note regarding the overlap with Arxiv submission \url{https://arxiv.org/abs/2407.21116}: the contents of that Arxiv submission from July 2024 are contained in the following sections of this thesis.
From Chapter \ref{chapter:background} on background related to the Grassmannian: Section \ref{sec:positroids}, most of Section \ref{sec:positroidvarieties}, Section \ref{sec:DelContr}, and Section \ref{section:theorder} are in the 2024 submission; and from the following chapter \ref{ch:combinbackgroundpipes} on combinatorial background related to pipe dreams, essentially everything other than \ref{subsec:affinepipesshapemotivation} appears in the previous Arxiv submission. Almost all of Chapter \ref{ch:posvarsmoothness} which deals primarily with the Main Theorem on Smoothness, with the exception of some more detail on the examples in Subsection \ref{subsec:AJSBilley} and the Combinatorial Proof \ref{proof:singDelsingcombinatorial} in Section \ref{sec:smoothnessordering}), is in 2407.21116. Chapter \ref{ch:maintheoremproof} on the proof of the Main Theorem and \ref{ch:atomicpairs} on Atomic Positroid Pairs are basically taken directly from the previous Arxiv submission. In the other direction, essentially everything in \url{https://arxiv.org/abs/2407.21116} is contained in this thesis, except the last section titled ``Skew Partitions and Schubert Varieties" and the proof that the ``diagonal diagrams" give the right deletion/contraction permutation in Appendix A.

A note on history of this thesis: this is a revised and cleaned up version of my thesis, written originally in August 2024. However, that version was a bit hastily thrown together before a deadline, so I immediately began reorganizing and filling in details. A good deal of this embellishment and refactoring was done in Fall of 2024, and an earlier version of this revised thesis was emailed to Thomas Lam in November 2024. Unfortunately I didn't have time to look at it again it until this summer (2026). Essentially all the ideas in the current version (August 2026) are in the earlier versions - the only completely new section is Chapter \ref{sec:compformula} on the Computational Formula for calculating multiplicity and the tables showing the multiplicities of all points for $Gr(k,n)$ up to $n=6$ (a bit ironic since this result was known before most of the other results in this thesis and was used to check these other results, but my notes on this weren't included in earlier versions). Nick Early had the excellent suggestion to include the tables cataloguing the smooth and singular points for all positroid varieties up to a certain $n$ - unfortunately the data generated from my code was too large so was just sitting there for years, until I asked Claude to generate the tables from my data this summer. This is the only place in this thesis where AI was used (other than asking Gemini to generate $\backslash$bibitem citations). Previously I felt that the thesis was too messy to post on Arxiv; although still far from perfect, hopefully this version is a bit better. I hope you enjoy reading it!

\part{Smooth and Singular Points on Positroid Varieties}

\chapter{Combinatorial and Geometric Background Related to the Grassmannian}\label{chapter:background}

\section{Grassmannian}\label{sec:Grassmannian}

\begin{Def}[Grassmannian]
Let $V$ be an n-dimensional vector space. For any integer $0\leq k\leq n$, the \textbf{$\bm{(k,n)}$ Grassmannian} or the \textbf{Grassmannian $\bm{Gr_k(V)}$ of $\bm{k}$-planes in $\bm{V}$} is the set of all $k$-dimensional linear subspaces (through the origin) of $V$. In this thesis, all our vector spaces will be over $\C$, and we will be concerned with $Gr_k(\C^n)$, denoted $Gr(k,n)$ for simplicity. 
\end{Def}

Remark: $Gr_k(V)$ can be given the structure of a smooth manifold of dimension $k(n-k)$.

\begin{Exa}

(1) $Gr(0,n)= \{0\}$ the origin

(2) $Gr(n,n)=\{ \C^n \}$ the whole space (is a point in the Grassmannian)

We have a natural duality: $Gr(k,n) \cong Gr(n-k,n)$ by considering the ``orthogonal complement" $L^\perp$ of $L\in Gr(k,n)$. 

(3) $Gr(1,n)\cong \PP^{n-1}$ projective space

(4) Thus by Grassmannian duality, $Gr(n-1,n)\cong \PP^{n-1}$
  
\end{Exa}
  
We tend to write $Gr(k,n)$ as $GL(k) \backslash [\text{full rank }k\times n \text{ matrices}]$, and think about it as $k$ row vectors spanning a plane up to invertible basis transformation inside the plane.

Let $[n]$ denote the set $\{1,2,...,n\}$ and $\binom{[n]}{k}$ the set of all $k$-element subsets of $[n]$. Each full rank $k\times n$ matrix $M$ has a point in $Gr(k,n)$ given by $GL(k)\backslash M$. Since $M$ is full rank, it has a set of columns $\lambda \in \binom{[n]}{k}$ whose determinant is nonzero, so we can use $GL(k)$ to set $\lambda$ to be the identity; thus, $M \sim \mqty(*&1 & * & 0 & & 0 & * \\ * & 0 & * &  1 & & 0 & * \\ \vdots & \vdots & \vdots & \vdots & \ddots & \vdots & \vdots \\ * & 0 &* &  0 & & 1 & *)$ where the columns that look like $(00\cdots1\cdots 0)^T$ are the ones in $\lambda$.

\subsection{The Standard Open Cover}\label{subsection:standardopencover}

These $\binom{n}{k}$ choices of $\lambda$ enumerate all the fixed points of the Grassmannian under the right torus action, that is, the action of $n\times n$ diagonal matrices by right multiplication. Although an abuse of notation, in the following, we will use $\lambda$ to denote the torus-fixed point on the Grassmannian corresponding to setting the columns $\lambda$ to be the identity (and all other entries zero), in the notation of the previous paragraph. Furthermore, if we let the other columns of $*$'s be anything, so that these denote open sets, then these open sets are isomorphic to $\C^{k(n-k)}$.

Furthermore, these open sets cover the Grassmannian. In other words, $Gr(k,n)$ has an open cover by affine patches (not just affine varieties but things isomorphic to $\A^N$ itself), and since any two patches have nonempty intersection, $Gr(k,n)$ is an irreducible variety of dimension $k(n-k)$.

We will work on these affine patches in the following, and we denote them by $U_\lambda$ where $\lambda$ refers to the $k$ columns we have set equal to the identity.

\begin{Exa}
For example, if we are in the $(k,n)=(2,4)$ Grassmannian, and we let $\lambda=\{1,3\}$, then $\lambda$ denotes the point: $\lambda = \mqty(1&0&0&0\\0&0&1&0)$, and $U_\lambda$ is the open set:

$U_\lambda = GL_k \bigg\backslash \left\{\mqty(1&*&0&* \\0&*&1&*) \right\}$
\end{Exa}

\section{Matroids}

Matroids were created to capture the combinatorial essence of linear dependence. There are many characterizations/equivalent definitions. Here's one:

\begin{Def}[Matroid]
A matroid $M \subseteq {[n] \choose k}$ satisfies the {\em exchange axiom}:

For all $I,J \in M$ and all $i\in I$, there exists a $j\in J$ such that $(I\bs \{i\}) \cup \{j \} \in M$. 

\end{Def}

\begin{Exa}

Take $n=3, k=2$. The following is a matroid $M = \{1,2\}, \{1,3\}$.

For example, if $I = {1,2}$ and $J = {1,3}$:

$i=2: I\bs i=\{1\}$, and let $j=3\in J$, then $\{1\}\cup \{3\}=\{1,3\}\in M$

$i=1: I\bs i=\{2\}$, and let $j=1\in J$, then $\{2\}\cup \{1\}=\{1,2\} \in M$

\end{Exa}

\subsection{Realizable Matroids}
  
We are specifically interested in realizable matroids.

\begin{Def}
From a matrix $A$ we can make a matroid $M$: the elements of $M$ are all the sets of columns of $A$ that give a nonzero maximal minor. We can describe such matroids as $M=(E,B)$, a ground set $E$ and bases $B$ making up the matroid. If a matroid arises this way, it is called \textbf{realizable}.
\end{Def}

\begin{Exa}

We can get a realizable matroid from the matrix 
\[
\begin{blockarray}{ccccc}
1 & 2 & 3 & 4 & 5 \\
\begin{block}{(ccccc)}
  1 & 0 & 2 & 0 & 1  \\
  0 & 0 & 0 & 1 & 1 \\
\end{block}
\end{blockarray}
 \]
 
 The ground set $E=\{1,2,3,4,5\}$ for the columns. This realizable matroid has the bases: $\{ \{1,4\}, \{1,5\}, \{3,4\}, \{3,5\}, \{4,5\} \}$.
 
 \end{Exa}

\begin{Thm}
Realizable matroids are matroids (satisfy the exchange axiom in the definition above). 
\end{Thm}

\begin{Proof}
We can take a matrix $A\in Gr(k,n)$, and we need to show that the columns of nonzero maximal minors do indeed satisfy the exchange relation given above.

$I,J\in M \Leftrightarrow \Delta_I(A)\neq 0$ and $\Delta_J(A)\neq 0 \Leftrightarrow \Delta_I \cdot \Delta_J \neq 0$.

But we have $\Delta_I \cdot \Delta_J = \sum_{j\in J}\pm \Delta_{(J\bs j)\cup i} \cdot \Delta_{(I\bs i)\cup j}$ by the Plucker relations with one element $i$. Thus this sum must have at least one nonzero element which means there is a $j$ such that $(I\bs i)\cup j\in M$.
\end{Proof}

\subsection{Matroids and Grassmannians}

\begin{Def}[Stratification]
A stratification is a collection of subvarieties such that the intersection of any two is a union of others. Then we can talk about stratification generated by some subvarieties: if we repeatedly intersect and decompose them, then we get smaller subvarieties in the stratification. 
\end{Def}

\begin{Def}[Matroid ``Strata"]
For a given realizable matroid $M$, we can define its \textit{matroid stratum} $S_M$ as all points in $Gr(k,n)$ such that a matrix representative has nonvanishing maximal minor precisely for the bases in $M$.
\end{Def}

Note that the matroid strata $S_M$ of a realizable matroid $M$ does not actually satisfy the definition of stratification above. However, the positroid strata \textit{do} form an actual stratification, which we will see later.

\section{Positroids}\label{sec:positroids}

\begin{Def}[Positroid \cite{Post}]
For a given $Gr(k,n)$, if we look at all matroid strata, and intersect with the totally nonnegative (TNN) real Grassmannian (all maximal minors real and nonnegative), then the matroids whose strata are nonempty are called \textit{positroids}.
\end{Def}

\begin{Def}[Bounded Affine Permutation]
A permutation $f:\Z \to \Z$ is called \textbf{affine} if it satisfies the periodicity condition $f(i+n)=f(i)+n$, $\forall i \in \Z$ where $n$ is the period.

If the affine permutation also satisfies $i\leq f(i)\leq i+n$, it is called \textbf{bounded}.
\end{Def}

Jugging patterns are a way to draw bounded affine permutations. The following gives an example:

\begin{figure}[htbp] \centering
	\includegraphics[scale=0.8,clip=true]{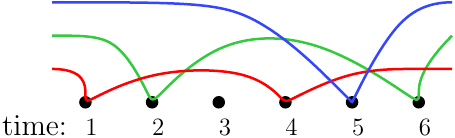}
	\caption{}
\end{figure}

Each dot is a moment in time: the first dot is the time at the start, let's call it $t=1$, the second dot is at $t=2$. This juggling pattern repeats every 6 seconds. There are 3 balls, one for each color. We can think of the dots as the hand; note that this is a simplified model in which there is only one hand.

Note that this looks somewhat like a permutation. We could describe it with the following list of numbers:

\begin{figure}[htbp] \centering
	\includegraphics[scale=0.8,clip=true]{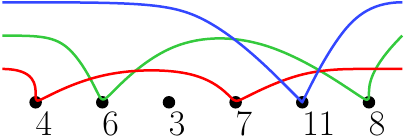}
	\caption{}
\end{figure}

We have written below each point the moment in time when the ball thrown at that moment lands. So the first point has a $4$ because the (red) ball thrown at $t=1$ lands at time $t=4$.  

We can write it like: $\mqty(1 & 2 & 3 & 4& 5 & 6 \\ 4 &6 &3 & 7 &11 &8)$. We call this a bounded affine permutation.

\begin{Exa}
In the juggling pattern given above, $g=(4,6,3,7,11,8)$ defines a bounded affine permutation of period $6$ sending $g(1)=4, g(2)=3,...,g(6)=8$ and repeating, so $g(7)=g(1)+6=10$, etc.
\end{Exa}

The information contained in a bounded affine permutation $f$ can be given a different, but equivalent, labelling called its \textit{siteswap}.

\begin{Def}[Siteswap]
The siteswap $f$ corresponding to a bounded affine permutation $g$ is $f(i):=g(i)-i, 1\leq 1 \leq n$. Notice that a siteswap $f(i)$ is only defined for $i\in[1,n]$.
\end{Def}

\begin{Exa}
The siteswap corresponded to $g=(4,6,3,7,11,8)$ is  $f=(3,4,0,3,6,2)$.
\end{Exa}

Remark: The data of juggling can also be equivalently formulated in terms of Grassmann necklaces (see \cite{Post}), which are $k$-element subsets $I_i$, and can be considered as a moving ``snapshots" which give the (mod $n$) times at which the balls in the air at time $i$ land back in the hand. Using the language of Grassmann necklaces, it is not difficult to map between siteswaps and ``positroids," which are matroids (collections of full-rank nonvanishing columns) corresponding to matrices where all minors are nonnegative (the ``positive Grassmannian"). This is where the terminology ``positroids" comes from. We will say a few more words about this at the end of the  section on positroid varieties.

\section{Positroid Varieties}\label{sec:positroidvarieties}

Now let's connect positroids to matrices and Grassmannians.

\begin{Def}\label{matrixtositeswap}
We can define a map from full rank $k$-by-$n$ matrices $M$ to bounded affine permutations: we take each column of a matrix, and assign the earliest next column for which it's in the span of the following: 
\[ \tilde{g}_M (i) := \min \{ j\geq i\colon \vec c_{i \bmod n} \in span(\vec c_{i+1 \bmod n},\ldots,\vec c_{j\bmod n}) \}.\]
 We only have to make a slight modification for siteswaps, in which case the map is: 
 \[ \tilde{f}_M (i) := \min \{ j\geq i\colon \vec c_{i \bmod n} \in span(\vec c_{i+1 \bmod n},\ldots,\vec c_{j\bmod n}) \}-i \] 
\end{Def}

\begin{Exa}

As an example, consider the following matrix:

$M=\mqty(
1 & 0 & 0 & 1 & 0 & 0 \\
0 & 1 & 0 & 0 & 0 & 1 \\
0 & 0 & 0 & 0 & 1 & 0 \\
)$

Its associated bounded affine permutation is: $g=4,6,3,7,11,8$. We obtain, for example, $\tilde{g}_M(1)=4$ because the first column is not in the span of the 2nd and 3rd columns, but only in the span of columns 2 through 4. 

\end{Exa}

Clearly, this map is many-to-one with a very large family of matrices mapping to a given bounded affine permutation or sitewap in general.

Positroid varieties denoted $\Pi_f$, which we will precisely define momentarily, are subsets of Grassmannia defined by the cyclic rank conditions imposed by siteswaps $f$, or equivalently, bounded affine permutations $g$ ($f$ might also denote the bounded affine permutations, but should be clear from context). It is essentially the reverse of the map we just defined.

\begin{Exa}

As an example, we can take the siteswap considered above: $f=340362$ (or bounded affine permutation $g=4,6,3,7,11,8$). The following matrix gives a point within $\Pi_f$: 

$\mqty(
1 & 0 & 0 & 1 & 0 & 0 \\
0 & 1 & 0 & 0 & 0 & 1 \\
0 & 0 & 0 & 0 & 1 & 0 \\
)$.
\end{Exa}

Consider the following figure. 

\begin{tikzcd}
& \{(y,Y^0) : y\in Y^0\}  \ar[dl, leftrightarrow, "bijective"'] \ar[dr]  &\\
X=Gr(k,n) \ar[rr] & &  \{ \text{siteswaps } f | avgf(i)=k \}  
\end{tikzcd}

\vspace{0.5cm}

Although we are primarily interested in the closed strata because only these have singularities, the figure shows the locally closed strata $Y^0$ since these map bijectively onto the whole Grassmannian $X$: they stratify $X$. At the top we have an incidence variety with some point $y$ chosen on each open positroid variety $Y^0$. These $y$'s are in bijection with $f$'s: the arrow to the southeast maps each $ \{(y,Y^0) \}$ to the unique siteswap $f$ corresponding to the $Y^0$. 

We now discuss the horizontal map along the bottom of the diagram. Going rightwards (the direction of the arrow), the map consists of taking any matrix representative $M$ of a point in $Gr(k,n)$ and doing the computation in Definition \ref{matrixtositeswap}. This produces the siteswap $f$ corresponding to the open positroid variety $Y^0$ in which the point $M$ lives. But suppose we want to go in the opposite direction and take the pullback. Then, the fiber over a given $f$ of this map defines an open positroid variety $\Pi_f^0$ (which is equal to one of the $Y^0$) for a given siteswap $f$ (or bounded affine permutation or positroid). 

\begin{Def}[Positroid Variety]
For a given siteswap $f$, we can define the positroid variety $\Pi_f$ using the map $\tilde{f}_M(i)$ from matrices to siteswaps we gave above:
\[\Pi_f = \{M\in Gr(k,n) \ | \ \tilde{f}_M(i) \leq f(i) \}\]
(The open positroid variety is $\Pi_f^\circ = \{M\in Gr(k,n) \ | \ \tilde{f}_M(i) = f(i) \}$, and these correspond to the $Y^0$ in the diagram above.)
\end{Def}
This completes the discussion of the fiber of the bottom horizontal map in the commutative diagram above. 
\\

We can give a second definition of positroid varieties using the ranks of all cyclic columns. Suppose given such a matrix $M$, we extend it periodically so that $[i,j]$ can be any (but at most quantity $n$) consecutive columns, since dealing with a bounded affine permutation $g$ (rather than a siteswap) is easier in the discussion below. In other words:

Let $[i,j] := \{i, i+1,...,j\}(\bmod n)$, so $| [i,j] | =\begin{cases} j-i+1 & j\geq i \\ (n-i+1)+j & i>j\end{cases}$. 
\\

Define $a_{g,i,j}:=\#\{[i,j]\backslash \{g(i+\mathbb{N}) \} \cap \{\leq j\} \}$. This can be visualized as follows. Draw a permutation matrix for the bounded affine permutation $g$. For say $g(\tilde{i})=\tilde{j}$, the $\tilde{i}$ axis runs vertically downwards in the figure below, while $\tilde{j}$ axis runs horizontally to the right. 

\begin{figure}[htbp] \centering
	\includegraphics[scale=0.6,clip=true]{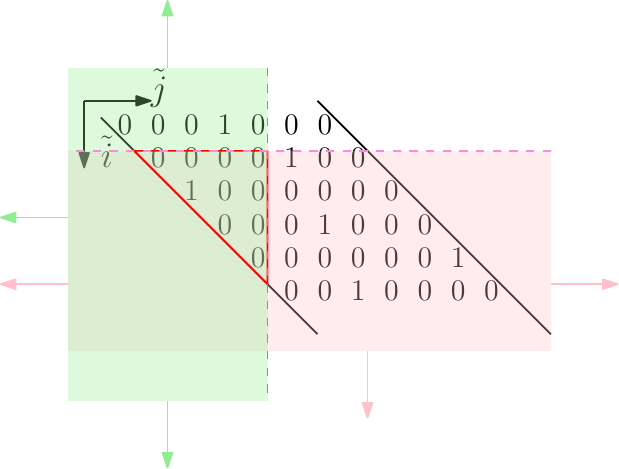}
	\caption{}
\end{figure}

Since we are dealing in particular with bounded affine permutations, the permutation matrix will be supported entirely within a strip of width $n+1$ running from northwest to southeast. The example just shown depicts the bounded affine permutation $g=(4,6,3,7,11,8)$ from before for $[i,j]=[2,5]$. The part shaded pink (which continues going south) depicts $i+\mathbb{N}$, so $\{g(i+\mathbb{N}) \}$ consists of the values of $\tilde{j}$ for the part in pink: here, $\{g(i+\mathbb{N}) \}=\{6,3,7,11,8,...\}$. The part shaded green (which continues going west) depicts $\{\leq j\}$. In this example, $j=5$, so only $\{3\}$ is in the intersection. Thus, $a_{g,i,j}=\#\{3\}=1$. An easier way to state the value of $a_{g,i,j}$ is that it is the number of $1$'s inside the red triangle. If the intuition is helpful carrying over the juggling analogy, the 1's inside this triangle are the number of balls that are thrown inside the time range $[2,5]$ that fall back into the hand and are re-thrown within this time interval. Mathematically, these are the columns (within columns 2 through 5) that are linearly dependent on previous (earlier) columns (in 2 through 5). 

Thus, we have: $rank_g[i,j]=\abs{[i,j]}-\{\# 1's \text{ to SW} \}=(j-i+1)-a_{g,i,j}$. 

``SW" means ``southwest." Now, we can define the positroid variety corresponding to a siteswap $f$ (we apply $\bmod n$ to the definitions using $g$ above):
\\

\begin{Def}[Positroid Variety]
\[\Pi_f := GL_k \backslash \{M\subseteq M_{k,n} \big\vert rankM=k \text{ and } rank_M([i,j]) \leq rank_f[i,j] , \  i\leq j \leq i+n\} \]
where $rank_M([i,j])$ is defined cyclically.
\end{Def}

Note that this is the second definition of the closed positroid variety. If we simply imposed equality ($rank_M([i,j]) = rank_f[i,j]$), then we would obtain the stratum $\Pi_f^0$. 
\\

We will not prove this next theorem, but will state it for completeness of terminology:

\begin{Thm}\label{thm:densestratum} [KLS]
Positroids are in bijection with siteswaps. More specifically, the bijection is that within the positroid variety $\Pi_f$ of each siteswap $f$, there is a single dense matroid (positroid) stratum $S_M$ contained within. They are not equal (in fact, the matroid strata do not form a stratification). However, their closures are equal. Hence, we can equivalently index the positroid varieties with positroids rather than siteswaps. 
\end{Thm}

\section{Positive (Real) Part of Grassmannian}\label{sec:positivereal}

 Given a point in the Grassmannian, one can write down the column-sets (say, $\{m_1,...,m_j\}$ whose maximal minors are nonzero: this gives us a collection of $k$-element subsets: this is a construction of a (realizable matroids) - those matroids that come from a point in the Grassmannian. The matroid strata consists of all matrices in $Gr(k,n)$ with this set $\{m_1,...,m_j\}$ of maximal minors nonzero. Note that this is a combinatorial invariant of the Grassmannian because it is invariant under action of GL(k): this action changes the determinants, but not which are zero or nonzero. Note that here are non-realizable matroids that fit the combinatorial definition but are bad from the linear point of view. We ignore those, focusing on realizable matroids. The map we just described of taking the $k$-element column subsets giving nonzero minors, gives us a (discontinuous constructible map) from the Grassmannian to a finite set of realizable matroids, which therefore gives us a decomposition of the Grassmannian into locally closed pieces which are the fibers of that map. 

A theorem of \cite{Mnev} showed that these fibers are not nice at all: the schemes we get are locally isomorphic to any other scheme one can define with integer coefficients, so in order to understand the matroid decomposition, one must understand all schemes (Mnev's Universality Theorem).
The proof uses rather simple compass and straightedge constructions to do multiplication, addition, and square roots, etc.; in this case, one can implement all the algebraic geometry one needs to do with just point and line configurations and thereby construct these very difficult-to-deal-with sets in the Grassmannian. 

The idea of \cite{Post} was to restrict our consideration to the matroid strata in the Grassmannian that meets the positive Grassmannian: that is, on which of these strata can one actually find positive real points: the matroids that have positive real points are called positroids. We can define the positive (or nonnegative) Grassmannian $Gr(k,n)$ as being the subset of the Grassmannian given by full-rank $k\times n$ matrices whose (maximal) minors are all positive or zero. Note that the introduction of nonnegativity (inequalities) takes one outside pure algebraic geometry since semi-algebraic sets are allowed.

A priori, it's not obvious that one can find real positive points, or even real points. For example, the scheme $\R[x]/(x^2+1)$ has no real points; it only has 2 complex-valued points. Thus, in this example, one can find a matroid that has no points in the real Grassmannian, only in the complex Grassmannian. However, \cite{Post} not only restricts to the real Grassmannian, but further restricts to the nonnegative Grassmannian where all minors are $0$ or positive real numbers. This turns out to be very restrictive. For example, if we consider the first Plucker relation: $p_{13}p_{24}=p_{12}p_{34}+p_{14}p_{23}$, over the complex numbers, if we knew one of each of the terms in the sums on the right is zero, then one of $p_{13}p_{24}$ is also zero - some terms being equal to zero imply other terms being zero. 
Now, let's set $p_{13}=0$ on the left side of the equation. In the complex case, this could happen because the terms on the right cancel each other. However, in the positroid case, because every nonzero term is positive, we must have both terms on the right being zero. This is like a converse to the previous statement: for the matroid, having one $p$ equal zero in each of the terms on the right implied one of the $p$'s on the left is zero; now for the positroids, one of the $p$'s on the left being zero implies that one $p$ equals zero in each of the terms on the right. 
 
When one considers this restriction to the nonnegative Grassmannian, many matroids do not work (they are not positroids). The ones that do have real positive points turn out to be very nice: they are parametrizable, isomorphic to $\R^n$ for some $n$ \cite{Post}. This is very nontrivial, since even a smooth curve of genus 0 cannot be algebraically parameterized. In fact, the equations defining positroid varieties are surprisingly simple:  just given by certain Plucker coordinates vanishing (with a minor vanishing if it’s not one of the fixed points on the positroid variety \cite{KLS} (one can write the cyclic rank conditions in terms of the vanishing of a set of Plucker coordinates). Thus, if one wants visualize the positroid variety by mapping the Grassmannian to Plucker space with the Plucker embedding, the equations defining positroid varieties are linear, given by intersecting the image of the Grassmannian in Plucker space with a linear subspace. Note that a priori, one might have thought that if $P$ is a positroid, the matroid stratum corresponding to $P$ has multiple components, but this turns out to not be the cause: positroid varieties are irreducible already in the complex Grassmannian. In other words, if we take the matroid stratum corresponding to a positroid $P$ and take its Zariski closure inside the complex Grassmannian, it is irreducible so it has a single component - it is a Zariski dense set (this is the statement that ``there is a single dense matroid (positroid) stratum $S_M$" in the previous theorem).

Here is the picture I have in mind for the preceding theorem (\ref{thm:densestratum}):

 \begin{figure}[htbp]  \centering 
	\includegraphics[scale=0.9,clip=true]{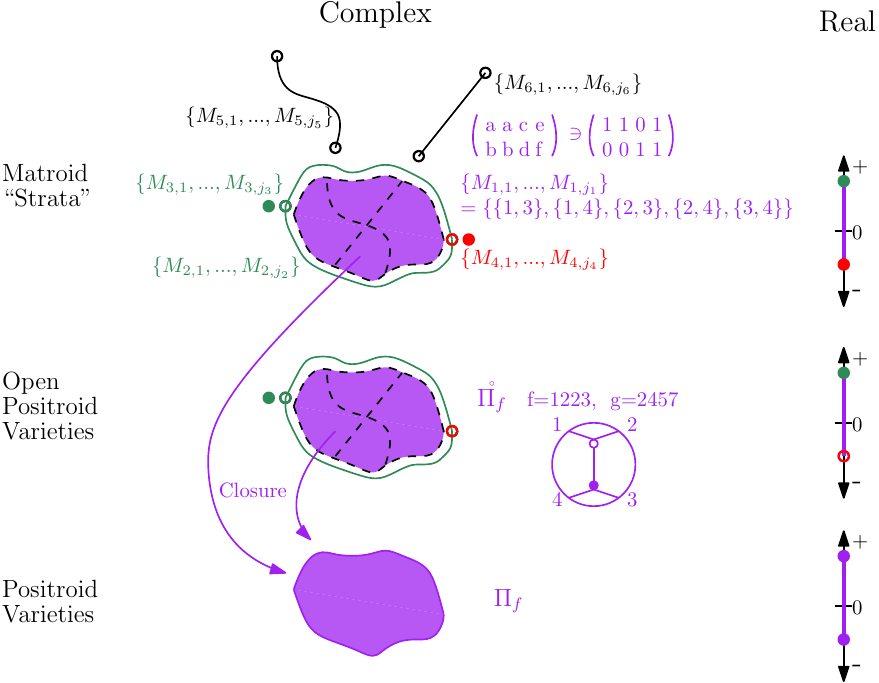}
	\caption{}
\end{figure}

This figure should be taken as suggestive, rather than precise (we are dealing with the Zariski topology after all, so the sets don't look like this), but it captures a few essential features of the matroid vs. positroid relationship. 

The positroids are the sets of columns with nonvanishing determinant when one considers only $k\times n$ matrices representing points in $Gr(k,n)$ whose maximal minors are all nonzero. One key point from this figure is that the matroid decomposition is much finer. We are considering the nonvanishing minors $M_i=\{M_{i,1},...,M_{i,j} \}$, where each $M_{i,l} \in {[n]\choose k}$ is a $k$-element subset of $[n]$, and this matroid has $j$ such $k$-element subsets. In general there will be more such matroids $M_i$ for matrices with complex entries allowed, than if we restricted ourselves to only looking at matrices with all minors nonnegative (in the latter case, this implies that all entries have to be real, since we can use a $GL(k)$ transformation to set a $k\times k$ block equal to the identity). If a (realizable) matroid $M_i$ contains another matroid $M_j$ ($M_i \supset M_j$), this implies the matrices/points in $M_j$ are a subset of those in $M_j$ and $M_j$ lives on a lower dimensional boundary of $M_i$; for example in $Gr(2,4)$ if we have the matroid $M_j = \{ \{1,2\} \}$, this has dense representation given by matrices of the form $\mqty(* & * & 0 & 0\\ *  & * & 0 & 0)$, and if $M_i = \{ \{1,2\},\{1,3\}, \{2,3\} \}$, this has matrices of the form $\mqty(* & * & * & 0\\ * & * & * & 0)$. 

For the previous figure, we take a positroid variety $\Pi_f$ for a siteswap like $f=1223$ (bounded affine permutation $g=2457$) in $Gr(2,4)$. This is the bottom row. At the top, we consider the matroid ``strata," consisting of all matroids in the complex Grassmannian that intersect $\Pi_f$. The point is that by Theorem \ref{thm:densestratum}, there is a unique dense matroid stratum, whose closure is equal to $\Pi_f$: we have drawn this in purple as the 2D piece corresponding to $M_1$. This clearly has real positive points: we give the example of the matrix $\mqty(1&1&0&1\\0&0&1&1)$ which has all maximal minors nonnegative. $M_4$, $M_5$, and $M_6$ are matroid strata that do not have positive real points. ($M_1$ is irreducible, so the picture might be misleading here since it seems to show $M_5$ and $M_6$ disconnecting the space, which does not happen in the Zariski topology) 

In green, we have drawn $M_2$ and $M_3$ to suggest lower dimensional positroid varieties living on the boundary of $\overset{\circ}{\Pi_f}$. This is a place where the picture is not perfect - in the case of $\Pi_f, f=1223$, there are multiple positroid varieties in the next dimension lower positroid varieties on the boundary of $\Pi_f$. We have tried to include two open points in the green and red circles to suggest that there may be strata for this lower-dimensional positroid variety as well, some of which might corresponds to positroids and some which do not (in this case, the green circle does correspond to a positroid but the red does not). (The picture is also imperfect here since we have drawn a disconnected set in $\R$ coming from restricting a connected set). As such, when we restrict the open positroid varieties to the real Grassmannian, which we have tried to draw in the right column, there may be points on the boundary that are missing (here, given by the red dot) since they correspond to non-positroid matroids. Of course, by theorem \ref{thm:densestratum}, all these points are there, and everything is equal when one considers the closures, as in the bottom row.

\section{Smoothness}

We are interested in determining whether a positroid variety is singular or smooth. Intuitively, it is the notion that tangent spaces are not ``too big":

 \begin{figure}[htbp]  \centering 
	\includegraphics[scale=0.5,clip=true]{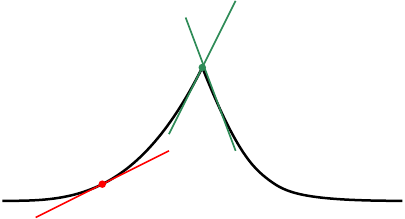}
	\caption{}
\end{figure}

In this picture, the manifold drawn in black is 1-dimensional, and the red point is a smooth point since the tangent space is the right dimension, but the tangent space at the green point is spanned by the two lines, so it is too big.

For completeness, we give the following definition of smoothness:

\begin{Def}
Suppose we are on an affine variety $X=Spec R$ of a ring $R$, and we have a point $a \in Spec R $ defined by a maximal ideal $I(a)$. Then $X$ is smooth at $a$ when the maximal ideal has a number of generators equal to the dimension of the ring, that is, $\dim I(a)=\dim O_{X,a}$
So a point $a\in X$ is smooth if its ring of local functions $A(X)_{I(a)}$ is regular, that is, if $\dim T_aX=\dim Hom_K(T_aX,K)=\dim I(a)/I(a)^2=\text{codim}_X\{a\}=\dim X=\dim A(X)_{I(a)}$. Equivalently, $X$ is smooth at $a$ if its tangent cone is the same as its tangent space at $a$.
\end{Def}

We often test smoothness using some sort of Jacobi condition, e.g.:

\begin{Thm} (Affine Jacobi Criterion)
Let $a\in X$ be a point on an affine variety $X\subseteq \mathbb{A}^n$ and let $I(x)=\langle f_1,...,f_n\rangle$. Then $X$ is smooth at $a$ $\Leftrightarrow$ the rank of the $r\times n$ Jacobian matrix $(\frac{\del f_i}{\del x_j}(a))_{ij}$ is at least $n-codim_X\{a\}$ (so it's equal to $n-codim_X\{a\}$)
\end{Thm}

\section{Deletion and Contraction (of Positroid Varieties)}\label{sec:DelContr}

The Grassmannian is represented by full-rank $k$-by-$n$ matrices (up to $GL(k)$ action on the left).

Then there is a circle action on the right via matrix multiplication by 
\[S_i=\mqty(\dmat{1,1,\ddots,1,z,1,\ddots,1})\]
 such that this scales $i$-th columns of $\{M\}$ by $z$. 
 \\
 
There is a ``Pascal recurrence" on Grassmannians:

\begin{Prop}

Let $S_i$ act on $\{M\}=Gr(k,n)$ on the right. 

Let

\[\Pi_{contr,i} := rowspan \left \{ \mqty(&&&0&&&\\  &*&&0&&*&\\ 0&\cdots&0&1&0&\cdots&0) \right \} \cong Gr(k-1,n-1) \] 

and

\[\Pi_{del,i} := rowspan \left \{ \mqty(&&&0&&&\\  &*&&0&&*&\\ &&&0&&&) \right \} \cong Gr(k,n-1) \]

where the rest of the matrix is full-rank. 

Then the fixed point set under the $S_i$ action is: 

\[Gr(k,n)^{S_i}= \Pi_{contr,i} \sqcup \Pi_{del,i} \cong Gr(k-1,n-1) \sqcup Gr(k,n-1).\]

\end{Prop}

Because of the extensive use of deletion and contraction in this thesis, we next explain all the details of this Pascal recurrence, including the proof of the above proposition, but the familiar reader can skip it without loss of continuity in reading this thesis. 
\\ 
 
If we consider the elements of the Grassmannian that are invariant under the $S_i$ action, there are two types:
\\

(1) The first type consists of all points in the Grassmannian that are planes which contain the $i$ coordinate. In other words, this is all $k$-by-$n$ matrices of the form 
\[\mqty(&&&*&&&\\  &&&&&&\\ 0&\cdots&0&1&0&\cdots&0)\] 
(or can be brought to such a form by a $GL(k)$ action on the left; furthermore, by an additional $k$-element permutation matrix, we can move this row to the bottom row). In addition, since all the other entries of the bottom row are zero, we can do further elementary matrix operations such that the matrix representatives all have the form: 
\[\mqty(&&&0&&&\\  &&&0&&&\\ 0&\cdots&0&1&0&\cdots&0),\]
 that is, the $i$-th column also has nothing but zeros above the bottom row. What this says is that if our point in the Grassmannian is a plane containing the $i$-th coordinate vector, then by a change of basis, we can pick a basis of $k$ vectors spanning the plane such that one of the vectors lies along the $i$ coordinate and the remaining $k-1$ vectors lie in the orthogonal complement of the $i$-coordinate. Thus, we can then project to the smaller Grassmannian $Gr(k-1,n-1)$ by removing the $i$ coordinate: this can be seen as a two-step process where we first remove the bottom row (or last vector) of the matrix, leaving us with a $k-1$-plane in $n$ dimensions, and then we project to the orthogonal complement of the $i$-th coordinate, which leaves a $k-1$-plane in the remaining $n-1$ dimensions. We denote the set of matrices 
 \[\Pi_{contr,i} := \mqty(&&&0&&&\\  &&&0&&&\\ 0&\cdots&0&1&0&\cdots&0)\] 
which corresponds to the siteswap $(k-1)^{n-k}k^{k-1}n$ rotated so that the $n$ is in the $i$th position. Note that in this matrix there are no restrictions on the rest of the $k-1$-by-$n-1$ matrix, as long as it is full rank, so $\Pi_{contr,i}\cong Gr(k-1,n-1)$. 
 \\
 
The process called ``contraction" has two steps. First we first intersect with $\Pi_{contr,i}: Contr_i(\Pi_f)=\Pi_f\cap \Pi_{contr,i}$. This restricts $\Pi_f$ to the subset that can be written in the form 
  \[\mqty(&&&0&&&\\  &&&0&&&\\ 0&\cdots&0&1&0&\cdots&0).\]
 Secondly, we do the projection that removes the last row and $i$-th column to bring it to a subset of $Gr(k-1,n-1)$. This intersection with $Contr_i(\Pi_f)$ turns out to be the largest positroid variety with $f(i)=i+n$ contained in $\Pi_f$ (the proof is below). As long as a given positroid variety $\Pi_f$ contains at least one point whose $k$-plane does not use the $i$-th coordinate, the contraction of this positroid variety will be nonempty. 
 \\
 
\emph{Remark}: there is sometimes an abuse of terminology where ``contraction" can also refer to only performing the first step, without the projection to $Gr(k-1,n-1)$. We define contraction as performing both steps, and when only the first step is involved, we call it \textit{projectionless contraction}. Similarly for \textit{projectionless deletion} below. 
 \\
 
\emph{Remark}: note that contraction (and similarly for deletion, which we describe next) is an operation that is performed on the entire Grassmannian $Gr(k,n)$ (in fact $Gr(k,n)=\Pi_f'$ with $f'=kk\cdots kk$), but since we will be working with individual positroid varieties throughout this paper, our notation here emphasizes the effect of contraction (or deletion) on a specific positroid variety $\Pi_f$. 
\\

(2) The second type consists of all points in the Grassmannian that are planes that are already contained within the $n-1$-dimensional complement of the $i$th-coordinate. Here, the matrices are of the form: 

\[\Pi_{del,i} := \mqty(&&&0&&&\\  &&&0&&&\\ &&&0&&&) \]
 since the $i$-the coordinate is not used at all. For the subset of the $Gr(k,n)$ Grassmannian that can be represented by matrices of this form, we can project to an isomorphic Grassmannian by removing the $i$-th column, and this will simply be $Gr(k,n-1)$ since there are no restrictions on the other entries except that they be full rank. Note that a difference with the previous case of contraction is that here, the entire $k$-plane already lives in the complement (it does not contain the $i$-coordinate vector at all), so the elements are still $k$-planes after projection to the complement of the $i$-th coordinate, not $k-1$-planes. 
\\
 
Similarly, the process of ``deletion" has two steps: we first intersect with $\Pi_{del,i}: Del_i(\Pi_f)=\Pi_f\cap \Pi_{del,i}$ where  $\Pi_{del,i}$. This restricts $\Pi_f$ to the subset that can be written in the form 
 
 \[\mqty(&&&0&&&\\  &&&0&&&\\ &&&0&&&).\] 
 Then we do the projection that removes the $i$-th column to bring it to a subset of $Gr(k,n-1)$. This intersection with $Del_i(\Pi_f)$ turns out to be the largest positroid variety with $f(i)=i$ contained in $\Pi_f$. As long as a given positroid variety $\Pi_f$ contains at least one point whose k-plane contains the $i$-th coordinate, the deletion of this positroid variety will be nonempty. 
\\

There is another description of the deletion and contraction of positroid varieties.
For positroid varieties of $Gr(k,n)$, deletion and contraction of the $i$th column of a positroid variety $\Pi_f$ can be described as:

\begin{Prop}\label{prop:del2def}
The projectionless deletion $Del_i(\Pi_f)$ of the $i$-th column of a positroid variety $\Pi_f$ is equal to $\Pi_{f'}^{d(i)}:=$ the largest positroid variety $\Pi_{f'}$ contained in $\Pi_f$ such that $f'(i)=0$. (By largest, we mean that $\Pi_{f'}$ is not properly contained inside another positroid variety $\Pi_{f''}$ properly contained in $\Pi_f$.) 

The projectionless contraction $Contr_i(\Pi_f)$ of the $i$-th column of a positroid variety $\Pi_f$ is equal to $\Pi_{f'}^{c(i)}:=$ the largest positroid variety $\Pi_{f'}$ contained in $\Pi_f$ such that $f'(i)=n$.
\end{Prop}

Recall that the projectionless contraction was defined as the intersection: $Contr_i(\Pi_f)=\Pi_f \cap \Pi_{contr,i}$, and similarly the projectionless deletion is  $Del_i(\Pi_f)=\Pi_f \cap \Pi_{del,i}$.

Before we prove Prop \ref{prop:del2def}, we prove another proposition that we will need in the proof. 

\begin{Prop}\label{prop:delirred}
The intersections $Contr_i(\Pi_f)=\Pi_f \cap \Pi_{contr,i}$ and $Del_i(\Pi_f)=\Pi_f \cap \Pi_{del,i}$ are irreducible.
\end{Prop}

\begin{proof}
Because the positroid varieties form a stratification, intersections of positroid varieties are in general unions of positroid varieties, but we show now that in this case of $Contr_i(\Pi_f)=\Pi_f \cap \Pi_{contr,i}$, this intersection is irreducible (when it's nonempty), which means that $Contr_i(\Pi_f)$ cannot split into multiple positroid varieties - it must be a single positroid variety. The intersection $\Pi_f \cap \Pi_{contr,i}$ is irreducible because we have 
\[\Pi_f \hookleftarrow \Pi_f \backslash \Pi_{del,i} \twoheadrightarrow \Pi_f \cap \Pi_{contr,i}\] 
The latter map is a surjection because we have a right inverse given by taking a plane (point) in $\Pi_f \cap \Pi_{contr,i}$ and taking its direct sum $\oplus \C$, a copy of $\C$ lying along the $i$th coordinate. It may be even easier to argue that since $\Pi_{contr,i}$ consists of planes containing the i-th coordinate, and $\Pi_{del,i}$ consists of planes in the orthogonal complement of the i-th coordinate, they are disjoint, so removing $\Pi_{del,i}$ from $\Pi_f$ doesn't remove any of the planes that contain the $i$-th coordinate (which are the ones in $\Pi_{contr,i}$).
Recall that irreducibility means the coordinate ring has no zerodivisors (for example, modding a polynomial ring by a prime ideal). Positroid vareities are irreducible, so $\Pi_f$ is irreducible. If we look at $\Pi_f \backslash \Pi_{del,i}$, ripping out a subvariety enlargens the coordinate ring by introducing denominators (localization), not zerodivisors; therefore, $\Pi_f \backslash \Pi_{del,i}$ is irreducible. Finally, if we consider $\Pi_f \backslash \Pi_{del,i} \twoheadrightarrow \Pi_f \cap \Pi_{contr,i}$, we use the fact that the image of an irreducible variety is an irreducible variety to conclude irreducibility of $\Pi_f \cap \Pi_{contr,i}$ (here, the image is everything due to surjectivity, so it is irreducible). This fact is true simply because the map of a variety onto its image corresponds to an inclusion of the coordinate rings in the contravariant direction, and a subring of a ring with no zerodivisors does not have zerodivisors either.  

The analogous statement for $\Pi_f \hookleftarrow \Pi_f \backslash \Pi_{contr,i} \twoheadrightarrow \Pi_f \cap \Pi_{del,i}$ follows from Grassmannian duality and proves the irreducibility of $\Pi_f \cap \Pi_{del,i}$.

\end{proof}

\begin{proof}[Proof of Prop \ref{prop:del2def}]

We make use of the fact that $f'(i)=0$ for a positroid variety $\Pi_{f'}$ is equivalent to $\Pi_{f'}$ being contained in $\Pi_{del,i}$; similarly, $f'(i)=n$ for a positroid variety $\Pi_{f'}$ is equivalent to $\Pi_{f'}$ being contained in $\Pi_{contr,i}$. 

$(\supseteq)$ We first show that $Del_i(\Pi_f) \supseteq \Pi_{f'}^{d(i)}$, that is, $\Pi_f \cap \Pi_{del,i} \supseteq \Pi_{f'}^{d(i)}$. This is true because (a) $\Pi_{f'}^{d(i)} \subseteq \Pi_f$ by definition (as the largest positroid variety \textit{contained within} $\Pi_f$) , and similarly $\Pi_{f'}^{d(i)} \subseteq \Pi_{del,i}$ follows from its definition as having $f'(i)=0$, in accordance with the fact that we just stated. The case of contraction is the same.

$(\subseteq)$ Finally we show that $Del_i(\Pi_f) \subseteq \Pi_{f'}^{d(i)}$, that is, $\Pi_f \cap \Pi_{del,i} \subseteq \Pi_{f'}^{d(i)}$. 
For contradiction, assume this is not true; that is, $\Pi_{f'}^{d(i)}$ is properly contained within $\Pi_f \cap \Pi_{del,i}$: $\Pi_{f'}^{d(i)} \varsubsetneqq \Pi_f \cap \Pi_{del,i}$. We will show the contradiction by making use of the irreducibility noted in Prop \ref{prop:delirred} (namely, that while an intersection of positroid varieties like $\Pi_f \cap \Pi_{del,i}$ can in general split into a union of positroid varieties: $ \Pi_{f_1} \cap \Pi_{f_2} = \Pi_{f_3} \cup \Pi_{f_4} \cup \cdots \cup \Pi_{f_j}$, in this special case of $\Pi_f \cap \Pi_{del,i}$, the irreducibility implies that it is a single positroid variety: $\Pi_f \cap \Pi_{del,i} = \Pi_g$). In general, $\Pi_g \varsubsetneqq \Pi_f$, so we would have $\Pi_{f'}^{d(i)} \varsubsetneqq \Pi_g \varsubsetneqq \Pi_f$, and this contradicts the defintion of $\Pi_{f'}^{d(i)}$ as the largest positroid variety contained in $\Pi_f$ such that $f'(i)=0$ (since $\Pi_g = \Pi_f \cap \Pi_{del,i}$ also satisfies $f'(i)=0$). Of course, in the special case when $\Pi_f \subseteq \Pi_{del,i}$, then $\Pi_f \cap \Pi_{del,i} = \Pi_f$, and the largest positroid variety $\Pi_{f'}$ contained within $\Pi_f$ satsifying $f'(i)=0$ is just $\Pi_f$ itself, so $\Pi_{f'}^{d(i)} = \Pi_f = Del_i(\Pi_f)$ and so everything holds automatically.

The proof for the contraction case is analogous.
\end{proof}

\section{The Order and Atomic Pairs}\label{section:theorder}

We describe the poset we will be using in this thesis. As explained in more detail in the previous section, we have the following two subGrassmannians inside $Gr(k,n)$:

\noindent $\Pi_{del,i} := \{V\in Gr(k,n) \mid  V\subseteq Proj_i(\C^n)\cong \C^{n-1} \text{, planes not using the $i$th coordinate line} \}$ \\
$\cong Gr(k,n-1)$ 

\noindent $\Pi_{contr,i} := \{V\in Gr(k,n) \mid  V\supseteq \C_i \text{, planes containing the $i$-th coordinate line} \}$ \\ $\cong Gr(k-1,n-1)$

Given a $T$-fixed point $\lambda$ in $Gr(k,n)$ and a coordinate $x_i$, then either that fixed point $\lambda$ (which corresponds to a coordinate $k$-plane) uses $x_i$ or it does not. This can be rephrased as: for each $T$-fixed point $\lambda$ and each $i\in \{1,...,n\}$, either $\lambda \in \Pi_{del,i}$ or $\lambda \in \Pi_{contr,i}$, but not both. We can intersect a positroid variety with one of these two subGrassmannians, and what we obtain is a smaller positroid variety $\Pi_{f'}$ and a point $\lambda'$ on it, where is $\Pi_{f'}$ the deletion or contraction. In other words: 
\\

(1) If we have a $T$-fixed point $\lambda$ on a positroid variety $\Pi_f$ and we have one of the coordinates $x_i$, then we can either delete or contract $x_i$, and we will get to another point $\lambda'$ on a smaller positroid variety $\Pi_{f'}$.
\\

(2) Furthermore, if the point $\lambda$ we started out with is smooth on $\Pi_f$, then the point $\lambda'$ after deletion or contraction will again be smooth on $\Pi_f \cap \Pi_{del,i}$ or $\Pi_f\cap \Pi_{contr,i}$ (this will be proved as Proposition \ref{prop:delsmooth}). 
\\

Because of (1), we can make an (infinite) poset on the set of pairs $(\Pi_f,\lambda)$ with deletion and contraction as the covering relations. By (2), the set of singular pairs is closed under going upwards. We can visualize the poset as:

 \begin{figure}[htbp]  \centering 
	\includegraphics[scale=0.5,clip=true]{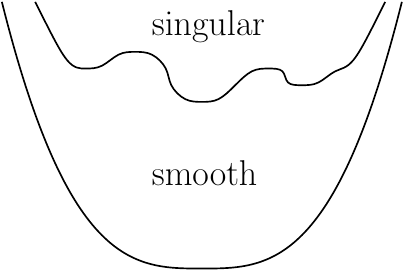}
	\caption{}
\end{figure}

We can give another name to going down in this ordering, a concept which can also be applied to \textit{(affine) pipe dreams} which are defined and discussed in the next chapter.

\begin{Def}\label{def:dcreduce}
(1) We say that $(\Pi_{f_1},\lambda_1)$ \textbf{d/c-reduces} (for deletion/contraction) to $(\Pi_{f_2},\lambda_2)$ if a series of deletions and contractions of $(\Pi_{f_1},\lambda_1)$ results in $(\Pi_{f_2},\lambda_2)$.

(2) While we have not yet defined (affine) pipe dreams, we note here that the notion of ``d/c-reduces" in particular will apply to pipe dreams. It is sufficient to note here that a given $(\Pi_{f_1},\lambda_1)$ corresponds to a set of mathematical objects called pipe dreams, and this set will be denoted $PD_{(\Pi_f,\lambda)}$. This becomes a bijection if we consider the elements in $(\Pi_{f_1},\lambda_1)$ up to equivalence by the ``moves" defined in Definition \ref{def:moves}. Furthermore, we have a notion of deletion and contraction for individual pipe dreams, that is, applying deletion and contraction to elements $\delta \in PD_{(\Pi_f,\lambda)}$ (see \S \ref{section:delcontrpipedream}).
 Thus, we can consider, more specifically, whether a given $\delta_1 \in PD_{(\Pi_{f_1},\lambda_1)}$ goes to $\delta_2 \in PD_{(\Pi_{f_2},\lambda_2)}$ under a series of deletions and contractions on affine pipe dreams described in \S \ref{section:delcontrpipedream}. Obviously, if we consider these affine pipe dreams only up to equivalence under moves, then the relations will be the same as (1) above. However, since it will be easier to work with individual $\delta \in PD_{(\Pi_f,\lambda)}$ in many of the following proofs, in the rest of the thesis when we use the term ``d/c-reduces," we will use it to refer to this latter notion (2). 
\end{Def}

Given this poset, we can consider the points along the boundary (the minimal singular pairs):

 \begin{figure}[htbp] \centering
	\includegraphics[scale=0.5,clip=true]{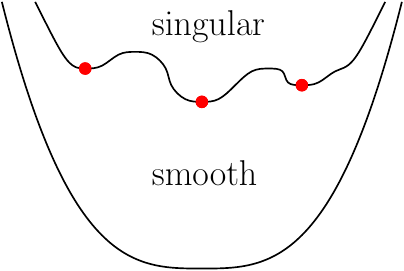}
	\caption{}\label{fig:atomicbound}
\end{figure}

Since $n$ decreases as one goes down in the ordering, we can only go down in the order a finite number of steps from any given point. By this finiteness of going down (a sort of Noetherianess of this poset), a pair is singular if and only if it is greater than one of these minimal singular pairs in the ordering. We give these minimal singular pairs a name:

\begin{Def}[Atomic Positroid Pair]
A positroid variety $\Pi_f$ is \textbf{atomic} on an open patch $U_\lambda$ if it is singular at $\lambda$, but any deletion or contraction (of any column) causes this positroid variety to be smooth (or empty) at $\lambda$. We call $(\Pi_f, \lambda)$ an \textbf{atomic positroid pair}.
\end{Def}

The important thing is that deletion and contraction give us a poset on the set of positroid varieties: a positroid variety $\Pi_{f_2}$ than can be obtained by either deleting or contracting one of the columns of another positroid variety $\Pi_{f_1}$ gives a covering relation, with $\Pi_{f_1}$ covering $\Pi_{f_2}$ (so $\Pi_{f_2}$ is lower in the partial order). The goal is not to understand deletion and contraction of positroid varieties in itself, but to use these operations (and the covering relation they define) to study and understand the smoothness and singularity properties of positroid varieties.

\section{Deletion and Contraction on Matrices and Matroids}\label{sec:delcontrmatrices}

In this section, we prove that deletion on matrices requires setting a column equal to zero. This will be used when we do explicit calculations and start working with actual matrices. 

Deletion of column $j$ means that we intersect with the positroid variety $Del_j$. Let us think on the matroid level. Deletion of column $j$ means that that we impose that every set of $k$ columns $\{v_{i_1},...,v_{i_k}\}$ where $j\in \{i_1,...,i_k\}$, satisfies $Det(v_{i_1},...,v_{i_k})=0$. 

Contraction of column $j$ means that we impose that every set of $k$ columns $\{v_{i_1},...,v_{i_k}\}$ where $j\not\in \{i_1,...,i_k\}$, satisfies $Det(v_{i_1},...,v_{i_k})=0$. 

Suppose we have some set of column vectors $\{v_{i_1},...,v_{i_k}\}$ such that $Det(v_{i_1},...,v_{i_k})\neq 0$. Clearly, such a set $\{v_{i_1},...,v_{i_k}\}$ must exist, because we are in the Grassmannian, so the matrix must be full-rank (of rank $k$). Now, we consider column $j$, which is not in $\{i_1,...,i_k\}$ (otherwise $Det(v_{i_1},...,v_{i_k}) = 0$). The condition (of deletion) is that for every set of $k-1$ columns $\{v_{h_1},...,v_{h_{k-1}}\}$, $Det(v_{h_1},...,v_{h_{k-1}},v_j) = 0$. We want to show that this must imply that $v_j$ is the zero vector. Let's go back and look at $\{v_{i_1},...,v_{i_k}\}$: if we remove any of the $v_i$ and replace it with $v_j$, the determinant vanishes. 
Because $\{v_{i_1},...,v_{i_k}\}$ is full rank for a $k$-dimensional space, it forms a basis. Thus, $v_j =  a_{i_1} v_{i_1} + ... + a_{i_k} v_{i_k}$. And if $v_j\neq0$, then at least one of the $a$'s is not equal to zero. Without loss of generality, say the first $p\leq k$ of them have nonzero coefficient: $v_j =  a_{i_1} v_{i_1} + ... + a_{i_p} v_{i_p}$; then, if we consider the set $\{v_{i_1},...,v_j,...v_{i_k}\}$ where we have replaced any one of $\{v_{i_1},...,v_{i_p}\}$ by $v_j =  a_{i_1} v_{i_1} + ... + a_{i_p} v_{i_p}$, then $Det(v_{i_1},...,v_j,...v_{i_k}) \neq 0$, and we have a contradiction. Therefore, $v_j$ must be the zero vector. 

Similarly, we want to prove that there exists a $Gl(k)$ transformation such that, for the contraction of column $j$, that column $j$ is an identity column (all zeros except a single $1$), and if that $1$ is in row $h$, then all other entires in row $h$ are zero. Since all minors that do not contain column $j$ vanish, but the matrix as an element of the Grassmannian is full rank, column vector $j$ must exist in a dimension orthogonal to the vector space spanned by the other $n-1$ column vectors. In particular, it must be nonzero, so we can always use a GL(k) transformation to set column $j$ to be of the form $(00\cdots 00 100 \cdots 0)^T$ (this would just be rotating/reflecting the axes so that one of the axes lies along $v_j$, and then we stretch or shrink so that $v_j$ has length $1$. (in fact, we can do a further GL(k) to bring it to the last row so column $j$ is $(00\cdots 00 1)^T$. From here, it is automatic that, for the row that contains the $1$, all other entries are $0$: this is because the contraction condition is that all other minors not involving column $v_j$ have determinant zero. Furthermore, we just noted that column vector $j$ must exist in a dimension orthogonal to the vector space spanned by the other $n-1$ column vectors, so all other entries in that row must be zero. We leave the details to the reader. 

We remark this it is nontrivial that the deletion and contraction of a positroid variety $\Pi_f$ is the unique positroid variety lowest in Bruhat order obtained from $f$ via Bruhat moves. This comes from the proof of irreducibility above - that the deletion (or contraction) is the largest positroid variety $\Pi_{f'}$ contained in $\Pi_f$ such that $f'(i)=0$ and due to irreducibility, there is a unique positroid variety of largest dimension, which corresponds to a minimum number of Bruhat moves.

\section{Deletion for Grassmann Necklaces}

We start by showing how one obtains the Grassmann necklace from a siteswap. We'll do this by example. Take a siteswap, say 340362 in Gr(3,6). Do the +1,2,3,4,5,6 to get the affine permutation 4,6,3,7,11,8. Since the $3\to3$ (its siteswap number is zero), $3$ will not show up in any Grassmann necklace set. Since $5\to 5+n=5+6$, $5$ will show up in every Grassmann necklace set. Start by finding the first $k=3$ numbers that satisfy the following 2 conditions: (1) they are not zeros (like $3$) and (2) they aren't dependent on a previous column. For example, for $I_1$, we start by checking the first column, which is a $4$; from here on out, let's write out $1\choose 4$ to denote that the permutation sends $1$ to $4$. This is $4-1=3$ for the siteswap number, which is not zero, so we satisfy that condition. In terms of the condition that it cannot be dependent on a previous column, because this is the first column we are considering, this requirement is automatically satisfied. Thus, let's add column $(1)$ (using parentheses for clarity/disambiguation) to $I_1$: $I_1=\{1\}$ Next, the second number is $2\choose 6$. This satisfies both requirements as well, so we add $(2)$: $I_1=\{1,2\}$. Now let's check $(3)$, but we already saw this was a zero, so we can't add that to $I_1$. Next we look at $4\choose 7$; but we already had $1\choose 4$, so this is linearly dependent. (This is the rule: if the number (here: 4) already showed up previously, then it's linearly dependent so we can't add it). So we move onto the next one: $5\choose 11$; this one satisfies both conditions (in fact, because it is $5\to 5+n=5+6$, we could have started by putting $5$ in every Grassmann necklace set: $I_1=\{5,...\},I_2=\{5,...\},I_3=\{5,...\}$; in other words, we could have stopped once we had two numbers (other than $5$), and this would have saved us some work. Therefore, we have $I_1=\{1,2,5\}$. The following ones are easier. We just have to shift over by one. In other words, we know that $I_2=\{2,5,\sigma(1)\}$, if we denote the permutation by $\sigma$. And of course, ${1\choose \sigma(1)}={1\choose 4}$, so $I_2=\{2,5,4\}$. $I_3$ stays the same $I_3=I_2$ since the $3$ has siteswap number zero. For $I_4$, we drop the $4$ and put in $\sigma(4)$, which is $7=1(mod6)$. 

Thus we can now consider the question of what deletion and contraction do to the Grassmann necklace. Unfortunately, there doesn't seems to be a super clean or efficient way to do this. We still have to look at the transpositions. In other words, the best way is still to look at the diagonal diagrams of zeros and ones, and consider the little staircase block of 0's and 1's that get changed under deletion or contraction. This will involve a series of transpositions. We can look at these transpositions, and use them to modify the Grassmann necklace. The following is the most general form of the sorts of transpositions we can get:

\begin{figure}[htbp]
\includegraphics[scale=0.4,clip=true]{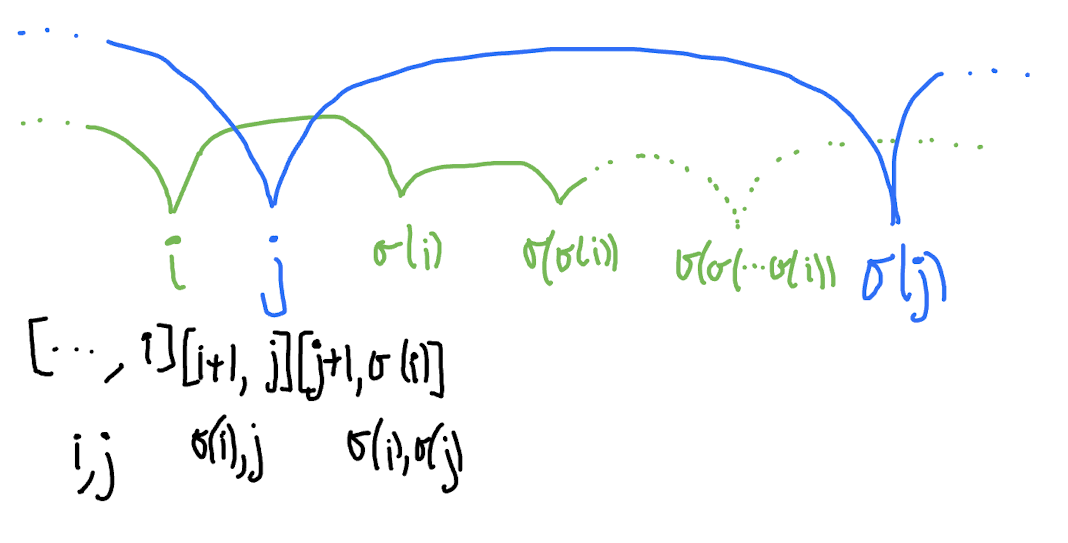}
\end{figure}

That is, we have $i<j$ and $\sigma(i)<\sigma (j)$, and we transpose them. Note that in the general form, the juggling ball corresponding to the throw at $i$ (drawn in green), can potentially land a bunch more times before $j$ lands (hence, the $\sigma(\sigma(\cdots(\sigma(i))\cdots)$. The Grassmann necklace is when the balls next land (so since we've depicted them in different colors, the soonest that the different colors land). We see that if we are, say, right before $i$, say $i-1$, supposing that $\sigma(i-1)\neq i,j$, then the green ball and blue ball next land on $i$ and $j$, so $I_{i-1}=\{i,j,...\}$, so nothing changes. If we are on $i$, then it's the same thing: we have $I_i=\{i,j,...\}$ so nothing changes. If we are at $\sigma(i)$, then $I_{\sigma(i)}=\{\sigma(i),\sigma(j),...\}$, so both numbers transposed are contained, so nothing changes. This is true for the region between $\sigma(i)$ and $\sigma(j)$ as well, as can be easily shown. The only place where something changes is the region $[i+1,j]$. Here, for $I_{i+1},I_{i+2},...,I_{j}$, these all contain $\sigma(i)$ (which is where the green ball next lands), and $j$, where the blue ball next lands. But once we do the transposition (see the figure below):

\begin{figure}[htbp]
\includegraphics[scale=0.4,clip=true]{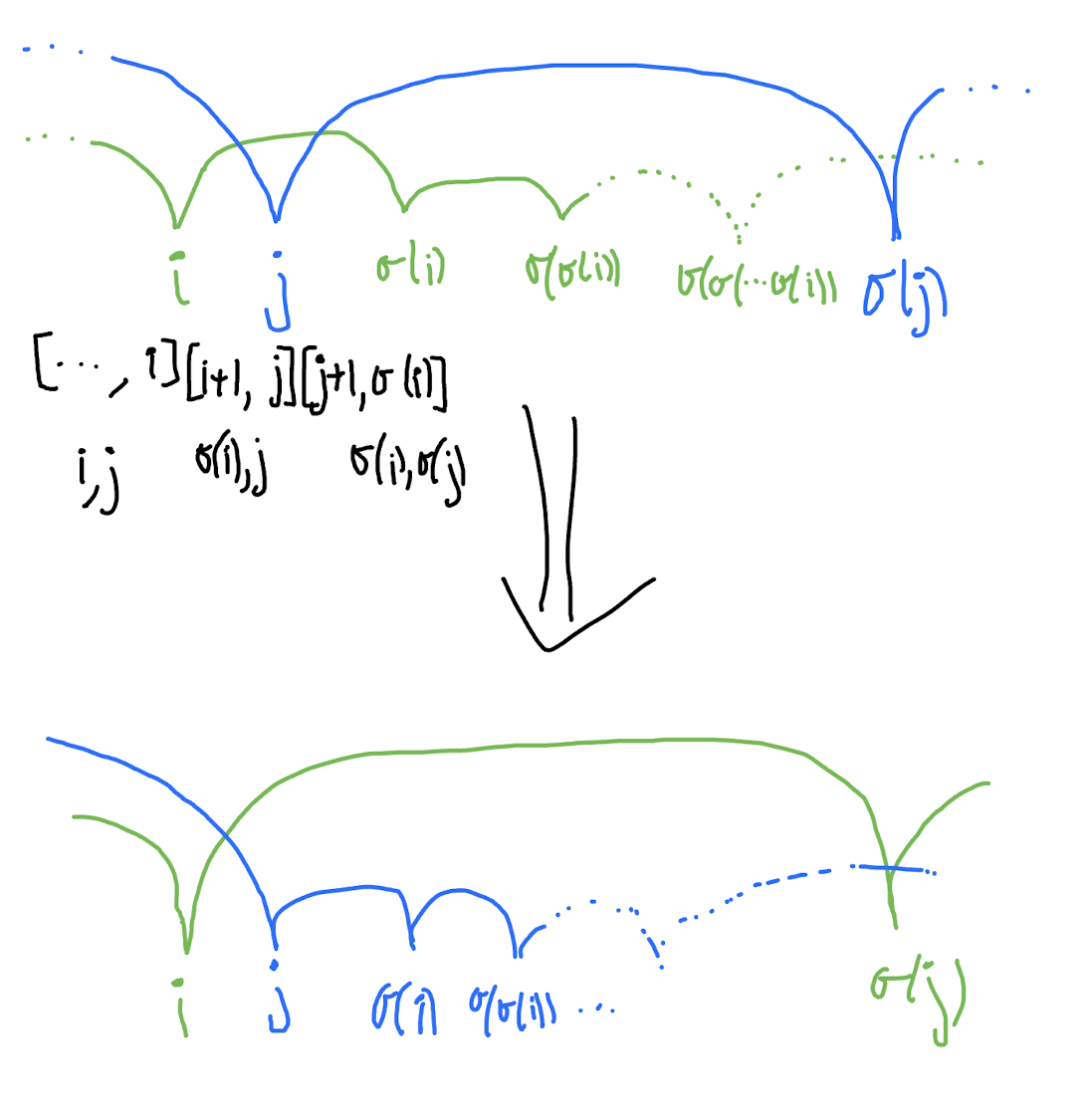}
\end{figure}

the $\sigma(i)$ changes to $\sigma(j)$. The $j$ stays the same. So this is the rule: doing the transposition for $i \choose \sigma(i)$ and $j \choose \sigma(j)$ in the deletion process (it will be a dual thing for contraction) causes the following change on the Grassmann necklace. Take $I_{i+1},I_{i+2},...,I_{j}$ and change $\sigma(i)$ to $\sigma(j)$ for these. 

Let's do a quick example. Let's use the example in my paper, with siteswap 340362, and with the process for deletion and contraction depicted in the following 2 figures (specifically, they show that deleting the $6$-th column produces siteswap 350460, or 4,7,3,8,11,6 as a bounded affine permutation):

\begin{figure}[htbp]
\includegraphics[scale=0.5,clip=true]{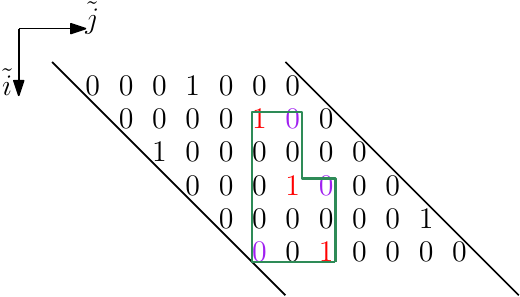}
\end{figure}

\begin{figure}[htbp]
\includegraphics[scale=0.5,clip=true]{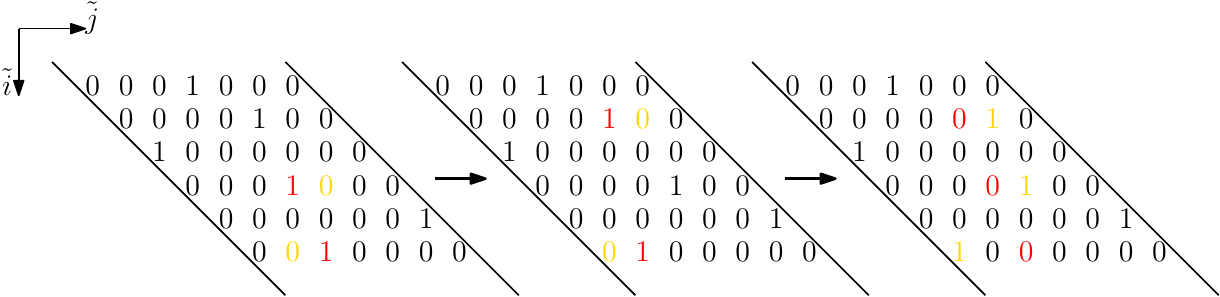}
\end{figure}

The Grassmann necklace for 340362 (=4,6,3,7,11,8 as a bounded affine permutation) is $\{I_1,I_2,...,I_6\}=\{125,245,456,456,156,156\}$. The figures show that the process of deleting the $6$-th column involves two transpositions: switching the 4th and 6th entries, and then switching the 2nd and 6th entries. Specifically, switching ${4\choose3} \leftrightarrow {6\choose2}$, or as a permutation (which let's use for the rest of this example) ${i \choose \sigma(i)} = {4\choose 7} \leftrightarrow {j\choose \sigma(j)} ={6\choose 8}$, our work above shows that we must replace $\sigma(i)=7=1(mod6)$ with $\sigma(j)=8=2(mod6)$ in $I_{i+1},...I_j$, which in this case is $I_5,I_6$. Thus $\{125,245,456,456,156,156\} \mapsto \{125,245,456,456,\rd{2}56,\rd{2}56\} $. Next, we have ${i \choose \sigma(i)} = {2\choose 6} \leftrightarrow {j\choose \sigma(j)}= {6\choose 7}$ means we $\sigma(i)=6$ with $\sigma(j)=7=1(mod6)$ in $I_{i+1},...I_j = I_3,I_4,I_5,I_6$. Thus, $\{125,245,456,456,256,256\} \mapsto \{125,245,45\rd{1},45\rd{1},25\rd{1},25\rd{1}\} $, and this is indeed the correct Grassmann necklace for siteswap 350460, or 4,7,3,8,11,6 as a bounded affine permutation. 

It is not hard to show that the final result is to take every set $I_i$ containing $6$ (or the column being deleted) and replace it with smallest number (in Gale order) in $I_6$ not in $I_i$. This reproduces (and gives an alternate proof of) the result in \cite{SOh}.

\section{Positroid vs. Siteswap/Permutation Bijection}

In this section, we describe how to map back and forth between the siteswap or bounded affine permutation and the positroid (the matroid).

First, note that positroids (for the complex Grassmannian, as well as the real), are the intersected of rotated Schubert matroids. The Schubert matroids  come from looking at the Schubert cells and stating the nonzero minors. For example, for Gr(2,4), one of the Schubert cells, the one with pivots in columns 1 and 4, looks like:
$\mqty(1&*&*&0\\0&0&0&1)$, which has Schubert matroid: $\{14, 24, 34\}$. These are the nonvanishing minors. The positroids would include all 4 rotations: $\{ 21,31,41\}, \{32,42,21 \}, \{43,13,32 \}$. We do this for all of the Schubert cells in Gr(2,4).

We are also allowed intersections. For example $\mqty(1&*&0&*\\0&0&1&*)$ corresponds to the condition $det(1,2)=0$, or $span(1,2)\leq1$. This has positroid/matroid $\{13,14,23,24,34\}$. We are allowed for positroids to impose any combination of cyclic rank conditions. So if we also impose that $det(3,4)=0$, then this corresponds to this positroid shifted by 2: $\{31,32,41,42,12\}$. So the intersection corresponds to the positroid $\{31,32,41,42\}$. Note that it can be the case that the intersection is reducible. The prototypical case is in Gr(2,3) where you impose $det(1,2)=0, det(2,3)=0$: you can satisfy this with either of the 2 varieties: $\mqty(a&a&a\\b&b&b)\cup \mqty(*&0&*\\ *&0&*)$.

To go between the siteswap $f$ and positroid you need to use Grassmann necklaces. To go between the siteswap and Grassmann necklace, do this: for the siteswap, draw the juggling pattern or diagonal strip of 0's and 1's, and for each row, find out where the soonest balls $k$ drop (for Gr(k,n)) - equivalently, the first $k$ linearly independent columns or unique lines. In the reverse direction, a Grassmann necklace consists of $(I_1,I_2,...,I_n)$ where each $I_i$ has $k$ elements. For each successive $I_i, I_{i+1}$, if the numbers are equal, then $i\not\in I_i$, so column $i$ is a column of zeros, so the siteswap has $f(i)=0$, otherwise, $i\in I_i, i\not\in I_{i+1}$, so there is a single new number in $I_{i+1}$, and this new number is precisely $f(i)$ (or rather, the bounded affine permutation mod n)

Then, to go between the Grassmann necklace and the positroid, use the following result stated in \cite{ARW}. In other words, look at the rotated (considering $i$ as the smallest number) Schubert matroid corresponding to each $I_i$ (each $I_i$ is a set of pivots columns, so you can consider the Schubert matroid like we did above), and intersect them. The map from positroids to Grassmann necklaces is given by the map denoted $\mathcal{M}$ map in Theorem 4.3 of of \cite{ARW} (coming from \cite{Post} and \cite{Oh08})

\begin{Thm}[\cite{Post}, \cite{Oh08}] Let $\mathcal{I}=(I_1,I_2,...,I_n)$ be a Grassmann necklace of type $(k,n)$. Then the collection 
\[ \mathcal{B}(\mathcal{I}):= \bigg\{ B\in {[n]\choose k} \vert B\geq_j I_j \text{ for all } j\in [n] \bigg\}\]
is the collection of bases of a rank $k$ positroid $\mathcal{M}(\mathcal{I}):= ([n],\mathcal{B}(\mathcal{I}))$. Moreover, for any positroid $M$ we have $\mathcal{M}(\mathcal{I}(M))=M$.
\end{Thm}

\chapter{Combinatorial Background Related to Pipe Dreams}\label{ch:combinbackgroundpipes}

\section{Square Pipe Dreams}

\begin{Def}[Pipe Dream \cite{BB93}]
A pipe dream in $S_n$ is a diagram in an $n$-by-$n$ square where each box can be filled by one of two types of tiles, elbows \elbow and crosses \+, such that all crosses occur above the southwest-northeast diagonal. Thus, we can think of the grid as a set of pipes that begin on the north and east edges and end on the west and south edges (with the south to east tiles being all elbows)

We called the pipe dream \textbf{reduced} if no two of the pipes cross twice. We will focus only on the reduced case in the rest of the thesis.

\end{Def}

\begin{Exa}
The following figure shows an example for all pipe dreams corresponding to the permutation $\pi=2143$ in $S_4$. The pipes running from the south edge to the east edge are uninteresting so are not shown. 
\\

\end{Exa}

\begin{figure}[htbp] \centering
\includegraphics[scale=0.5,clip=true]{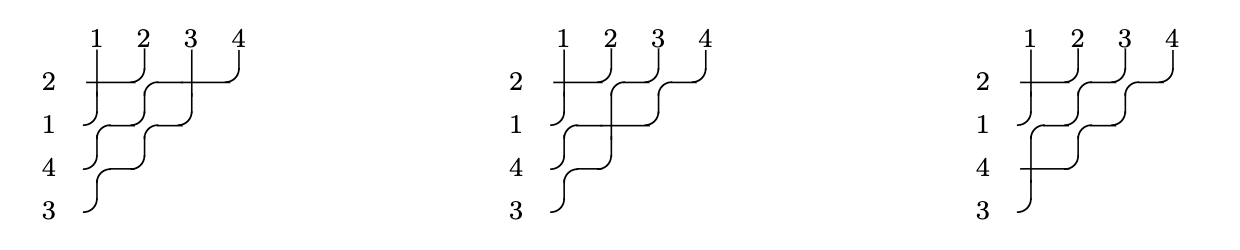}
\caption{}
\end{figure}

\begin{Def}\label{def:moves}
A \textbf{cross-elbow move} (in the following: ``\textbf{move}," for short) in a pipe dream is an interchange of a cross tile and an elbow tile that preserves the connectivity of the pipes - in other words, that leaves the permutation invariant. We think of the cross as ``moving" to the position where the elbow it is replacing was located. If we perform all possible moves in a pipe dream, and keep repeating, taking all the new pipe dreams obtained after performing moves and doing moves on these new ones, then we will obtain the set of all possible pipe dreams for that particular permutation (see \cite{BB93}). If the pipe dream has no possible moves, so is the unique pipe dream for a particular permutation, then we call the pipe dream \textbf{rigid}. The notion of moves carries over without change to the case of bounded affine permutations which we define next (the chute and ladder moves of \cite{BB93} are a subset of our moves), and will be used extensively below. 
\end{Def}

\begin{Exa}
In the previous figure, for the pipe dream on the very left, we see that the upper cross has no moves, but the lower cross has 2 different places it could potentially move, which does not change the connectivity/permutation. Therefore, this gives us the 3 different pipe dreams. 
\end{Exa}

Next, we describe a couple of different and important ways to think about pipe dreams, focusing on examples. 

\section{Matrix Schubert Varieties and Rothe Diagrams}

\begin{Def}
Given a permutation $\pi$ on $n$ elements, we can write the permutation matrix $N_\pi$ corresponding to $\pi$. Then we can define a \textbf{matrix Schubert variety} ${\overline X}_\pi$ as all $n$-by$n$ matrices $M$ satisfying the condition that the rank of the $p\times q$ northwest submatrix, denoted $r_{pq}(M)$ is less than or equal to that of $N_\pi$: $r_{pq}(M)\leq r_{pq}(N_\pi)$
\end{Def}

\begin{Exa}
Let's take a permutation of 4 elements, say: $1423$. 

The matrix corresponding to this permutation is:

$N_\pi=N_{1423}=\mqty(1 & & & \\ & & & 1 \\ & 1 & & \\ & & 1 &)$. 

Then ${\overline X}_\pi$ is the set of all 4-by-4 matrices $M$, such that the rank of the $p\times q$ northwest submatrix, denoted $r_{pq}(M)$, is less than or equal to that of $N_\pi=N_{1423}$: $r_{pq}(M)\leq r_{pq}(N_{1423})$.

\end{Exa}

\begin{Def}
The \textbf{Rothe diagram} associated to a permutation $\pi$ is obtained by taking the permutation matrix $N_\pi$ of $\pi$ and drawing lines south and east from each $1$ in $N_\pi$. All squares/boxes within the $n$-by-$n$ permutation matrix that have a line going through them are considered eliminated. 
\end{Def}

\begin{Exa}
Let us draw the Rothe diagram associated with $\pi = {1423}$, given by drawing lines south and east from each $1$ in the permutation matrix:

\begin{figure}[htbp] \centering
	\includegraphics[scale=0.8,clip=true]{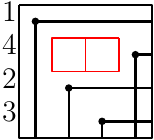}
\caption{}
\end{figure}
\end{Exa}

Note that $\text{codim} \overline X_\pi$ equals the number of boxes in $\pi$'s Rothe diagram. In other words, the number of boxes is the number of conditions imposed, equivalently the number of inversions or length of the shortest word for $\pi$ (see below), each which drops the dimension by one. 

If we turn the boxes into crosses and push them to the left, we get the `bottom pipe dream"; to the top is the ``top pipe dream." The red dots in the figure below show a potential ``move": a crossing and then near-miss (important for what follows).

\begin{figure}[htbp] \centering
	\includegraphics[scale=0.6,clip=true]{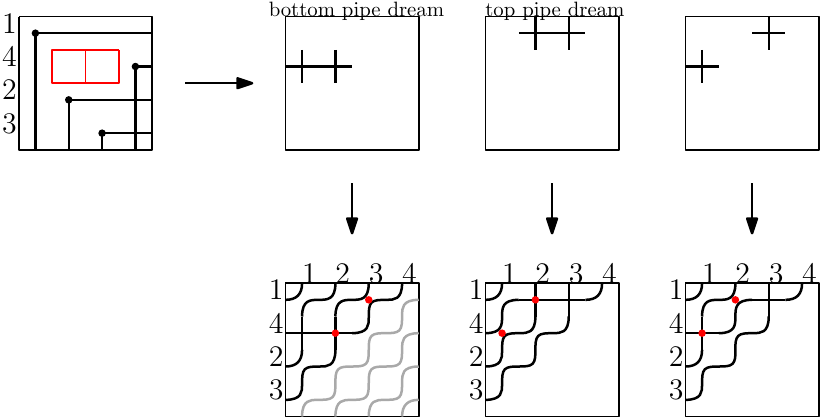}
	\caption{}
\end{figure}

\newpage

\section{Wiring Diagrams}

Below, we describe another way of thinking about pipe dreams, as wiring diagrams, but this will not be used after this section, so can be skipped without loss of continuity. 

\begin{Def}[\cite{KM03}]
Let $\Pi$ denote the set of permutations on $n$ elements, and let $\Sigma$ denote the set of transpositions $s_1,...,s_n$ where $s_i$ switches $i$ and $i+1$. A word of size $m$ is an ordered sequence  $Q=(\sigma_1,...,\sigma_m)$ of elements of $\Sigma$. An ordered subsequence $P$ of $Q$ is called a \textbf{subword} of $Q$. 

 $P$ \textbf{represents} $\pi\in \Pi$ if the ordered product of the simple reflections in $P$ is a reduced decomposition for $\pi$.

$P$\textbf{contains} $\pi\in\Pi$ if some subsequence of $P$ represents $\pi$.

The subword complex $\Delta(Q,\pi)$ is the set of subwords $Q\backslash P$ whose complements $P$ contain $\pi$.
\end{Def}

Square pipe dreams correspond bijectively to reduced words and embeddings into square words, or equivalently, reduced words and compatible sequences. Wiring diagrams correspond bijectively to reduced words. The map from square pipe dreams is injective if we only consider the top (or just the bottom) pipe dream (so that any given reduced word has at most one pipe dream representative); however, the map is not surjective since there are reduced words that do not satisfy the compatible sequence condition. 

\begin{Exa}(Non-example)
If we consider $S_3$ with long word 321,432,543, we can find the subword 212 in it, but not 121. Now, both 121 and 212 give the permutation $123 \choose 321$, so we do obtain all permutations (just not all reduced words)
\end{Exa}

The following is an example from Postnikov of a wiring diagram (thus corresponding to a reduced word) that does have a square pipe dream mapping to it. We describe explicitly how to 

\begin{figure}[htbp] \centering
	\includegraphics[scale=0.4,clip=true]{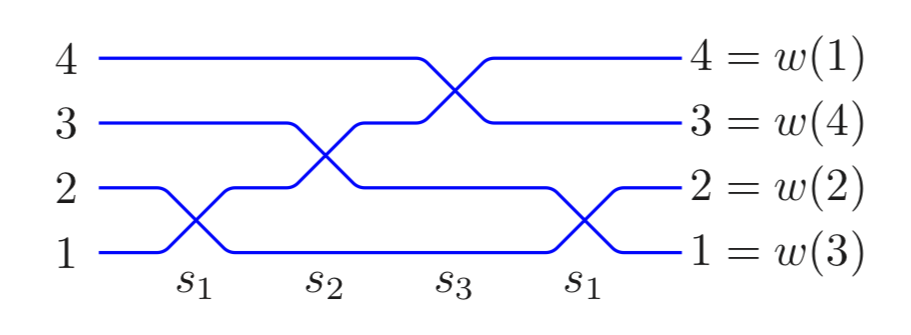}
	\caption{}
\end{figure}

In this context, reduced means that no two lines cross twice. In this example, if we trace the pipes from left to right, we get the permutation $4213$. With transpositions $(i,i+1)$ denoted by $s_i$, we see that the wiring diagram gives us a word $s_1 s_2 s_3 s_1$ which produces this permutation: 

$s_1 s_2 s_3 s_1(1234)=s_1 s_2 s_3 (2134)=s_1 s_2 (2143) = s_1 (2413) = (4213)$.

For these words, reduced corresponds to it being the shortest word producing a given permutation.

Now let's see what the transpositions are in the square from before:

\begin{figure}[htbp] \centering
	\includegraphics[scale=0.4,clip=true]{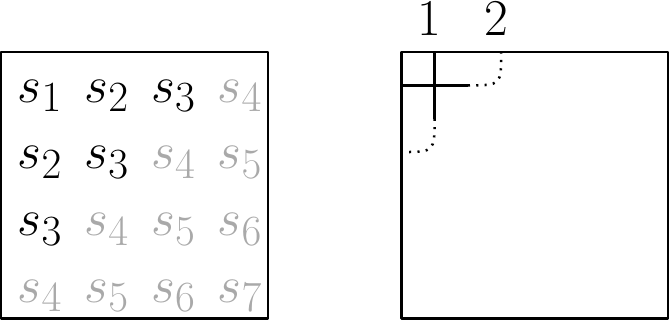}
	\caption{}
\end{figure}

We read from top to bottom and left to right. The second figure shows that if we just had the northwesternmost cross, it would switch pipes 1 and 2, so it makes sense for it to be $s_1$. For the example above of $1423$ it's clear that all pipe dreams give the same word $s_3 s_2$. We can thus view the pipe dream as a rotated wiring diagram producing a reduced word:

\begin{figure}[htbp] \centering
	\includegraphics[scale=0.5,clip=true]{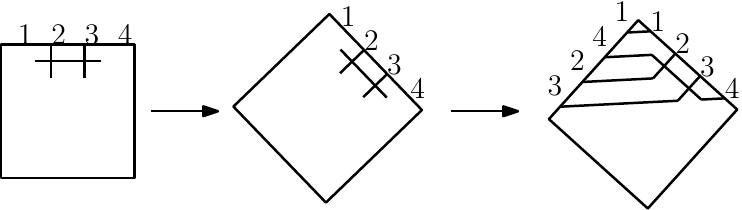}
	\caption{}
\end{figure}

Here we are reading right to left, so in the order opposite to Postnikov, we can verify: $s_3 s_2 (1234)= s_2(1243) = 1423$.

\section{Affine Pipe Dreams}\label{affinepipedreams}

Recall that we are studying positroid varieties on the open affine patches $U_\lambda$, these open neighborhoods of fixed points $\lambda$. We have an analogue of pipe dreams for bounded affine permutations called \textbf{affine pipe dreams} invented in \cite{Sni}. 

\begin{Def}
\textbf{Affine pipe dreams} are diagrams constructed by picking a \textbf{distinguished path} corresponding to $\lambda$ and using it to cut a matrix into two pieces, which are then shifted northeast and southwest infinitely. The details are in Chapter 4 of \cite{Sni}, and are illustrated in rest of this section via an example. 
\end{Def}

We do an example to show how this works.

Let's take an example of a siteswap $2225377$ (bounded affine permutation $34598,13,14$) inside Gr(3,7), and let's test it on $\lambda=\{2,5,7\}$. 

We first write the $U_\lambda$ corresponding to the point $\lambda=\{2,5,7\}$:

 \begin{figure}[htbp] \centering
	\includegraphics[scale=0.8,clip=true]{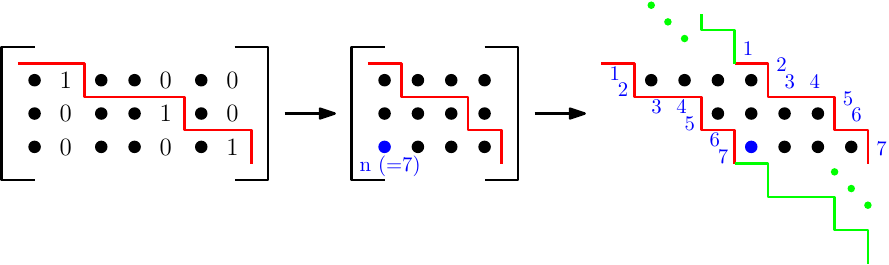}
	\caption{}
\end{figure}

(1) In the first figure on the very left, the distinguished path (illustrated in red) is formed by moving horizontally starting from the northwest corner of the matrix until a $1$ in one of the identity columns is encountered, at which point the distinguished path moves down. Thus, the distinguished borders the north and east of each $1$ and otherwise is horizontal. \\
(2) In the middle figure, we collapse the identity columns to obtain the distinguished path. The tiles will be filled by transpositions $s_i$. Here, we have highlighted in blue the southwesternmost square. This will be $s_n$ (in the depicted example, $s_n = s_7$). \\
(3) Finally, we take one copy of the distinguished path and move it to the right (depicted also in red in the figure on the very right). We can continue adjoining distinguished paths (depicted in green) to the northwest and southeast of these two red distinguished paths, forming an infinite strip, which is the boundary of the affine pipe dream. Note the labelling along the boundary: it simply starts with ``$1$" at the top line segment of the distinguished path and increases moving downwards. Thus, we see why the blue dot will be the transposition $s_7$ (or just ``$7$" in the next figure): a cross where the blue dot is will switch the $7$-th and $8=1(\bmod 7)$-th pipes coming from the bottom boundary.

Filling in the rest of the tiles, we get: 

 \begin{figure}[htbp] \centering
	\includegraphics[scale=0.9,clip=true]{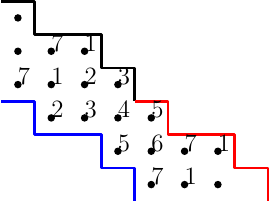}
	\caption{}
\end{figure}

If we just look at one block, denoting the transpositions by $s_i$ again, it is:
 \begin{figure}[htbp] \centering
	\includegraphics[scale=0.7,clip=true]{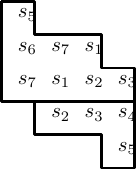}
	\caption{}
\end{figure}
To get the bottom pipe dream, we have to pick a way of reading this. There are equivalent ways of reading it, given by: 
 \begin{figure}[htbp] \centering
	\includegraphics[scale=0.7,clip=true]{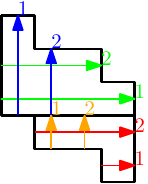}
	\caption{}
\end{figure}

In other words, if we choose the direction of red followed by green, then the ordering will be $s_5, s_2 s_3 s_4, s_7 s_1 s_2 s_3, s_6 s_7 s_1, s_5$. If we apply this word to $(k+1,k+2,...k+n)$, then we get the bounded affine permutation correspond to $+n$ for the $\lambda$ columns $2,5,7$ and $+0$ for all other columns. This shows that this pipe dream shape represents the neighborhood $U_\lambda$.

Here's how reading it this way produces a word for siteswap $2225377$ (bounded affine permutation $34598,13,14$):

We start with:

\begin{center}
\begin{tabular}{c|c|c|c|c|c|c|c}
1 & 2 & 3 & 4 & 5 & 6 & 7 & 8(=1) \\ \hline
3 & 4 & 5 & 9 & 8 & 13 & 14
\end{tabular}
\end{center}

Since the big word is $s_5, s_2 s_3 s_4, s_7 s_1 s_2 s_3, s_6 s_7 s_1, s_5$, we find the first $s_i$ (starting from the left with $s_5$) where there is a transposition with $f(i)>f(i+1)$. Then we can apply $s_i$ to switch $f(i)$ and $f(i+1)$. Here, the first one is $s_4$ which switches the $98$ to $89$:
\\

\begin{center}
\begin{tabular}{c|c|c|c|c|c|c|c}
1 & 2 & 3 & 4 & 5 & 6 & 7 & 8(=1) \\ \hline
3 & 4 & 5 & 9 & 8 & 13 & 14 & \\ \hline
3 & 4 & 5 & 8 & 9 & 13 & 14 & 10 \\
\end{tabular}
\end{center}

We continue on until we get the most general permutation:
\\

\begin{center}
\begin{tabular}{c|c|c|c|c|c|c|c}
1 & 2 & 3 & 4 & 5 & 6 & 7 & 8(=1) \\ \hline
3 & 4 & 5 & 9 & 8 & 13 & 14 & \\ \hline
3 & 4 & 5 & 8 & 9 & 13 & 14 & 10 \\ \hline
7 & 4 & 5 & 8 & 9 & 13 & 10 & 14 \\ \hline
4 & 7 & 5 & 8 & 9 & 13 & 10 & 11 \\ \hline
4 & 5 & 7 & 8 & 9 & 13 & 10 & 11 \\ \hline
4 & 5 & 7 & 8 & 9 & 10 & 13 & 11 \\ \hline
6 & 5 & 7 & 8 & 9 & 10 & 11 & 13 \\ \hline
5 & 6 & 7 & 8 & 9 & 10 & 11 & 12 \\
\end{tabular}
\end{center}
Keeping track of which transpositions we used, we get:
 \begin{figure}[htbp] \centering
	\includegraphics[scale=0.7,clip=true]{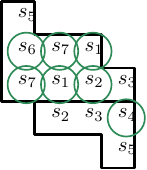}
	\caption{}
\end{figure}
If we put crosses where the circles are, and elbows elsewhere, then we get the bottom pipe dream.

In contrast, the top pipe dream corresponds to the pipe dream that we would obtain if we applied the same process as above, but where we make a few changes: \\
(1) We must use the inverse permutation. We write this as the numbers $g^{-1}(i) +n $. One way to obtain these is to draw out the juggling pattern and for each $i$, look at the ball that lands on it, and consider the number of step backwards this is; these will all be negative numbers. If we then add $n$ to all entries, we obtain the inverse siteswap.\\
(2a) We could change change the $s_i$'s such that they correspond to the labelling along the top boundary of the pipe dream shape rather than the bottom (this is depicted in the left column of the following figure). We will then read the word in the opposite ordering to the one above; in fact, any ordering that only considers a square $x$ when all squares to the northeast of $x$ have all been read, will work. If we want to read in a consistent direction within each block (above or below the distinguished path), the choices are shown in the colored arrows in the next figure. \\
(2b) Alternatively to (2a) (but essentially equivalent), we could use Grassmannian duality to reflect the pipe dream shape across the northwest-southeast line, before going through with the same process of finding the bottom pipe dream; obviously, the changes of the $s_i$'s to match the (originally) top boundary happens automatically here. Finally, after doing this, we have to reflect back across the northwest-southeast line to get the top pipe dream.  

In our example the inverse permutation of 2225377 is 4255500; explicitly, we must start with:

 \begin{figure}[htbp] \centering
	\includegraphics[scale=0.7,clip=true]{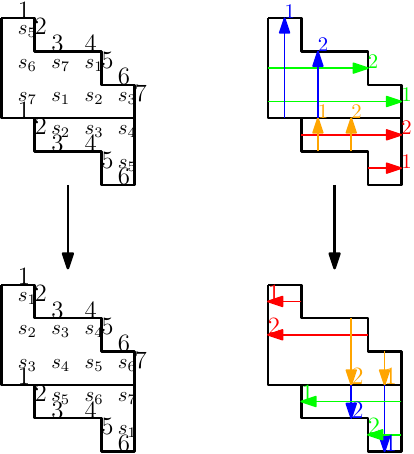}
	\caption{}
\end{figure}

\subsection{Motivational Proof for Affine Pipe Dream Shape}\label{subsec:affinepipesshapemotivation}

This section can be skipped without loss of continuity for understanding this thesis.

One question is why one chooses this particular shape.

It is sufficient to consider the main theorem proven in \cite{Sni}, that the intersection of a positroid variety with $U_\lambda$ is isomorphic to an affine Kazhdan-Lusztig variety. \cite{Sni} furthermore shows that there exists a correspondence between the stratifications, generated by an intersect-decompose algorithm. In the next section, we will use this to show that this implies that the equivariant cohomology class of a positroid variety restricted to a point is therefore given by a sum over these affine pipe dreams (for the shape of the point being restricted to). We will then show that by restricting our torus action to a circle acting by dilation, we obtain a corollary that the multiplicity of that point on the positroid variety is equal to the number of pipe dreams for that positroid variety in that shape. The details are below. In any case, because these affine pipe dreams have already been proven to give the isomorphism is used to prove the theorem on smoothness stated in this paper, we can use them for the results we prove below. However, we also include for motivational purposes (which the author found helpful) another argument in analogy with matrix Schubert varieties. 
\\

It is helpful to consider the case of matrix Schubert varieties. Let us go back to the case of $S_4$ (where a specific case of $\pi=2143$ was given in the example above). We note that if the square is empty, by following the pipes we get the identity permutation: 1234. By having crosses in the entire NW region where crosses are allowed, we get the largest permutation in Bruhat order: 4321. These permutations and pipe dreams are discrete things; there are a finite number of them. But each one represents an algebro-geometric object: namely, the matrix Schubert variety given by imposing rank conditions on these square matrices. In fact, each pipe dream for a given matrix Schubert variety $X_\pi$ has length($\pi$), equivalently codimension($X_\pi$) many crosses. 
\\

By analogy with the matrix Schubert variety case, the affine pipe dreams are diagrams where filling the shape with crosses gives us the largest permutation in Bruhat order: these are the points of the form $00nn0n\cdots 0n$ (a string of $n$ total numbers, of which $k$ of them are $n$'s, and $n-k$ of them are $0$'s; as a check, the average is then $n\cdot k/n=k$ as it should be). Recall that these are the fixed points we are considering in our partial order. In fact, we have a proposition:

\begin{Prop}
Let $X$ be an affine (periodic, infinite) pipe dream shape with squares in which crosses and elbows can be placed. Suppose that this pipe dream shape represents a neighborhood around given fixed point $p$; that is, filling in the shape completely with crosses produces the siteswap corresponding to $p$, whereas configurations that contain elbow(s) give siteswaps corresponding to positroid varieties that contain $p$. Then, the shape must be of the form introduced above, with the north-south and east-west parallel-shifted labeling on the two boundary lines being the same. 

\end{Prop}

\begin{proof}
First, we make a few observations. The fact that this is a pipe dream implies that there is a labeling on two different line segments, the labeling on one line segment telling us where the pipes start, and the labeling on the other one giving the permutation. The numbers on both sides must be labelled consecutively $...1,2,3,...,n,...$, as in the matrix Schubert variety case (although here in the affine case the numbers extend infinitely). Furthermore, the lines can touch, but they can never cross; this is because the pipe dream is an oriented object, with the elbow and cross tiles moving in one direction only. Also, the labeling must run ``in the same direction." In other words, since the infinite lines are periodic and never cross, they both must on average move in 180-degree directions, both oriented in the same way, and therefore, since they are ``parallel", if the labelling was in opposite directions, then (viewing the labelling as that of a bounded affine permutation) the pipe dream would give a non-periodic permutation of increasing larger and larger jumps, which would be neither periodic nor bounded. 
\\

We put the word ``parallel" in quotations since the non-intersection implies that they are parallel on average, but we have not shown that the two lines must be parallel in the strong sense: that is, one line is just a shifted copy of the other line. This will be proved in the following. 
\\

Consider the following picture:  

\begin{figure}[htbp] \centering
\includegraphics[scale=0.7,clip=true]{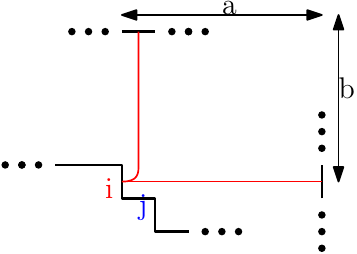}
\caption{}
\end{figure}

Without loss of generality there is a periodic line segment on the bottom where the pipes begin. We have randomly labelled a pair of consecutive positions along this representative bottom line segment $i$ and $j$ (so here $j=i+1$ [although we could label the pipes going northwest instead, but this would just involve a reflection]). There is a maximal window or range of places where the pipe starting from position $i$ can end. This is shown by the red lines. The farthest left it can go is the ending position of the horizontal line at the top, the farthest right it can go is the vertical segment on the right. Thus, this maximal range is the width plus the height, which we have labelled as $a$ and $b$. Since this is a bounded affine permutation, this gives us a condition on these numbers: $a+b\leq n+1$. This implies that the two ``parallel" lines/boundaries are bounded in terms of how far they can diverge from each other. In particular, the line at the top must be $\geq i$ and the line on the right must be $\leq i+n$. 
\\

However, we had an additional assumption, namely that the lowest dimensional variety, the point $p$, is the configuration given entirely with crosses. In order to consider this configuration, we look at the following figure:    

\begin{figure}[htbp] \centering
\includegraphics[scale=0.7,clip=true]{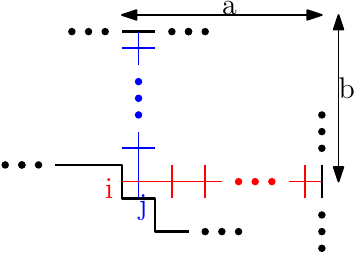}
\caption{}
\end{figure}

The lines along the bottom that are vertical (like $i$) move to the maximal position along horizontal strings of crosses (shown in red). The point $p$ is given by a siteswap like $00nn0n\cdots 0n$ consisting entirely of $0$'s and $n$'s. Thus, the vertical lines must correspond to the $n$'s. Similarly, the lines along the bottom that are horizontal (like $j$) move to the maximal position along vertical strings of crosses (shown in white). These then, must correspond to $0$'s in the siteswap. Thus, if we label the boundary of the pipe dream complex with the bounded affine permutation numbers, $i$ moves along the red string of crosses and ends at $i+n$, whereas $j$ moves along the blue string of crosses and ends at $j$. Since the line segments are continuous, this fixes the shape of the other (northeast) boundary of the pipe dream complex entirely: it must be a parallel shift of the other boundary. 
\\

(Remark: We can do this constructively: consider the position $i+2$ after $j=i+1$. If it is again horizontal like $j$ then clearly we must also extend the line on top horizontally since by what we have already shown, if $i+2$ is horizontal along the bottom boundary then it must end at at horizontal position $i+2$ as well. However, if $i+2$ is vertical along the bottom, then along the top boundary, we should attach a vertical piece of the right edge of the horizontal segment to the north of $j$; otherwise, if we made the next line a horizontal $j+1=i+2$, then at some point along the bottom boundary when it moves horizontally we would have a position greater than $i+2$ capable of moving to the $i+2$ position along the top, which would be a permutation going negative number of steps, so not a bounded affine permutation.) The only thing that remains to be fixed is the amount the parallel northeast boundary of the complex is to be shifted, and this is easily fixed by the fact that the configuration with all elbows has siteswap $kk\cdots kk$; in particular, since we showed that the $n$'s in the point $p$'s siteswap are vertical line segments along the boundary and there are $k$ $n$'s, the two boundaries should be separated by height $b=k$ (and width $a=n-k$, which is fixed by everything else).  

\end{proof}

This concludes the explanation of why our affine pipe dreams have their particular shape. We already described what the numbers within the pipe dream complex mean in terms of  rotated wiring diagrams, using the connection to words in Coxeter groups, specifically transpositions defining permutations.

\subsection{Moves on Affine Pipe Dreams}

\begin{Def}[see \cite{BB93}]
A \textbf{move} on a pipe dream or affine pipe dream is an exchange of a cross tile and a near-miss tile that does not change the connectivity of the pipe dream (so all pipes begin and end at the same places). 
\end{Def}

Now we show that for these affine pipe dreams, it is still the case that if a particular affine pipe dream is \textbf{rigid} (it has no possible moves), then it is the only one that produces its corresponding bounded affine permutation.

\begin{Prop}\label{prop:BAPmoves}
For any affine pipe dream $\Omega$ with bounded affine permutation $g$, every other affine pipe dream $\Omega'$ that also has bounded affine permutation $g$ can be reached via $\Omega$ by performing a sequence of moves on $\Omega$.
\end{Prop}

In the section \ref{affinepipedreams} above, we did an example to illustrate how a subword $P_{f,\lambda}$ for the siteswap $f$ at the point $\lambda$ is obtained from the word $Q_\lambda$ of the affine patch around the point $\lambda$. Here, we state the following facts in a more general context.

\begin{Prop}[Corollary 1 in \cite{Sni}]\label{prop:subwdcomp} 
The \textbf{affine pipe dream complex} $\Delta(\lambda,\pi)$ is the simplicial complex with vertices labeled by entries $(i, j)$ in the periodic strip and faces labeled by the elbow sets in the affine pipe dreams for $\pi$ with shape defined by $\lambda$. It is isomorphic to the subword complex $\Delta(Q_\lambda,f)$.
\end{Prop}

\begin{Thm}[KM04, Theorem 3.7]\label{thm:ballsphere} The topological realization of any subword complex $\Delta(Q,\pi)$ homeomorphic to a ball or sphere; in particular, every ridge (codimension 1 facet) is contained in one or two facets. The dual graph (whose vertices are the facets of $\Delta(Q,\pi)$ and whose edges are the ridges of $\Delta(Q,\pi)$) is connected.
\end{Thm}

\begin{Def}[Definition \# 2 of pipe dreams moves]
We defined pipe dream moves above in terms of exchanges of tiles on a particular pipe dream. In the context of subword complexes, moves correspond to going to a ridge and then going to another codimension-1 facet. This follows from the fact the notion give above of a move in terms of exchanging a cross and a near-miss can actually be thought of as adding a cross in the new location before removing the original one, and the ridges in the complex are where we have nonreduced pipe dreams.
\end{Def}

\begin{proof}[Proof of Prop \ref{prop:BAPmoves}]
By Proposition \ref{prop:subwdcomp}, there exists a subword complex $\Delta(Q_\lambda,f)$ for which the set of all pipe dreams for a siteswap $f$ on a patch $\lambda$ are the codimension-1 facets of this complex. By Theorem \ref{thm:ballsphere}, since subword complexes are homeomorphic to balls or spheres, they are connected in codimension-1. Therefore, any two pipe dreams for the same siteswap on the same patch must be connected by a sequence of moves. 
\end{proof}

\section{Deletion and Contraction on Affine Pipe Dreams}\label{section:delcontrpipedream}

We present a picture of what $\Pi_{del,i}$ looks like on an open patch. Recall that this is the positroid variety which has siteswap $(k+1)^k0k^{n-k-1}$ rotated so the $0$ is in the $i$ slot. It is the largest positroid variety whose siteswap has a $0$ in the $i$-th slot. It will be easy to see that the pipe dream must look like the following: 

\begin{figure}[htbp] \centering
\includegraphics[scale=0.45,clip=true]{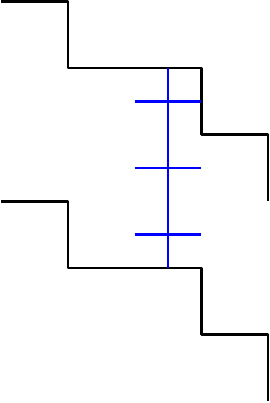}
\caption{}
\end{figure}

Here is the proof. The siteswap is $(k+1)^k0k^{n-k-1}$. In the pipe dream above, the crosses are placed along the $i$-th slot, so that the only difference from the pipe dream with all elbows (producing the siteswap $kk\cdots kk$), is that the $k$ crosses along that vertical line cause the $k$ slots immediately prior to $i$ to move one slot ahead of what they would be in the all-elbow pipe dream. This is precisely the $(k+1)^k$ part immediately prior to the $0$ in the siteswap. Thus, this pipe dream produces the siteswap we wanted, and in fact it is the unique such pipe dream since there are no moves in the pipe dream. 

Another way to prove this is that, by the proof above, not only must the column to be deleted lie along a horizontal segment, but the entire column must be filled with crosses so that it ends in the minimal spot (in this case: itself, given by the siteswap notation of $0$ steps forward). Any additional crosses lowers the dimension of the variety. $\Pi_{del,i}$ was defined as the maximal positroid variety with a $0$ in the $i$-th slot, so it must be depicted by this pipe dream. 
\\

The case of  $\Pi_{contr,i}$ (which has siteswap $(k-1)^{n-k}k^{k-1}n$ rotated so the $n$ is in the $i$ slot) is analogous, with the crosses proceeding from the vertical line segment labeled $i$ in a horizontal line from the left to the right. 
\\

The next question to be asked is then what the deletion or contraction of a positroid variety looks like. We know what $\Pi_{del,i}$ and $\Pi_{contr,i}$ look like, from what we have just done. Previously, we defined the deletion or contraction of a positroid variety $\Pi_f$ as the intersection of $\Pi_f$ with $\Pi_{del,i}$ or $\Pi_{contr,i}$. So the question now is what this intersection looks like on an affine patch. We have the following proposition:

\begin{Prop}\label{prop:delpipedream}

The pipe dreams of the $i$-th deletion of a positroid variety $\Pi_f$ (that is, $del_i(\Pi_f)=\Pi_f \cap \Pi_{del,i}$) on the open patch $U_\lambda$ are given by taking the pipe dreams for $\Pi_f$ on $U_\lambda$ with the maximal number of crosses along the $i$th column, and then filling in this $i$th column with crosses.
\end{Prop}

\begin{Exa}
Here is a quick pictorial representation of $Del_4(\Pi_f)$ (assuming the blank parts of the pipe dream are all elbows): 

 \begin{figure}[htbp] \centering
	\includegraphics[scale=0.5,clip=true]{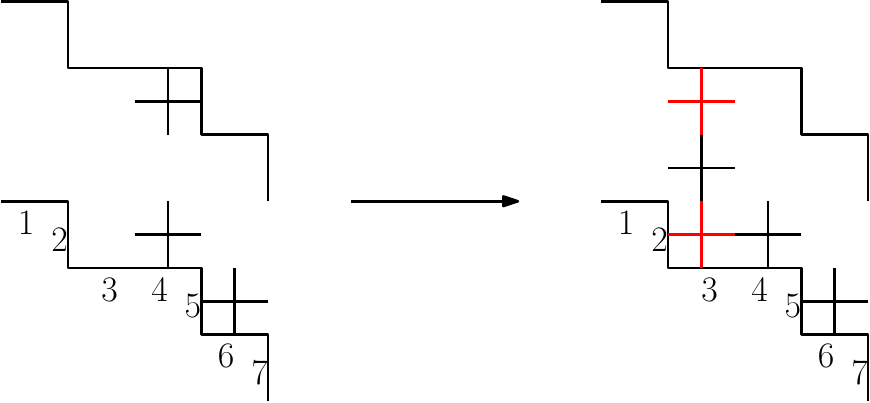}
	\caption{}
\end{figure}

That is, we take the pipe dream(s) with the largest number of crosses along that column, and fill in the rest of this column with crosses (denoted by the red crosses in this example here).
\end{Exa}

We can read out the new (deleted) affine permutation by following the pipes along this new pipe dream; however, note that the affine permutation $del_i(f)$ is not dependent on which patch $U_\lambda$ we draw the pipe dream on. Thus, if the goal was simply to compute $del_i(f)$, rather than drawing the pipe dream (on a particular patch $U_\lambda$), we could do this more efficiently using an algorithm that is presented in Appendix A of \cite{Flu}.  

The case of contraction is analogous with the crosses filled in horizontally (in a row rather than a column of crosses).

\begin{proof}[Proof of Prop \ref{prop:delpipedream}]

We take $Del_i(\Pi_f)=\Pi_f \cap \Pi_{del,i}$ and intersect everything with the vector space $U_\lambda$:

\[Del_i(\Pi_f)\cap U_\lambda=(\Pi_f \cap U_\lambda) \cap (\Pi_{del,i} \cap U_\lambda).\] 

In other words, we get the components of $(\Pi_f \cap U_\lambda) \cap (\Pi_{del,i} \cap U_\lambda)$ by taking the components of the pipe dreams for $\Pi_f$ and the components of the pipe dreams for $\Pi_{del,i}$ and intersecting them and looking at the maximal element. Indeed, as we showed in the proof of Prop \ref{prop:del2def}, the intersection $Del_i(\Pi_f)=\Pi_f \cap \Pi_{del,i}$ is irreducible so is a single maximal positroid variety. We then intersect this positroid variety with $U_\lambda$, which is just looking at the pipe dreams for $Del_i(\Pi_f)$ on the shape defined by $\lambda$. 

Furthermore, since we know that the maximal intersection will have the fewest number of crosses (since each cross represents a variable that is set to zero - thus an extra condition that lowers the dimension of the variety), and since we know that deletion involves filling in the column entirely with crosses, we must choose to intersect with the pipe dream(s) for $(\Pi_f \cap U_\lambda)$ that already had the greatest number of crosses along this column so that filling in this column involves adding the minimal number of additional crosses.

\end{proof}

\chapter{Positroid Varieties and Smoothness}\label{ch:posvarsmoothness}

We have already seen how to study positroid varieties on an open cover using affine pipe dreams. We have also seen how easy it is perform deletion and contraction on positroid varieties by using these affine pipe dreams. The main goal of this section is to prove Theorem \ref{thm:smoothone}, which will provide a means of studying smoothness or singularity of positroid varieties on this open cover using pipe dreams.

\section{Singularity on Each Open Patch}

Recall from \S \ref{subsection:standardopencover} that the $T$-fixed points on $Gr(k,n)$ are given by the span of $k$ coordinates lines, which we have denoted by $\lambda$, so there are $n\choose k$ such fixed points.

\begin{Thm}\label{thm:testonnchoosek}
One can check whether a positroid variety $\Pi_f$ is singular (or smooth) just by testing singularity at the fixed points $\lambda$.
\end{Thm}

\begin{proof}

Fact (1): Borel's Fixed Point Theorem: If a solvable group (such as a torus) acts on some nonempty proper (such as projective) variety, then the fixed point set is also nonempty.
\\

Fact (2): The torus $T$ acts on the singular locus.

Proof: Whenever a group $G$ acts on a variety, the singular locus will be $G$-invariant (i.e. a group of symmetries will leave the singular locus alone). The ``proof" is simply that the notion of singularity does not require coordinates. 
\\

Fact (3): Positroid varieties are proper. It is a standard fact that projective implies proper. 
Positroid varieties are closed in Grassmannians which are closed in projective space (in the Plücker embedding); thus positroid varieties are proper. 
\\

Fact (4): The singular locus of any variety is closed. This follows from the fact that a smooth locus is given by a nonvanishing (open) determinantal Jacobian condition.
\\

Suppose $\Pi_f$ has singular point(s), i.e. the singular locus $(\Pi_f)_{sing}$ is nonempty. Then we know that $(\Pi_f)_{sing}$ is closed in $\Pi_f$ by (4). Since $\Pi_f$ is proper by (3) and $(\Pi_f)_{sing}$ is closed in $\Pi_f$, we conclude that $(\Pi_f)_{sing}$ is also proper. By (2), the torus acts on the singular locus, so we can apply (1) Borel's Fixed point Theorem to it, to conclude that the singular locus includes at least one fixed point (under the torus action).  

Denoting the fixed points of the torus action with superscript $T$, we summarize:
\[ (\Pi_f)_{sing} \neq \emptyset \Leftrightarrow (\Pi_f)_{sing}^T \neq \emptyset\]
\end{proof}

Thus, we have reduced checking singularity (or smoothness) to a finite set of points $\{\lambda\}={[n]\choose k}$!

In what follows, we will check these by checking on the affine open cover of the $U_\lambda$ described earlier. Clearly, this will include all the torus-fixed points $ \{\lambda\}$.

\section{Statement of Main Theorem on Smoothness}\label{sec:smoothness}

\begin{Thm}\label{thm:smoothone}
A positroid variety $\Pi_f$ is smooth at the point $\lambda$ (it meets $\lambda$ and is smooth there) if and only if there is a single affine pipe dream representative \cite{Sni} of the affine permutation $f$ on the open set $U_\lambda$. Specifically $\Pi_f$ meets the point $\lambda$ if and only if there exists (at least) one affine pipe dream, and it is smooth if and only if that affine pipe dream is unique. Thus, we can take $\Pi_f$ and look at pipe dreams for all choices of $\lambda$ to test for smoothness. 
\end{Thm}

We will spend most of the rest of this section proving this theorem. We begin by developing some of the background facts.

\subsection{Proof of Theorem - Preliminary Facts}

\begin{Def} 
The \textbf{multiplicity} $Mult(x\in X)$ is the degree of the tangent cone.
\end{Def}

\begin{Thm}
[\cite{Hartshorne} Exercise 5.3a] $Mult(x\in X)=1$ if and only if $X$ is regular (or smooth) at $x$. 
\end{Thm}

Thus, proving that ``smoothness means one pipe dream" has become ``multiplicity 1 means one pipe dream"; the next theorem shows that alternatively we have ``coefficient 1 in equivariant cohomology means one pipe dream." For an introduction to equivariant cohomology, see for example sections 1.2, 2.1 of \cite{KnTao01}.  

\begin{Thm}\label{thm:Rossmann}[\cite{Rossmann}]
Let $V$ be a vector space, and $X \subseteq V$ be a cone, i.e. a subvariety that is invariant under the dilation action of $\C^\times$ on $V$ (that is, $\C^\times$ acts with all weights $1$ on $V$). Denote the equivariant cohomology $H_{\C^\times}^*(V)$ by $\Z[h]$. Then the multiplicity of $X$ at $0\in V$ is the coefficient of $h^{\codim X}$: 
\[[X]_{\C^\times} = (\deg X)h^{\codim X} \in H^*_{\C^\times}(V)\]
\end{Thm}

Remark: We use the notation $[X]_{\C^\times}$ to disambiguate from $[X]_T \in H^*_T(V)$, which we will see later. 

In our case, this means that we want to show that $[\Pi_f \cap U_\lambda] = 1 \cdot h^{codim X} \in H_{\C^\times}^*(U_\lambda)$ is equivalent to $\Pi_f$ having a unique pipe dream on $U_\lambda$. More generally, we will show that the multiplicity counts the number of affine pipe dreams. 

The proof goes by showing the following sequence of equalities:

\begin{equation}\label{equivcoheq} 
[\Pi_f \cap U_\lambda \subseteq U_\lambda] \overset{(1)}{=} [X_{v(f)} \cap X_\circ^{v(\lambda)} \subseteq X_\circ^{v(\lambda)}] \overset{(2)}{=} [X_{v(f)}]|_{v(\lambda)} \overset{(3)}{=} \sum_R \prod_{r\in R}  \hat{\beta_r} \in H_T^*(U_\lambda) \cong \Z[y_1,...,y_n] 
\end{equation}

Note on notation: $f$ is a siteswap and $v(f)$ its corresponding bounded affine permutation, which indexes a $T$-fixed point on the affine flag variety $\widehat{GL_n}/\widehat{B}$ (note: $\widehat{GL_n}$ denotes $n\times n$ invertible matrices where the entries are Laurent series in one variable, and this is infinite-dimensional). Then $X_\circ^{v(\lambda)}$ is the $B$-orbit through the point $v(\lambda)B/B$; $X_\circ^{v(\lambda)} := B v(\lambda) B/B$. And  $X_{v(f)}$ is the $B_-$-orbit closure through $v(f)$; $X_{v(f)} := \overline{B_- v(f) B}/B$.

\subsection{Equation \ref{equivcoheq} (1) - Snider's Results}

Equality (1) follows from the work of \cite{Sni} Section 4.2, which says we have a T-equivariant commutative diagram:

$$
 \begin{tikzcd}
U_\lambda \ar[r,leftarrow,"\cong"]  & X_\circ^{v(\lambda)} \\
\Pi_f \cap U_\lambda \ar[r,leftarrow,"\cong"] \ar[u,hookrightarrow] & X_{v(f)} \cap X_\circ^{v(\lambda)} \arrow[u,hookrightarrow] 
\end{tikzcd}
$$

On the level of T-equivariant cohomology, this gives rise to:
\\

$$
\begin{tabular}{ccc}
$g^*: H_T^*(U_\lambda)$ & $\ra{\cong}$ & $H_T^*(X_\circ^{v(\lambda)})$ \\
$\rotatebox[origin=c]{90}{$\in$}$ & & $\rotatebox[origin=c]{90}{$\in$}$ \\
$[\Pi_f \cap U_\lambda ]$ & $\mapsto$ & $[X_{v(f)} \cap X_\circ^{v(\lambda)} ]$ 
\end{tabular}
$$

In words, \cite{Sni} says that there is a $T$-equivariant isomorphism (this map is not simply an isomorphism of subvarieties, but as a map on the whole space, an isomorphism of their cohomologies) $g: X_\circ^{v(\lambda)} \cong U_\lambda $ such that each $\Pi_f \cap U_\lambda$ maps to $X_{v(f)} \cap X_\circ^{v(\lambda)}$. Thus, once we descend to $T$-equivariant cohomology, the map called ``$g^*$" in the diagram above takes 
$[\Pi_f \cap U_\lambda ] \mapsto [X_{v(f)} \cap X_\circ^{v(\lambda)} ]$.
\\

Since $U_\lambda$ and $X_\circ^{v(\lambda)}$ are contractible, we have canonical (because there is a unique map from a space $X\to \{pt\}$) isomorphisms: $H_T^*(U_\lambda) \cong H_T^*(pt) \cong H_T^*(X_\circ^{v(\lambda)})$. 
$$\begin{tikzcd} & H_T^*(pt) \ar[dl,leftrightarrow,"\cong"'] \ar[dr, leftrightarrow,"\cong"] & \\ H_T^*(U_\lambda)  \ar[rr,"g^*"] & & H_T^*(X_\circ^{v(\lambda)}) \end{tikzcd}$$ 
Using these isomorphisms to identify each with $H_T^*(pt)$, the polynomials $[\Pi_f \cap U_\lambda]$ and $[X_{v(f)} \cap X_\circ^{v(\lambda)} ]$ are equal. 

This completes the justification of equality (1). Next we prove equality (2): $[X_{v(f)} \cap X_\circ^{v(\lambda)} \subseteq X_\circ^{v(\lambda)}] = [X_{v(f)}]|_{v(\lambda)}$.  

\subsection{Equation \ref{equivcoheq} (2) - Kleiman Transversality}

The right side of this equality is the restriction of the class $[X_{v(f)}]$ to the point $v(\lambda)B/B$ (in the following, we will denote this as simply $v$ or $v(\lambda)$). Naively, we might consider this via the diagram: 

$$\begin{tikzcd}
 & X_{v(f)} \ar[d, hookrightarrow] \\
v  \ar[r, hookrightarrow, "\iota"]  & G/B 
\end{tikzcd}
$$

We factor this inclusion as follows:

$$\begin{tikzcd}
& X_{v(f)} \cap X_\circ^{v(\lambda)} \ar[d, hookrightarrow]  \ar[r] & X_{v(f)} \ar[d, hookrightarrow] \\
v  \ar[r, hookrightarrow] & X_\circ^{v(\lambda)} \ar[r, hookrightarrow]  & G/B 
\end{tikzcd}
$$

The bottom row induces maps in cohomology:
\[H^*_T(v)  \leftarrow  H^*_T(X_\circ^{v(\lambda)})  \leftarrow H^*_T(G/B) \]
The goal is to show that under this induced map:
\[ [X_{v(f)}] |_{v} \mapsfrom  [X_{v(f)} \cap X_\circ^{v(\lambda)}]  \mapsfrom [X_{v(f)}] \]

We know that the composition of the two maps will take $[X_{v(f)}]$ to $[X_{v(f)}] |_{v}$ because cohomology is functorial so the pullback on cohomology of the composite map $v \to X_\circ^{v(\lambda)} \to G/B$ is the composite of the pullback; in other words, once we know where $v$ maps to, given by the map $v\hookrightarrow G/B$, this determines the map $H^*_T(G/B) \to H^*_T(v)$ regardless of how we might factor the map in the middle, and this defines $[X_{v(f)}] |_{v}$ as the image of this map. 

The harder thing is we need to show that this map hits $[X_{v(f)} \cap X_\circ^{v(\lambda)}]$ in the factorization through the middle space $X_\circ^{v(\lambda)}$. This will follow from transversality of $X_{v(f)}$ and $X^{v(\lambda)}_\circ$ \cite{Kleiman}. We can apply Kleiman transversality because $X_\circ^{v(\lambda)}$ is a $B$-orbit, and and $X_{v(f)}$ is $B_-$ invariant, so then we conclude that under the map $g: X_\circ^v \hookrightarrow G/B$ going to $g^*: H_T^*(G/B) \to H_T^*(X_\circ^v)$, $g^*([X_{v(f)} ]) = [X_{v(f)} \cap X_\circ^{v(\lambda)} ]$.

Finally, to prove equality (2), we need to show that  $[X_{v(f)} \cap X_\circ^{v(\lambda)}]$ doesn't just map to $[X_{v(f)}] |_{v}$, but that they are equal: we could see this already from the discussion preceding the triangle commutative diagram above, since $X_\circ^{v(\lambda)}$ is equivariantly contractible to $v$, so $H_T^*(v) \cong H_T^*(X_\circ^{v(\lambda)})$. 

This completes the proof of equality (2). 

Remark. In analogy by the same sort of arguments, writing \[\lambda \hookrightarrow U_\lambda \hookrightarrow Gr(k,n) \ni \Pi_f\] and taking cohomology to get \[ [\Pi_f] \to [\Pi_f \cap U_\lambda] \to [\Pi_f]|_\lambda , \] we can show that $[\Pi_f \cap U_\lambda]= [\Pi_f]|_\lambda$. 
\\

\subsection{Equation \ref{equivcoheq} (3) - AJS/Billey Formula}\label{subsec:AJSBilley}
 
Equality (3) is the AJS/Billey formula \cite{AJSBilley}, which gives the restriction of a Schubert class to a fixed point. We state it here:

\begin{Def}
Let $Q$ be a word in the simple roots of a Kac-Moody group (see \cite{Kumar} Prop 11.1.11 for the Kac-Moody case). We define $\beta_i$ as follows:
\[\beta_i := (\prod_{j=1}^{i-1} r_{Q(j)})\cdot (x_{Q(i)} - x_{Q(i)+1}) =  r_{Q(1)}\cdots r_{Q(i-1)} \cdot (x_{Q(i)} - x_{Q(i)+1}), i\in [m] \]
\end{Def}
Let $\prod Q$ be a word in the letters $Q(i)$ so $\prod Q$ is an element of a Coxeter group (in our case, this will be a product of transpositions representing a bounded affine permutation), and $[X_w]|_{\prod Q}$ denote the restriction of the Schubert class $[X_w]$ to the point in the (affine) flag variety represented by $\prod Q$. Then the \textbf{AJS/Billey formula} says:

\[[X_w]|_{\Pi Q}= \sum_{\substack{R\subseteq Q \\ R\text{ reduced}}} \prod_{r\in R}  \hat{\beta_r}  \]
\\

Our work so far has proven:
\[ [\Pi_f \cap U_\lambda]_T = \sum_{\substack{R\subseteq Q \\ R\text{ reduced}}} \prod_{r\in R}  \hat{\beta_r} \]
We state a proposition for what the right hand side looks like.

\begin{Prop}\label{prop:AJSrowcolform}
\[ [\Pi_f \cap U_\lambda]_T = \sum_{\substack{R\subseteq v(\lambda) \\ \prod R=v(f) \\ R \text{ reduced}}} \prod_{r\in R}  [(y_i)_r - (y_j)_r] \]

The AJS/Billey formula produces a sum of products of $(y_i-y_j)$: the terms multiplied together to form a product are indexed by $r\in R$, and each $r$ corresponds to a cross at $(j,i)$ in a pipe dream representative for the given positroid variety on the given patch. Specifically, the $i$ corresponds to the column and the $j$ to the row in the pipe dream. The labelling is given by the numbering of the distinguished path explained in \S \ref{affinepipedreams}; thus, the rows $j$ will come from the vertical boundary line segments, which correspond to $\lambda$, where the columns $i$ will be from the horizontal boundary line segments found in $[n]\backslash \lambda$. 
\end{Prop}





\begin{proof}

There are 2 cases: (1) the case when the cross of interest is above the distinguished path, and (2) when this cross is below the distinguished path. Case (1) is easier, so we start with this. The general picture is the following:

\begin{figure}[htbp] \centering
	\includegraphics[scale=1.4,clip=true]{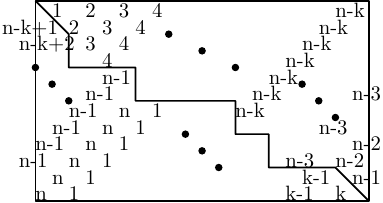}
	\caption{}
\end{figure}

It may be helpful to keep a specific example in mind during the proof. Below, we provide an example in $Gr(4,10)$:

  \begin{figure}[htbp] \centering
	\includegraphics[scale=0.6,clip=true]{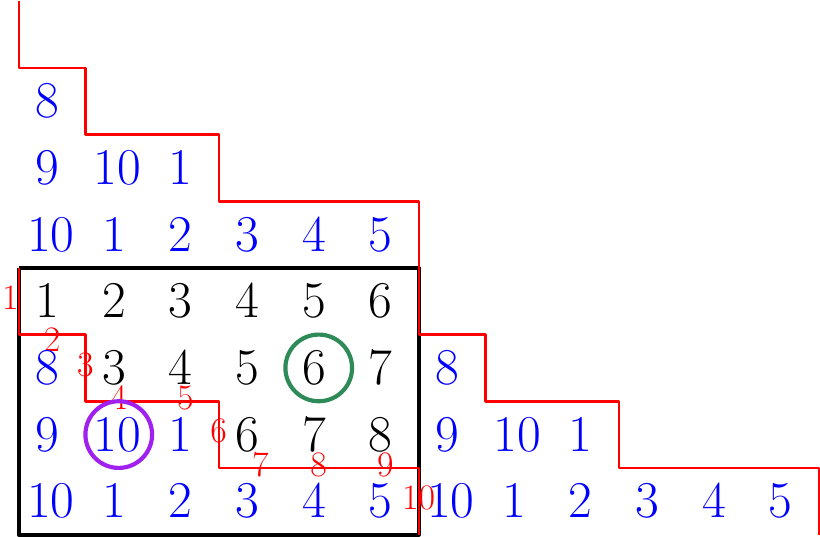}
	\caption{}
\end{figure}

Our goal is to calculate $\beta_i=\prod_{j=1}^{i-1} r_{Q(j)} (x_{Q(i)} - x_{Q(i)+1})$; specifically, we want to show that it equals $x_{row}-x_{col}$. We recall the convention in \cite{Sni} that the words are read from the lower rows (starting above the distinguished path) to the upper rows, and within each row, from left to right. For example, using our example in $Gr(4,10)$, the order would be as follows:

\begin{figure}[htbp] \centering
	\includegraphics[scale=0.6,clip=true]{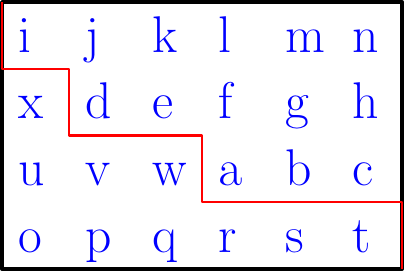}
	\caption{}
\end{figure}

Thus, the word would be: $Q(a)Q(b)Q(c)Q(d)Q(e)Q(f)Q(g)Q(h)Q(i)Q(j)Q(k)$
$Q(l)Q(m)Q(n)Q(o)Q(p)Q(q)Q(r)Q(s)Q(t)Q(u)Q(v)Q(w)Q(x)$. However, it is important to note that these letters act on the roots in the opposite order: in this example, $Q(x)$ would act on the root first, and then secondly $Q(w)$, etc.

Thus, in case (1) where $Q(i)$ is above the distinguished path, we would have (using the letter $b$ to denote ``beginning"): 
\begin{align*}
\beta_i &=\prod_{j=1}^{i-1} r_{Q(j)} (x_{Q(i)} - x_{Q(i)+1}) \\
&= r_{Q(1)} r_{Q(2)} r_{Q(3)} \cdots r_{Q(i-1)} (x_{Q(i)} - x_{Q(i)+1}) \\
&= r_{Q(1)} r_{Q(2)} r_{Q(3)} \cdots r_{Q(i-1)} (x_b - x_{b+1}) 
\end{align*}

In our example, we have circled $Q(i)$ in green, and highlighted in orange $r_{Q(1)} r_{Q(2)} r_{Q(3)} \cdots r_{Q(i-1)}$. 

\begin{figure}[htbp] \centering
	\includegraphics[scale=0.7,clip=true]{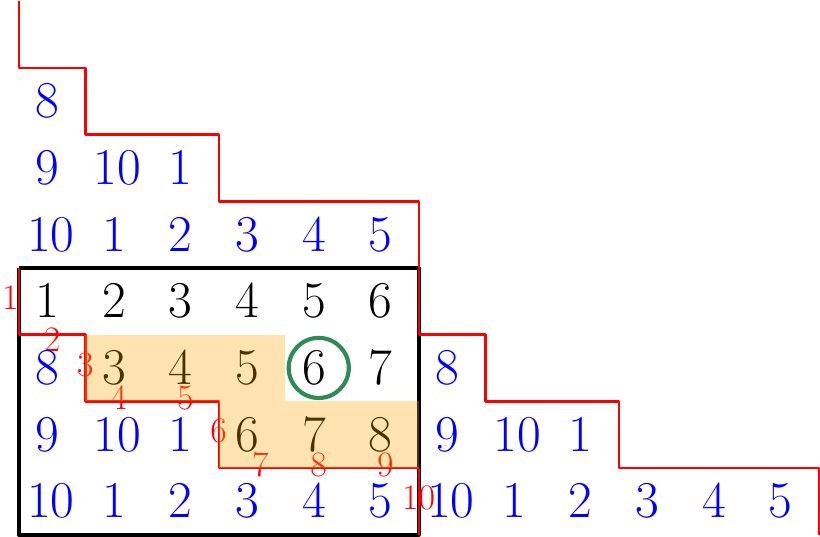}
	\caption{}
\end{figure}

Here, $r_{Q(1)} r_{Q(2)} r_{Q(3)} \cdots r_{Q(i-1)} Q(i) = r_{Q(1)} r_{Q(2)} r_{Q(3)} r_{Q(4)} r_{Q(5)} r_{Q(6)} Q(7) = r_6 r_7 r_8 r_3 r_4 r_5 (x_6 - x_7)$

Since the word is read from left to right within each row, and it acts on the root in the opposite order, if $Q(i)$ is not the leftmost letter in its row (if it is not adjacent to the distinguished path), then $r_{Q(i-1)}$ transposes $b-1$ and $b$, producing $x_{b-1} - x_{b+1}$; then, if $Q(i-1)$ is not the leftmost letter in its row, $r_{Q(i-2)}$ transposes $b-2$ and $b-1$, producing $x_{b-2} - x_{b+1}$, and so on. This will of course continue until we have $x_{l} - x_{b+1}$, where $l$ is precisely the number along the vertical line segment of the distinguished path to the left of $Q(i)$. 

Note that now $x_l$ can only be changed by $Q(m)=l$ or $Q(m)=l-1$, but none of these will appear again in the word as we read left. This is clear from the way the numbers are laid out in the grid, with numbers running diagonally from southwest to northeast; thus, since we have $l$ to the very left in its row, all other $l$'s occur to the northeast in the grid, while all $l-1$ occur in the diagonal that is adjacent to the northwest. However, the remaining letters to be read are all in rows below. As such, we might as well write the polynomial as $x_{row} - x_{b+1}$, where $x_{row}$ is fixed. 

We now consider that the remaining letters $Q(m)$ will be coming from rows below the first row, and acting starting from the right. The $x_{b+1}$ can be changed by either a $Q(m)=b+1$ or $Q(m)=b$. Since we started with $Q(i)=s_b$, the number directly beneath it is $s_{b+1}$ (unless $r_{Q(1)} r_{Q(2)} r_{Q(3)} \cdots r_{Q(i-1)}$ only involved a single row, in which case we are already done), which changes the polynomial to $x_{row} - x_{b+2}$. Every other transposition to the left in the row is less than $b+1$, so does not effect any change. Therefore, we see that each additional row that is lower (but still above the distinguished path) increases the index of the second $x$ by one. Therefore, the final result will be $x_{row} - x_{b+j+1}$, where $j$ is the number of squares between $Q(i)$ and the distinguished path, if we go down vertically. But we also see that the labelling of the column, which is the labelling on the distinguished path for horizontal line segments, is precisely $Q(m)+1$, where $Q(m)$ is the letter directly above any horizontal line segment along the distinguished path. Therefore, $x_{row} - x_{col}$.

Now we turn to case (2). 

Here, we can write 
\begin{align*}
\beta_i &=\prod_{j=1}^{i-1} r_{Q(j)} (x_{Q(i)} - x_{Q(i)+1}) \\
&= r_{Q(1)} r_{Q(2)} r_{Q(3)} \cdots r_{Q(i-1)} (x_{Q(i)} - x_{Q(i)+1}) \\
&= r_{Q(1)} r_{Q(2)} r_{Q(3)} \cdots r_{Q(p-1)} r_{Q(p)} r_{Q(p+1)} r_{Q(p+2)} \cdots r_{Q(i-1)} (x_{Q(i)} - x_{Q(i)+1}) \\
&= r_{Q(1)} r_{Q(2)} r_{Q(3)} \cdots r_{Q(p-1)} r_{n-k} r_{n} r_{Q(p+2)} \cdots r_{Q(i-1)} (x_{Q(i)} - x_{Q(i)+1})
\end{align*}

where $r_{Q(p)} = r_{n-k}$ is the last number in the section above the distinguished path before going into the letters in the section below the distinguished path, and $r_{Q(p+1)} = r_{n}$ is the first (reading the word from left to right) letter in the section below the distinguished path. By the same arguments as above, we see that $(x_{Q(i)} - x_{Q(i)+1}) = x_b - x_{b+1}$ will first go to $x_{\tilde{i}} - x_{\tilde{j}}$, where $\tilde{i}$ is the leftmost square in the same row as $Q(i)$, and $\tilde{j}$ is the bottommost square directly beneath $Q(i)$ in the section below the distinguished path. 

We can consider a different example within $Gr(4,10)$, where we have circled an $s_{10}$:

\begin{figure}[htbp] \centering
	\includegraphics[scale=0.6,clip=true]{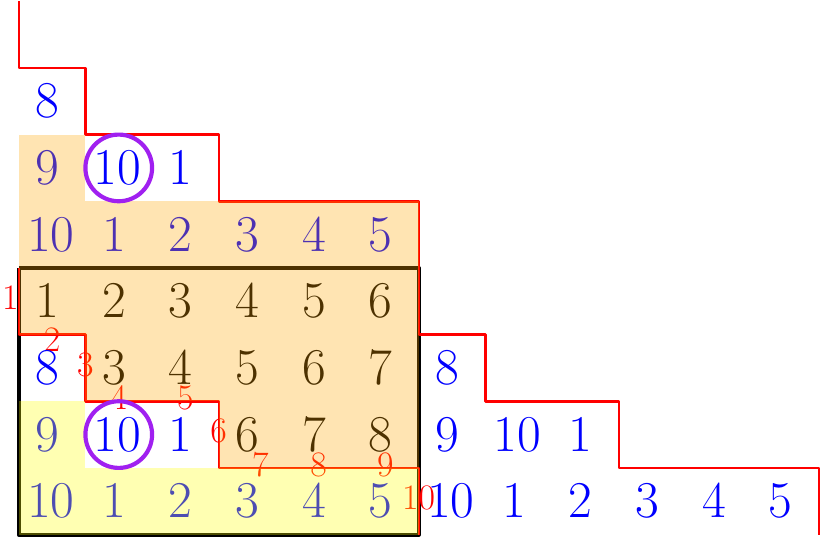}
	\caption{}
\end{figure}

Here, we would have: 
\begin{align*}
\beta_i &= r_{Q(1)} r_{Q(2)} r_{Q(3)} \cdots r_{Q(p-1)} r_{n-k} r_{n} r_{Q(p+2)} \cdots r_{Q(i-1)} (x_{Q(i)} - x_{Q(i)+1}) \\
&= r_6 r_7 r_8 \cdots r_5 r_6 r_{10} r_1 \cdots r_5 r_9 (x_{10} - x_1) \\
&= r_6 r_7 r_8 \cdots r_5 r_6 r_{10} r_1 \cdots r_5 (x_9 - x_1) \\
&= r_6 r_7 r_8 \cdots r_5 r_6 (x_9 - x_2)
\end{align*}

Next, we will be applying all the letters in the section above the distinguished path. This will first move $x_{\tilde{i}} - x_{\tilde{j}}$ to $x_{\tilde{i}} - x_{col}$. Then, it will move $x_{\tilde{i}} - x_{col}$ to $x_{row} - x_{col}$. This is perhaps easier to see by considering the region below the distinguished path as lying above the region above the distinguished path, like in the top portion highlighted in orange in this example. This shows that we will always hit the transposition increasing $\tilde{j}$ successively as we go down in rows, until we hit a bottom horizontal line segment along the distinguished path. Before we hit this horizontal line segment, the $x_{\tilde{i}}$ cannot change; there will be no ambiguity - it will always change $x_{\tilde{j}}$ to $x_{col}$ first, while leaving $x_{\tilde{i}}$ untouched, before changing $x_{\tilde{i}}$ to $x_{row}$. The proof is essentially the following diagram:

\begin{figure}[htbp] \centering
	\includegraphics[scale=0.6,clip=true]{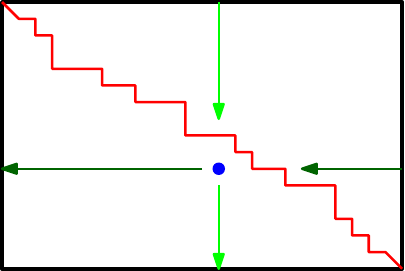}
	\caption{}
\end{figure}

In other words, since the starting point is a spot below the distinguished path, and since, by the arguments above, the letters will move the indices on the $x_b - x_{b+1}$ straight down or left, we see that once we go to the section above the distinguished path, the vertical portion will hit the distinguished path, producing $x_{col}$ before getting to the row that the $x_{\tilde{i}}$ is in. 

It is easy to see that the same arguments apply to the case when there are no squares beneath the distinguished path (or the case when there are no squares above the distinguished path).

\end{proof}

\begin{Exa}
Here, we complete the examples in Gr(4,10) that we introduced in the proof above. Below, we will calculate $\beta_i=\prod_{j=1}^{i-1} r_{Q(j)} (x_{Q(i)} - x_{Q(i)+1})$, showing that it equals $x_{row}-x_{col}$. 

For case (1), we recall that we chose $Q(i)$ to be the letter circled in green, where we have highlighted in orange $\prod_{j=1}^{i-1} r_{Q(j)} = r_{Q(1)} r_{Q(2)} r_{Q(3)} \cdots r_{Q(i-1)}$: 

\begin{figure}[htbp] \centering
	\includegraphics[scale=0.5,clip=true]{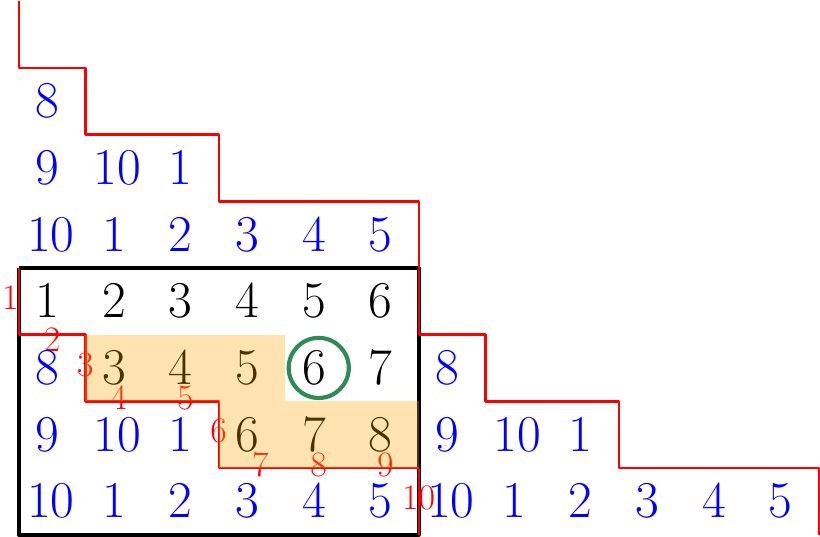}
\end{figure}

We obtain:
\begin{align*}
\beta_i &= r_{Q(1)} r_{Q(2)} r_{Q(3)} \cdots r_{Q(i-1)} (x_{Q(i)} - x_{Q(i)+1}) \\
&= r_{Q(1)} r_{Q(2)} r_{Q(3)} r_{Q(4)} r_{Q(5)} r_{Q(6)} (x_{Q(7)} - x_{Q(7)+1}) \\
&= r_6 r_7 r_8 r_3 r_4 r_5 (x_6 - x_7) \\
&= r_6 r_7 r_8 r_3 r_4  (x_5 - x_7) \\
&= r_6 r_7 r_8 r_3  (x_4 - x_7) \\
&= r_6 r_7 r_8 (x_3 - x_7) \\
&= r_6 r_7  (x_3 - x_7) \\
&= r_6 (x_3 - x_8) \\
&= x_3 - x_8
\end{align*}

and this is precisely $x_{row}-x_{col}$ if we draw lines west and south in the figure, as we wanted to show.

For case (2) of $Q(i)$ below the distinguished path, given the way the word is read in \cite{Sni}, it is more natural to consider the lower section being placed above the distinguished path to form the word we have highlighted in orange below:

\begin{figure}[htbp] \centering
	\includegraphics[scale=0.5,clip=true]{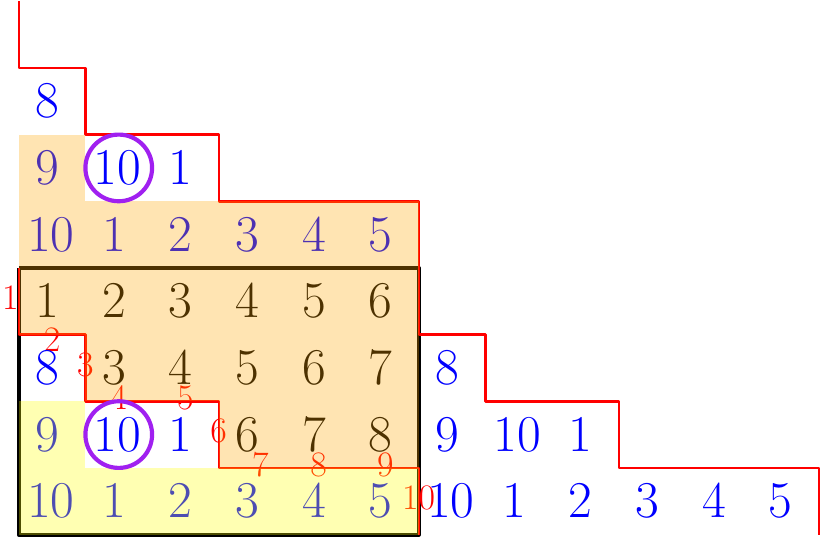}
\end{figure}

Here, we would have: 
\begin{align*}
\beta_i &= r_{Q(1)} r_{Q(2)} r_{Q(3)} \cdots r_{Q(p-1)} r_{n-k} r_{n} r_{Q(p+2)} \cdots r_{Q(i-1)} (x_{Q(i)} - x_{Q(i)+1}) \\
&= r_6 r_7 r_8 r_3 \cdots r_7 r_1 \cdots r_6 r_{10} r_1 \cdots r_5 r_9 (x_{10} - x_1) \\
&= r_6 r_7 r_8 r_3 \cdots r_7 r_1 \cdots r_6 r_{10} r_1 \cdots r_5 (x_9 - x_1) \\
&= r_6 r_7 r_8 r_3 \cdots r_7 r_1 \cdots r_6 (x_9 - x_2) \\
&= r_6 r_7 r_8 r_3 \cdots r_7 (x_9 - x_3) \\
&= r_6 r_7 r_8 (x_9 - x_4) \\
&= x_6 - x_4
\end{align*}
and this also equals $x_{row}-x_{col}$, as we wanted to show.

 \end{Exa}


We just need to show one more thing. Note that the AJS/Billey formula produces an element of the $n$-dimensional torus $T$-equivariant cohomology $H_T^*(V)\cong \Z[y_1,...,y_n]$. However, the theorem from \cite{Rossmann} involves $H_{\C^\times}^*(V)\cong \Z[h]$. Thus, our final step will be to produce a map $H_T^*(V) \to H_{\C^\times}^*(V)$. We will show that this map takes the rightmost term 
$\sum \prod_{r\in R}  \hat{\beta_r}$ and maps it to $(\text{\# pipe dreams on } U_\lambda)\cdot h \in H_{\C^\times}^*(V)\cong \Z[h]$. Then the proof will be complete.

\subsection{Equation \ref{equivcoheq} - Restricting Torus to Circle}\label{subsec:restrictcircle}

Finally, it remains to be shown, under the map from $H_T^*(U_\lambda)$ to $H_{\C^\times}^*(U_\lambda)$, that the right side of the AJS/Billey formula corresponds to the number of pipe dreams.

Recall that we showed above the isomorphisms $H_T^*(U_\lambda) \cong H_T^*(pt) \cong H_T^*(X_\circ^{v(\lambda)})$. However, to use Rossmann's Theorem \ref{thm:Rossmann}, we need to have $H_T^*(V)$ for a vector space $V$. We will get this via an isomorphism of $U_\lambda$ with the vector space $M_{k,n-k}$, the space of $k$-by-($n-k$) matrices. This is depicted in the following diagram, which shows $M_{k,n-k}$ being identified with another space of matrices with the identity in columns $\lambda$, which we can call $(M_{\lambda})_{k,n-k}$; the row span of $(M_{\lambda})_{k,n-k}$ is equal to $U_\lambda \subseteq Gr(k,n)$. 

\begin{tikzcd}
M_{k,n-k} \ar[r,"\sim"] & \left\{ \mqty(&1&0&\\ *&0&1&*\\&\vdots&\vdots&) \right\} \ar[rr,"\sim"]  \ar[drr,shorten >=3ex,"row span"'] & & U_\lambda \curvearrowleft T \hookleftarrow \C^\times \ar[d,hookrightarrow, shift right=5ex]  \\
\left[ v_1\cdots v_{n-k}\right] \ar[u,"\rin"] \ar[r,mapsto] & \mqty(&&1&0&\\ v_1 & v_2 &0&1& \ar[u,"\rin"] \hdots \\&&\vdots&\vdots&) & & Gr(k,n)
\end{tikzcd}

We also recall that the right side of AJS/Billey is: $\sum_{R\subseteq Q, R\text{ reduced}} \prod_{r\in R}  \hat{\beta_r}$, where $Q$ corresponds to the word representing $v(\lambda)$, and then the sum is over all $R$, which are words representing the siteswap $f$. In other words, there is precisely one summand in this formula for each pipe dream representative of $\Pi_f \cap U_\lambda$. We will show that our choice of map $H_T^*(V) \to H_{\C^\times}^*(V)$ takes each of these summands to $1\cdot h^{codimX} \in H_{\C^\times}^*(V)\cong \Z[h]$, so the result will be $(\text{\# pipe dreams of }\Pi_f \text{ on } U_\lambda)\cdot h^{codimX}$. Then we will be done.

The map $H_T^*(V)$ to $H_{\C^\times}^*(V)$ is entirely determined by the map $\C^\times \to T$. We describe the relevant maps explicitly in a moment, but first we give the general theory. 

Suppose we have a map of tori $S\to T$. Then we get the corresponding Lie algebra map $\mf{s} \to \mf{t}$ as well as the dual map $\mf{s}^* \leftarrow \mf{t}^*$. Then we have the following commutative diagram:

$$\begin{tikzcd}
\mf{t}^* \ar[d,hookrightarrow]  \ar[r] & \mf{s}^* \ar[d,hookrightarrow] \\
Sym(\mf{t}^*) \ar[d,leftrightarrow,"\cong"]\  \ar[r] &  Sym(\mf{s}^*) \ar[d,leftrightarrow,"\cong"]\\
H_T^* \ar[r] & H_S^*
\end{tikzcd}
$$

Note that if we pick a basis for $\mf{t}^*$, then $Sym(\mf{t}^*)$ becomes a polynomial algebra, and these basis elements are degree-2 generators (the cohomology is in even degree). Also note that $H_T^* := H_T^*(pt) = H_T^*(V)$ since a vector space $V$ is contractible. Recall that an application of Theorem \ref{thm:Rossmann} required that $\C^\times$ acts with all weights 1. 

\begin{Lem}
We satisfy the condition of the circle $\C^\times \subseteq T$ acting on $U_\lambda$ with all weights 1 due to our choice of map $\C^\times \to T$, which is: \[z\mapsto (z,1,1,z,...,1)\] with the entries $z$ in the coordinates of $[n]\backslash \lambda$.
\end{Lem}

\begin{Exa}
Suppose we are in $Gr(2,4)$ with $\lambda=\{1,3\}$. Then we can represent $U_\lambda$, and the action of $\C^\times$ on the right as:
\[\mqty(1 & * & 0 & * \\ 0 & * & 1 & *) \mqty(\dmat{1,z,1,z}).\]
This scales all the $*$'s, which give coordinates of $U_\lambda \cong \mathbb{A}^{k(n-k)}=\mathbb{A}^4$, by $z^1$. Thus, $\C^\times$ acts via dilation on $U_\lambda$, as required.

Alternatively, we can consider the full action of $T$ on $U_\lambda$ and find the action of $\C^\times$ inside $T$ using:
\[\mqty(\dmat{z_1^{-1},z_3^{-1}}) \mqty(1 & * & 0 & * \\ 0 & * & 1 & *) \mqty(\dmat{z_1,z_2,z_3,z_4}) \] \[ = \mqty(\dmat{z_1^{-1},z_3^{-1}}) \mqty(z_1 & z_2* & 0 & z_4* \\ 0 & z_2* & z_3 & z_4*) = \mqty(1 & \frac{z_2}{z_1}* & 0 & \frac{z_4}{z_1}* \\ 0 & \frac{z_2}{z_3}* & 1 & \frac{z_4}{z_3}*).\]

We see that in order to maintain our matrix representative in the form $\mqty(1 & * & 0 & * \\ 0 & * & 1 & *)$ we must ``unscale" on the left, which results in the coordinates having factors $\frac{z_{col}}{z_{row}}$, where the number $row$ refers to: within this number's row, find where the $1$ is, and record the column of this $1$.
\end{Exa}

\begin{proof}
Note that for a general embedding of $\C^\times $ in $T$ acting with arbitrary weight, we have $z_i=z^{y_i}$. Then $\frac{z_{col}}{z_{row}}= z^{y_{col}-y_{row}}$. Since we want the action to be by dilation, we need all these $y_{col}-y_{row}=1$. We just saw that we could do this via the map $z\mapsto diag(z,1,1,z,...,1)=diag(z^1,z^0,z^0,z^1,...,z^0)$ with the $z$'s in the non-$\lambda$ spots. In other words, by setting \\
$y_i= \begin{cases} 0 & i\in \lambda \\ 1 & i\not\in \lambda \end{cases}$, \\
then $y_{col}-y_{row}=1-0=1$ since the $col$ corresponds to the columns not in $\lambda$ and the $row$ correspond to the unit vectors in the matrix which are where $\lambda$ is. 
\end{proof}

If the $T$-action on the tangent spaces to the fixed points on a flag manifold $G/P$ contains dilation, $G/P$ is called cominuscule. We have just recovered the fact that Grassmannians are cominuscule (something not true for full flag manifolds).

Let's now write the remaining maps explicitly.
$\C^\times \to T, z\mapsto diag(1,z,z,1,....,1,z)$ where the $z$'s exist precisely in the places in the diagonal where $\lambda$ is not gives us:
$\mf{s} \to \mf{t}, 1 \mapsto [1,0,0,1,...,0]$ the map of Lie algebras, where we have just taken the derivative of the map above and extracted the diagonal.

$ \mf{s}^* \leftarrow \mf{t}^* : \mqty(1\\0\\0\\1\\ \vdots \\ 0)$ 

Finally, once we go to $Sym(\mf{t}^*) \to Sym(\mf{s}^*)$ this is just a map of polynomial rings, where this map $\mqty(1 & 0 & 0 & 1& \cdots & 0)^T$ sends all variables $y_i$ in $Sym(\mf{t}^*)$ corresponding to the numbers in $\lambda$ to $0$, and all other variables to $h$. 

Recall Proposition \ref{prop:AJSrowcolform}, which noted that 
\[ [\Pi_f \cap U_\lambda]_T = \sum_{\substack{R\subseteq v(\lambda) \\ \prod R=v(f) \\ R \text{ reduced}}} \prod_{r\in R}  [(y_i)_r - (y_j)_r] \] where $i$ is the column and $j$ is the row of a cross in a pipe dream.

Now, we can apply to this our choice of map $\tilde{\lambda}:\C^\times \to T$ sets all $y_{row}=0, y_{col}=h\in \Z[h]$, so since each summand in the AJS/Billey formula is a product of $codim(\Pi_f)$ many $\beta$'s, each summand will be $1\cdot h^{codim \Pi_f}$. As we already showed above, there is a single summand for each pipe dream representative of $\Pi_f$ on $U_\lambda$, so the AJS/Billey formula is equal to $(\text{\# pipe dreams of }\Pi_f \text{ on } U_\lambda)\cdot h^{codim \Pi_f}$. This completes the proof of the uniqueness (smoothness) part of Theorem \ref{thm:smoothone}.

Finally for the proof of the existence part, that $\lambda \in \Pi_f$ is equivalent to there existing at least one affine pipe dream, this follows from the fact that each affine pipe dream for $\Pi_f$ on $U_\lambda$ is a subword for $\Pi_f$ in the word for $\lambda$ in accordance with what we have already described in Section \ref{affinepipedreams}. In addition, it is a general fact of words in Bruhat order that a point $\lambda$ being contained in a Schubert variety defined by $f$ is equivalent to the word $v(\lambda) > v(f)$ in Bruhat order, and furthermore, this is equivalent to $v(f)$ being a subword in a reduced for for $v(\lambda)$.

\section{Smoothness in Relation to the Ordering under Deletion/Contraction}\label{sec:smoothnessordering}

Smooth positroid varieties remain smooth under deletion and contraction, while singular positroid varieties if obtained from another positroid variety via deletion or contraction must have come from a singular positroid variety. We can state this in the following proposition:

\begin{Prop}\label{prop:delsmooth}
If $\lambda$ is smooth on $\Pi_f$, then $del_i(\lambda)$ is smooth on $\Pi_{del_i(f)}$ and $contr_i(\lambda)$ is smooth on $\Pi_{contr_i(f)}$. Or this can be stated contrapositively: if $del_i(\lambda)$ is singular on $\Pi_{del_i(f)}$ or $contr_i(\lambda)$ is singular on $\Pi_{contr_i(f)}$, then $\lambda$ is singular on $\Pi_f$.
\end{Prop}

This means that in the ordering (given by deletion and contraction) of positroid varieties and points, the smooth versus singular pairs can be divided in the following way:

 \begin{figure}[htbp] \centering
	\includegraphics[scale=0.5,clip=true]{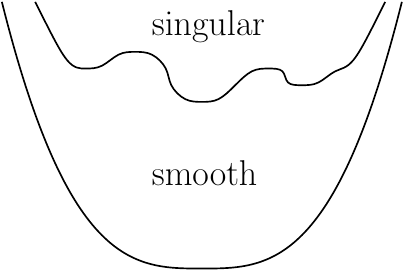}
	\caption{}
\end{figure}

\begin{proof}[Geometric Proof]

By e.g. \cite{BilBra03}, if the group $G$ is reductive, then whenever $G$ acts on a regular set $X$, the fixed point set $X^G$ is also regular. In our case, $G=T$ is the torus, so the conclusion follows. 

In particular, if $\Pi_f$ is smooth at $p$, then $\Pi_f^{\C^\times}$ is smooth at $p$, and $\Pi_f^{\C^\times}=(\Pi_f\cap \Pi_{del,i}) \cup (\Pi_f\cap \Pi_{contr,i}) = Del_i(\Pi_f) \cup Contr_i(\Pi_f)$. 
\end{proof}

We can give a combinatorial proof of this proposition using the machinery of affine pipe dreams:

\begin{proof}\label{proof:singDelsingcombinatorial}[Combinatorial Proof]

Without loss of generality, let us do deletion on the $i$th column in a rigid pipe dream (the case of contraction will be analogous but for a row). 

The column will alternate between consecutive sets of crosses and consecutive sets of elbows. It will be helpful for the proof to look at the following figure which zooms in on the $i$th column, where the consecutive sets, for simplicity, have been drawn as groups of 3:

 \begin{figure}[htbp] \centering
	\includegraphics[scale=0.7,clip=true]{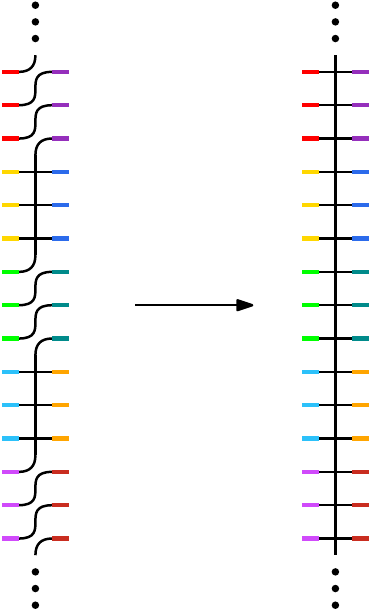}
	\caption{}
\end{figure}

Note that the colors on the two pipe dreams in the figure are the same; for each one, reading from top to bottom, the left side (of column $i$) goes: light red - light orange - light green - light blue - light purple; and on the right side: dark purple - dark blue - dark green - dark orange - dark red. The fact that we have flipped the order of the coloring (while doing dark instead of light to distinguish them) is intentional: the pipe dream has a northwest-southeast symmetry, coming from the fact that the pipes (and the moves) go southwest to northeast. Thus, it is sufficient to make arguments on the just left side (light) colors, so for ease of reference, we call the left/light colors just: red, orange, green, blue, purple; whereas if we need to distinguish the colors on the right side, we will call them: dark purple, dark blue, cyan (or dark green), dark orange, dark red. The green to cyan pipes in the middle are where we will focus. After all, only the elbow tiles are changed in the deletion. 

Therefore, what we are doing here is picking some random block of consecutive elbow tiles in the column, and showing that, even though those tiles are changed (turned into crosses), no new moves can be created. If we can show this for the green to cyan tiles in this diagram, then this applies to any other elbow tile in this column that gets deleted, and the proof will be complete.

If we need to further distinguish the pipes within each block, we can call them, for example: red 1, red 2, red 3, going from top to bottom so that red 1 is the highest. 

We first make several observations. 

(1) There can be no crossing among any of the green pipes anywhere because they are elbows/near-misses along the $i$th column; any crossing would then provide a move (specifically, the cross could then be moved to the $i$th column), contradicting rigidity. For the same reason, there can be no crossing among the cyan pipes. (1a) Note, however, that the orange pipes have no such restriction among themselves.

(2) The green pipes cannot cross any of the orange pipes. The reason is that, by (1), in order for any of the green pipes to cross any of the orange pipes, green 1 would have to cross (first). However, since green 1 crosses all of the orange pipes along column $i$, this would create a double crossing. 

(3) The green pipes cannot cross any of the red pipes. The reason is that, by (1), in order for any of the green pipes to cross any of the red pipes, green 1 would have to cross red 3 (first). However, since green 1 and red 3 have a near miss, this would create a move, contradicting rigidity. The same argument show that the purple pipes cannot cross any of the green pipes. 

(4) The green pipes \textit{can} cross the blue pipes. 

(5) By symmetry, all these points are valid for the right side of column $i$ as well: just replace ``color" by ``dark color."

We can recap all of this succinctly. The only crossings among pipes that are allowed are: (I) \textit{within} each of block pipes coming from consecutive crosses along the $i$th column, and (II) between a block of consecutive elbows and the block of consecutive crosses directly below them (on the left side of column $i$; or, on the right side, by symmetry, between consecutive elbows and consecutive crosses directly above)

We now show that doing the deletion does not create any new moves. 

The deletion only changes the connectivity of the elbow tiles. In other words, the orange and blue pipes do not change, so there cannot be any new moves coming from the orange or blue pipes. Let us focus in on the green pipes, and ask the question of whether there can be any new moves coming from the change in the connectivity of the green pipes.

A move consists of 2 pipes having a cross and a near-miss. Obviously, many new crosses are created along column $i$; however, these crosses cannot participate in any move because the $i$th column pipe runs straight up and down vertically, crossing every pipe along the $i$th column, and thus has no near-misses. Therefore, a new move must come from a crossing of 2 pipes on the left side of column $i$, for which these 2 same pipes have a near-miss on the right side of column $i$ (or vice-versa switching left and right, but it is sufficient to just consider the former case). Since this is a new move, it did not exist before.  The cross and near-miss cannot have both been on the same side (left or right); otherwise, this move would have existed prior to the deletion and thus not only not be a new move, but the original pipe dream would have had a move, contradicting rigidity.

Thus, the question now becomes, (a) what possible crosses involving the green pipes on the left are there, and (b) what near-misses involving these pipes are there on the right. The work we did above now comes in handy. We know that the green pipes cannot cross among themselves by (1), and in fact cannot cross any pipes except the blue ones by (II). This answers (a). Now we consider the question of (b). Since the green pipes can only potentially cross the blue pipes, the question is whether it is possible for a potential green-blue crossing on the left, after we do the deletion which causes these pipes to turn into dark green - dark orange pipes on the right, to involve a dark green - dark orange near-miss on the right.

This is not possible for the following reason. Pipe dark green 3 crosses all three of the dark orange pipes. Since the pipe dream is rigid, that means that dark green 3 cannot have a near miss with any of the dark orange pipes. By (1) above, none of the dark green pipes can cross each other, which mean that none of the dark green pipes can have a near miss with any of the dark orange pipes. This completes the proof.

\end{proof}

However, our work above allows us to say something further about this last point. In order for dark green 3 to not have a near miss with any of the dark orange pipes, we must have another color pipe come between the dark green and dark orange pipes. We saw above in (4) or (I) that, on the left, the elbow pipes could cross with the cross pipes in the block below, which, on the right, corresponds to elbow pipes crossing with the cross pipes above. This means that one, or any (since the dark blue can cross among themselves by (1a)), of the dark blue pipes can come down by crossing the dark green pipes to create separation from the dark orange pipes; or, dark red pipe 1 would have to come up (the other dark red could also come up, but dark red pipe 1 would have to come up first by (1)).

Given a pipe dream with a maximum number of crossings along column i, we can view the column $i$ to be deleted, once we find a pipe dream that has a maximal number of crosses along column $i$, as consisting of alternating sets of consecutive crosses followed by consecutive elbows. The pipes entering crosses along column $i$ do not change after deletion. The pipes entering elbows along column $i$ do change, but note that because we have moved the maximum number of crosses to this column, this implies that there are no crossings among any of the pipes entering a set of consecutive elbows along column $i$. In these cases, deletion and contraction do not create a new moves that did not exist prior to the deletion or contraction. The only case in which a new move can be created is at the border between a set of crosses and a set of elbows along column $i$ (the details are a bit long and tedious, so are not included here), and is depicted in the following figure:
 \begin{figure}[htbp] \centering	\includegraphics[scale=0.7,clip=true]{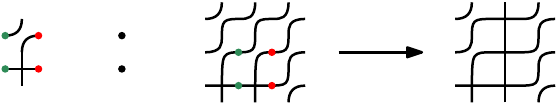}	\caption{}  \end{figure}
In more detail, if we have a cross followed by an elbow directly above it in column $i$, then if the pipes emerging from the green dots on the left cross at some point, and the pipes emerging from the red dots on the right have a near-miss at some point, then a new move is created. The figure gives an example. (If the cross is directly above the elbow, then we have to switch the ``left" and ``right" in the previous statement). This is the only case where deletion or contraction can create a new move, but this required that a move already existed, so the positroid variety was singular at the point already. If the positroid variety was smooth, then this exception cannot occur, so no new moves could be created.

Thus, if there were no moves prior to the deletion, then there are no moves after the deletion. The contrapositive statement is that a singular positroid variety $\Pi_{f'}$ coming from another positroid variety $\Pi_f$ via deletion or contraction implies that $\Pi_f$ is itself singular, since a smooth positroid variety being deleted to a singular positroid variety is impossible.

\chapter{Proof of the Main Theorem}\label{ch:maintheoremproof}

\section{Preliminary Notions}

Since the following parts involve the notion of a maximal rectangle, we include a figure to display what the maximal rectangles look like in our pipe dream shapes:

 \begin{figure}[htbp] \centering
	\includegraphics[scale=0.5,clip=true]{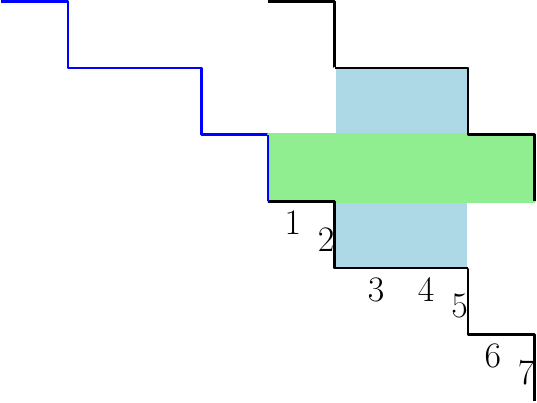}
	\caption{}
\end{figure}

\begin{Def}\label{def:partition}
A pipe dream in a rectangular shape \textbf{reduces to NW/SE partitions} if, after deleting all entire rows and columns of crosses, the result is a partition in the southeast and northwest corners. 

This is equivalent to the condition that all crosses must be contained in: (1) entire rows of crosses, (2) entire columns of crosses, (3) a partition in the northwest, (4) a partition in the southeast. In particular, a rectangular pipe dream that reduces to NW/SE partitions is one that d/c reduces (see Definition \ref{def:dcreduce}) to partitions (in the NW and SE corners), by removing specifically those (1) full rows of crosses and (2) full columns of crosses. The notion of ``d/c reduces" is more general, since it can be applied to nonrectangular pipe dreams and deletion/contraction applied to rows and columns that include elbows.
\end{Def}

\begin{Exa}\label{exa:SENWpart}
The following is an example of a rectangular pipe dream that reduces to NW/SE partitions:
 \begin{figure}[htbp] \centering
	\includegraphics[scale=0.5,clip=true]{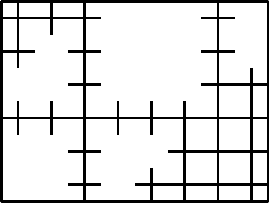}
	\caption{}
\end{figure}
\end{Exa}

\section{Statement and Proof of the Main Theorem}

\begin{Thm}\label{thm:mainthm} [Main Theorem]
Let $\delta \in PD_{(\Pi_f,\lambda)}$ be an affine pipe dream for $\Pi_f \cap U_\lambda$. If each maximal rectangle in $\delta$ reduces to NW/SE partitions, then $\Pi_f$ is smooth at $\lambda$; otherwise, it is singular. 
\end{Thm}

\begin{proof}

This theorem consists of showing two parts:\\ 
(1) Showing that smoothness can be determined just by looking at the maximal rectangles in an affine pipe dream. This is proven in Proposition \ref{prop:insidemaxrect} below. \\
(2) Showing that the thing that needs to be looked at on each maximal rectangle is that each maximal rectangle is a NW/SE partition. This is proven in Theorem \ref{thm:NW/SE} below. 

However, first we give the full proof, referencing both Proposition \ref{prop:insidemaxrect} and Theorem \ref{thm:NW/SE} as needed. 
\\

We want to show that smoothness of $(\Pi_f,\lambda)$ is equivalent to $\delta$ having no moves \textit{in any maximal rectangle}. But we already have Theorem \ref{thm:smoothone}, which says that smoothness of $(\Pi_f,\lambda)$ is equivalent to its pipe dream having no moves (anywhere), and, by Proposition \ref{prop:insidemaxrect}, all moves exist within maximal rectangles. Thus, if all maximal rectangles have no moves, then the pipe dream for $(\Pi_f,\lambda)$ has no moves. The other direction, that the pipe dream for $(\Pi_f,\lambda)$ having no moves implies that all its maximal rectangles have no moves, is automatic. 

Therefore, smoothness of $(\Pi_f,\lambda)$ is equivalent to $\delta$ having no moves in any maximal rectangle.
Thus, we just have to show that rigid rectangles reduce to NW/SE partitions and the proof will be complete, and this is shown in Theorem \ref{thm:NW/SE}.
\end{proof}

\section{First Main Idea Used in Proof}

\begin{Prop}\label{prop:insidemaxrect}
Every cross-elbow move exists within some maximal rectangle. More precisely, every pair consisting of a cross and a near-miss where this cross can move is contained within a maximal rectangle. 

Trivially, this implies that the portions of the two pipes involved in the move, starting from the cross and running all the way to the near-miss, are entirely contained in this maximal rectangle.
\end{Prop}
In fact, the proof below shows the stronger statement that if we consider a cross and the set of all places where it can move, these will all be contained in a maximal rectangle.  
\begin{proof}
The very notion itself of nonrigidity is that there is a cross to the southwest and a near-miss of the pipes in that cross somewhere to the northeast (or vice-versa with southwest and northeast switched). In the sorts of shapes given by affine pipe dreams we are considering, such a cross and then near-miss combination must exist within a full rectangle because all the jagged (or convex) parts of the shape are in the northeast and southwest.

This is depicted in the red rectangle in the following figure: 

\centerline{
	\includegraphics[scale=0.5,clip=true]{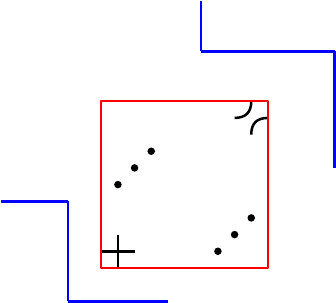}
}
\end{proof}

\section{Second Main Idea Used in Proof}


\begin{Lem}\label{lem:fillinrect}
In a reduced pipe dream, any time there is an elbow followed by $m$ consecutive crosses to its east and $n$ to its south,  then the $n\times m$ box formed by having the easternmost and southernmost crosses as corners must be entirely filled in with crosses. The same is true with ``east and south" replaced by ``west and north." That is: 

 \begin{figure}[htbp] \centering
	\includegraphics[scale=0.5,clip=true]{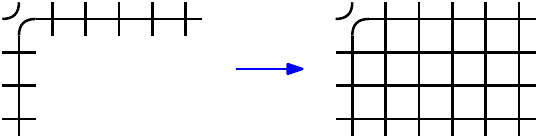}
	\caption{}
\end{figure}
\end{Lem}

\begin{proof}
The following figure provides the proof of this lemma.
 \begin{figure}[htbp] \centering
	\includegraphics[scale=0.5,clip=true]{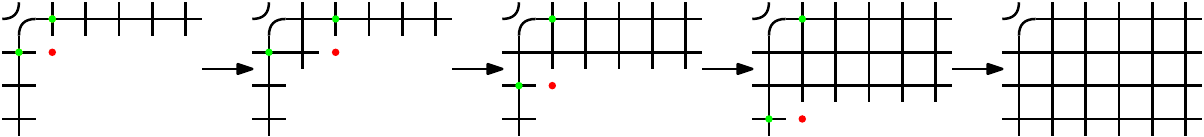}
	\caption{}
\end{figure}
Consider the red dot in the first figure on the very left. It cannot be an elbow, or else the two green dots would be a double crossing (non-reduced); so the red dot must be a cross. So once we have filled that tile in with a cross, move the red dot to the right, where it again must be a cross, by a similar argument. We continue inductively until the entire first row is filled in with crosses. Then we move on to the 2nd row, and apply the same argument to fill it in with crosses. We continue on, row by row, until the $n\times m$ rectangle is filled in with crosses. 
\end{proof}

\begin{Prop}\label{prop:Le}
A rectangular pipe dream is rigid if and only if: (1) no cross has both an elbow to its north as well as to its east, and (2) no cross has both an elbow to its south as well as to its west.

\end{Prop}

This is essentially a ``double-Le condition," and the proof is similar to statements in Section 19 of \cite{Post} or Section 4.3.3 of \cite{Sni}. As hinted in \cite{Sni}, one of these conditions basically corresponds to the bottom pipe dream (see Appendix B), so the other would correspond to the top pipe dream.

\begin{proof}

First we prove that if a pipe dream satisfies this double-Le condition, then it is rigid. We show the contrapositive, that if a pipe dream is not rigid, then it must violate this condition. The proof is essentially the following picture:
\\

 \begin{figure}[htbp] \centering
	\includegraphics[scale=0.5,clip=true]{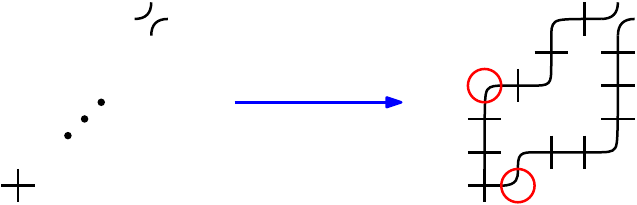}
	\caption{}
\end{figure}

Let's assume that this pipe dream is not rigid. This means that, as shown at the top of the figure, there exists two pipes that cross somewhere and later have a near-miss, either to the northeast of the cross as in the figure, or to the southwest (not depicted, but analogous). This requires the two pipes to each have a ``first bend," circled in red in the representative example bottom left figure, because if they never bent, then they would continue on indefinitely north and east and never cross. Thus, a pipe dream with a move either contains a cross with elbows both to its north and east (as shown in this figure), thus violating the north/east Le condition, or it has a move with the figure reflected across the $x=-y$ line; that is, with a violation of the south/west Le condition. 
\\

\minip{.77}{12mm}{ Now we prove the other direction: if a pipe dream violates this double-Le condition, then it must not be rigid. Without loss of generality, let's again assume that it violates the north/east Le condition (the south/west case is analogous). This means that there is a cross with elbows somewhere to its north and east, and we can then assume the elbows to be the closest in these directions. That is, we have:
}
\hfill
\minip{.20}{0mm}{
	\includegraphics[scale=0.4,clip=true]{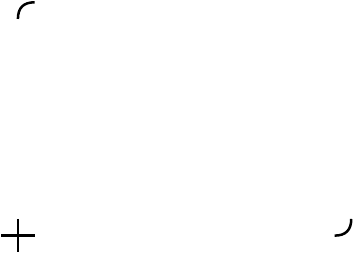}
}
\vskip .05in

\minip{.77}{12mm}{ Furthermore, once they bend at an elbow, they cannot keep moving along straight lines, as depicted in the dotted lines in the following figure, or else they would cross again, which is violation of reducedness:
}
\hfill
\minip{.20}{0mm}{
	\includegraphics[scale=0.53,clip=true]{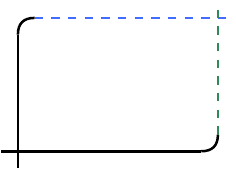}
}
\vskip .05in

\minip{.77}{12mm}{Because of this, at least one of the pipes must have another bend. In fact, we must have each pipe bending a 2nd time because otherwise we would be able to apply Lemma \ref{lem:fillinrect} to this rectangular region of interest to arrive at a nonreduced pipe dream, like the following figure: 
}
\hfill
\minip{.20}{0mm}{
	\includegraphics[scale=0.53,clip=true]{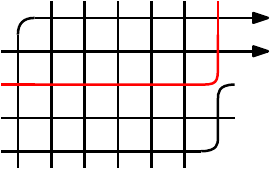}
}
\vskip .05in

\minip{.77}{12mm}{ Thus, putting in these additional bends, and filling in the parts to the west and south of the bends with crosses via an application of Lemma \ref{lem:fillinrect}), the pipe dream looks something like:
}
\hfill
\minip{.20}{0mm}{
	\includegraphics[scale=0.5,clip=true]{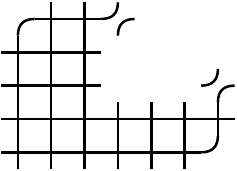}
}
\vskip .05in

In other words, by imposing a violation of the double Le condition and reducedness (simply having a valid pipe dream), we have concluded that we have two pipes that cross, that each bend, and then bend again, and everything to the southwest of the rectangle formed by the places where their second bends are must be filled in with crosses.
\\

\minip{.77}{12mm}{ Next, we apply the assumption of rigidity we made for the sake of contradiction. In the following figure, if the elbow on the right were to bend inwards, then we would have the red elbow depicted, and there would be the possible move depicted by the two red dots, which is not allowed by rigidity:
}
\hfill
\minip{.20}{0mm}{
	\includegraphics[scale=0.5,clip=true]{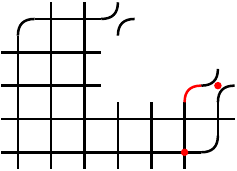}
}
\vskip .05in

\minip{.77}{12mm}{ So we must have the green cross below. But then, inductively, we can make the same argument for the spot next to it, for if it was an elbow, then we would have the possible move with the two red dots:
}
\hfill
\minip{.20}{0mm}{
	\includegraphics[scale=0.5,clip=true]{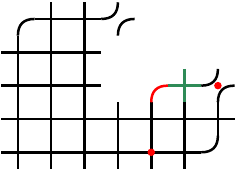}
}
\vskip .05in

\minip{.77}{12mm}{ Thus, it is clear that by this argument, we must fill in everything along the same lines as the elbows with crosses, to get the following:
}
\hfill
\minip{.20}{0mm}{
	\includegraphics[scale=0.5,clip=true]{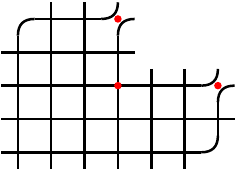}
}
\vskip .05in

\minip{.77}{12mm}{ If we look at the red dots, we see that the configuration that we started with (with two northwest southeast elbows and a cross at one of the corners) is here, only smaller. We can keep doing this inductively, but due to the finite size, we cannot go on forever. We must reach a point where we have a move, giving us the contradiction:
}
\hfill
\minip{.20}{0mm}{
	\includegraphics[scale=0.5,clip=true]{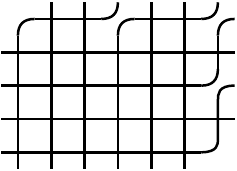}
}

\end{proof}

\begin{Cor}
If a pipe dream is rigid within a rectangular shape, then every time there is a cross with an elbow to its south, then it must have crosses extending all the way west. We can visualize it as follows:

 \begin{figure}[htbp] \centering
	\includegraphics[scale=0.35,clip=true]{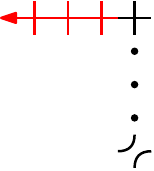}
	\caption{}
\end{figure}

The same statement is true with ``west" and ``south" reversed, as well as replaced with ``north" and ``east". 
\end{Cor}



We note that the previous Prop \ref{prop:Le} was only valid within rectangular shapes. In nonrectangular shapes, a violation of the Le condition is not sufficient to conclude nonrigidity. A counterexample is provided by the following figure, where the two red dots represent a Le-violating condition, yet the portion of this pipe dream depicted is perfectly rigid. 
 \begin{figure}[htbp] \centering
	\includegraphics[scale=0.5,clip=true]{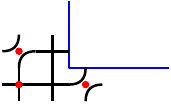}
	\caption{}
\end{figure}

The following corollary, however, gives a condition for rigidity in any pipe dream shape (not just rectangular ones).

\begin{Cor}\label{cor:box} [of the proof]

An affine pipe dream is rigid if and only the following two types of rectangular pipe crossings do not appear anywhere:

 \begin{figure}[htbp] \centering
	\includegraphics[scale=0.5,clip=true]{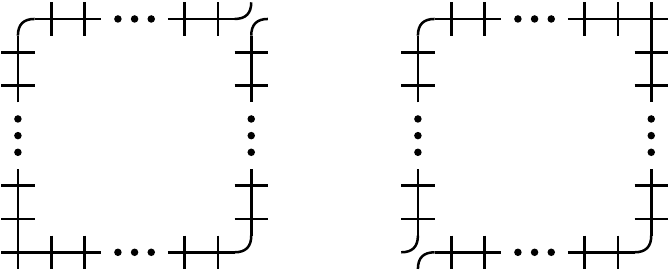}
	\caption{}
\end{figure}
\end{Cor}

{\em Remark 1.} This corollary does not state that in a nonrigid \textit{pipe dream}, that every nonrigid \textit{pair of pipes} must be of this form. In fact, two pipes that cross and bend inwards at some point before having a near-miss is indeed possible; however, even in this case, the nonrigid pipe dream must contain nonrigid pairs of pipes of the form depicted in this corollary. The following figure gives an example, where the red and blue types have a move that is not of the type stated in the corollary, but there are 3 pairs of pipes that do have a move of the stated type, which are depicted in the 3 pairs of (different shades of) green dots: 

 \begin{figure}[htbp] \centering
	\includegraphics[scale=0.6,clip=true]{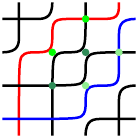}
	\caption{}
\end{figure}

{\em Remark 2.} First, we recall that the main theorem gives a criterion for testing smoothness (equivalent to rigidity) which is checkable on maximal rectangles. The fact that we only need to check on maximal rectangles is evident from this corollary, since this sort of pipe dream configuration can only exist within rectangular shapes. Of course, this is just a special case of Proposition \ref{prop:insidemaxrect}.

\begin{proof}
The proof of Corollary \ref{cor:box} is already in the proof of the Proposition \ref{prop:Le}. One direction is obvious: a box configuration as depicted implies nonrigidity. Conversely, if we have something that is nonrigid, then we have a cross and then a near-miss for a pair of pipes, like in the previous figure. Then these pipes must bend somewhere. The south and west portions of these bends must be filled in with crosses via the Lemma, as explained in the previous proof. Then, if there are any elbows in the remaining portion, they will form a rectangular configuration like that depicted. If not, then we arrive at a smaller box, as explained previously. 
\end{proof}

Note that this rule, that neither of the two box configurations above occur, encompasses the sorts of moves considered in \cite{BB93}. Specifically, they considered moves of the following form:

 \begin{figure}[htbp] \centering
	\includegraphics[scale=0.5,clip=true]{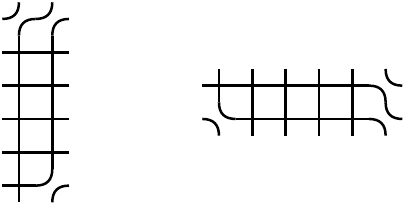}
	\caption{}
\end{figure}

These can be viewed as size $(1,j)$ or $(j,1)$ boxes of the form in the previous corollary. 

\begin{Thm}\label{thm:NW/SE}
All rigid rectangles reduce to NW/SE partitions.
\end{Thm}

\begin{proof}

What we will prove below is that every cross (schematically shown by circling the crosses in red) is contained in one of the four configurations below.

\begin{figure}[htbp]\centering
	\includegraphics[scale=0.5,clip=true]{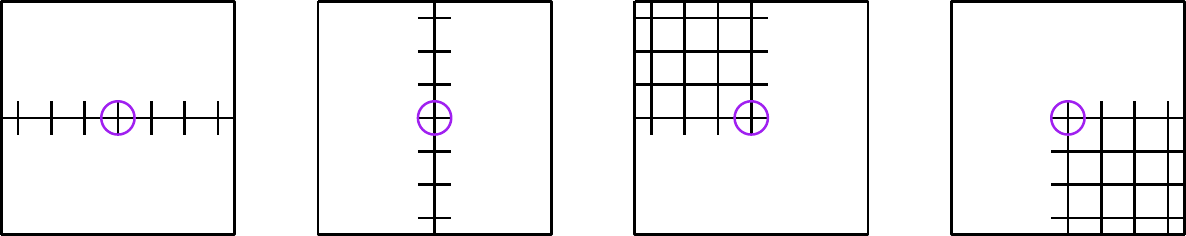}
	\caption{}
\end{figure}

The proof is mostly constructive, and proceeds by going through the following diagram:

  \begin{figure}[htbp]\label{thegrid}\centering
	\includegraphics[scale=0.45,clip=true]{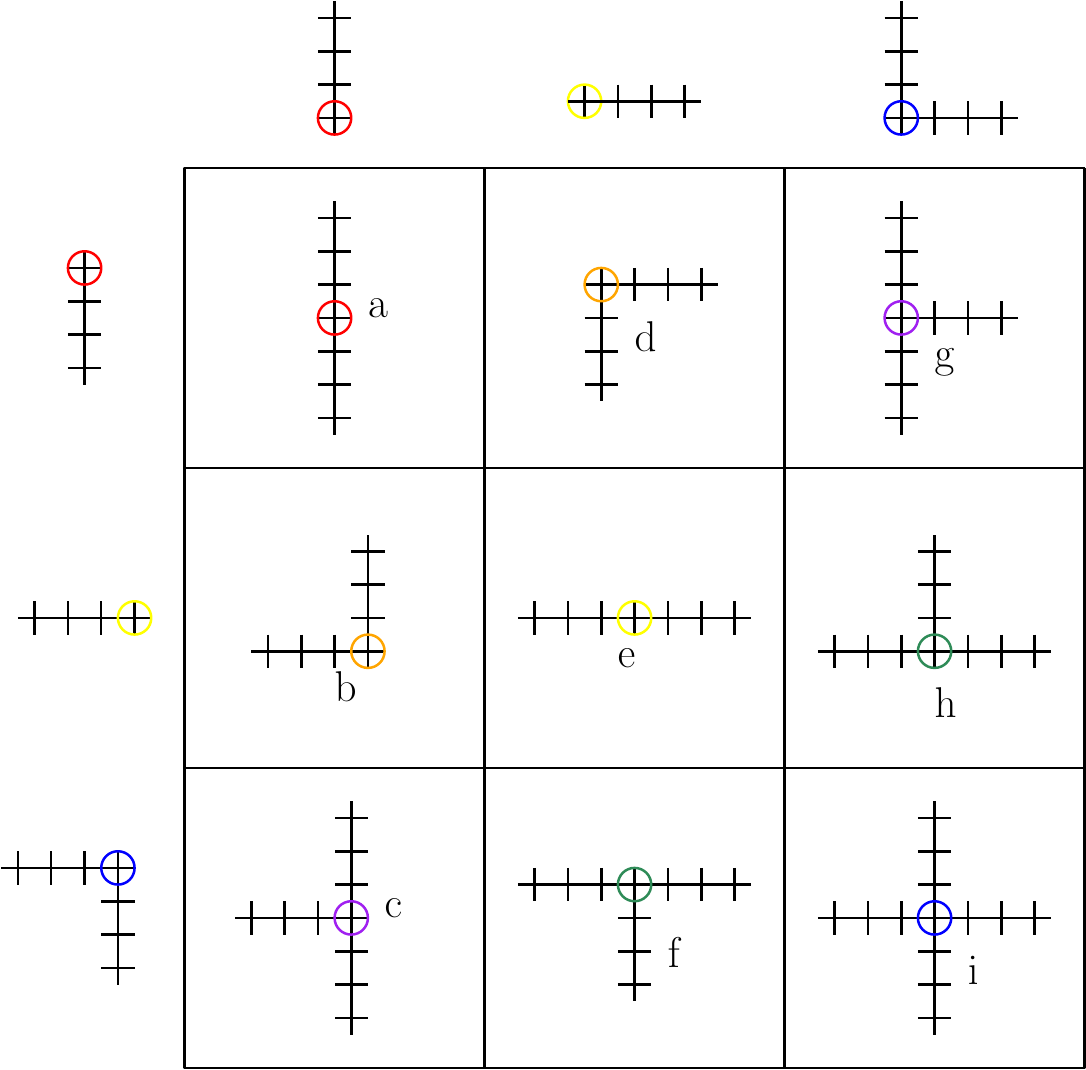}
	\caption{Cases of the theorem}
\end{figure}

We can think of these as two steps, which we can call NE and SW. These come from the previous theorem which shows that any cross within a rigid diagram cannot have elbows both to its north and east, or both to its south and west. Thus, to ensure that both the NE condition and SW condition are satisfied, we can think of doing 2 steps: ensuring that the NE condition is satisfied, and ensuring that the SW condition is satisfied. 
\\

Looking at the diagram, we start with a single cross, which we mark by giving it a circle. The top (separating the columns) then shows the result of applying the NE condition to this cross. Any time we have a series of three crosses coming out from the circled cross, that pictorially depicts an unbroken line of crosses stretching all the way to the wall. Thus, for example, the leftmost configuration at the top shows the cross (circled in red), with a line of crosses stretching north all the way to the wall. This satisfies the NE condition, since it can no longer have elbows \textit{both} to its north and east. If we then apply the SW step to this configuration, then we obtain diagrams (a), (b) and (c). Diagram (a) is obtained by satisfying the SW rule by extending south, (b) by ending west, and (c) by extending both south and west. Diagrams (d)-(f), and (g)-(i) follow the same pattern of elongation starting from their diagram satisfying the SW condition.
\\

Now, each (a), (e), and (i) is already contained in an entire row or column of crosses (or union thereof), so they already satisfy the conditions of the theorem. Thus, we have to deal with the rest. We do the argument for two of these: (c) and (d). The others follow by analogous arguments, and are left to the reader.
\\

\minip{.77}{12mm}{ (c) looks like the following:
}
\hfill
\minip{.20}{0mm}{
  \includegraphics[scale=0.5,clip=true]{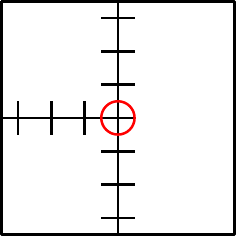}
}
\vskip .05in

\minip{.77}{12mm}{ The crosses along the vertical line all satisfy both NE and SW conditions. However, note that the empty space is unknown: at the moment, it can be either crosses or elbows. Thus, if we look at the crosses boxed in grey in the next figure, they do not satisfy the NE condition if the empty space is all elbows. 
}
\hfill
\minip{.20}{0mm}{
  \includegraphics[scale=0.5,clip=true]{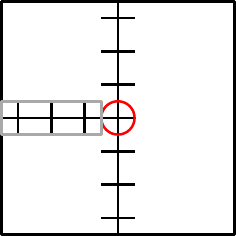}
}
\vskip .05in

\minip{.77}{12mm}{ In order to have these crosses satisfy the NE condition, we can extend them east, or extend them north. If we extend them east, then this extra set of crosses is sufficient for satisfying the NE condition for all the crosses boxed in grey:
}
\hfill
\minip{.20}{0mm}{
  \includegraphics[scale=0.5,clip=true]{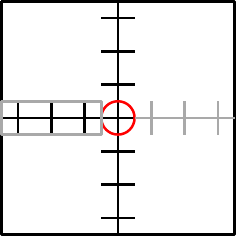}
}
\vskip .05in

\minip{.77}{12mm}{ However, if we cannot extend east (that is, if there is any elbow to the east of the cross circled red), then we must extend north individually for each of the crosses boxed in grey:
}
\hfill
\minip{.20}{0mm}{
  \includegraphics[scale=0.5,clip=true]{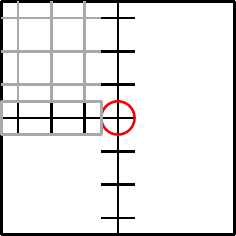}
}
\vskip .05in

In fact, we have shown something slightly stronger than the theorem. The configuration we started with here already had the cross circled in red contained within an entire column of crosses. Thus, what we showed is that if a cross is contained within an entire column or row of crosses, and the perpendicular direction has a line of crosses extending all the way to the wall in only one of the two possible directions, then it must be part of a northwest or southeast partition. 
\\

Now let us take a look at (d):

\begin{figure}[htbp] \centering
	\includegraphics[scale=0.5,clip=true]{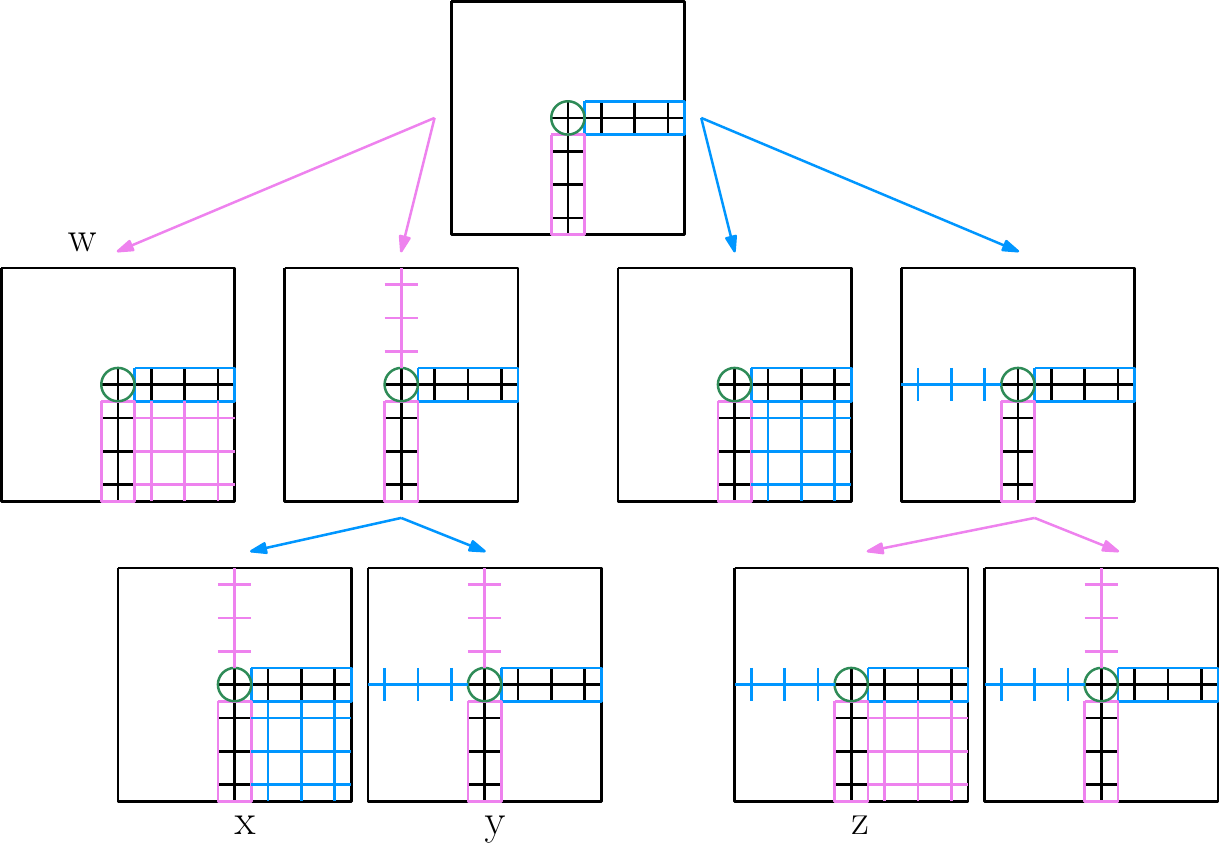}
	\caption{}
\end{figure}

Here we have boxed the crosses in pink and light blue. The pink crosses satisfy the SW condition but not the NE. The light blue crosses satisfy the NE condition but not the SW. Following the arrows through the figure, where the pink arrows represent resolving the SW condition, and analogously for the light blue arrows (the details are the same as the description for the case of configuration (c) above), we obtain 6 final configurations. Some are redundant, so the unique ones are labelled (w), (x), (y), and (z). 
\\

It is clear that by symmetry, all other cases in (a)-(i) will just be certain reflections/rotations of what we have already shown. This proves the theorem. 

\end{proof}

\begin{proof}[Third proof of Proposition \ref{prop:delsmooth}]

We are trying to show that if the point $\lambda$ is smooth on $\Pi_f$, then $\lambda$ is also smooth at any deletion or contraction of $\Pi_f$. 

Since we assume that $\lambda$ is smooth on $\Pi_f$, we can apply Theorem \ref{thm:smoothone}, which says that the pipe dream for $\Pi_f$ on $U_\lambda$ is rigid. (This implies that any subset of the pipe dream (for example, any maximal rectangle) is also rigid, since no two pipes containing any cross and near-miss combo in the whole pipe dream implies the same thing for any subset; but we will not need this fact for this proof, and this can already by deduced from the Main Theorem since rectangles that reduce to NW/SE partitions are rigid). This also implies that if we do deletion or contraction by filling in entire rows or columns of crosses, as explained in \S \ref{section:delcontrpipedream}, we do not need to (since in fact we cannot) do any moves prior to the filling in - rather, in this smooth/rigid case, deletion (contraction) of column (row) $i$ simply involves taking the $ith$ column (row) and brute-force changing every elbow tile in this column (row) into a cross tile. In the second (combinatorial) proof of Proposition \ref{prop:delsmooth} above, we did a detailed argument to show that this ``brute-force filling in" did not create any new moves, but now we show that we can bypass that argument by applying the Main Theorem.

Recall the definition of ``reducing to NW/SE partitions": this means having a rectangular shape where the northwest and southeast were partitions, and additionally, any number of rows or columns could be filled in with crosses; everything else is elbows. The Main Theorem says that (in this smooth case) all maximal rectangles reduce to NW/SE partitions; now, taking any one of these maximal rectangles (which reduce to NW/SE partitions) and filling in a row or column with crosses is still precisely a shape that reduces to NW/SE partitions. Therefore, it is still smooth.  


\end{proof}

\chapter{Atomic Positroid Pairs}\label{ch:atomicpairs}

Next, we describe the singular pairs that are lowest in the order explained in \S \ref{section:theorder}: recall that we called these ``atomic positroid pairs." By definition, any deletion or contraction of an atomic positroid pair results in a smooth positroid pair (see Figure \ref{fig:atomicbound}). The order we defined on positroid pairs has the property that for any given pair, as we go down in the order, we can only delete or contract a finite number of times. Because of this, any singular pair lies over an atomic pair, i.e. these minimal elements generate the order ideal.

We will use the atomic positroid pairs to give a slightly modified proof of the Main Theorem. 

\begin{Thm}\label{thm:atomicpd}

The atomic positroid pairs have affine pipe dreams that look like the following:

 \begin{figure}[htbp] \centering
	\includegraphics[scale=0.6,clip=true]{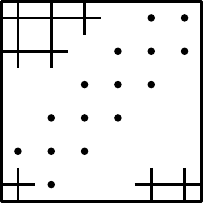}
	\caption{}
\end{figure}

More precisely, the affine pipe dreams looks like the following, with each square labeled ``A" being the same, and taking the form of the previous figure:

 \begin{figure}[htbp] \centering
	\includegraphics[scale=0.6,clip=true]{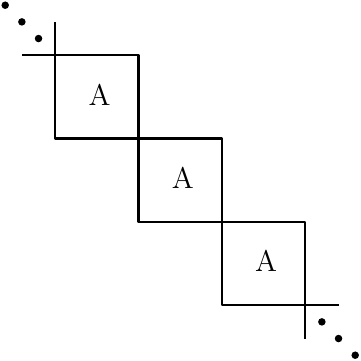}
	\caption{}
\end{figure}

Specifically:

(1) The atomic positroid pairs have pipe dreams in a square shape: so they are in $Gr(k,2k)$ and on points where $\lambda$ consists of $k$ consecutive columns. 

(2) In this pipe dream, there is a single cross along the southwest to northeast (longest) diagonal, with the diagonals directly adjacent to this longest diagonal being free of any crosses.

(3) Besides this one cross, there are partitions of crosses in the northwest and southeast corners, and no others.

\end{Thm}

\begin{proof}

We begin by proving that for an atomic positroid pair, in each of its pipe dreams there must be a single cross that can move. Suppose for contradiction that there exist two crosses that we call cross A and cross B, both with moves. Recall that by the definition of atomicity, this implies that any deletion or contraction must eliminate both crosses. We show that in general we can always eliminate one of the crosses while preserving the other one through a deletion or contraction - unless the pipe dream stretches infinitely in both the northeast and southwest directions, which gives us the contradiction. 
\\

In detail: suppose that cross A is at $(i_1,j_1)$ and cross B at $(i_2,j_2)$. The argument below will use the fact that, if a deletion
eliminates both crosses, then they have to be in the same column (the case of contraction: they have to be in the same row). Without loss of generality, at least the first or second coordinate must differ, or else cross A and cross B are the same cross, so let's take $j_1=j_2$: they are along the same column, and let's take $i_1$ to be above $i_2$ in the section of the pipe dream we are looking at. Now, clearly, deleting along this column eliminates both crosses, but since they are on different rows with cross A above cross B, we could contract cross A, which would leave cross B unchanged - that is, unless cross B has a move northeast to a new position $(i_1,j_3)$ along the same row as A but to its east. Now, we have solved the issue of contraction along $i_1$ eliminating both crosses, but the problem is that B is in a new easternmost column $j_3$ so we can delete $j_3$ to only eliminate B - that is, unless  A has a move northeast to a new position $(i_3,j_3)$ that is in the same column as B but north of B. But now we run into the issue that cross A has a new possible contraction at $i_3$. Thus, continuing on, we would have to do this infinitely in the northeast (as well as analogously in the southwest ) direction in order to preserve atomicity, which is impossible since our pipe dream shapes are bounded. 
\\

Now that we have proven that there is a single cross with a move, let us consider the possible positions of this cross: it must exist along the main diagonal of a square-shaped pipe dream; this is because, if this cross's moves missed any row or column, then we could delete this skipped row or column while leaving the cross (and its moves) intact, contradicting atomicity. Furthermore, this means that we have a clear diagonal along which this cross moves; in other words, the diagonals above and below this diagonal are free of any crosses, otherwise such a cross would hinder the free movement of this cross along the main diagonal (this relies on the fact just proven that no other cross can move): this proves (2). Since the pipe dream shape includes every row and column along a southwest to northeast diagonal, and no other rows/columns, and because in the shapes we are considering the non-convex portions only occur along the southwest and northeast corners, the boundaries of the pipe dreams shape must extend straight from the southeast and northeast corners: thus, the shape must be a rectangle. But since we already showed that every row and column must intersect the main diagonal, the shape in fact must be square. This proves (1). Since we already proved that this is the only cross that can move, everything else must be rigid, and this proves (3). 
\end{proof}

\begin{Prop}\label{prop:atom=sing}
A point $\lambda$ being singular on a positroid variety $\Pi_f$  is equivalent to the pipe dream for $\Pi_f$ on $U_\lambda$ being able to reach an atomic configuration via a series of deletions and contractions (that is, d/c-reduces to an atomic configuration). 
\end{Prop}

\begin{proof}

If $(\Pi_f,\lambda)$ d/c reduces to an atomic configuration, since atomic positroid pairs are singular by definition (minimally singular), and since we already proved in Prop \ref{prop:delsmooth} that the order on deletion and contraction is such that, for a singular element, everything higher up in the order is also singular, $(\Pi_f,\lambda)$ must be singular. 

In the other direction, $\Pi_f$ is finite dimensional, so a finite number of deletions and contractions will take it to the point $\lambda$, which is smooth in itself. Thus, every path going down in the order from a pair $(\Pi_f, \lambda)$ becomes smooth. Suppose $\Pi_f\in Gr(k,n)$. Then there are $k$ contractions (which can be applied to the columns in $\lambda$) and $n-k$ deletions (for the other columns, which cannot be contracted since $U_\lambda$ contains only the $T$-fixed point $\lambda$). Furthermore, suppose $\Pi_f$ is singular. Then these $n$ deletions and contractions are either (1) all smooth, in which case $(\Pi_f, \lambda)$ is atomic, or (2) at least one is singular, which let's call $(\Pi_{f'},\lambda')$. We can apply the same argument to $(\Pi_{f'},\lambda')$, which will have $n-1$ possible deletions and contractions - either it is atomic, or it has a deletion or contraction that remains singular. In this way, we can view each point in the ordering as spawning a tree of children below. By Prop \ref{prop:delsmooth}, every singular point only has singular parents, and so singular points are connected. Because all children become smooth eventually, these connected subtrees of singular points eventually end, with the last point being an atomic pair. See the next figure for an example.

\end{proof}

\begin{Exa}
Here, we provide an example of a part of the ordering. In particular, we show the positroid variety given by siteswap $f=342333$ containing point $\lambda=\{1,2,5\}$. This point is singular, which we have denoted by putting this positroid pair inside a blue rectangle. 

\begin{figure}[htbp] \centering
	\includegraphics[scale=0.65,clip=true]{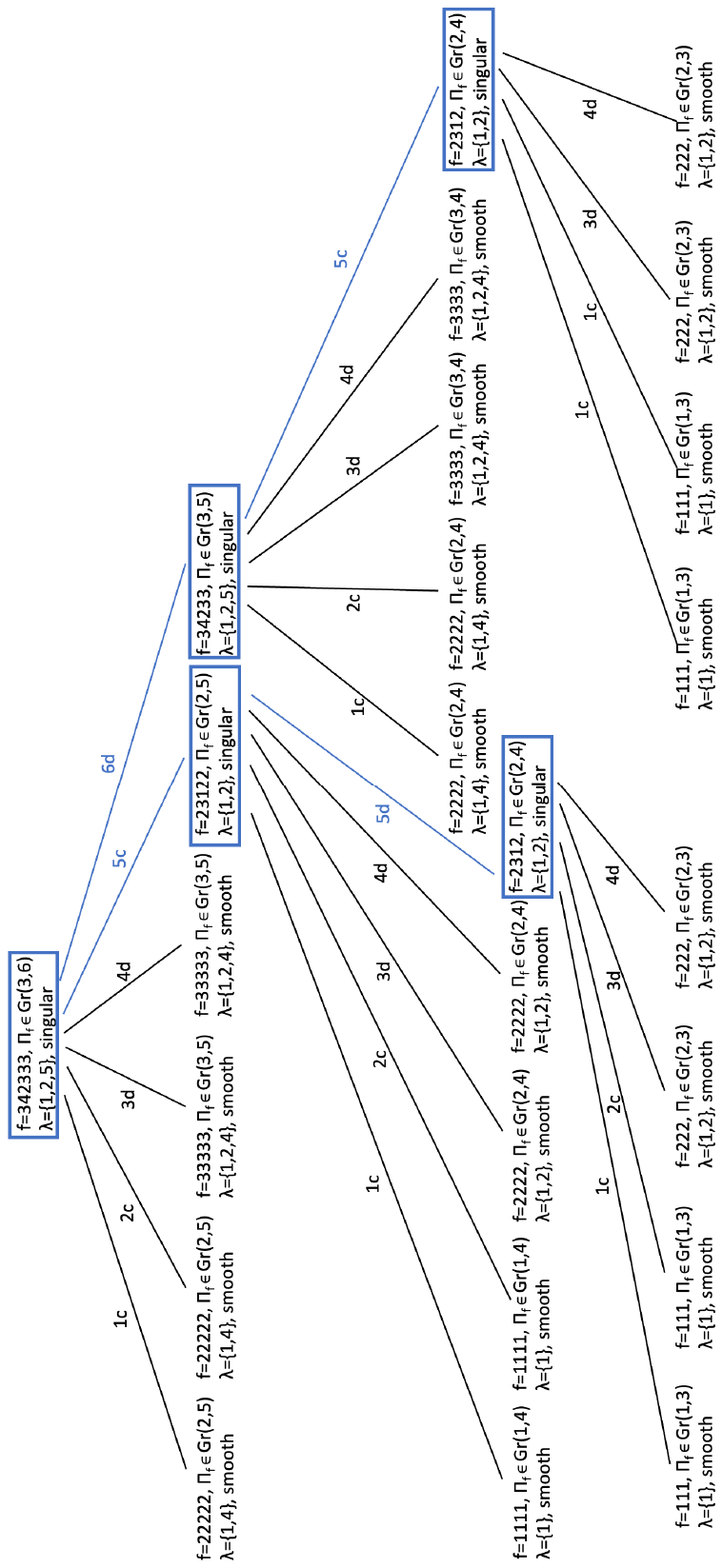}
\caption{}\label{fig:atomicorderexa}
\end{figure}

\newpage

We consider the set of all its deletions and contractions, denoted by the lines, and labeled for example, by ``2c" meaning that we contract the 2nd column or ``5d" meaning that we delete the 5th column. We stop doing further deletions and contractions once we reach a smooth point, since by Proposition \ref{prop:delsmooth}, anything lower in the order is still smooth. We see that the pair $(342333,\{1,2,5\})$ has exactly 2 deletions or contractions that are still singular, which are contraction of the 5th column or deletion of the 6th column. We can continue to delete or contract and obtain something singular until we get to $f=2312, \lambda=\{1,2\}$. Any further deletions or contractions result in a smooth point. Thus, by definition, this is an atomic pair (in fact, the only atomic pair that d\c reduces from $(342333,\{1,2,5\})$). And indeed, if we draw out the pipe dream, we see that $(2312, \lambda=\{1,2\})$ has a pipe dream that matches the form depicted in Theorem \ref{thm:atomicpd}:

 \begin{figure}[htbp]
	\includegraphics[scale=1,clip=true]{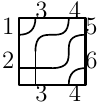}
\end{figure} 

\end{Exa}

We prove one more proposition in this section before making some final remarks.

\begin{Prop}\label{prop:atomrect}
Every atomic positroid pair comes from within a maximal rectangle. More precisely, if we have a pipe dream $\delta$ and there exists a series of deletions and contractions that take it (d/c reduces) to an atomic configuration that we call $\delta_{atomic}$, the process of deletion and contraction will have removed every tile outside of some maximal rectangle $R$ contained in $\delta$ (and most likely, some tiles within $R$ as well). 
\end{Prop}

\begin{proof}
This can be seen via an argument based on the following figure:

 \begin{figure}[htbp] \centering
	\includegraphics[scale=0.6,clip=true]{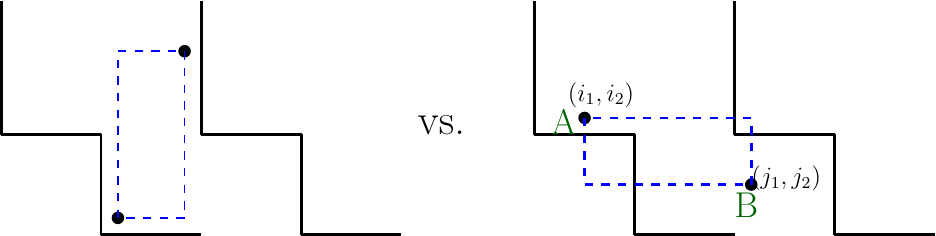}
\end{figure}

On the left side, we see that any pair of tiles that are southwest/northeast of each other must be contained within a maximal rectangle: this follows purely from the shapes we are considering. However, if they are southeast/northwest of each other, this might not be the case. In this case, if northwest tile A has coordinates $(i_1,i_2)$ and southeast tile B has coordinates $(j_1,j_2)$, then they are not contained in a maximal rectangle if either $(i_1,j_2)$ or $(i_2,j_1)$ is not contained in the shape. In such a case, there can be no series of deletions or contractions that put A and B into the same atomic pipe dream. The reason is that we cannot delete row $i_1$ or contract column $i_2$ without eliminating tile A, and we cannot delete row $j_1$ or contract column $j_2$ without eliminating tile B. Therefore, we cannot eliminate the spaces $(i_1,j_2)$ or $(i_2,j_1)$: they will always exist, so there will always be a convex hole, meaning that tiles A and B will never be able to fit into a square shape. This completes the proof of the fact that every atomic positroid pair comes from within a maximal rectangle.
\end{proof}

Let $\delta \in PD_{(\Pi_f,\lambda)}$ be an affine pipe dream for the positroid variety $\Pi_f$ on $U_\lambda$. We recap some of the statements that have been proven in this thesis.

\begin{equation*}
\begin{split} 
&\text{$\delta$ has a move within some maximal rectangle} \\
&\iff \text{$\delta$ has a move (is not rigid)} \\
&\iff \text{$\Pi_f$ is singular} \\
&\iff \text{$\delta$ d/c-reduces to an atomic pair} \\
&\iff \text{$\delta$ d/c-reduces to an atomic pair in some maximal rectangle}
\end{split} 
\end{equation*}

The first equivalence comes from Proposition \ref{prop:insidemaxrect} (one direction is trivial) or alternatively Corollary \ref{cor:box}, the second equivalence comes from Theorem \ref{thm:smoothone} (as well as Prop \ref{prop:BAPmoves}), the third equivalence is the content of Proposition \ref{prop:atom=sing}, and the last equivalence was just shown in Proposition \ref{prop:atomrect}. It might be hoped that, due to the apparent similarity of the first and last statements (both involving a condition within a maximal rectangle), that an equivalence between these two statements could be proved directly, which would produce a different proof of Proposition \ref{prop:insidemaxrect} using atomicity, and thereby yield an alternate proof of the Main Theorem \ref{thm:mainthm}. Unfortunately, the proof is more difficult than it seems. In one direction, going from $\delta$ having a move within some maximal rectangle to showing that $\delta$ d/c reduces to an atomic pair in that same maximal rectangle can be done by arguments involving deleting and contracting all rows and columns except for the positions where two pipes have a cross and near-misses (the argument could be made a bit easier by making use of Corollary \ref{cor:box}). The other direction is more difficult because the process of deletion or contraction occurs by doing all possible moves (in order to make sure that all pipe dreams for a given bounded permutation with the maximum number of crosses along the column/row to be deleted/contracted are obtained). This means that immediately upon invoking the notion of d/c-reduction, we are no longer just considering the single affine pipe dream $\delta$ but rather a possibly very large set of affine pipe dreams. In particular, we can no longer focus in on a single maximal rectangle (e.g. the maximal rectangle inside which we hope the atomic pair arises or the maximal rectangle where we hope to prove there exists a move) because we have to consider a set that includes pipe dreams where a cross-elbow move has shifted to being outside that single maximal rectangle we originally intended to focus on.

\chapter{A Computational Formula for Smoothness}\label{sec:compformula}

The following theorem provides a formula that tells us the multiplicity of a positroid variety at a torus-fixed point.
\begin{Thm}
Let $\Pi_f$ be a positroid variety with siteswap $f$ and $\lambda$ a T-fixed point. Then:
\[
[\Pi_f]|_\lambda = 
\sum_{u\mapsto \lambda} \frac{[X^w_v]|_u}{\Pi_{i<j\leq k} (y_{u(i)}-y_{u(j)}) \Pi_{k<i< j} (y_{u(i)}-y_{u(j)})} 
\]
\end{Thm}

The simple version of the proof is as follows: 

\begin{align*}
[\Pi_f]|_\lambda &= \pi_*([X_v^w])|_\lambda \\
&=\pi_*(\sum_u \frac{[X_v^w]|_u}{[u]|_u} [u])|_\lambda \\
&= \sum_u \frac{[X_v^w]|_u}{[u]|_u} (\pi_*(u))|_\lambda \\
&= \sum_{u\mapsto \lambda} \frac{[X_v^w]|_u}{[u]|_u} [\lambda]|_\lambda \\
&=\sum_{u\mapsto \lambda} \frac{[X^w_v]|_u}{(\Pi_{i<j\leq k} (y_{u(i)}-y_{u(j)}) \Pi_{k<i< j} (y_{u(i)}-y_{u(j)})) }
\end{align*}

$X_v^w$ denotes a Richardson variety, which is projected to the Grassmannian to obtain the positroid variety $\Pi_f$. Richardson varieties are intersection of Schubert varieties (denoted by subscripts e.g. $X_v$) and opposite Schubert varieties (denoted by superscripts e.g. $X^w$)

For $G$ a reductive group, $B$ a Borel subgroup and $P$ a parabolic subgroup, we have a projection map $\pi: G/B \to G/P$. If $G=GL_n(\C)$, we can take $B$ as the upper triangular matrices and $P$ as the matrices with nonzero values for everything north and east of some point along the diagonal; then $G/B$ is the flags in $\C^n$, and $G/P=Gr(k,n)$ for some $k\leq n$. 
\\

Then the intersection of the Schubert variety $x_v$ with the opposite Schubert variety $x^w$ is the Richardson variety denoted as $x_v^w:= x_v\cap x^w$. It is always possible to choose $w$ so that it is a Grassmannian permutation (the argument is in KLS 3.4 showing a bijection with k-Burhat intervals which are equivalence classes); this means that is has only one descent, so $w=a_1,...,a_k, a_{k+1},...,a_n$ where $a_1<\cdots < a_k > a_{k+1} < \cdots < a_n$. 
\\

We can apply the projection map $\pi$ to this Richardson variety to obtain a positroid variety corresponding to the siteswap $f$ that we denote $\Pi_f$: $\pi(X_v^w)=\Pi_f$. If we just consider this map on cohomology, we get: $ [\Pi_f]=\pi_* [X_v^w] $. Now, we restrict this to the point $\lambda$: $ [\Pi_f]|_\lambda=\pi_* [X_v^w]|_\lambda $. This lives in the torus-equivariant cohomology of a point: $H_T^*(\lambda)=\Z[y_1,...,y_n]$. 

For calculating the numerator $[x_v^w]|_u$: 

$[X^w_v]|_u = [X_v]|_u [X^w]|_u$ because they intersect transversely. This follows from Kleiman transversality since one is a $B$ orbit, and the other is a $B_-$ orbit. Restricting to a point is pullback along inclusion of a point, and pullback is a ring homomorphism so it takes products to products. 

$X^w=w_0\cdot x_{w_0w}$ (this is just left-multiplying by the matrix $w_0$ and is effectively taking Schubert varieties to opposite Schubert varieties, but one can also talk about it doing things to individual points, taking say $u$ to $w_0 u$)

$[X^w]|_u=w_0\cdot([X_{w_0 w}]|_{w_0 u}$) where $w_0 \cdot$performs: $y_i\mapsto y_{n+1-i}$. 

Putting it all together: $[X_v^w]|_u = [X_v]|_u * [X^w]|_u = [X_v]|_u * w_0\cdot([X_{w_0*w}]|_{w_0*u})$. We will do an example of this in a later section.

We expand on this proof in more detail next.

Denoting $\iota_\lambda: \{\lambda\} \to G/P$ an inclusion of a torus fixed point $\lambda \in G/P$ and $\gamma: G/B \to G/P$, 

\begin{align}\label{eq:compproof}
\begin{split}
[\Pi_f]|_\lambda &= \iota_\lambda^* [\Pi_f] \\
&= \iota_\lambda^*(\gamma_*[X_v^w]) \\
&= \iota_\lambda^* \gamma_* \bigg( \sum_u \frac{[X_v^w]|_u}{[uB/B]|_u} [uB/B] \bigg) \\
&= \iota_\lambda^*  \bigg( \sum_u \frac{[X_v^w]|_u}{[uB/B]|_u} [uP/P] \bigg) \\
&= \sum_u \frac{[X_v^w]|_u}{[uB/B]|_u} \iota_\lambda^*([uP/P])  \\
&= \sum_{u\in S_n, u\mapsto \lambda} \frac{[X_v^w]|_u}{[uB/B]|_u} [\lambda]|_\lambda \\
&= \sum_{u\in S_n, u\mapsto \lambda} [X_v^w]|_u \frac{[\lambda]|_\lambda}{[uB/B]|_u} 
\end{split}
\end{align}

where the last equality is clearly just rearranging the fraction, and for the fraction in the last term: \\
$\frac{[\lambda]|_\lambda}{[uB/B]|_u} = \frac{\text{product of T-weights on }T_\lambda G/P}{\text{product of T-weights on }T_{uB/B} G/B} = (\text{product of T-weights in } ker(T_{uB/B} \gamma))^{-1} = u\cdot \prod_{i<j\leq k} (y_i-y_j)^{-1} \prod_{k<i<j} (y_i-y_j)^{-1}$.

The first equality in Equation \ref{eq:compproof} is essentially by definition, and we have already explained the second equality above. The main step is the third equality, so we will spend most of the remainder of this section unpacking it. The other equalities are not difficult. 

Let's expand on what's going on by considering the maps of spaces:

\begin{tikzcd}
W \ar[r,hookrightarrow] \ar[d,twoheadrightarrow] & G/B \ar[d,twoheadrightarrow,"\gamma"] \\
W/W_P \ar[r,hookrightarrow] & G/P
\end{tikzcd}

where $W$ is the Weyl group $W = N_G(T)/T$ of the algebraic group $G$ with respect to a maximal torus $T$ contained in the Borel subgroup $B$.

We can take the T-equivariant cohomology of the diagram above by applying the contravariant functor $H_T^*$:

\begin{tikzcd}
H_T^*(W) & \ar[l] H_T^*(G/B)  \\
H_T^*(W/W_P) \ar[u] & H_T^*(G/P) \ar[l] \ar[u]
\end{tikzcd}

We will be working primarily with this diagram of equivariant cohomology for the equalities above, but first let us consider the first diagram in our specific case:

\begin{tikzcd}
u \ar[r,dotted,"\in"] \ar[d,mapsto] & W = S_n \ar[r,hookrightarrow] \ar[d,twoheadrightarrow] & G/B \ar[d,twoheadrightarrow]  \ar[r,hookleftarrow] & X_w^v \ar[d] \\
uW_P \ar[r,dotted,"\in"] & W/W_P \cong {[n]\choose k} \ar[r,hookrightarrow] & G/P \ar[r,hookleftarrow] & \Pi_f 
\end{tikzcd}

In our case, $G=GL_n$ and $B$ is the upper triangular matrices, and $P$ is the parabolic subgroup giving us $G/P = Gr(k,n)$ the Grassmannian. Matrix Schubert varieties give us a basis of the cohomology of $G/B$. The Weyl group is represented by permutation matrices $S_n$, so the fixed points under the torus action are labeled by permutations. Although to be more precise, we could write the points as:

\begin{tikzcd}
uB/B \ar[r,hookrightarrow] \ar[d] & G/B \ar[d,"\gamma"] \\
uP/P \ar[r] & G/P
\end{tikzcd}

with derivative map on tangent spaces: $T_{uB/B} \gamma: T_{uB/B} G/B \twoheadrightarrow T_\lambda G/P$

Let $u,v,w\in W=S_n$ be permutations. Writing $B_+$ for the upper triangular matrices and $B_-$ the lower triangular matrices, the matrix Schubert variety corresponding to $u$ is $X_u = \overline{B_+ u B_+/B_+}$, and the opposite Schubert variety (given equivalently by taking the flag in the opposite direction) would be $X^u = \overline{B_- u B_+/B_+}$. We obtain the Richardson variety $X^u_w$ as the intersection $X^w_v :=  X^w\cap X_v$ \cite{KLS2}. The Richardson variety will project down to a positroid variety in $G/P$ that we call $\Pi_f$. We can go in the other direction, where we choose $\Pi_f$ first, and then try to obtain $v$ and $w$ such that $X^u_w$ projects to it - obviously the choice is not unique (but if we choose $u$ to be a Grassmannian permutation, then it is) - we will into more detail about how to make this choice below (c.f. \cite{KLS}). The fixed point set $W$ of $G/B$ projects to the fixed point set of $G/P$, which we can denote $W_P$ ($W_P$ is obviously the kernel of the projection) $W_P=(1,2,3,...,k,*)$ where $*$ can be any permutation of $\{k+1,k+2,...,n\}$. That is, the projection map literally just takes the first $k$ elements of the permutation, and forgets the order (hence just a subset).

\begin{Def}

$\tilde{H}_T^*(X) := frac(H_T^*) \otimes_{H_T^*} H_T^*(X)$. Thus, $\tilde{H}_T^*(X)$ is the extension of $H_T^*(X)$ to the field of fraction, which allows us to not just multiply, but also divide by polynomials. 
\end{Def}

\begin{tikzcd}
\oplus \tilde{H}_T^*  & \tilde{H}_T^* (G/B) \ar[l,"\cong(Thm C)","\rho^*"'] & \\
H_T^*(W) = \oplus_{u\in W} H_T^* \ar[u,hookrightarrow] \ar[r, hookleftarrow,"Thm B"',"\rho^*"]  & H_T^*(G/B) \ar[d,"\gamma_*"] \ar[u, hookrightarrow,"Thm A"'] \ar[r,dotted,leftarrow,"\ni"] & \left[X_v^w\right] \ar[d,mapsto]  \\
H_T^*(W/W_p) = \oplus_{W/W_P} H_T^* \ar[u] & H_T^*(G/P) \ar[l] \ar[u,shift left,"\gamma^*"] \ar[r,dotted,leftarrow,"\ni"] & \left[ \Pi_f\right]
\end{tikzcd}

\begin{Thm}[A]
The statement here is that there is no torsion when we extend the coefficients to the fraction field; in other words, $H^*_T(G/B)$ is a free module over $H^*_T(pt)$. \end{Thm}
\begin{Proof}
It is a free module because there's a basis given by the Schubert classes (this is not the only basis, but at least it has a basis, whereas most modules do not have bases).
\end{Proof}

\begin{Thm}[B]
 $H_T^*(G/B) \hookrightarrow H_T^*(W) = \oplus_{u\in W} H_T^*$ in injective
This is a very powerful theorem. It says that you can look at the cohomology literally just by restricting to points. We call this map $\rho^*$.
\end{Thm}

\begin{Proof}Proposition 10 of \cite{KnutQuiver}
\end{Proof}

\begin{Thm}[C]
Now, when we extend to scalars, we actually get an isomorphism: \[ \oplus_{u\in W} \tilde{H}^*_T(pt)\cong \tilde{H}^*_T (G/B) \]
\end{Thm}
\begin{Proof}
This is proved in \cite{AtiBott84}
\end{Proof}

The next step is to take the Schubert classes of $H_T^*(G/B)$ and map via $\rho^*$, and see what that is in terms of the plain basis $\delta_u$.  

Cohomology of a (topological) disjoint union is a direct sum of the thing being unioned. In our case, they are just (disjoint) points. That's why $H_T^*(W/W_p) = \oplus_{W/W_P} H_T^*$.

For $H_T^*(W) = \bigoplus_{u\in W} H_T^*$, we have a basis consisting of $\delta_u$. If we denote by $e$ the map that extends the coefficients to the fraction field, then in the formula we have: $e(\delta_u)=\rho^*(\frac{e([uB/B])}{[uB/B]|_u})$. We can write $\delta_u=(0,0,1,0,0,...,0)$, with a 1 in the slot where $u$ is. These are orthogonal idempotents: they square to the identity and their product with others is equal to zero.  Now, $\rho^*([uB/B])=(0,0,[uB/B]|_u,0,0,...,0)$ is almost $\delta_u$, but we have to divide by $[uB/B]|_u$. $\rho^*$ is doing the restriction to points; for a class $[X_\pi]\in H^*_T(G/B)$:

\[\rho^*([X_\pi])=\sum_u [X_\pi]|_u \delta_u \xmapsto{e} \sum_u [X_\pi]|_u \frac{\rho^*([uB/B])}{[uB/B]|_u}= \rho^*( \sum_u \frac{[X_\pi]|_u}{[uB/B]|_u} [uB/B])\] 
(technically, every $u$ in this equation should be $uB/B$. And the last line follows since $\rho^*$ is $H_T^*$- and $\tilde{H}_T^*$-linear)

We should check that the $[uB/B]|_u$ is not equal to zero (otherwise we would be dividing by zero): $\rho^*([uB/B]) = (0,0,[uB/B]|_u\neq 0,0,0,...)$. These are not zero since they are the product of the weights on the tangent space to the point $u$, and because $u$ is an isolated fixed point, none of those weights are zero, so their product is not zero. Note that $[uB/B]$ is not a zerodivisor because it lives in a polynomial ring where there are no zerodivisors.

\subsection{Calculating Double Schubert Polynomials}\label{sec:doubSchub}

We calculate the double Schubert polynomials recursively using the formula below:

\[S_\rho = (x_i-y_j)S_{\rho'} + \sum_{\rho''\gtrdot \rho'} S_{\rho''} \]

where $\rho''=\rho'\circ (a i), a<i$, $i =$ position of last descent of $\rho$, $k=max\{l>i: \rho(l),\rho(i)\}$, $j=\rho(k)$, $\rho'=\rho\circ (ik)$. Note that both $\circ (ai)$ and $\circ (ik)$ are acting on the right so are acting on positions. 
\\

Example: Take $\rho=2143$. The last descent is in the 3rd position, and the furthest to the right we can go while staying less than $\rho(i)=\rho(3)=4$ is just the $j=3$ in the $k=4$-th slot. Thus, $\rho'=\rho\circ (ik)=(2143)\circ (34)=(2134)$, and $S_{2143}=(x_3-y_3)S_{2134}+\sum_{\rho''\gtrdot \rho'} S_{\rho''}$.
\\

Thus, we have to calculate the $\rho''$ pieces. The covering relation means that $\rho''$ is related to $\rho'$ by an inversion that is a weak Bruhat move. Since $\rho''=\rho'\circ (a i), a<i$, and we have $i=3$, we have to check $a=1,2$. In this case, both work: $(2134) \circ (13)=(3124)$ and $(2134) \circ (23)=(2314)$. 
\\

This completes the first step of the recursion: $S_{2143}=(x_3-y_3)S_{2134}+S_{3124}+S_{2314}$. For the second step, we first calculate $S_{2134}$. It has $i=1$, so clearly $\rho''=\emptyset$ since there can be no $a<i$. Since $k=2, j=1$, we get $S_{2143}=(x_1-y_1)S_{1234}=(x_1-y_1)$ since $S_{1234}=1$. For $S_{3124}$, again the last descent is at $a=1$ so $\rho''=\emptyset$, and here $k=3$, $j=2$, so its double Schubert polynomial is $(x_1-y_2)S_{2134}$. For $S_{2314}$, we have $i=2$, $k=3$, $j=1$, but after obtaining the new $\rho'=\rho\circ(ik)=(2314)\circ(23)=(2134)$, choosing $a=1$ would not give a covering relation since applying $\circ (ai)=\circ(12)$ would undo an already existing inversion, so here as well $\rho''=\emptyset$. Thus, $S_{2314}=(x_2-y_1)S_{2134}$.
\\

Putting it all together, we end up with:

\begin{align*}S_{2143}&=(x_3-y_3)S_{2134}+(x_1-y_2)S_{2134}+(x_2-y_1)S_{2134} \\
&=(x_1+x_2+x_3-y_1-y_2-y_3)S_{2134} \\
&=[x_1+x_2+x_3-(y_1+y_2+y_3)](x_1-y_1)
\end{align*}

\subsection{Obtaining $v$ and $w$ in $[X_v^w]$}\label{sec:getvw}

Now, we discuss how to obtain $v$ and $w$ in the formula for $\Pi_f=\pi(X_v^w)$, given a bounded affine permutation $f$, where $v$ and $w$ are elements of the symmetric group $S_n$. 

From section 3.4 of \cite{KLS}, $f_{u,w} = v t_{\omega_k} w^{-1}$. $[v,w]_k \in \mathcal{Q}(k, n)$ the equivalence classes of k-Bruhat intervals. In particular, the length $l(w)>l(v)$. The formula above gives us a way to go from $u$ and $w$ to $f_{u,w}$. We want to go in the reverse direction. The decomposition is not unique but it is unique if we choose $w$ to be a Grassmannian permutation (so it has a single descent).

$t$ is a translation element, and $\omega_k =(1,1...,1,0,0,...0)$ consist of $k$ one's followed by $n-k$ zero's. 

From equations (2) and (1) of Section 3.4 of \cite{KLS}: $a t_{\omega} a^{-1} = t_{a\cdot \omega}$ where $a\cdot (\omega_1,...,\omega_n) = (\omega_{a^{-1}(1)} ,...,\omega_{a^{-1}(n)})$

Denoting $f_m=f\mod n$, $f= t_\mu f_m =  (a^{-1} a t_\mu a^{-1}) a f(mod n) = a^{-1} t_{\omega_k} a f_m$.

Thus, in the formula $f_{v,w} = v t_{\omega_k} w^{-1}$, our candidate for $w^{-1}$ is $a f_m$, so let's call $w:= (a f_m)^{-1}$. However, our choice is not unique as we just noted. We can make it unique by letting the last part be a Grassmannian permutation. Let $w_g$ be the unique permutation such that the $[v,w]_k \sim [v,w_g] \in \mathcal{Q}(k, n)$ are the same class, but $w_g$ is Grassmannian. Then there exists a permutation $g$ such that $w=w_g \cdot g$, so $a\cdot f_m = w^{-1}= g^{-1} \cdot w_g^{-1}$, and thus $w_g^{-1} = g\cdot a \cdot f_m$. 

Using, this, we have $f = a^{-1} t_{\omega_k} a f_m =  a^{-1} t_{\omega_k} (g^{-1} g) a  f_m = a^{-1} t_{\omega_k} g^{-1} (g a  f_m) = a^{-1} t_{\omega_k} g^{-1} w_g^{-1}$. 

To finally get the expression for $f$ in the form we want, we need to let flip $g^{-1}$ to the left side of $t_{\omega_k}$. That is, we want to get the permutation $x$ such that $t_{\omega_k}  g^{-1} = x t_{\omega_k}$. Thus, $x= t_{\omega_k}  g^{-1} t_{\omega_k}^{-1}$. 

Putting all this together, $f= u t_{\omega_k} w_g^{-1}$ where $u= a^{-1} t_{\omega_k}  g^{-1} t_{\omega_k}^{-1}$ and $w_g^{-1} = g a f_m$. 

This is all done in the code, which we explain how to use in the next section.

\subsection{Example of Computational Formula}\label{sec:CompExa}

Here, we explicitly work out the formula:

\[
[\Pi_f]|_\lambda =
\sum_{u\in S_n, u\mapsto \lambda} \frac{[X_v]|_u  \,\, w_0\cdot([X_{w_0*w}]|_{w_0*u})}{(\Pi_{i<j\leq k} (y_{u(i)}-y_{u(j)}) \Pi_{k<i< j} (y_{u(i)}-y_{u(j)})) }
\]

Let's just do a very simple example: let $f=\{1,5,3\}$ and $\lambda=\{2\}$ in Gr(1,3).

In the formula, $u$ will consist of all permutations whose first $k=1$ element is $\lambda$. There are two such permutations $u$: $213$ and $231$.

We can use the formulas in the previous section to obtain $v$ and $w$ where $w$ is Grassmannian (easier just to use my code, described in more detail in the next section). Here, in this simple example $v=w=213$. We can do a quick check that $f=ut_\omega v$: 

In one-line notation, we start with $(213)$, and $t_\omega = t_{(100)}$ takes this to $(243)$. Then applying $(213)$ again, for example, sends the number that is 1 (mod 3) to $2$ in its window, i.e. $4\mapsto 5$, so the result is $(153)$ as we wanted. 

We have two terms, one corresponding to (a) $u=213$ and teh second to (b) $u=231$:

\[[X^w_v]|_u = [X_v]|_u w_0\cdot ([X_{w_0 w}]|_{w_0 u})\]

$\Z[y_1,...,y_n] \cong H_T^* \ni [X_v]|_u = S_v(x_1,...,x_n,y_1,...,y_n)|_{x_i\mapsto y_{u(i)}}$, where $S_v$ denotes the double Schubert polynomial associated with $v$. We explained how to calculate this in Section \ref{sec:doubSchub} above. 

For (a) $u=213$, $[X_v]=(x_1-y_1)|_{213}$ sends $x_i\mapsto y_{u(i)}$ so $x_1 \mapsto (y_{u(1)} = y_{2})$, so $[X_v]=y_2-y_1$.

For (b) $u=231$, $[X_v]=(x_1-y_1)|_{231}$, $u(1)$ is the same as (a), so this one also equals $y_2-y_1$. 

Now let's look at $w_0\cdot ([X_{w_0 w}]|_{w_0 u})$. 

$w_0 w = (321)\cdot (213)$: we have (213) and then (321) takes the 1 to 3, and 3 to 1, so the result is (231). Calculating the double Schubert polynomial $[X_{w_0 w}] = [X_{231}] = x_1x_2-x_1y_1-x_2 y_1+y_1^2$

For (a), $w_0 u = (321)\cdot (213)=(231)$, and for (b) $w_0 u = (321)\cdot (231)=(213)$. Thud, taking $[X_{w_0 w}] = x_1 x_2 - x_1 y_1 - x_2 y_1 + y_1^2$, $x_1\mapsto y_2$ for both (a) and (b), but $x_2\mapsto y_3$ for (a) versus $x_2\mapsto y_1$ for (b):

(a) $[X_{w_0 w}]|_{w_0u} = y_2 y_3 - y_2 y_1 - y_3 y_1 + y_1^2$ and (b) $[X_{w_0 w}]|_{w_0u} = y_2 y_1 - y_2 y_1 - y_1 y_1 + y_1^2=0$.

This means the term corresponding to (b) $u=(231)$ is zero.

Finally for the numerator in (a), $w_0 * ([X_{w_0 w}]|_{w_0u})$, the $w_0 *$ performs: $y_i\mapsto y_{n+1-i}$, so $y_1\mapsto y_3, y_2 \mapsto y_2, y_3\mapsto y_1$:
$w_0 * ([X_{w_0 w}]|_{w_0u}) = y_2 y_1 - y_2 y_3 - y_3 y_1 + y_3^2$, and

$[X_v]|_u \,\, w_0 * ([X_{w_0 w}]|_{w_0u}) = (y_2-y_1)(y_1 y_2 - y_1 y_3 - y_2 y_3 + y_3^2)$.

For the denominator $u\cdot \prod_{i<j\leq k} (y_i-y_j)^{-1} \prod_{k<i<j} (y_i-y_j)^{-1} = (\Pi_{i<j\leq k} (y_{u(i)}-y_{u(j)}) \Pi_{k<i< j} (y_{u(i)}-y_{u(j)}))$, 
 
since $k=1$, $\Pi_{i<j\leq k} (y_{u(i)}-y_{u(j)})$ has no terms since $i<j\leq 1$ is not possible, and $\Pi_{k<i< j} (y_{u(i)}-y_{u(j)})) = \Pi_{1<2< 3} (y_{u(2)}-y_{u(3)})) = y_1-y_3$.

Putting everythign together, and noting that $y_1 y_2 - y_1 y_3 - y_2 y_3 + y_3^2 = y_3(y_3-y_2)+y_1(y_2-y_3)=(y_3-y_1)(y_3-y_2)$:

\[\sum_{u\in S_n, u\mapsto \lambda} \frac{[X_v]|_u  \,\, w_0\cdot([X_{w_0*w}]|_{w_0*u})}{(\Pi_{i<j\leq k} (y_{u(i)}-y_{u(j)}) \Pi_{k<i< j} (y_{u(i)}-y_{u(j)})) }
= \frac{(y_2-y_1)(y_3-y_1)(y_3-y_2)}{y_1-y_3} = -(y_2-y_1)(y_3-y_2)
\]

In our code (see next section), we have written the final answer in terms of $z$'s, so it will show up as: $(z-2-z_3)(z_1-z_2)(-1)$.

It is a nice check to see that the result we get matches the other approach we described in Subsection \ref{subsec:AJSBilley}. There, we noted that we can calculate $[\Pi_f]|_\lambda$, we will get a sum of terms, one term for each pipe dream for $f$ in the shape determined by $\lambda$. If there are multiple terms (that is, pipe dreams), $\Pi_f$ is singular at $\lambda$, if there is a single pipe dreams, then it is smooth, and if there are no pipe dreams for $f$ on $\lambda$, then $\lambda \not\in \Pi_f$. In particular, Subsection \ref{subsec:AJSBilley} described how, each term in $[\Pi_f]|_\lambda$ consists of a product of $x_{row}-x_{col}$ for each cross in the pipe dream for that term. Let's verify that this matches what we've found in our example in this section.

Our shape is that given by $\lambda=\{2\}$. We form the distinguished path as before and form the shape on the right:

\begin{figure}[htbp] \centering
\includegraphics[scale=0.6,clip=true]{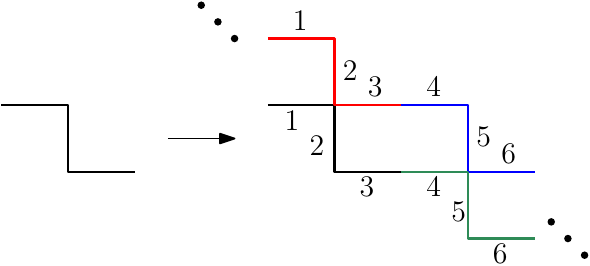}
\caption{}
\end{figure}

It is obvious in this example that there is a single pipe dream given by filling in the diagram entirely with crosses: the bounded affine permutation is $g=153$, so the siteswap is just $030$, consisting entirely of $0$'s and $n=3$'s in Gr(1,3):

\begin{figure}[htbp] \centering
\includegraphics[scale=0.6,clip=true]{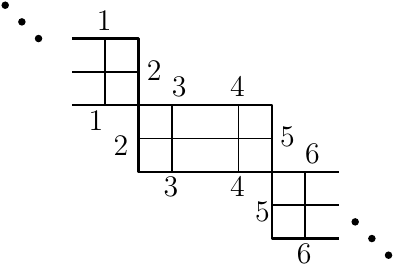}
\caption{}
\end{figure}

In particular, here the point $\lambda$ is smooth on $\Pi_f$, with a single term in the equivariant cohomology corresponding to the single pipe dream. 

If we just take one period:

\begin{figure}[htbp] \centering
\includegraphics[scale=0.6,clip=true]{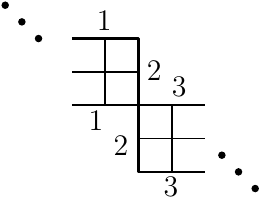}
\caption{}
\end{figure}

and consider the product $x_{row}-x_{col}$, or $y_{row}-y_{col}$ to be consistent with the notation in this section, this diagram corresponds to the product of $y_2-y_1$ for the upper cross and $y_2-y_3$ for the lower cross: $(y_2-y_1)(y_2-y_3)$ exactly matches our expression $-(y_2-y_1)(y_3-y_2)$ from our example in this section using our computational formula, as well as $(z-2-z_3)(z_1-z_2)(-1)$ from the code that we noted a few lines above.

\section{Using the Code for the Computational Formula}\label{sec:usingcode}

This formula is easy to run on a computer. Please see the website \url{josephflueg.github.io} for the code. 

The code is written for Macaulay2. Once you open on Macaulay2, say in Terminal on a Mac, load the scripts:

load ".../ClassFlambdaNew.m2"  \\
load ".../RunAllkn.m2"  \\
load ".../OrderBy.m2"  \\
load ".../SortOrder.m2"  \\
load ".../DoubleSchubert.m2"

where the ... stands for the location on your computer where these files are saved.

Then, to obtain the multiplicity of every pair $(f,\lambda)$ with $f$ a siteswap and $\lambda$ a point in $Gr(k,n)$, simply run:

RunAllkn(k,n)

This will generate 3 lists:

(1) smooth: a list of pairs $\{f,\lambda \}$ for all smooth torus-fixed points $\lambda$ on positroid varieties $\Pi_f$

(2) nonpoint: a list of pairs $\{f,\lambda \}$ for points $\lambda$ that are not contained in $f$, and

(3) singular: a list of 3-tuples $\{f,\lambda, m \}$ where m is the multiplicity of the point $\lambda$ on $\Pi_f$

If you would like to just calculate the class and multiplicity for a single choice of $f$ and $\lambda$, rather than all of them, after loading the scripts as above, run:

kk=2;  \\
nn=4;  \\
k=lift(kk,ZZ);  \\
n=lift(nn,ZZ);  \\
FinalFinalRing=QQ[h];  \\
$FinalRing=QQ[z_1..z_n];$  \\
$SchubRing=QQ[x_1..x_n,y_1..y_n];$  \\
DS = new MutableHashTable;

Note: change 2 and 4 in kk=2, nn=4 above to the actual k and n you want. Then, put in the bounded affine permutation $g=\{g_1,...,g_n\}$ (not siteswap) and point $lambda=\{l_1,...,l_k\}$ before running ClassFlambdaNew and run:

$g=\{g_1,...,g_n\};$ \\
$lambda=\{l_1,...,l_k\};$\\
ClassFlambdaNew(g,lambda)

The result will consist of pairs $(p(z),(deg X) h^{codim X})$ for $X$ the positroid variety, where $p(z)= \sum_{u\in S_n, u\mapsto \lambda} \frac{[X_v]|_u  \,\, w_0\cdot([X_{w_0*w}]|_{w_0*u})}{(\Pi_{i<j\leq k} (y_{u(i)}-y_{u(j)}) \Pi_{k<i< j} (y_{u(i)}-y_{u(j)}))} \in \Z[z_1,...,z_n]$, which is $H_T^*(pt)$. Finally, we restrict to the circle acting by dilation to get an element of $H_{\C^\times}^*(pt)=\Z[h]$ in the way we already described in Subsection \ref{subsec:restrictcircle}: we send each $z_i\mapsto h$ if $i\in \lambda$, otherwise $z_i \mapsto 0$. The important thing is the coefficient $deg X$ in front of $h^{codim X}$ in the 2nd term; $degX$ is the multiplicity - if it is $0$, then $\lambda\not \in \Pi_f$, if $degX=1$ then $\lambda$ is a smooth point of $\Pi_f$, and if $degX>1$, then $\lambda$ is a singular point of $\Pi_f$. 

Referring to the notation above in section \ref{sec:getvw}, in the code, we have called $a \mapsto A$, $g\mapsto grassperm$, $w_g \mapsto wgrass$, $f_m=switch$, $x\mapsto X$.

\section{Tables Cataloguing all Smooth and Singular Positroid Varieties up to $n=6$}\label{ch:catalogue}

In this section, we include tables showing all positroid varieties up to dimension $n=6$. We thank Nick Early for suggesting that we create these tables. 

The tables show the siteswap $f$ along the left, and the torus-fixed points $\lambda$ along the top. An entry of the table in the row $f$ and column $\lambda$ is equal to $0$ if $\lambda \not\in \Pi_f$ is not a point on the positroid variety. It is equal to $1$ if $\Pi_f$ is smooth at $\lambda$ (and contains $\lambda$). If $\Pi_f$ is singular at $\lambda$, then we include the multiplicity (the degree of the tangent cone) of $\Pi_f$ at $\lambda$, calculated using the formula we gave in the previous section. Finally, for easy visualization, a last column is included at the end that says whether the positroid variety is smooth. 

\begin{table}[h]
\centering
\begin{tabular}{|c|cccc|c|}
\hline
$Gr(1,4)$ & $\{1\}$ & $\{2\}$ & $\{3\}$ & $\{4\}$ & Smooth \\
\hline
$\{4,0,0,0\}$ & 1 & 0 & 0 & 0 & $\vee$ \\
\hline
$\{1,3,0,0\}$ & 1 & 1 & 0 & 0 & $\vee$ \\
\hline
$\{1,1,2,0\}$ & 1 & 1 & 1 & 0 & $\vee$ \\
\hline
$\{1,1,1,1\}$ & 1 & 1 & 1 & 1 & $\vee$ \\
\hline
$\{2,0,2,0\}$ & 1 & 0 & 1 & 0 & $\vee$ \\
\hline
\end{tabular}
\end{table}

The Grassmannian $Gr(1,4)$ is simply projective space $\mathbb{P}^3$, which is smooth, so all points are smooth. This same will be true of $Gr(1,n)\cong \mathbb{P}^{n-1}$ for any $n$, so we omit the tables $Gr(1,n)$ for $n=5,n=6$, since the tables just contain information about whether the positroid variety contains a point or not. 

\begin{table}[h]
\centering
\begin{tabular}{|c|cccc|c|}
\hline
$Gr(3,4)$ & $\{1,2,3\}$ & $\{1,2,4\}$ & $\{1,3,4\}$ & $\{2,3,4\}$ & Smooth \\
\hline
$\{4,4,4,0\}$ & 1 & 0 & 0 & 0 & $\vee$ \\
\hline
$\{4,4,1,3\}$ & 1 & 1 & 0 & 0 & $\vee$ \\
\hline
$\{4,2,3,3\}$ & 1 & 1 & 1 & 0 & $\vee$ \\
\hline
$\{4,2,4,2\}$ & 1 & 0 & 1 & 0 & $\vee$ \\
\hline
$\{3,3,3,3\}$ & 1 & 1 & 1 & 1 & $\vee$ \\
\hline
\end{tabular}
\end{table}

As expected, these points are also all smooth, since there is Grassmannian duality $Gr(k,n)\cong Gr(n-k,n)$. In particular, $Gr(1,n)\cong Gr(n-1,n)$.

The correspondence between $Gr(k,n)\cong Gr(n-k,n)$ is as follows. 

For points $\lambda_1$ on $Gr(k,n)$, the corresponding point $\lambda_2$ on $Gr(n-k,n)$ is given by $\lambda_2 = \{1,2,3,...,n-1,n\}\backslash \lambda_1$. For example

For the permutations, we have to take the inverse permutation. So for a siteswap $f$ with bounded affine permutation $g$ giving $\Pi_g\subseteq Gr(k,n)$, this corresponds to $g^{-1}$ for $\Pi_{g^{-1}}\subseteq Gr(n-k,n)$. We form $g^{-1}$ by looking, for each $i$, seeing which entry $h$ landed on $i$; for example if $g(h)=i$, then $g^{-1}(i)=h$. The number of steps backwards $h$ is from $i$ is $i-h$. Take the list of numbers $-(i-h)$ and add $n$ to all these numbers: this gives $f^{-1}$.

For example, let's take $\{1,1,2,0\}$. If we put two of these together (consider two periods), we have $\{1,1,2,0,1,1,2,0\}$. For small permutations, it might be easiest to draw the juggling pattern and see which line enters each spot. The number that lands on the first $1$ (in red) is the $2$ (in blue): $\{1,1,\blue{2}\black{,}0,\rd{1},1,2,0\}$; thus, we start with $-2$. Then, considering the second number, which is also $1=f(2)$, it's the previous $1$ (again in blue) which lands on it: $\{1,1,2,0,\blue{1},\rd{1}\black{,}2,0\}$; so now we have $(-2,-1)$, Continuing on, we get $(-2,-1,-1,0)$. Finally, we add $+n=+4$ to all these numbers, giving us $(2,3,3,4)$; notice that this is a siteswap in $Gr(3,4)$. 

Thus, we can apply the Grassmannian duality to $f=\{1,1,2,0\}$ which contained the points (and was smooth at the points) $\{1\},\{2\},\{3\}$ and did not contain the point $\{4\}$. Under the duality, $f=\{2,3,3,4\}$, which should contain the points $\{2,3,4\},\{1,3,4\},\{1,2,4\}$ and not contain the point $\{1,2,3\}$. This is exactly what we see in the second table for $Gr(3,4)$ (which shows $\{4,2,3,3\}$, so we need to rotate by $1$. As such, below we include just the table for $Gr(2,6)$, and not $Gr(4,6)$ to save space (noting that it is easy to calculate the dual using the process just described).

\begin{table}[h]
\centering
\begin{tabular}{|c|cccccc|c|}
\hline
$Gr(2,4)$ & $\{1,2\}$ & $\{1,3\}$ & $\{1,4\}$ & $\{2,3\}$ & $\{2,4\}$ & $\{3,4\}$ & Smooth \\
\hline
$\{4,4,0,0\}$ & 1 & 0 & 0 & 0 & 0 & 0 & $\vee$ \\
\hline
$\{4,1,3,0\}$ & 1 & 1 & 0 & 0 & 0 & 0 & $\vee$ \\
\hline
$\{4,1,1,2\}$ & 1 & 1 & 1 & 0 & 0 & 0 & $\vee$ \\
\hline
$\{4,2,0,2\}$ & 1 & 0 & 1 & 0 & 0 & 0 & $\vee$ \\
\hline
$\{2,3,3,0\}$ & 1 & 1 & 0 & 1 & 0 & 0 & $\vee$ \\
\hline
$\{2,3,1,2\}$ & 2 & 1 & 1 & 1 & 1 & 0 &  \\
\hline
$\{2,2,2,2\}$ & 1 & 1 & 1 & 1 & 1 & 1 & $\vee$ \\
\hline
$\{3,4,0,1\}$ & 1 & 0 & 0 & 0 & 1 & 0 & $\vee$ \\
\hline
$\{3,1,3,1\}$ & 1 & 1 & 0 & 0 & 1 & 1 & $\vee$ \\
\hline
$\{4,0,4,0\}$ & 0 & 1 & 0 & 0 & 0 & 0 & $\vee$ \\
\hline
\end{tabular}
\end{table}

$Gr(2,4)$ has 10 positroid varieties, of which just one - the one corresponding to $\{2,3,1,2\}$ - is singular. 

\newpage

Next, let's look at the tables corresponding to $n=5$.

Gr(2,5) has 7 singular positroid varieties, out of a total of 27 positroid varieties.

\begin{table}[h]
\hspace*{-2cm}
\small
\begin{tabular}{|c|cccccccccc|c|}
\hline
$Gr(2,5)$ & $\{1,2\}$ & $\{1,3\}$ & $\{1,4\}$ & $\{1,5\}$ & $\{2,3\}$ & $\{2,4\}$ & $\{2,5\}$ & $\{3,4\}$ & $\{3,5\}$ & $\{4,5\}$ & Smooth \\
\hline
$\{5,5,0,0,0\}$ & 1 & 0 & 0 & 0 & 0 & 0 & 0 & 0 & 0 & 0 & $\vee$ \\
\hline
$\{5,1,4,0,0\}$ & 1 & 1 & 0 & 0 & 0 & 0 & 0 & 0 & 0 & 0 & $\vee$ \\
\hline
$\{5,1,1,3,0\}$ & 1 & 1 & 1 & 0 & 0 & 0 & 0 & 0 & 0 & 0 & $\vee$ \\
\hline
$\{5,1,1,1,2\}$ & 1 & 1 & 1 & 1 & 0 & 0 & 0 & 0 & 0 & 0 & $\vee$ \\
\hline
$\{5,1,2,0,2\}$ & 1 & 1 & 0 & 1 & 0 & 0 & 0 & 0 & 0 & 0 & $\vee$ \\
\hline
$\{5,2,0,3,0\}$ & 1 & 0 & 1 & 0 & 0 & 0 & 0 & 0 & 0 & 0 & $\vee$ \\
\hline
$\{5,2,0,1,2\}$ & 1 & 0 & 1 & 1 & 0 & 0 & 0 & 0 & 0 & 0 & $\vee$ \\
\hline
$\{5,3,0,0,2\}$ & 1 & 0 & 0 & 1 & 0 & 0 & 0 & 0 & 0 & 0 & $\vee$ \\
\hline
$\{2,4,4,0,0\}$ & 1 & 1 & 0 & 0 & 1 & 0 & 0 & 0 & 0 & 0 & $\vee$ \\
\hline
$\{2,4,1,3,0\}$ & 2 & 1 & 1 & 0 & 1 & 1 & 0 & 0 & 0 & 0 &  \\
\hline
$\{2,4,1,1,2\}$ & 3 & 1 & 1 & 1 & 1 & 1 & 1 & 0 & 0 & 0 &  \\
\hline
$\{2,4,2,0,2\}$ & 2 & 1 & 0 & 1 & 1 & 0 & 1 & 0 & 0 & 0 &  \\
\hline
$\{2,2,3,3,0\}$ & 1 & 1 & 1 & 0 & 1 & 1 & 0 & 1 & 0 & 0 & $\vee$ \\
\hline
$\{2,2,3,1,2\}$ & 2 & 2 & 1 & 1 & 2 & 1 & 1 & 1 & 1 & 0 &  \\
\hline
$\{2,2,2,2,2\}$ & 1 & 1 & 1 & 1 & 1 & 1 & 1 & 1 & 1 & 1 & $\vee$ \\
\hline
$\{2,3,4,0,1\}$ & 1 & 1 & 0 & 0 & 2 & 0 & 1 & 0 & 1 & 0 &  \\
\hline
$\{2,3,1,3,1\}$ & 2 & 1 & 1 & 0 & 2 & 2 & 2 & 0 & 1 & 1 &  \\
\hline
$\{3,4,0,3,0\}$ & 1 & 0 & 1 & 0 & 0 & 1 & 0 & 0 & 0 & 0 & $\vee$ \\
\hline
$\{3,5,0,2,0\}$ & 1 & 0 & 0 & 0 & 0 & 1 & 0 & 0 & 0 & 0 & $\vee$ \\
\hline
$\{3,5,0,1,1\}$ & 1 & 0 & 0 & 0 & 0 & 1 & 1 & 0 & 0 & 0 & $\vee$ \\
\hline
$\{3,1,3,3,0\}$ & 1 & 1 & 2 & 0 & 0 & 1 & 0 & 1 & 0 & 0 &  \\
\hline
$\{3,1,4,2,0\}$ & 1 & 1 & 0 & 0 & 0 & 1 & 0 & 1 & 0 & 0 & $\vee$ \\
\hline
$\{3,1,4,1,1\}$ & 1 & 1 & 0 & 0 & 0 & 1 & 1 & 1 & 1 & 0 & $\vee$ \\
\hline
$\{4,5,0,0,1\}$ & 1 & 0 & 0 & 0 & 0 & 0 & 1 & 0 & 0 & 0 & $\vee$ \\
\hline
$\{4,1,4,0,1\}$ & 1 & 1 & 0 & 0 & 0 & 0 & 1 & 0 & 1 & 0 & $\vee$ \\
\hline
$\{5,0,5,0,0\}$ & 0 & 1 & 0 & 0 & 0 & 0 & 0 & 0 & 0 & 0 & $\vee$ \\
\hline
$\{5,0,1,4,0\}$ & 0 & 1 & 1 & 0 & 0 & 0 & 0 & 0 & 0 & 0 & $\vee$ \\
\hline
\end{tabular}
\end{table}

\newpage

\begin{table}[h]
\hspace*{-2.5cm}
\small
\begin{tabular}{|c|cccccccccc|c|}
\hline

$Gr(3,5)$ & \footnotesize$\{1,2,3\}$ & \footnotesize$\{1,2,4\}$ & \footnotesize$\{1,2,5\}$ & \footnotesize$\{1,3,4\}$ & \footnotesize$\{1,3,5\}$ & \footnotesize$\{1,4,5\}$ & \footnotesize$\{2,3,4\}$ & \footnotesize$\{2,3,5\}$ & \footnotesize$\{2,4,5\}$ & \footnotesize$\{3,4,5\}$ &\footnotesize Smooth \\
\hline
$\{5,5,5,0,0\}$ & 1 & 0 & 0 & 0 & 0 & 0 & 0 & 0 & 0 & 0 & $\vee$ \\
\hline
$\{5,5,1,4,0\}$ & 1 & 1 & 0 & 0 & 0 & 0 & 0 & 0 & 0 & 0 & $\vee$ \\
\hline
$\{5,5,1,1,3\}$ & 1 & 1 & 1 & 0 & 0 & 0 & 0 & 0 & 0 & 0 & $\vee$ \\
\hline
$\{5,5,2,0,3\}$ & 1 & 0 & 1 & 0 & 0 & 0 & 0 & 0 & 0 & 0 & $\vee$ \\
\hline
$\{5,2,4,4,0\}$ & 1 & 1 & 0 & 1 & 0 & 0 & 0 & 0 & 0 & 0 & $\vee$ \\
\hline
$\{5,2,4,1,3\}$ & 2 & 1 & 1 & 1 & 1 & 0 & 0 & 0 & 0 & 0 &  \\
\hline
$\{5,2,5,3,0\}$ & 1 & 0 & 0 & 1 & 0 & 0 & 0 & 0 & 0 & 0 & $\vee$ \\
\hline
$\{5,2,5,1,2\}$ & 1 & 0 & 0 & 1 & 1 & 0 & 0 & 0 & 0 & 0 & $\vee$ \\
\hline
$\{5,2,2,3,3\}$ & 1 & 1 & 1 & 1 & 1 & 1 & 0 & 0 & 0 & 0 & $\vee$ \\
\hline
$\{5,2,2,4,2\}$ & 1 & 1 & 0 & 2 & 1 & 1 & 0 & 0 & 0 & 0 &  \\
\hline
$\{5,3,4,0,3\}$ & 1 & 0 & 1 & 0 & 1 & 0 & 0 & 0 & 0 & 0 & $\vee$ \\
\hline
$\{5,3,5,0,2\}$ & 1 & 0 & 0 & 0 & 1 & 0 & 0 & 0 & 0 & 0 & $\vee$ \\
\hline
$\{5,3,1,3,3\}$ & 1 & 1 & 2 & 0 & 1 & 1 & 0 & 0 & 0 & 0 &  \\
\hline
$\{5,3,1,4,2\}$ & 1 & 1 & 0 & 0 & 1 & 1 & 0 & 0 & 0 & 0 & $\vee$ \\
\hline
$\{3,4,4,4,0\}$ & 1 & 1 & 0 & 1 & 0 & 0 & 1 & 0 & 0 & 0 & $\vee$ \\
\hline
$\{3,4,4,1,3\}$ & 3 & 1 & 1 & 1 & 1 & 0 & 1 & 1 & 0 & 0 &  \\
\hline
$\{3,4,5,3,0\}$ & 1 & 0 & 0 & 1 & 0 & 0 & 1 & 0 & 0 & 0 & $\vee$ \\
\hline
$\{3,4,5,1,2\}$ & 2 & 0 & 0 & 1 & 1 & 0 & 1 & 1 & 0 & 0 &  \\
\hline
$\{3,4,2,3,3\}$ & 2 & 2 & 2 & 1 & 1 & 1 & 1 & 1 & 1 & 0 &  \\
\hline
$\{3,4,2,4,2\}$ & 2 & 2 & 0 & 2 & 1 & 1 & 2 & 1 & 1 & 0 &  \\
\hline
$\{3,3,3,3,3\}$ & 1 & 1 & 1 & 1 & 1 & 1 & 1 & 1 & 1 & 1 & $\vee$ \\
\hline
$\{4,4,5,0,2\}$ & 1 & 0 & 0 & 0 & 1 & 0 & 0 & 1 & 0 & 0 & $\vee$ \\
\hline
$\{4,4,1,4,2\}$ & 1 & 1 & 0 & 0 & 1 & 1 & 0 & 1 & 1 & 0 & $\vee$ \\
\hline
$\{4,5,5,0,1\}$ & 1 & 0 & 0 & 0 & 0 & 0 & 0 & 1 & 0 & 0 & $\vee$ \\
\hline
$\{4,5,1,4,1\}$ & 1 & 1 & 0 & 0 & 0 & 0 & 0 & 1 & 1 & 0 & $\vee$ \\
\hline
$\{5,5,0,5,0\}$ & 0 & 1 & 0 & 0 & 0 & 0 & 0 & 0 & 0 & 0 & $\vee$ \\
\hline
$\{5,1,4,5,0\}$ & 0 & 1 & 0 & 1 & 0 & 0 & 0 & 0 & 0 & 0 & $\vee$ \\
\hline
\end{tabular}
\end{table}

Gr(3,5) is a bit redundant, given what we said above about Grassmann duality, but we included it here since the table is not large, and it allows us to check the duality for the multiplicity. For example, $\{2,4,1,1,2\}$ in $Gr(2,5)$ goes to $\{1,3,3,4,4\}$ in $Gr(3,5)$. $\{2,4,1,1,2\}$ had one singular point $\{1,2\}$ (which has multiplicity 3) so this should correspond to point $\{3,4,5\}$ for $\{1,3,3,4,4\}$; $\{2,4,1,1,2\}$ does not contain three points $\lambda = \{3,4\}, \{3,5\},\{4,5\}$, so these should correspond to point $\{1,2,5\}, \{1,2,4\}, \{1,2,3\}$ for $\{1,3,3,4,4\}$, and this is exactly what we see in the table for $Gr(3,5)$.

Next, we show the table for $Gr(2,6)$. This has 86 distinct positroid varieties, with 33 of them being singular, the rest (53) are smooth. The positroid variety for $f=\{3,1,3,1,3,1\}$ is special since it only contains torus-fixed points that are singular. 

We also notice that for each of these tables, the only positroid variety in $Gr(k,n)$ that contains every point and where every point is smooth is the top dimensional one, corresponding to $f=\{k,k,...,k,k\}$

{\newgeometry{top=1cm, bottom=1cm, left=1.5cm, right=1.5cm}
\setlength{\tabcolsep}{2pt}
\renewcommand{\arraystretch}{1.1}
\footnotesize
\begin{longtable}{|c|ccccccccccccccc|c|}
\hline
\text{ }$Gr(2,6)$ \text{  } & \scriptsize$\{1,2\}$ & \scriptsize$\{1,3\}$ & \scriptsize$\{1,4\}$ & \scriptsize$\{1,5\}$ & \scriptsize$\{1,6\}$ & \scriptsize$\{2,3\}$ & \scriptsize$\{2,4\}$ & \scriptsize$\{2,5\}$ & \scriptsize$\{2,6\}$ & \scriptsize$\{3,4\}$ & \scriptsize$\{3,5\}$ & \scriptsize$\{3,6\}$ & \scriptsize$\{4,5\}$ & \scriptsize$\{4,6\}$ & \scriptsize$\{5,6\}$ & Smooth \\
\hline
\endhead
\scriptsize$\{6,6,0,0,0,0\}$ & 1 & 0 & 0 & 0 & 0 & 0 & 0 & 0 & 0 & 0 & 0 & 0 & 0 & 0 & 0 & $\vee$ \\
\hline
\scriptsize$\{6,1,5,0,0,0\}$ & 1 & 1 & 0 & 0 & 0 & 0 & 0 & 0 & 0 & 0 & 0 & 0 & 0 & 0 & 0 & $\vee$ \\
\hline
\scriptsize$\{6,1,1,4,0,0\}$ & 1 & 1 & 1 & 0 & 0 & 0 & 0 & 0 & 0 & 0 & 0 & 0 & 0 & 0 & 0 & $\vee$ \\
\hline
\scriptsize$\{6,1,1,1,3,0\}$ & 1 & 1 & 1 & 1 & 0 & 0 & 0 & 0 & 0 & 0 & 0 & 0 & 0 & 0 & 0 & $\vee$ \\
\hline
\scriptsize$\{6,1,1,1,1,2\}$ & 1 & 1 & 1 & 1 & 1 & 0 & 0 & 0 & 0 & 0 & 0 & 0 & 0 & 0 & 0 & $\vee$ \\
\hline
\scriptsize$\{6,1,1,2,0,2\}$ & 1 & 1 & 1 & 0 & 1 & 0 & 0 & 0 & 0 & 0 & 0 & 0 & 0 & 0 & 0 & $\vee$ \\
\hline
\scriptsize$\{6,1,2,0,3,0\}$ & 1 & 1 & 0 & 1 & 0 & 0 & 0 & 0 & 0 & 0 & 0 & 0 & 0 & 0 & 0 & $\vee$ \\
\hline
\scriptsize$\{6,1,2,0,1,2\}$ & 1 & 1 & 0 & 1 & 1 & 0 & 0 & 0 & 0 & 0 & 0 & 0 & 0 & 0 & 0 & $\vee$ \\
\hline
\scriptsize$\{6,1,3,0,0,2\}$ & 1 & 1 & 0 & 0 & 1 & 0 & 0 & 0 & 0 & 0 & 0 & 0 & 0 & 0 & 0 & $\vee$ \\
\hline
\scriptsize$\{6,2,0,4,0,0\}$ & 1 & 0 & 1 & 0 & 0 & 0 & 0 & 0 & 0 & 0 & 0 & 0 & 0 & 0 & 0 & $\vee$ \\
\hline
\scriptsize$\{6,2,0,1,3,0\}$ & 1 & 0 & 1 & 1 & 0 & 0 & 0 & 0 & 0 & 0 & 0 & 0 & 0 & 0 & 0 & $\vee$ \\
\hline
\scriptsize$\{6,2,0,1,1,2\}$ & 1 & 0 & 1 & 1 & 1 & 0 & 0 & 0 & 0 & 0 & 0 & 0 & 0 & 0 & 0 & $\vee$ \\
\hline
\scriptsize$\{6,2,0,2,0,2\}$ & 1 & 0 & 1 & 0 & 1 & 0 & 0 & 0 & 0 & 0 & 0 & 0 & 0 & 0 & 0 & $\vee$ \\
\hline
\scriptsize$\{6,3,0,0,3,0\}$ & 1 & 0 & 0 & 1 & 0 & 0 & 0 & 0 & 0 & 0 & 0 & 0 & 0 & 0 & 0 & $\vee$ \\
\hline
\scriptsize$\{6,3,0,0,1,2\}$ & 1 & 0 & 0 & 1 & 1 & 0 & 0 & 0 & 0 & 0 & 0 & 0 & 0 & 0 & 0 & $\vee$ \\
\hline
\scriptsize$\{6,4,0,0,0,2\}$ & 1 & 0 & 0 & 0 & 1 & 0 & 0 & 0 & 0 & 0 & 0 & 0 & 0 & 0 & 0 & $\vee$ \\
\hline
\scriptsize$\{2,5,5,0,0,0\}$ & 1 & 1 & 0 & 0 & 0 & 1 & 0 & 0 & 0 & 0 & 0 & 0 & 0 & 0 & 0 & $\vee$ \\
\hline
\scriptsize$\{2,5,1,4,0,0\}$ & 2 & 1 & 1 & 0 & 0 & 1 & 1 & 0 & 0 & 0 & 0 & 0 & 0 & 0 & 0 &  \\
\hline
\scriptsize$\{2,5,1,1,3,0\}$ & 3 & 1 & 1 & 1 & 0 & 1 & 1 & 1 & 0 & 0 & 0 & 0 & 0 & 0 & 0 &  \\
\hline
\scriptsize$\{2,5,1,1,1,2\}$ & 4 & 1 & 1 & 1 & 1 & 1 & 1 & 1 & 1 & 0 & 0 & 0 & 0 & 0 & 0 &  \\
\hline
\scriptsize$\{2,5,1,2,0,2\}$ & 3 & 1 & 1 & 0 & 1 & 1 & 1 & 0 & 1 & 0 & 0 & 0 & 0 & 0 & 0 &  \\
\hline
\scriptsize$\{2,5,2,0,3,0\}$ & 2 & 1 & 0 & 1 & 0 & 1 & 0 & 1 & 0 & 0 & 0 & 0 & 0 & 0 & 0 &  \\
\hline
\scriptsize$\{2,5,2,0,1,2\}$ & 3 & 1 & 0 & 1 & 1 & 1 & 0 & 1 & 1 & 0 & 0 & 0 & 0 & 0 & 0 &  \\
\hline
\scriptsize$\{2,5,3,0,0,2\}$ & 2 & 1 & 0 & 0 & 1 & 1 & 0 & 0 & 1 & 0 & 0 & 0 & 0 & 0 & 0 &  \\
\hline
\scriptsize$\{2,2,4,4,0,0\}$ & 1 & 1 & 1 & 0 & 0 & 1 & 1 & 0 & 0 & 1 & 0 & 0 & 0 & 0 & 0 & $\vee$ \\
\hline
\scriptsize$\{2,2,4,1,3,0\}$ & 2 & 2 & 1 & 1 & 0 & 2 & 1 & 1 & 0 & 1 & 1 & 0 & 0 & 0 & 0 &  \\
\hline
\scriptsize$\{2,2,4,1,1,2\}$ & 3 & 3 & 1 & 1 & 1 & 3 & 1 & 1 & 1 & 1 & 1 & 1 & 0 & 0 & 0 &  \\
\hline
\scriptsize$\{2,2,4,2,0,2\}$ & 2 & 2 & 1 & 0 & 1 & 2 & 1 & 0 & 1 & 1 & 0 & 1 & 0 & 0 & 0 &  \\
\hline
\scriptsize$\{2,2,2,3,3,0\}$ & 1 & 1 & 1 & 1 & 0 & 1 & 1 & 1 & 0 & 1 & 1 & 0 & 1 & 0 & 0 & $\vee$ \\
\hline
\scriptsize$\{2,2,2,3,1,2\}$ & 2 & 2 & 2 & 1 & 1 & 2 & 2 & 1 & 1 & 2 & 1 & 1 & 1 & 1 & 0 &  \\
\hline
\scriptsize$\{2,2,2,2,2,2\}$ & 1 & 1 & 1 & 1 & 1 & 1 & 1 & 1 & 1 & 1 & 1 & 1 & 1 & 1 & 1 & $\vee$ \\
\hline
\scriptsize$\{2,2,3,4,0,1\}$ & 1 & 1 & 1 & 0 & 0 & 2 & 2 & 0 & 1 & 2 & 0 & 1 & 0 & 1 & 0 &  \\
\hline
\scriptsize$\{2,2,3,1,3,1\}$ & 2 & 2 & 1 & 1 & 0 & 4 & 2 & 2 & 2 & 2 & 2 & 2 & 0 & 1 & 1 &  \\
\hline
\scriptsize$\{2,3,4,0,3,0\}$ & 1 & 1 & 0 & 1 & 0 & 1 & 0 & 1 & 0 & 0 & 1 & 0 & 0 & 0 & 0 & $\vee$ \\
\hline
\scriptsize$\{2,3,5,0,2,0\}$ & 1 & 1 & 0 & 0 & 0 & 2 & 0 & 1 & 0 & 0 & 1 & 0 & 0 & 0 & 0 &  \\
\hline
\scriptsize$\{2,3,5,0,1,1\}$ & 1 & 1 & 0 & 0 & 0 & 3 & 0 & 1 & 1 & 0 & 1 & 1 & 0 & 0 & 0 &  \\
\hline
\scriptsize$\{2,3,1,3,3,0\}$ & 2 & 1 & 1 & 2 & 0 & 1 & 1 & 2 & 0 & 0 & 1 & 0 & 1 & 0 & 0 &  \\
\hline
\scriptsize$\{2,3,1,4,2,0\}$ & 2 & 1 & 1 & 0 & 0 & 2 & 2 & 2 & 0 & 0 & 1 & 0 & 1 & 0 & 0 &  \\
\hline
\scriptsize$\{2,3,1,4,1,1\}$ & 2 & 1 & 1 & 0 & 0 & 3 & 3 & 2 & 2 & 0 & 1 & 1 & 1 & 1 & 0 &  \\
\hline
\scriptsize$\{2,3,1,2,3,1\}$ & 2 & 1 & 1 & 2 & 0 & 2 & 2 & 4 & 2 & 0 & 2 & 1 & 2 & 1 & 2 &  \\
\hline
\scriptsize$\{2,3,3,0,3,1\}$ & 1 & 1 & 0 & 1 & 0 & 2 & 0 & 2 & 1 & 0 & 2 & 1 & 0 & 0 & 1 &  \\
\hline
\scriptsize$\{2,4,5,0,0,1\}$ & 1 & 1 & 0 & 0 & 0 & 2 & 0 & 0 & 1 & 0 & 0 & 1 & 0 & 0 & 0 &  \\
\hline
\scriptsize$\{2,4,1,4,0,1\}$ & 2 & 1 & 1 & 0 & 0 & 2 & 2 & 0 & 2 & 0 & 0 & 1 & 0 & 1 & 0 &  \\
\hline
\scriptsize$\{2,4,1,1,3,1\}$ & 3 & 1 & 1 & 1 & 0 & 2 & 2 & 2 & 3 & 0 & 0 & 1 & 0 & 1 & 1 &  \\
\hline
\scriptsize$\{2,4,2,0,3,1\}$ & 2 & 1 & 0 & 1 & 0 & 2 & 0 & 2 & 2 & 0 & 0 & 1 & 0 & 0 & 1 &  \\
\hline
\scriptsize$\{3,5,0,4,0,0\}$ & 1 & 0 & 1 & 0 & 0 & 0 & 1 & 0 & 0 & 0 & 0 & 0 & 0 & 0 & 0 & $\vee$ \\
\hline
\scriptsize$\{3,5,0,1,3,0\}$ & 2 & 0 & 1 & 1 & 0 & 0 & 1 & 1 & 0 & 0 & 0 & 0 & 0 & 0 & 0 &  \\
\hline
\scriptsize$\{3,6,0,3,0,0\}$ & 1 & 0 & 0 & 0 & 0 & 0 & 1 & 0 & 0 & 0 & 0 & 0 & 0 & 0 & 0 & $\vee$ \\
\hline
\scriptsize$\{3,6,0,1,2,0\}$ & 1 & 0 & 0 & 0 & 0 & 0 & 1 & 1 & 0 & 0 & 0 & 0 & 0 & 0 & 0 & $\vee$ \\
\hline
\scriptsize$\{3,6,0,1,1,1\}$ & 1 & 0 & 0 & 0 & 0 & 0 & 1 & 1 & 1 & 0 & 0 & 0 & 0 & 0 & 0 & $\vee$ \\
\hline
\scriptsize$\{3,6,0,2,0,1\}$ & 1 & 0 & 0 & 0 & 0 & 0 & 1 & 0 & 1 & 0 & 0 & 0 & 0 & 0 & 0 & $\vee$ \\
\hline
\scriptsize$\{3,1,4,4,0,0\}$ & 1 & 1 & 2 & 0 & 0 & 0 & 1 & 0 & 0 & 1 & 0 & 0 & 0 & 0 & 0 &  \\
\hline
\scriptsize$\{3,1,4,1,3,0\}$ & 2 & 2 & 2 & 2 & 0 & 0 & 1 & 1 & 0 & 1 & 1 & 0 & 0 & 0 & 0 &  \\
\hline
\scriptsize$\{3,1,5,3,0,0\}$ & 1 & 1 & 0 & 0 & 0 & 0 & 1 & 0 & 0 & 1 & 0 & 0 & 0 & 0 & 0 & $\vee$ \\
\hline
\scriptsize$\{3,1,5,1,2,0\}$ & 1 & 1 & 0 & 0 & 0 & 0 & 1 & 1 & 0 & 1 & 1 & 0 & 0 & 0 & 0 & $\vee$ \\
\hline
\scriptsize$\{3,1,5,1,1,1\}$ & 1 & 1 & 0 & 0 & 0 & 0 & 1 & 1 & 1 & 1 & 1 & 1 & 0 & 0 & 0 & $\vee$ \\
\hline
\scriptsize$\{3,1,5,2,0,1\}$ & 1 & 1 & 0 & 0 & 0 & 0 & 1 & 0 & 1 & 1 & 0 & 1 & 0 & 0 & 0 & $\vee$ \\
\hline
\scriptsize$\{3,1,3,4,0,1\}$ & 1 & 1 & 2 & 0 & 0 & 0 & 2 & 0 & 1 & 2 & 0 & 1 & 0 & 2 & 0 &  \\
\hline
\scriptsize$\{3,3,0,3,3,0\}$ & 1 & 0 & 1 & 1 & 0 & 0 & 1 & 1 & 0 & 0 & 0 & 0 & 1 & 0 & 0 & $\vee$ \\
\hline
\scriptsize$\{3,3,0,4,2,0\}$ & 1 & 0 & 1 & 0 & 0 & 0 & 2 & 1 & 0 & 0 & 0 & 0 & 1 & 0 & 0 &  \\
\hline
\scriptsize$\{3,3,0,4,1,1\}$ & 1 & 0 & 1 & 0 & 0 & 0 & 3 & 1 & 1 & 0 & 0 & 0 & 1 & 1 & 0 &  \\
\hline
\scriptsize$\{3,4,0,4,0,1\}$ & 1 & 0 & 1 & 0 & 0 & 0 & 2 & 0 & 1 & 0 & 0 & 0 & 0 & 1 & 0 &  \\
\hline
\scriptsize$\{4,5,0,0,3,0\}$ & 1 & 0 & 0 & 1 & 0 & 0 & 0 & 1 & 0 & 0 & 0 & 0 & 0 & 0 & 0 & $\vee$ \\
\hline
\scriptsize$\{4,6,0,0,2,0\}$ & 1 & 0 & 0 & 0 & 0 & 0 & 0 & 1 & 0 & 0 & 0 & 0 & 0 & 0 & 0 & $\vee$ \\
\hline
\scriptsize$\{4,6,0,0,1,1\}$ & 1 & 0 & 0 & 0 & 0 & 0 & 0 & 1 & 1 & 0 & 0 & 0 & 0 & 0 & 0 & $\vee$ \\
\hline
\scriptsize$\{4,1,4,0,3,0\}$ & 1 & 1 & 0 & 2 & 0 & 0 & 0 & 1 & 0 & 0 & 1 & 0 & 0 & 0 & 0 &  \\
\hline
\scriptsize$\{4,1,5,0,2,0\}$ & 1 & 1 & 0 & 0 & 0 & 0 & 0 & 1 & 0 & 0 & 1 & 0 & 0 & 0 & 0 & $\vee$ \\
\hline
\scriptsize$\{4,1,5,0,1,1\}$ & 1 & 1 & 0 & 0 & 0 & 0 & 0 & 1 & 1 & 0 & 1 & 1 & 0 & 0 & 0 & $\vee$ \\
\hline
\scriptsize$\{4,1,1,4,2,0\}$ & 1 & 1 & 1 & 0 & 0 & 0 & 0 & 1 & 0 & 0 & 1 & 0 & 1 & 0 & 0 & $\vee$ \\
\hline
\scriptsize$\{4,1,1,4,1,1\}$ & 1 & 1 & 1 & 0 & 0 & 0 & 0 & 1 & 1 & 0 & 1 & 1 & 1 & 1 & 0 & $\vee$ \\
\hline
\scriptsize$\{4,2,0,4,2,0\}$ & 1 & 0 & 1 & 0 & 0 & 0 & 0 & 1 & 0 & 0 & 0 & 0 & 1 & 0 & 0 & $\vee$ \\
\hline
\scriptsize$\{5,6,0,0,0,1\}$ & 1 & 0 & 0 & 0 & 0 & 0 & 0 & 0 & 1 & 0 & 0 & 0 & 0 & 0 & 0 & $\vee$ \\
\hline
\scriptsize$\{5,1,5,0,0,1\}$ & 1 & 1 & 0 & 0 & 0 & 0 & 0 & 0 & 1 & 0 & 0 & 1 & 0 & 0 & 0 & $\vee$ \\
\hline
\scriptsize$\{5,1,1,4,0,1\}$ & 1 & 1 & 1 & 0 & 0 & 0 & 0 & 0 & 1 & 0 & 0 & 1 & 0 & 1 & 0 & $\vee$ \\
\hline
\scriptsize$\{5,2,0,4,0,1\}$ & 1 & 0 & 1 & 0 & 0 & 0 & 0 & 0 & 1 & 0 & 0 & 0 & 0 & 1 & 0 & $\vee$ \\
\hline
\scriptsize$\{6,0,6,0,0,0\}$ & 0 & 1 & 0 & 0 & 0 & 0 & 0 & 0 & 0 & 0 & 0 & 0 & 0 & 0 & 0 & $\vee$ \\
\hline
\scriptsize$\{6,0,1,5,0,0\}$ & 0 & 1 & 1 & 0 & 0 & 0 & 0 & 0 & 0 & 0 & 0 & 0 & 0 & 0 & 0 & $\vee$ \\
\hline
\scriptsize$\{6,0,1,1,4,0\}$ & 0 & 1 & 1 & 1 & 0 & 0 & 0 & 0 & 0 & 0 & 0 & 0 & 0 & 0 & 0 & $\vee$ \\
\hline
\scriptsize$\{6,0,2,0,4,0\}$ & 0 & 1 & 0 & 1 & 0 & 0 & 0 & 0 & 0 & 0 & 0 & 0 & 0 & 0 & 0 & $\vee$ \\
\hline
\scriptsize$\{4,0,4,0,4,0\}$ & 0 & 1 & 0 & 1 & 0 & 0 & 0 & 0 & 0 & 0 & 1 & 0 & 0 & 0 & 0 & $\vee$ \\
\hline
\scriptsize$\{5,0,6,0,0,1\}$ & 0 & 1 & 0 & 0 & 0 & 0 & 0 & 0 & 0 & 0 & 0 & 1 & 0 & 0 & 0 & $\vee$ \\
\hline
\scriptsize$\{5,0,1,5,0,1\}$ & 0 & 1 & 1 & 0 & 0 & 0 & 0 & 0 & 0 & 0 & 0 & 1 & 0 & 1 & 0 & $\vee$ \\
\hline
\scriptsize$\{6,0,0,6,0,0\}$ & 0 & 0 & 1 & 0 & 0 & 0 & 0 & 0 & 0 & 0 & 0 & 0 & 0 & 0 & 0 & $\vee$ \\
\hline
\scriptsize$\{3,1,3,1,3,1\}$ & 2 & 2 & 2 & 2 & 0 & 0 & 2 & 2 & 2 & 2 & 2 & 2 & 0 & 2 & 2 &  \\
\hline
\end{longtable}
\restoregeometry}

Finally, Gr(3,6). This table is a bit large. There are 154 distinct siteswaps, of which 76 are singular (so 78 smooth).

{\newgeometry{top=1cm, bottom=1cm, left=1.5cm, right=1.5cm}
\begin{landscape}
\setlength{\tabcolsep}{2pt}
\renewcommand{\arraystretch}{1.1}
\footnotesize
\begin{longtable}{|c|cccccccccccccccccccc|c|}
\hline
\text{   } $Gr(3,6)$ \text{   } & \scriptsize$\{1,2,3\}$ & \scriptsize$\{1,2,4\}$ & \scriptsize$\{1,2,5\}$ & \scriptsize$\{1,2,6\}$ & \scriptsize$\{1,3,4\}$ & \scriptsize$\{1,3,5\}$ & \scriptsize$\{1,3,6\}$ & \scriptsize$\{1,4,5\}$ & \scriptsize$\{1,4,6\}$ & \scriptsize$\{1,5,6\}$ & \scriptsize$\{2,3,4\}$ & \scriptsize$\{2,3,5\}$ & \scriptsize$\{2,3,6\}$ & \scriptsize$\{2,4,5\}$ & \scriptsize$\{2,4,6\}$ & \scriptsize$\{2,5,6\}$ & \scriptsize$\{3,4,5\}$ & \scriptsize$\{3,4,6\}$ & \scriptsize$\{3,5,6\}$ & \scriptsize$\{4,5,6\}$ & \scriptsize Smooth \\
\hline
\endhead
\scriptsize$\{6,6,6,0,0,0\}$ & 1 & 0 & 0 & 0 & 0 & 0 & 0 & 0 & 0 & 0 & 0 & 0 & 0 & 0 & 0 & 0 & 0 & 0 & 0 & 0 & $\vee$ \\
\hline
\scriptsize$\{6,6,1,5,0,0\}$ & 1 & 1 & 0 & 0 & 0 & 0 & 0 & 0 & 0 & 0 & 0 & 0 & 0 & 0 & 0 & 0 & 0 & 0 & 0 & 0 & $\vee$ \\
\hline
\scriptsize$\{6,6,1,1,4,0\}$ & 1 & 1 & 1 & 0 & 0 & 0 & 0 & 0 & 0 & 0 & 0 & 0 & 0 & 0 & 0 & 0 & 0 & 0 & 0 & 0 & $\vee$ \\
\hline
\scriptsize$\{6,6,1,1,1,3\}$ & 1 & 1 & 1 & 1 & 0 & 0 & 0 & 0 & 0 & 0 & 0 & 0 & 0 & 0 & 0 & 0 & 0 & 0 & 0 & 0 & $\vee$ \\
\hline
\scriptsize$\{6,6,1,2,0,3\}$ & 1 & 1 & 0 & 1 & 0 & 0 & 0 & 0 & 0 & 0 & 0 & 0 & 0 & 0 & 0 & 0 & 0 & 0 & 0 & 0 & $\vee$ \\
\hline
\scriptsize$\{6,6,2,0,4,0\}$ & 1 & 0 & 1 & 0 & 0 & 0 & 0 & 0 & 0 & 0 & 0 & 0 & 0 & 0 & 0 & 0 & 0 & 0 & 0 & 0 & $\vee$ \\
\hline
\scriptsize$\{6,6,2,0,1,3\}$ & 1 & 0 & 1 & 1 & 0 & 0 & 0 & 0 & 0 & 0 & 0 & 0 & 0 & 0 & 0 & 0 & 0 & 0 & 0 & 0 & $\vee$ \\
\hline
\scriptsize$\{6,6,3,0,0,3\}$ & 1 & 0 & 0 & 1 & 0 & 0 & 0 & 0 & 0 & 0 & 0 & 0 & 0 & 0 & 0 & 0 & 0 & 0 & 0 & 0 & $\vee$ \\
\hline
\scriptsize$\{6,2,5,5,0,0\}$ & 1 & 1 & 0 & 0 & 1 & 0 & 0 & 0 & 0 & 0 & 0 & 0 & 0 & 0 & 0 & 0 & 0 & 0 & 0 & 0 & $\vee$ \\
\hline
\scriptsize$\{6,2,5,1,4,0\}$ & 2 & 1 & 1 & 0 & 1 & 1 & 0 & 0 & 0 & 0 & 0 & 0 & 0 & 0 & 0 & 0 & 0 & 0 & 0 & 0 &  \\
\hline
\scriptsize$\{6,2,5,1,1,3\}$ & 3 & 1 & 1 & 1 & 1 & 1 & 1 & 0 & 0 & 0 & 0 & 0 & 0 & 0 & 0 & 0 & 0 & 0 & 0 & 0 &  \\
\hline
\scriptsize$\{6,2,5,2,0,3\}$ & 2 & 1 & 0 & 1 & 1 & 0 & 1 & 0 & 0 & 0 & 0 & 0 & 0 & 0 & 0 & 0 & 0 & 0 & 0 & 0 &  \\
\hline
\scriptsize$\{6,2,6,4,0,0\}$ & 1 & 0 & 0 & 0 & 1 & 0 & 0 & 0 & 0 & 0 & 0 & 0 & 0 & 0 & 0 & 0 & 0 & 0 & 0 & 0 & $\vee$ \\
\hline
\scriptsize$\{6,2,6,1,3,0\}$ & 1 & 0 & 0 & 0 & 1 & 1 & 0 & 0 & 0 & 0 & 0 & 0 & 0 & 0 & 0 & 0 & 0 & 0 & 0 & 0 & $\vee$ \\
\hline
\scriptsize$\{6,2,6,1,1,2\}$ & 1 & 0 & 0 & 0 & 1 & 1 & 1 & 0 & 0 & 0 & 0 & 0 & 0 & 0 & 0 & 0 & 0 & 0 & 0 & 0 & $\vee$ \\
\hline
\scriptsize$\{6,2,6,2,0,2\}$ & 1 & 0 & 0 & 0 & 1 & 0 & 1 & 0 & 0 & 0 & 0 & 0 & 0 & 0 & 0 & 0 & 0 & 0 & 0 & 0 & $\vee$ \\
\hline
\scriptsize$\{6,2,2,4,4,0\}$ & 1 & 1 & 1 & 0 & 1 & 1 & 0 & 1 & 0 & 0 & 0 & 0 & 0 & 0 & 0 & 0 & 0 & 0 & 0 & 0 & $\vee$ \\
\hline
\scriptsize$\{6,2,2,4,1,3\}$ & 2 & 2 & 1 & 1 & 2 & 1 & 1 & 1 & 1 & 0 & 0 & 0 & 0 & 0 & 0 & 0 & 0 & 0 & 0 & 0 &  \\
\hline
\scriptsize$\{6,2,2,5,3,0\}$ & 1 & 1 & 0 & 0 & 2 & 1 & 0 & 1 & 0 & 0 & 0 & 0 & 0 & 0 & 0 & 0 & 0 & 0 & 0 & 0 &  \\
\hline
\scriptsize$\{6,2,2,5,1,2\}$ & 1 & 1 & 0 & 0 & 3 & 1 & 1 & 1 & 1 & 0 & 0 & 0 & 0 & 0 & 0 & 0 & 0 & 0 & 0 & 0 &  \\
\hline
\scriptsize$\{6,2,2,2,3,3\}$ & 1 & 1 & 1 & 1 & 1 & 1 & 1 & 1 & 1 & 1 & 0 & 0 & 0 & 0 & 0 & 0 & 0 & 0 & 0 & 0 & $\vee$ \\
\hline
\scriptsize$\{6,2,2,2,4,2\}$ & 1 & 1 & 1 & 0 & 2 & 2 & 1 & 2 & 1 & 1 & 0 & 0 & 0 & 0 & 0 & 0 & 0 & 0 & 0 & 0 &  \\
\hline
\scriptsize$\{6,2,3,4,0,3\}$ & 1 & 1 & 0 & 1 & 1 & 0 & 1 & 0 & 1 & 0 & 0 & 0 & 0 & 0 & 0 & 0 & 0 & 0 & 0 & 0 & $\vee$ \\
\hline
\scriptsize$\{6,2,3,5,0,2\}$ & 1 & 1 & 0 & 0 & 2 & 0 & 1 & 0 & 1 & 0 & 0 & 0 & 0 & 0 & 0 & 0 & 0 & 0 & 0 & 0 &  \\
\hline
\scriptsize$\{6,2,3,1,3,3\}$ & 2 & 1 & 1 & 2 & 1 & 1 & 2 & 0 & 1 & 1 & 0 & 0 & 0 & 0 & 0 & 0 & 0 & 0 & 0 & 0 &  \\
\hline
\scriptsize$\{6,2,3,1,4,2\}$ & 2 & 1 & 1 & 0 & 2 & 2 & 2 & 0 & 1 & 1 & 0 & 0 & 0 & 0 & 0 & 0 & 0 & 0 & 0 & 0 &  \\
\hline
\scriptsize$\{6,3,5,0,4,0\}$ & 1 & 0 & 1 & 0 & 0 & 1 & 0 & 0 & 0 & 0 & 0 & 0 & 0 & 0 & 0 & 0 & 0 & 0 & 0 & 0 & $\vee$ \\
\hline
\scriptsize$\{6,3,5,0,1,3\}$ & 2 & 0 & 1 & 1 & 0 & 1 & 1 & 0 & 0 & 0 & 0 & 0 & 0 & 0 & 0 & 0 & 0 & 0 & 0 & 0 &  \\
\hline
\scriptsize$\{6,3,6,0,3,0\}$ & 1 & 0 & 0 & 0 & 0 & 1 & 0 & 0 & 0 & 0 & 0 & 0 & 0 & 0 & 0 & 0 & 0 & 0 & 0 & 0 & $\vee$ \\
\hline
\scriptsize$\{6,3,6,0,1,2\}$ & 1 & 0 & 0 & 0 & 0 & 1 & 1 & 0 & 0 & 0 & 0 & 0 & 0 & 0 & 0 & 0 & 0 & 0 & 0 & 0 & $\vee$ \\
\hline
\scriptsize$\{6,3,1,4,4,0\}$ & 1 & 1 & 2 & 0 & 0 & 1 & 0 & 1 & 0 & 0 & 0 & 0 & 0 & 0 & 0 & 0 & 0 & 0 & 0 & 0 &  \\
\hline
\scriptsize$\{6,3,1,4,1,3\}$ & 2 & 2 & 2 & 2 & 0 & 1 & 1 & 1 & 1 & 0 & 0 & 0 & 0 & 0 & 0 & 0 & 0 & 0 & 0 & 0 &  \\
\hline
\scriptsize$\{6,3,1,5,3,0\}$ & 1 & 1 & 0 & 0 & 0 & 1 & 0 & 1 & 0 & 0 & 0 & 0 & 0 & 0 & 0 & 0 & 0 & 0 & 0 & 0 & $\vee$ \\
\hline
\scriptsize$\{6,3,1,5,1,2\}$ & 1 & 1 & 0 & 0 & 0 & 1 & 1 & 1 & 1 & 0 & 0 & 0 & 0 & 0 & 0 & 0 & 0 & 0 & 0 & 0 & $\vee$ \\
\hline
\scriptsize$\{6,3,1,2,3,3\}$ & 1 & 1 & 2 & 2 & 0 & 1 & 1 & 1 & 1 & 2 & 0 & 0 & 0 & 0 & 0 & 0 & 0 & 0 & 0 & 0 &  \\
\hline
\scriptsize$\{6,3,1,2,4,2\}$ & 1 & 1 & 2 & 0 & 0 & 2 & 1 & 2 & 1 & 2 & 0 & 0 & 0 & 0 & 0 & 0 & 0 & 0 & 0 & 0 &  \\
\hline
\scriptsize$\{6,3,3,0,3,3\}$ & 1 & 0 & 1 & 1 & 0 & 1 & 1 & 0 & 0 & 1 & 0 & 0 & 0 & 0 & 0 & 0 & 0 & 0 & 0 & 0 & $\vee$ \\
\hline
\scriptsize$\{6,3,3,0,4,2\}$ & 1 & 0 & 1 & 0 & 0 & 2 & 1 & 0 & 0 & 1 & 0 & 0 & 0 & 0 & 0 & 0 & 0 & 0 & 0 & 0 &  \\
\hline
\scriptsize$\{6,4,5,0,0,3\}$ & 1 & 0 & 0 & 1 & 0 & 0 & 1 & 0 & 0 & 0 & 0 & 0 & 0 & 0 & 0 & 0 & 0 & 0 & 0 & 0 & $\vee$ \\
\hline
\scriptsize$\{6,4,6,0,0,2\}$ & 1 & 0 & 0 & 0 & 0 & 0 & 1 & 0 & 0 & 0 & 0 & 0 & 0 & 0 & 0 & 0 & 0 & 0 & 0 & 0 & $\vee$ \\
\hline
\scriptsize$\{6,4,1,4,0,3\}$ & 1 & 1 & 0 & 2 & 0 & 0 & 1 & 0 & 1 & 0 & 0 & 0 & 0 & 0 & 0 & 0 & 0 & 0 & 0 & 0 &  \\
\hline
\scriptsize$\{6,4,1,5,0,2\}$ & 1 & 1 & 0 & 0 & 0 & 0 & 1 & 0 & 1 & 0 & 0 & 0 & 0 & 0 & 0 & 0 & 0 & 0 & 0 & 0 & $\vee$ \\
\hline
\scriptsize$\{6,4,1,1,3,3\}$ & 1 & 1 & 1 & 3 & 0 & 0 & 1 & 0 & 1 & 1 & 0 & 0 & 0 & 0 & 0 & 0 & 0 & 0 & 0 & 0 &  \\
\hline
\scriptsize$\{6,4,1,1,4,2\}$ & 1 & 1 & 1 & 0 & 0 & 0 & 1 & 0 & 1 & 1 & 0 & 0 & 0 & 0 & 0 & 0 & 0 & 0 & 0 & 0 & $\vee$ \\
\hline
\scriptsize$\{6,4,2,0,3,3\}$ & 1 & 0 & 1 & 2 & 0 & 0 & 1 & 0 & 0 & 1 & 0 & 0 & 0 & 0 & 0 & 0 & 0 & 0 & 0 & 0 &  \\
\hline
\scriptsize$\{6,4,2,0,4,2\}$ & 1 & 0 & 1 & 0 & 0 & 0 & 1 & 0 & 0 & 1 & 0 & 0 & 0 & 0 & 0 & 0 & 0 & 0 & 0 & 0 & $\vee$ \\
\hline
\scriptsize$\{3,5,5,5,0,0\}$ & 1 & 1 & 0 & 0 & 1 & 0 & 0 & 0 & 0 & 0 & 1 & 0 & 0 & 0 & 0 & 0 & 0 & 0 & 0 & 0 & $\vee$ \\
\hline
\scriptsize$\{3,5,5,1,4,0\}$ & 3 & 1 & 1 & 0 & 1 & 1 & 0 & 0 & 0 & 0 & 1 & 1 & 0 & 0 & 0 & 0 & 0 & 0 & 0 & 0 &  \\
\hline
\scriptsize$\{3,5,5,1,1,3\}$ & 6 & 1 & 1 & 1 & 1 & 1 & 1 & 0 & 0 & 0 & 1 & 1 & 1 & 0 & 0 & 0 & 0 & 0 & 0 & 0 &  \\
\hline
\scriptsize$\{3,5,5,2,0,3\}$ & 3 & 1 & 0 & 1 & 1 & 0 & 1 & 0 & 0 & 0 & 1 & 0 & 1 & 0 & 0 & 0 & 0 & 0 & 0 & 0 &  \\
\hline
\scriptsize$\{3,5,6,4,0,0\}$ & 1 & 0 & 0 & 0 & 1 & 0 & 0 & 0 & 0 & 0 & 1 & 0 & 0 & 0 & 0 & 0 & 0 & 0 & 0 & 0 & $\vee$ \\
\hline
\scriptsize$\{3,5,6,1,3,0\}$ & 2 & 0 & 0 & 0 & 1 & 1 & 0 & 0 & 0 & 0 & 1 & 1 & 0 & 0 & 0 & 0 & 0 & 0 & 0 & 0 &  \\
\hline
\scriptsize$\{3,5,6,1,1,2\}$ & 3 & 0 & 0 & 0 & 1 & 1 & 1 & 0 & 0 & 0 & 1 & 1 & 1 & 0 & 0 & 0 & 0 & 0 & 0 & 0 &  \\
\hline
\scriptsize$\{3,5,6,2,0,2\}$ & 2 & 0 & 0 & 0 & 1 & 0 & 1 & 0 & 0 & 0 & 1 & 0 & 1 & 0 & 0 & 0 & 0 & 0 & 0 & 0 &  \\
\hline
\scriptsize$\{3,5,2,4,4,0\}$ & 2 & 2 & 2 & 0 & 1 & 1 & 0 & 1 & 0 & 0 & 1 & 1 & 0 & 1 & 0 & 0 & 0 & 0 & 0 & 0 &  \\
\hline
\scriptsize$\{3,5,2,4,1,3\}$ & 5 & 5 & 2 & 2 & 2 & 1 & 1 & 1 & 1 & 0 & 2 & 1 & 1 & 1 & 1 & 0 & 0 & 0 & 0 & 0 &  \\
\hline
\scriptsize$\{3,5,2,5,3,0\}$ & 2 & 2 & 0 & 0 & 2 & 1 & 0 & 1 & 0 & 0 & 2 & 1 & 0 & 1 & 0 & 0 & 0 & 0 & 0 & 0 &  \\
\hline
\scriptsize$\{3,5,2,5,1,2\}$ & 3 & 3 & 0 & 0 & 3 & 1 & 1 & 1 & 1 & 0 & 3 & 1 & 1 & 1 & 1 & 0 & 0 & 0 & 0 & 0 &  \\
\hline
\scriptsize$\{3,5,2,2,3,3\}$ & 3 & 3 & 3 & 3 & 1 & 1 & 1 & 1 & 1 & 1 & 1 & 1 & 1 & 1 & 1 & 1 & 0 & 0 & 0 & 0 &  \\
\hline
\scriptsize$\{3,5,2,2,4,2\}$ & 3 & 3 & 3 & 0 & 2 & 2 & 1 & 2 & 1 & 1 & 2 & 2 & 1 & 2 & 1 & 1 & 0 & 0 & 0 & 0 &  \\
\hline
\scriptsize$\{3,5,3,4,0,3\}$ & 2 & 2 & 0 & 2 & 1 & 0 & 1 & 0 & 1 & 0 & 1 & 0 & 1 & 0 & 1 & 0 & 0 & 0 & 0 & 0 &  \\
\hline
\scriptsize$\{3,5,3,5,0,2\}$ & 2 & 2 & 0 & 0 & 2 & 0 & 1 & 0 & 1 & 0 & 2 & 0 & 1 & 0 & 1 & 0 & 0 & 0 & 0 & 0 &  \\
\hline
\scriptsize$\{3,5,3,1,3,3\}$ & 5 & 2 & 2 & 5 & 1 & 1 & 2 & 0 & 1 & 1 & 1 & 1 & 2 & 0 & 1 & 1 & 0 & 0 & 0 & 0 &  \\
\hline
\scriptsize$\{3,5,3,1,4,2\}$ & 5 & 2 & 2 & 0 & 2 & 2 & 2 & 0 & 1 & 1 & 2 & 2 & 2 & 0 & 1 & 1 & 0 & 0 & 0 & 0 &  \\
\hline
\scriptsize$\{3,3,4,4,4,0\}$ & 1 & 1 & 1 & 0 & 1 & 1 & 0 & 1 & 0 & 0 & 1 & 1 & 0 & 1 & 0 & 0 & 1 & 0 & 0 & 0 & $\vee$ \\
\hline
\scriptsize$\{3,3,4,4,1,3\}$ & 3 & 3 & 1 & 1 & 3 & 1 & 1 & 1 & 1 & 0 & 3 & 1 & 1 & 1 & 1 & 0 & 1 & 1 & 0 & 0 &  \\
\hline
\scriptsize$\{3,3,4,5,3,0\}$ & 1 & 1 & 0 & 0 & 2 & 1 & 0 & 1 & 0 & 0 & 2 & 1 & 0 & 1 & 0 & 0 & 2 & 0 & 0 & 0 &  \\
\hline
\scriptsize$\{3,3,4,5,1,2\}$ & 2 & 2 & 0 & 0 & 5 & 1 & 1 & 1 & 1 & 0 & 5 & 1 & 1 & 1 & 1 & 0 & 2 & 2 & 0 & 0 &  \\
\hline
\scriptsize$\{3,3,4,2,3,3\}$ & 3 & 2 & 2 & 2 & 2 & 2 & 2 & 1 & 1 & 1 & 2 & 2 & 2 & 1 & 1 & 1 & 1 & 1 & 1 & 0 &  \\
\hline
\scriptsize$\{3,3,4,2,4,2\}$ & 3 & 2 & 2 & 0 & 4 & 4 & 2 & 2 & 1 & 1 & 4 & 4 & 2 & 2 & 1 & 1 & 3 & 2 & 2 & 0 &  \\
\hline
\scriptsize$\{3,3,3,3,3,3\}$ & 1 & 1 & 1 & 1 & 1 & 1 & 1 & 1 & 1 & 1 & 1 & 1 & 1 & 1 & 1 & 1 & 1 & 1 & 1 & 1 & $\vee$ \\
\hline
\scriptsize$\{3,4,4,5,0,2\}$ & 1 & 1 & 0 & 0 & 2 & 0 & 1 & 0 & 1 & 0 & 2 & 0 & 1 & 0 & 1 & 0 & 0 & 2 & 0 & 0 &  \\
\hline
\scriptsize$\{3,4,4,1,4,2\}$ & 3 & 1 & 1 & 0 & 2 & 2 & 3 & 0 & 1 & 1 & 2 & 2 & 3 & 0 & 1 & 1 & 0 & 2 & 2 & 0 &  \\
\hline
\scriptsize$\{3,4,5,5,0,1\}$ & 1 & 1 & 0 & 0 & 1 & 0 & 0 & 0 & 0 & 0 & 3 & 0 & 1 & 0 & 1 & 0 & 0 & 1 & 0 & 0 &  \\
\hline
\scriptsize$\{3,4,5,1,4,1\}$ & 3 & 1 & 1 & 0 & 1 & 1 & 0 & 0 & 0 & 0 & 3 & 3 & 3 & 0 & 1 & 1 & 0 & 1 & 1 & 0 &  \\
\hline
\scriptsize$\{3,4,6,3,0,2\}$ & 1 & 0 & 0 & 0 & 1 & 0 & 1 & 0 & 0 & 0 & 1 & 0 & 1 & 0 & 0 & 0 & 0 & 1 & 0 & 0 & $\vee$ \\
\hline
\scriptsize$\{3,4,6,4,0,1\}$ & 1 & 0 & 0 & 0 & 1 & 0 & 0 & 0 & 0 & 0 & 2 & 0 & 1 & 0 & 0 & 0 & 0 & 1 & 0 & 0 &  \\
\hline
\scriptsize$\{3,4,6,1,2,2\}$ & 2 & 0 & 0 & 0 & 1 & 1 & 2 & 0 & 0 & 0 & 1 & 1 & 2 & 0 & 0 & 0 & 0 & 1 & 1 & 0 &  \\
\hline
\scriptsize$\{3,4,6,1,3,1\}$ & 2 & 0 & 0 & 0 & 1 & 1 & 0 & 0 & 0 & 0 & 2 & 2 & 2 & 0 & 0 & 0 & 0 & 1 & 1 & 0 &  \\
\hline
\scriptsize$\{3,4,2,4,4,1\}$ & 2 & 2 & 2 & 0 & 1 & 1 & 0 & 1 & 0 & 0 & 3 & 3 & 2 & 3 & 2 & 2 & 0 & 1 & 1 & 1 &  \\
\hline
\scriptsize$\{3,4,2,5,2,2\}$ & 2 & 2 & 0 & 0 & 3 & 1 & 2 & 1 & 2 & 0 & 3 & 1 & 2 & 1 & 2 & 0 & 0 & 3 & 1 & 1 &  \\
\hline
\scriptsize$\{3,4,2,5,3,1\}$ & 2 & 2 & 0 & 0 & 2 & 1 & 0 & 1 & 0 & 0 & 5 & 2 & 2 & 2 & 2 & 0 & 0 & 2 & 1 & 1 &  \\
\hline
\scriptsize$\{4,5,5,0,4,0\}$ & 1 & 0 & 1 & 0 & 0 & 1 & 0 & 0 & 0 & 0 & 0 & 1 & 0 & 0 & 0 & 0 & 0 & 0 & 0 & 0 & $\vee$ \\
\hline
\scriptsize$\{4,5,6,0,3,0\}$ & 1 & 0 & 0 & 0 & 0 & 1 & 0 & 0 & 0 & 0 & 0 & 1 & 0 & 0 & 0 & 0 & 0 & 0 & 0 & 0 & $\vee$ \\
\hline
\scriptsize$\{4,5,6,0,1,2\}$ & 2 & 0 & 0 & 0 & 0 & 1 & 1 & 0 & 0 & 0 & 0 & 1 & 1 & 0 & 0 & 0 & 0 & 0 & 0 & 0 &  \\
\hline
\scriptsize$\{4,5,1,4,4,0\}$ & 1 & 1 & 3 & 0 & 0 & 1 & 0 & 1 & 0 & 0 & 0 & 1 & 0 & 1 & 0 & 0 & 0 & 0 & 0 & 0 &  \\
\hline
\scriptsize$\{4,5,1,5,3,0\}$ & 1 & 1 & 0 & 0 & 0 & 1 & 0 & 1 & 0 & 0 & 0 & 1 & 0 & 1 & 0 & 0 & 0 & 0 & 0 & 0 & $\vee$ \\
\hline
\scriptsize$\{4,5,1,5,1,2\}$ & 2 & 2 & 0 & 0 & 0 & 1 & 1 & 1 & 1 & 0 & 0 & 1 & 1 & 1 & 1 & 0 & 0 & 0 & 0 & 0 &  \\
\hline
\scriptsize$\{4,5,1,2,4,2\}$ & 2 & 2 & 5 & 0 & 0 & 2 & 1 & 2 & 1 & 2 & 0 & 2 & 1 & 2 & 1 & 2 & 0 & 0 & 0 & 0 &  \\
\hline
\scriptsize$\{4,5,3,0,4,2\}$ & 2 & 0 & 2 & 0 & 0 & 2 & 1 & 0 & 0 & 1 & 0 & 2 & 1 & 0 & 0 & 1 & 0 & 0 & 0 & 0 &  \\
\hline
\scriptsize$\{4,6,4,0,4,0\}$ & 1 & 0 & 1 & 0 & 0 & 0 & 0 & 0 & 0 & 0 & 0 & 1 & 0 & 0 & 0 & 0 & 0 & 0 & 0 & 0 & $\vee$ \\
\hline
\scriptsize$\{4,6,6,0,2,0\}$ & 1 & 0 & 0 & 0 & 0 & 0 & 0 & 0 & 0 & 0 & 0 & 1 & 0 & 0 & 0 & 0 & 0 & 0 & 0 & 0 & $\vee$ \\
\hline
\scriptsize$\{4,6,6,0,1,1\}$ & 1 & 0 & 0 & 0 & 0 & 0 & 0 & 0 & 0 & 0 & 0 & 1 & 1 & 0 & 0 & 0 & 0 & 0 & 0 & 0 & $\vee$ \\
\hline
\scriptsize$\{4,6,1,3,4,0\}$ & 1 & 1 & 2 & 0 & 0 & 0 & 0 & 0 & 0 & 0 & 0 & 1 & 0 & 1 & 0 & 0 & 0 & 0 & 0 & 0 &  \\
\hline
\scriptsize$\{4,6,1,5,2,0\}$ & 1 & 1 & 0 & 0 & 0 & 0 & 0 & 0 & 0 & 0 & 0 & 1 & 0 & 1 & 0 & 0 & 0 & 0 & 0 & 0 & $\vee$ \\
\hline
\scriptsize$\{4,6,1,5,1,1\}$ & 1 & 1 & 0 & 0 & 0 & 0 & 0 & 0 & 0 & 0 & 0 & 1 & 1 & 1 & 1 & 0 & 0 & 0 & 0 & 0 & $\vee$ \\
\hline
\scriptsize$\{4,6,1,2,4,1\}$ & 1 & 1 & 2 & 0 & 0 & 0 & 0 & 0 & 0 & 0 & 0 & 2 & 1 & 2 & 1 & 2 & 0 & 0 & 0 & 0 &  \\
\hline
\scriptsize$\{4,6,3,0,4,1\}$ & 1 & 0 & 1 & 0 & 0 & 0 & 0 & 0 & 0 & 0 & 0 & 2 & 1 & 0 & 0 & 1 & 0 & 0 & 0 & 0 &  \\
\hline
\scriptsize$\{4,2,4,4,4,0\}$ & 1 & 1 & 2 & 0 & 1 & 2 & 0 & 2 & 0 & 0 & 0 & 1 & 0 & 1 & 0 & 0 & 1 & 0 & 0 & 0 &  \\
\hline
\scriptsize$\{4,2,4,2,4,2\}$ & 3 & 2 & 4 & 0 & 4 & 8 & 4 & 4 & 2 & 3 & 0 & 4 & 2 & 2 & 1 & 2 & 3 & 2 & 4 & 0 &  \\
\hline
\scriptsize$\{4,2,5,3,4,0\}$ & 2 & 1 & 2 & 0 & 1 & 2 & 0 & 0 & 0 & 0 & 0 & 2 & 0 & 1 & 0 & 0 & 1 & 0 & 0 & 0 &  \\
\hline
\scriptsize$\{4,2,5,5,2,0\}$ & 1 & 1 & 0 & 0 & 1 & 0 & 0 & 0 & 0 & 0 & 0 & 1 & 0 & 1 & 0 & 0 & 1 & 0 & 0 & 0 & $\vee$ \\
\hline
\scriptsize$\{4,2,5,5,1,1\}$ & 1 & 1 & 0 & 0 & 1 & 0 & 0 & 0 & 0 & 0 & 0 & 1 & 1 & 1 & 1 & 0 & 1 & 1 & 0 & 0 & $\vee$ \\
\hline
\scriptsize$\{4,2,5,2,4,1\}$ & 2 & 1 & 2 & 0 & 1 & 2 & 0 & 0 & 0 & 0 & 0 & 5 & 2 & 2 & 1 & 2 & 2 & 1 & 2 & 0 &  \\
\hline
\scriptsize$\{4,4,5,0,4,1\}$ & 1 & 0 & 1 & 0 & 0 & 1 & 0 & 0 & 0 & 0 & 0 & 3 & 1 & 0 & 0 & 1 & 0 & 0 & 1 & 0 &  \\
\hline
\scriptsize$\{4,4,6,0,2,2\}$ & 1 & 0 & 0 & 0 & 0 & 1 & 1 & 0 & 0 & 0 & 0 & 1 & 1 & 0 & 0 & 0 & 0 & 0 & 1 & 0 & $\vee$ \\
\hline
\scriptsize$\{4,4,6,0,3,1\}$ & 1 & 0 & 0 & 0 & 0 & 1 & 0 & 0 & 0 & 0 & 0 & 2 & 1 & 0 & 0 & 0 & 0 & 0 & 1 & 0 &  \\
\hline
\scriptsize$\{4,4,1,4,4,1\}$ & 1 & 1 & 3 & 0 & 0 & 1 & 0 & 1 & 0 & 0 & 0 & 3 & 1 & 3 & 1 & 3 & 0 & 0 & 1 & 1 &  \\
\hline
\scriptsize$\{4,4,1,5,2,2\}$ & 1 & 1 & 0 & 0 & 0 & 1 & 1 & 1 & 1 & 0 & 0 & 1 & 1 & 1 & 1 & 0 & 0 & 0 & 1 & 1 & $\vee$ \\
\hline
\scriptsize$\{4,4,1,5,3,1\}$ & 1 & 1 & 0 & 0 & 0 & 1 & 0 & 1 & 0 & 0 & 0 & 2 & 1 & 2 & 1 & 0 & 0 & 0 & 1 & 1 &  \\
\hline
\scriptsize$\{5,5,6,0,0,2\}$ & 1 & 0 & 0 & 0 & 0 & 0 & 1 & 0 & 0 & 0 & 0 & 0 & 1 & 0 & 0 & 0 & 0 & 0 & 0 & 0 & $\vee$ \\
\hline
\scriptsize$\{5,5,1,5,0,2\}$ & 1 & 1 & 0 & 0 & 0 & 0 & 1 & 0 & 1 & 0 & 0 & 0 & 1 & 0 & 1 & 0 & 0 & 0 & 0 & 0 & $\vee$ \\
\hline
\scriptsize$\{5,6,6,0,0,1\}$ & 1 & 0 & 0 & 0 & 0 & 0 & 0 & 0 & 0 & 0 & 0 & 0 & 1 & 0 & 0 & 0 & 0 & 0 & 0 & 0 & $\vee$ \\
\hline
\scriptsize$\{5,6,1,5,0,1\}$ & 1 & 1 & 0 & 0 & 0 & 0 & 0 & 0 & 0 & 0 & 0 & 0 & 1 & 0 & 1 & 0 & 0 & 0 & 0 & 0 & $\vee$ \\
\hline
\scriptsize$\{5,6,1,1,4,1\}$ & 1 & 1 & 1 & 0 & 0 & 0 & 0 & 0 & 0 & 0 & 0 & 0 & 1 & 0 & 1 & 1 & 0 & 0 & 0 & 0 & $\vee$ \\
\hline
\scriptsize$\{5,6,2,0,4,1\}$ & 1 & 0 & 1 & 0 & 0 & 0 & 0 & 0 & 0 & 0 & 0 & 0 & 1 & 0 & 0 & 1 & 0 & 0 & 0 & 0 & $\vee$ \\
\hline
\scriptsize$\{5,2,4,5,0,2\}$ & 1 & 1 & 0 & 0 & 2 & 0 & 2 & 0 & 2 & 0 & 0 & 0 & 1 & 0 & 1 & 0 & 0 & 2 & 0 & 0 &  \\
\hline
\scriptsize$\{5,2,5,5,0,1\}$ & 1 & 1 & 0 & 0 & 1 & 0 & 0 & 0 & 0 & 0 & 0 & 0 & 1 & 0 & 1 & 0 & 0 & 1 & 0 & 0 & $\vee$ \\
\hline
\scriptsize$\{5,2,5,1,4,1\}$ & 2 & 1 & 1 & 0 & 1 & 1 & 0 & 0 & 0 & 0 & 0 & 0 & 2 & 0 & 1 & 1 & 0 & 1 & 1 & 0 &  \\
\hline
\scriptsize$\{5,2,6,3,0,2\}$ & 1 & 0 & 0 & 0 & 1 & 0 & 2 & 0 & 0 & 0 & 0 & 0 & 1 & 0 & 0 & 0 & 0 & 1 & 0 & 0 &  \\
\hline
\scriptsize$\{5,2,6,4,0,1\}$ & 1 & 0 & 0 & 0 & 1 & 0 & 0 & 0 & 0 & 0 & 0 & 0 & 1 & 0 & 0 & 0 & 0 & 1 & 0 & 0 & $\vee$ \\
\hline
\scriptsize$\{5,2,6,1,2,2\}$ & 1 & 0 & 0 & 0 & 1 & 1 & 3 & 0 & 0 & 0 & 0 & 0 & 1 & 0 & 0 & 0 & 0 & 1 & 1 & 0 &  \\
\hline
\scriptsize$\{5,2,6,1,3,1\}$ & 1 & 0 & 0 & 0 & 1 & 1 & 0 & 0 & 0 & 0 & 0 & 0 & 1 & 0 & 0 & 0 & 0 & 1 & 1 & 0 & $\vee$ \\
\hline
\scriptsize$\{5,2,2,5,2,2\}$ & 1 & 1 & 0 & 0 & 3 & 1 & 3 & 1 & 3 & 0 & 0 & 0 & 1 & 0 & 1 & 0 & 0 & 3 & 1 & 1 &  \\
\hline
\scriptsize$\{5,2,2,5,3,1\}$ & 1 & 1 & 0 & 0 & 2 & 1 & 0 & 1 & 0 & 0 & 0 & 0 & 1 & 0 & 1 & 0 & 0 & 2 & 1 & 1 &  \\
\hline
\scriptsize$\{5,3,5,0,4,1\}$ & 1 & 0 & 1 & 0 & 0 & 1 & 0 & 0 & 0 & 0 & 0 & 0 & 1 & 0 & 0 & 1 & 0 & 0 & 1 & 0 & $\vee$ \\
\hline
\scriptsize$\{5,3,6,0,2,2\}$ & 1 & 0 & 0 & 0 & 0 & 1 & 2 & 0 & 0 & 0 & 0 & 0 & 1 & 0 & 0 & 0 & 0 & 0 & 1 & 0 &  \\
\hline
\scriptsize$\{5,3,6,0,3,1\}$ & 1 & 0 & 0 & 0 & 0 & 1 & 0 & 0 & 0 & 0 & 0 & 0 & 1 & 0 & 0 & 0 & 0 & 0 & 1 & 0 & $\vee$ \\
\hline
\scriptsize$\{5,3,1,5,3,1\}$ & 1 & 1 & 0 & 0 & 0 & 1 & 0 & 1 & 0 & 0 & 0 & 0 & 1 & 0 & 1 & 0 & 0 & 0 & 1 & 1 & $\vee$ \\
\hline
\scriptsize$\{6,6,0,6,0,0\}$ & 0 & 1 & 0 & 0 & 0 & 0 & 0 & 0 & 0 & 0 & 0 & 0 & 0 & 0 & 0 & 0 & 0 & 0 & 0 & 0 & $\vee$ \\
\hline
\scriptsize$\{6,6,0,1,5,0\}$ & 0 & 1 & 1 & 0 & 0 & 0 & 0 & 0 & 0 & 0 & 0 & 0 & 0 & 0 & 0 & 0 & 0 & 0 & 0 & 0 & $\vee$ \\
\hline
\scriptsize$\{6,1,5,6,0,0\}$ & 0 & 1 & 0 & 0 & 1 & 0 & 0 & 0 & 0 & 0 & 0 & 0 & 0 & 0 & 0 & 0 & 0 & 0 & 0 & 0 & $\vee$ \\
\hline
\scriptsize$\{6,1,5,1,5,0\}$ & 0 & 1 & 1 & 0 & 1 & 1 & 0 & 0 & 0 & 0 & 0 & 0 & 0 & 0 & 0 & 0 & 0 & 0 & 0 & 0 & $\vee$ \\
\hline
\scriptsize$\{6,1,2,4,5,0\}$ & 0 & 1 & 1 & 0 & 1 & 1 & 0 & 2 & 0 & 0 & 0 & 0 & 0 & 0 & 0 & 0 & 0 & 0 & 0 & 0 &  \\
\hline
\scriptsize$\{6,1,2,6,3,0\}$ & 0 & 1 & 0 & 0 & 1 & 0 & 0 & 1 & 0 & 0 & 0 & 0 & 0 & 0 & 0 & 0 & 0 & 0 & 0 & 0 & $\vee$ \\
\hline
\scriptsize$\{6,1,2,6,1,2\}$ & 0 & 1 & 0 & 0 & 1 & 0 & 0 & 1 & 1 & 0 & 0 & 0 & 0 & 0 & 0 & 0 & 0 & 0 & 0 & 0 & $\vee$ \\
\hline
\scriptsize$\{6,1,3,6,0,2\}$ & 0 & 1 & 0 & 0 & 1 & 0 & 0 & 0 & 1 & 0 & 0 & 0 & 0 & 0 & 0 & 0 & 0 & 0 & 0 & 0 & $\vee$ \\
\hline
\scriptsize$\{6,3,0,4,5,0\}$ & 0 & 1 & 1 & 0 & 0 & 0 & 0 & 1 & 0 & 0 & 0 & 0 & 0 & 0 & 0 & 0 & 0 & 0 & 0 & 0 & $\vee$ \\
\hline
\scriptsize$\{6,3,0,6,3,0\}$ & 0 & 1 & 0 & 0 & 0 & 0 & 0 & 1 & 0 & 0 & 0 & 0 & 0 & 0 & 0 & 0 & 0 & 0 & 0 & 0 & $\vee$ \\
\hline
\scriptsize$\{6,4,0,6,0,2\}$ & 0 & 1 & 0 & 0 & 0 & 0 & 0 & 0 & 1 & 0 & 0 & 0 & 0 & 0 & 0 & 0 & 0 & 0 & 0 & 0 & $\vee$ \\
\hline
\scriptsize$\{4,5,0,4,5,0\}$ & 0 & 1 & 1 & 0 & 0 & 0 & 0 & 1 & 0 & 0 & 0 & 0 & 0 & 1 & 0 & 0 & 0 & 0 & 0 & 0 & $\vee$ \\
\hline
\scriptsize$\{4,6,0,3,5,0\}$ & 0 & 1 & 1 & 0 & 0 & 0 & 0 & 0 & 0 & 0 & 0 & 0 & 0 & 1 & 0 & 0 & 0 & 0 & 0 & 0 & $\vee$ \\
\hline
\scriptsize$\{4,6,0,6,2,0\}$ & 0 & 1 & 0 & 0 & 0 & 0 & 0 & 0 & 0 & 0 & 0 & 0 & 0 & 1 & 0 & 0 & 0 & 0 & 0 & 0 & $\vee$ \\
\hline
\scriptsize$\{4,6,0,6,1,1\}$ & 0 & 1 & 0 & 0 & 0 & 0 & 0 & 0 & 0 & 0 & 0 & 0 & 0 & 1 & 1 & 0 & 0 & 0 & 0 & 0 & $\vee$ \\
\hline
\scriptsize$\{4,6,0,2,5,1\}$ & 0 & 1 & 1 & 0 & 0 & 0 & 0 & 0 & 0 & 0 & 0 & 0 & 0 & 2 & 1 & 1 & 0 & 0 & 0 & 0 &  \\
\hline
\scriptsize$\{5,5,0,6,0,2\}$ & 0 & 1 & 0 & 0 & 0 & 0 & 0 & 0 & 1 & 0 & 0 & 0 & 0 & 0 & 1 & 0 & 0 & 0 & 0 & 0 & $\vee$ \\
\hline
\scriptsize$\{5,6,0,6,0,1\}$ & 0 & 1 & 0 & 0 & 0 & 0 & 0 & 0 & 0 & 0 & 0 & 0 & 0 & 0 & 1 & 0 & 0 & 0 & 0 & 0 & $\vee$ \\
\hline
\scriptsize$\{5,6,0,1,5,1\}$ & 0 & 1 & 1 & 0 & 0 & 0 & 0 & 0 & 0 & 0 & 0 & 0 & 0 & 0 & 1 & 1 & 0 & 0 & 0 & 0 & $\vee$ \\
\hline
\scriptsize$\{5,1,5,1,5,1\}$ & 0 & 1 & 1 & 0 & 1 & 1 & 0 & 0 & 0 & 0 & 0 & 0 & 0 & 0 & 1 & 1 & 0 & 1 & 1 & 0 & $\vee$ \\
\hline
\scriptsize$\{6,0,6,6,0,0\}$ & 0 & 0 & 0 & 0 & 1 & 0 & 0 & 0 & 0 & 0 & 0 & 0 & 0 & 0 & 0 & 0 & 0 & 0 & 0 & 0 & $\vee$ \\
\hline
\scriptsize$\{6,0,6,1,5,0\}$ & 0 & 0 & 0 & 0 & 1 & 1 & 0 & 0 & 0 & 0 & 0 & 0 & 0 & 0 & 0 & 0 & 0 & 0 & 0 & 0 & $\vee$ \\
\hline
\scriptsize$\{6,0,3,6,0,3\}$ & 0 & 0 & 0 & 0 & 1 & 0 & 0 & 0 & 1 & 0 & 0 & 0 & 0 & 0 & 0 & 0 & 0 & 0 & 0 & 0 & $\vee$ \\
\hline
\scriptsize$\{6,0,6,0,6,0\}$ & 0 & 0 & 0 & 0 & 0 & 1 & 0 & 0 & 0 & 0 & 0 & 0 & 0 & 0 & 0 & 0 & 0 & 0 & 0 & 0 & $\vee$ \\
\hline
\scriptsize$\{3,4,2,3,4,2\}$ & 2 & 2 & 2 & 0 & 2 & 2 & 2 & 2 & 2 & 2 & 2 & 2 & 2 & 2 & 2 & 2 & 0 & 2 & 2 & 2 &  \\
\hline
\end{longtable}
\end{landscape}
\restoregeometry}

Here we find another positroid variety, with $f=\{3,4,2,3,4,2\}$, where all its (torus-fixed) points are singular. It is also interesting to note that $\{4,2,4,2,4,2\}$ has the point with the highest multiplicity: $\lambda=\{1,3,5\}$, which has multiplicity 8. If we draw the pipe dream for $\{4,2,4,2,4,2\}$ at $\{1,3,5\}$, it looks as follows:

 \begin{figure}[htbp]
 \centering
	\includegraphics[scale=0.8,clip=true]{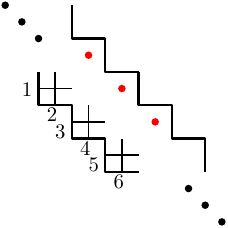}
\end{figure} 

We verify here that there are indeed 8 pipe dreams, since there are three crosses (all other tiles are elbows). Each of the three crosses can move one square to the northeast (to the tiles labeled by the red dots). So there are in total $2^3=8$ configurations.

\chapter{Affine Vs. Finite Pipe Dreams}\label{ch:affinefinitepipes}

We connect with the pipe dreams of \cite{Post} (or see \cite{ARW} for another nice explanation). These are just partitions or Young Diagrams that are filled with crosses and elbows. We distinguish them from affine pipe dreams by calling them finite pipe dreams, or simply nonaffine pipe dreams.


In this section, we will prove 3 propositions: \\
(1) Finite pipe dreams are embedded in ours by letting everything below the distinguished path be filled with crosses \\
(2) Our results about smoothness from earlier in the thesis (e.g. the main theorem, moves, deletion and contraction on pipe dreams etc.) apply to finite pipe dreams as well \\
(3) Finite pipe dreams are on a patch given by $\lambda=I_1$, the first Grassmann necklace element 

We first note, that one of the main differences between the affine and finite pipe dreams (besides the fact that the affine ones are infinite, so repeat indefinitely), is that the affine pipe dreams for a positroid variety $\Pi_f$ can take place on any patch $U_\lambda$ corresponding to the point $\lambda$ as long as $\lambda \in \Pi_f$; this patch is represented by the shape $S_\lambda$. More specifically, if we make a choice of $\lambda$, every permutation $f$ that can be obtained by filling the in $S_\lambda$ with crosses and elbow corresponds to a positroid variety $\Pi_f$ that contains the point $\lambda$. To give a simple example, both positroid varieties corresponding to $f_1=2222$ and $f_2=2231$ contain the point $\lambda=\{1,2\}$, so we can draw both $f_1$ and $f_2$ on $S_\lambda$; and $\mu$ is another point on $\Pi_{f_1}$, so we can draw $f_1$ on both $\lambda$ and $\mu$. 

In contrast, for finite pipe dreams, there is no such choice of patch or shape (let's denote these by $S'$ to disambiguate from the affine shape). Rather, the permutations $f$ that can be represented by each shape are disjoint. This can be proved easily by noting that the choice of $\lambda$ is the choice of vertical lines in the lower boundary of a finite pipe dream, and by construction, these will be the excedences (numbers in the bounded affine permutation that are greater than $n$). Therefore, a different choice of $\lambda$ (and shape $S'_\lambda$) will give rise to a permutation with different numbers. This of course assumes that we are fixing the starting point of the permutation (i.e. not considering cyclic shifts to be the same permutation). Given a permutation $f$, there is a unique shape $S'(f)$ that the finite pipe dream must exist on (and this shape must have $\lambda=I_1$, the first Grassmann necklace element - this is the content of the 3rd proposition of this section). In fact, the choice of patch or $\lambda$ for the finite pipe dreams is precisely the one such that (this will be relevant later in Section \ref{section:pipestoplabic} when we build a plabic diagram on a pipe dream using BCFW bridges).

\section{(1) The Embedding}

Here is a figure showing how the finite and affine pipe dreams are related to each other. The finite pipe dreams (left) are embedded in the affine ones (right)

\begin{figure}[htbp] \centering
\includegraphics[scale=0.6,clip=true]{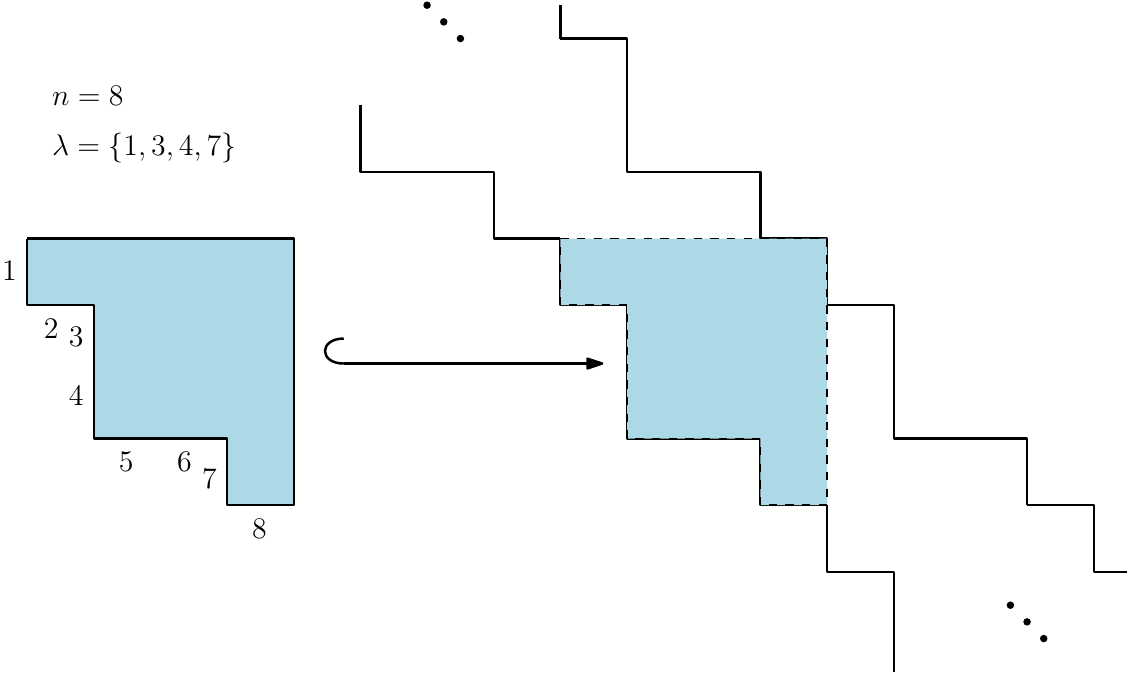}
\caption{}
\end{figure}

In the language of \cite{Post}, given $Gr(k,n)$, if we have a Young diagram $(\mu_1, \mu_2,...,\mu_k)$, the vertical lines (corresponding to the $\lambda$ in the notation in this thesis) is $\lambda= \{ n-\mu_1 +1, n-\mu_2+2,...,n-\mu_k+k\}$. 

The proof that this is the correct way to embed simply comes from the way the diagrams are read: $g(i)$ is read off with the values of $i$ being along the left (lower) boundary, and the values of $g(i)$ being the labels along the right (upper boundary). So we have to match up the labels on the lower paths; that is, the distinguished paths must match.

This provides the placement of the shapes. However, the affine pipe dreams contain more tiles. The correct embedding will fill in (every tile in) the piece below the distinguished path with crosses. The proof of this is trivial. Given that we already have an embedding of the finite pipe dreams given by letting the lower boundary coincide with the affine pipe dream, to represent the same permutation, the pipes must land on the same numbers $g(i)$ on the upper boundary. Given the way that these diagrams are constructed, this implies that once we hit a boundary of the finite pipe dream shape, within the affine pipe dream shape, we must continue going in the same direction (that is, continue going horizontally at a vertical boundary line segment, and continue going vertically at a horizontal boundary line segment). Continuing on with the previous figure example, to make things concrete, this would mean filling in the red portions entirely with crosses:

\begin{figure}[htbp] \centering
\includegraphics[scale=0.6,clip=true]{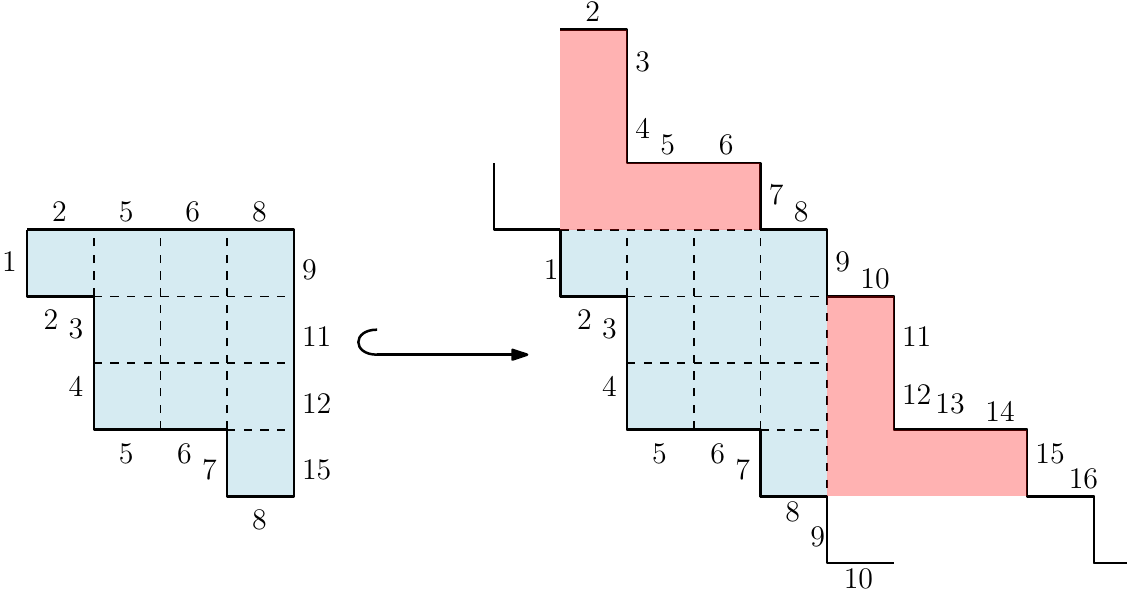}
\caption{}
\end{figure}

The pipe dream with (every tile in) the piece below the distinguished path filled in with crosses is the largest positroid variety that has Grassmann necklace $I_1=\lambda$; in other words, if any tile inside this piece is an elbow, the positroid variety is too big, so it has $I_1$ that is lower in value (Gale order) than $I_1=\lambda$. We can provide a proof simply by noting that the finite pipe dreams, embedded in the affine ones, only give the correct permutations (that they assign) when we follow the pipes if we fill in the portion below the distinguished path entirely with crosses. 

\section{(2) Smoothness on Finite Pipe Dreams}

Note that these crosses are entirely fixed - they have no moves because the entire northeast of the southwest corner is filled in and this corner touches a boundary, so it cannot move; thus, since all of the other crosses inside this block are separated by more than one pipe so cannot participate in a move and require this southwest corner to move first to create a move. This is displayed for example by the green crosses in the following diagram.
\begin{figure}[htbp] \centering
\includegraphics[scale=0.6,clip=true]{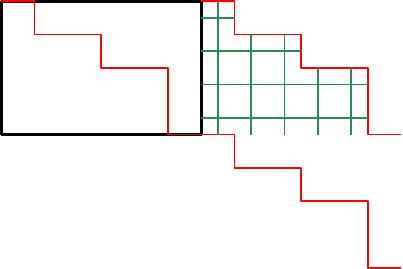}
\caption{}
\end{figure}

This proves that deletion and contraction on the smaller pipe dreams happens in the same way as in the affine pipe dreams: we have to move crosses and near-misses in this restricted portion (partition). Furthermore, since the Le condition of \cite{Post} corresponds to the bottom pipe dream as noted in \cite{Flu}, moves must consists of crosses moving northeast to switch with near-misses.  

Our Main Theorem in this case applies to maximal rectangles in the finite pipe dream as well. The reason/proof is as follows: any maximal rectangle in the affine pipe dream that includes boxes beyond the finite pipe dream will only contain entire rows and columns of crosses; these entire rows and columns can thus be eliminated without affecting the analysis of the moves (in the remaining diagram). Therefore, any move will come from a consideration of the boxes inside the finite pipe dream.(*) Thus, it is sufficient to consider all maximal rectangles inside the finite pipe dream. To give an example from the figure below: a maximal rectangle in the original affine pipe dream is boxed in green. However, the maximal rectangle of the finite pipe dream that is boxed in purple contains all the relevant information in the green rectangle (more precisely, any move in the green rectangle will show up in the purple rectangle, so it is sufficient to consider the purple rectangle). 

\begin{figure}[htbp] \centering
\includegraphics[scale=0.6,clip=true]{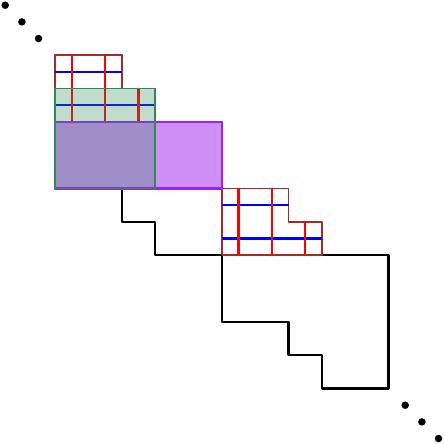}
\caption{}
\end{figure}

(*) We should comment here that, if the piece underneath the distinguished path (let's call this the \textit{lower partition} is filled in with crosses, then they are locked in place - they cannot participate in any moves. The proof is that the piece above the distinguished path (the "upper partition") touches at the NW and SE corners (if it touches at all). The pipes in the lower partition can be divided into two disjoint sets $A$ and $B$: every pipe in $A$ crosses all the pipes in $B$ and vice-versa (see the figure above). The red pipes correspond to say $A$, and the blue pipes to say $B$. They seem to exit into different parts of the pipe dream that have no contact with each other; however, we must remember that this diagram is actually an "unfurled" cylinder (see the last section in \cite{Flu}: since the affine pipe dream is actually periodic, it is indeed possible theoretically for two pipes that cross in the lower partition to move into two different regions in the affine pipe dream that are actually the same region on a cylinder. This connects to the issue of "double crossings" that is noted in \cite{Flu} - as we noted there as well, this issue does not arise in these affine pipe dreams because they are bounded of size $n-k$ in the horizontal direction and size $k$ in the vertical direction. 

Now, given this new Main Theorem, this implies that we can consider smoothness in this more restricted case (where for a positroid variety $\Pi_f$ and a point $\lambda$, we only consider points $\lambda$ on $\Pi_f$ such that $I_1$ in the Grassmann necklace for $f$ is greater than or equal to $\lambda$ (for example, if $I_1=\{1,3,5\}$, then $\lambda=\{1,3,4\}$ is ok, but $\lambda=\{1,4,5\}$ will not work). In this case, we can just consider a word for the finite pipe dream, and see if there are multiple subwords for $f$ (if there are, then it is singular, vs. if there is a single one, then it is smooth). This is really not much different than the more general case we have already considered, since we are essentially just noting that the letters in the word corresponding to the lower partition are not moveable (do not participate in moves) meaning that they are fixed, so we just remove them from the word, and look at subwords of this smaller word. 

As such, given that by subsection (1) the finite pipe dream, when embedded in the affine pipe dream, is just a copied over, with some additional extra tiles that cannot participate in any moves, it is essentially a corollary of this that our results involving moves, deletion/contraction, and smoothness on affine pipe dreams from the first part of this thesis carry over essentially unchanged to the finite pipe dream case. In particular, we can consider smoothness the same way on finite pipe dreams: we can also look at at maximal rectangles. Indeed, all maximal rectangles in an affine pipe dream are the same as maximal rectangles in finite pipe dreams, mod some extra entire rows and columns of crosses (as in the previous figure), so will produce no differences when looking for rectangles reducing to NW/SE partitions.

\section{(3) The Choice of Shape for Embedding}

Finally, we need to state what patch corresponding to a point $\lambda$ a finite pipe dream for a permutation $f$ is drawn on. Once we have a $\lambda$, subsection (1) has already shown us how to embed the finite pipe dream within the affine pipe dream and we arrive at subsection (2) as a corollary that the previous results on moves, deletion/contraction, and smoothness for affine pipe dreams apply in the finite pipe dream case as well. However, we have not yet shown which patch (that is, $\lambda$) should be chosen for a finite pipe dream of permutation $g$, so we discuss that here. 

Now, due to the way the diagrams are labeled, this is perhaps trivial. The excedences (numbers greater than $n$) label the boundary on the right, and the corresponding number $-n$ label the parallel (vertical) line on the left boundary. Since $\lambda$ is equal to the vertical lines on the left boundary, we just take the excedences $mod n$, and we have $\lambda$. However, we can show something a little more general. Suppose that we only knew that the upper and right boundary numbers $mod n$, so we didn't know which ones were excedences. Or suppose that we only knew the first Grassmann necklace element $I_1$ of $g$ (and we did not know $g$). We can prove that $I_1$ is equal to the vertical lines in the left boundary. This will show that the excedences ($mod n$) are equal to $I_1$. We have the following proposition:

\begin{Prop}
Suppose we are in $Gr(k,n)$ and we have a finite pipe dream shape $S'_\lambda$ corresponding to a distinguished path/lower boundary given by a $k$-element set $\lambda$ of numbers in $\{1,2,...,n\}$. (Note we do not know what pipe dream we have; that is, we do not know how the shape is filled in with crosses and elbows, so we do not know the permutation $g$ - we only know the outer boundary). Regardless of how the pipe dream is filled in, it must produce a permutation $g$ whose first Grassmann necklace element $I_1$ satisfies $I_1=\lambda$.
\end{Prop}

\begin{proof}

We want to show that the \textit{affine} pipe dream drawn on the patch corresponding to $\lambda$, with the portions below the distinguished path entirely consisting of crosses, corresponds to bounded affine permutations with the first Grassmann necklace element $I_1$ being equal to $\lambda$. 

First, we note that if we fill in the affine pipe dreams entirely with crosses, then we get exactly the point $\lambda$, given by a siteswap with $n$ where the $\lambda$ entires are, and $0$'s elsewhere. This is the siteswap that is highest in Bruhat order: all other siteswaps $f$ that exist within this affine pipe dream shape will be obtained via transpositions from the siteswaps with all $0$'s and $n$'s. Therefore, every siteswap $f$ coming from a set of pipes in this affine pipe dream shape will contain the point $\lambda$ (this is really just a re-statement of the idea of words/heaps from before), and this implies that $I_1$ cannot be greater than $\lambda$ in Gale order: since any $f$ in this pipe dream gives a positroid variety $\Pi_f$ that contains the points $\lambda$, the determinant of the $\lambda$ columns of a matrix representative of $\Pi_f$ must be nonvanishing (so $I_1$ must be less than or equal to $\lambda$ in Gale order).

We now show that $I_1$ must be equal to $\lambda$ on the nose. In other words, any $\mu < \lambda$ in Gale order must not be an element in the positroid for $f$; stated another way for a matrix representative $M$ of $\Pi_f$, $det_{\mu}(M)=0$, where we take the determinant of the $\mu$ columns of $M$. 

Let's denote the elements of $\lambda = \{\lambda_1,\lambda_2,....,\lambda_k\}$. These are the vertical line segments along the left boundary of the pipe dream (we are going with the orientation of the affine pipe dreams, so for the non-affine pipe dreams, we may have to reflect across the vertical axis to get it in the same orientation as in e.g. \cite{Post}). We can also consider the permutation elements, that is, for bounded affine permutation $g$, the set: $g(\lambda) = \{g(\lambda_1),g(\lambda_2),....,g(\lambda_k)\}$.

We use the fact that, column $i$ is linearly dependent on the columns $i+1,...,g(i)$, so any minor containing columns $i$ to $g(i)$ must vanish.

Suppose we have a $\mu<\lambda$. This means that, if we match up all $\mu_i \in \{\mu\}$ with the corresponding $\lambda_i \in \{ \lambda\}$, that there is some $\hat{i}$ with $\mu_{\hat{i}} < \lambda_{\hat{i}}$. Since we are considering $I_1$, so $1$ is smallest in the Gale order, and $\mu<\lambda$, this means that all $\mu_i \in \{\mu\}$ satisfy $1 \leq \mu_i \leq \lambda_k$ (and obviously $1\leq \lambda_1 \leq \lambda_i \leq \lambda_k$). Let's define another set: $\overline{\{\lambda_k\}}$ consisting of all numbers $\{1,2,...,\lambda_k\}$ that are \textit{not} in $\lambda$: these are the horizontal line segments in the range $\{1,2,...,\lambda_k\}$ along the bottom/left boundary. By what we have noted about there being some $\hat{i}$ with $\mu_{\hat{i}} < \lambda_{\hat{i}}$, there must be at least element of $\lambda$ that has become replaced by an element of $\overline{\{\lambda_k\}}$ in $\mu$: that is, at least one vertical line segment in $\lambda$ has become a horizontal line segment in $\mu$. Consider the least such horizontal line segment, which we can denote $\mu_{\overline{max}}$. Consider which pipe \textit{lands on} $\mu_{\overline{max}}$; i.e. which $i$ has $g(i)=\mu_{\overline{max}}$. It cannot be that $\mu_{\overline{max}}$ lands on itself, because this would imply that we have a vertical line of crosses, which would correspond to a zero column ($v_{\mu_{\overline{max}}} = \vec{0}$ a column in $M$), so the determinant $det_{\mu}(M)=0$, so $I_1$ (the minimal set of columns, start from 1, with nonvanishing determinant) cannot be equal to $\mu$. Therefore, we must have $i< \mu_{\overline{max}}$. 

$\mu_{\overline{min}}$ is the minimal horizontal line, and this implies that all elements of $\mu$ that are less than it are vertical lines. If $\mu_{\overline{min}}$ is the minimal element in $\mu$, that implies that there are no vertical lines smaller than $\mu_{\overline{min}}$ (since we can only replace vertical lines with horizontal lines that are less than it). This implies that $\mu_{\overline{min}}$ must be a part of a segment of $1,2,...,i$ consisting entirely of horizontal line segments. And for these shapes, this means that all of $1,2,...,i$ are zero columns. Thus, $det_\mu(M)=0$. However, if $\mu_{\overline{min}}$ is not the minimal element of $\mu$, then there exists at least one vertical line smaller than it (and, as we just noted, in fact all vertical lines less than $\mu_{\overline{min}}$ must be in $\mu$). $\mu_{\overline{min}}$ can't go to itself (otherwise it would be a zero column), so it must come from a pipe $i$ such that $1\leq i \leq \mu_{\overline{min}}$. If it comes from a pipe $i$ that is a vertical line segment, then we are done. If not, then it comes from a horizontal line segment - but we can ask the same question about this horizontal line segment $h$ and what pipe $j$ maps to it, and we must have $1\leq j\leq h$. And we can continue the same argument. 

In the juggling diagrams way of viewing things, you need the $k$ spots to be different colors - different lines. So if we can show that two of the $k$ numbers in $\mu$ are the same line, then we are done. This gives us an even more direct proof of the other direction: the vertical lines are independent. In other words, ``no vertical line can connect to another vertical line."

Let's do an example that demonstrates this proof. 

Suppose we are in $Gr(3,6)$ with $\lambda=\{1,3,5\}$ (see the following figure). Let $\mu=\{ \mu_1, \mu_2, \mu_3\}$ be a 3-number set less than $\lambda$ in Gale order. Since $I_1=\{i_1,i_2,i_3\}$ is the first (in Gale order) set of 3 columns of a matrix representative with nonvanishing determinant $Det(v_{i_1},v_{i_2},v_{i_3})\neq 0$, and since a permutation $g$ encodes the linear dependencies, we cannot have $g(i_1)=i_2,i_3$ or $g(i_2)=i_3$; otherwise we would have a linear dependency so that $Det(v_{i_1},v_{i_2},v_{i_3})= 0$. It's easy to see that, in the diagram above this condition is met automatically from the shape: it is impossible to have $g(1)=3,5$ or $g(3)=5$.

\begin{figure}[htbp] \centering
\includegraphics[scale=0.6,clip=true]{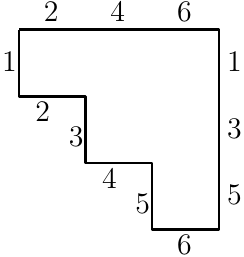}
\caption{}
\end{figure}

We take $\mu$ that is smaller in Gale order than $\lambda$, so it must have at least one horizontal line (since is it not equal to $\lambda$). Let's take the last such number. For example, suppose $\mu=\{1,3,4\}$ (so we have replaced $5$ in $\lambda$ by $4$). We ask the question of what pipe $i$ lands on $4$ (what $i$ satisfies $g(i)=4$). Diagrammatically, the placement of $4$ at the top implies that the $i$ landing on $4$ must be a smaller number: $1$, $2$, or $3$. By the comment in the previous paragraph about vanishing determinants, $i$ must not be in $\mu$ ($i\not\in \mu$), otherwise $Det(v_{\mu_1},v_{\mu_2},v_{\mu_3})= 0$. Now, there are two possibilities for $i$: it is either (1) a vertical or (2) a horizontal line along the bottom boundary. If it is a (1) horizontal line, then we can ask the same question: which $1\leq j<i$ satisfies $g(j)=i$? If (2), a vertical line, then by what we just said, namely that $i\not\in \mu$, since $i$ is a vertical line, the fact that we are smaller in Gale order implies that there is another $1\leq \hat{i}<i$, with $\hat{i}$ along a horizontal line such that $\hat{i}\in \mu$ replaces $i\in \lambda$. Then, we can ask the same question about $\hat{i}$: which $1\leq \hat{}j<\hat{i}$ satisfies $g(\hat{j})=\hat{i}$? We see that we cannot continue forever since this is a finite set: eventually, we will find that we must have a horizontal line $h$ such that $h$ is either a zero column $g(h)=h$ or one of the vertical lines $l\in \lambda,\mu$ has $g(l)=h$, in which case   $Det(v_{\mu_1},...,v_{\mu_k})= 0$.

In our example here, let's suppose that the largest horizontal line is $4$. There's only one number in $\lambda$ that is greater than $4$, so clearly the $4\in \mu$ is replacing the $5\in \lambda$. We ask which line goes to $4$. \\
Case (1) Suppose $g(2)=4$. We then ask the same question about $2$: which $i$ satisfies $g(i)=2$. Either $g(2)=2$, which contradicts that $g(2)=4$, or $g(1)=2$ which means that $Det(v_{1},v_{\mu_2},v_{4})= 0$ (and $1\in \mu$ necessarily since $1\in \lambda$ and $\mu <\lambda$ in Gale order). \\
Case (2) Suppose $g(3)=4$. Then $Det(v_{\mu_1},v_{\mu_2},v_{4})= 0$ if $3\in \mu$ so we must have $3\not\in\mu$. But since $3\in \lambda$, since $\mu$ is smaller is Gale order, there must be an $i<3$ with $i\in \mu$ and $i$ horizontal. Here, the only choice is $2$. So we must have $2\in \mu$. But then we can ask the same question about which $j$ satisfies $g(j)=2$. Then the only possibility is $g(1)=2$, so we end up with the same conclusion as case (1). \\
Here $g(1)=4$ is also a vertical line going to $4$, but then $1\in \mu$ so $1$ must be replaced by some other number, and there are only larger numbers to replace $1$, so $\mu$ cannot be smaller in Gale order, which is a contradiction. 

\end{proof}

\part{Planar $N=4$ Supersymmetric Yang-Mills Theory}

\chapter{Quantum Field Theory}\label{ch:QFT}

\section{Basic Quantum Field Theory}

Perhaps the fundamental quantity that is calculated in quantum field theory is the amplitude. 

In traditional quantum field theory, these amplitudes are calculated using Feynman diagrams. A (quantum field) ``theory" is a Lagrangian: an expression that tells you what particles there are and what kinds of interactions are possible in this theory. Quantum field theory seeks to combine quantum mechanics and special relativity. Since the language of special relativity uses 4-dimensional vectors in flat $(t,\vec{x})=(1,3)$ spacetime that have a symmetry under Lorentz transformations (rotations and boosts), the objects (fields) in the Lagrangian have a symmetry under Lorentz transformations (such that the Lagrangian is invariant (up to a total derivative)); in particular, the fields are representations of the Lorentz group. For example, a very simple Lagrangian involving a single scalar field $\phi$ would be: $\mathcal{L}=   \frac{1}{2}(\del_\mu \phi)^2- \frac{1}{2} m^2 \phi^2 + g\phi^3$. In this theory, there is a ``kinetic energy" term $\frac{1}{2}(\del_\mu \phi)^2$, a mass term for the field $\frac{1}{2} m^2 \phi^2$, and a single interaction term $g\phi^3$ for the interaction of 3 particles. 

In the early days of quantum field theory, researchers came up with tricks to take a Lagrangian and compute the amplitude; this was called old-fashioned perturbation theory. But then Richard Feynman discovered that an intermediate step made things much easier: one could take the Lagrangian, write down some rules and then draw all diagrams using those rules; those rules tell you what algebraic expression to write down for each diagram that it contributes to the total amplitude. This approach is ``perturbative": to good approximation, the amplitude is calculated by drawing all diagrams up to some given order (say, given by the number of nodes in the diagrams, or the number of loops in the diagrams - which keeps track of the order of one of the constant coefficients in the Lagrangian expression). One problem with this approach is that the manifest Lorentz-invariance of the Lagrangian implies that the expressions usually put the information into tensor products of 4-vectors, which gives rise to gauge symmetries/redundant information (there are these extra symmetries of the fields under which the Lagrangian is invariant up to total derivative, so when integrating, one needs to make a particular choice (gauge-fix)). This can cause problems when applying the Feynman rules in the same way for simpler theories, so (in say, nonabelian gauge theory) additional rules involving ghost fields are needed to get the rules to work out the right way. Because the on-shell approach to amplitudes does not allow for virtual particles, there is no issue with gauge symmetry; instead, one loses the manifest Lorentz invariance, which can show up in what are called spurious poles.  

Feynman diagrams give us: $\langle f |T|i\rangle = (2\pi)^2 \delta(\sum p) M$. We square this and put it into $d\sigma =\frac{1}{(2E_1)(2E_2) \abs{\vec{v_1}-\vec{v_2}}} \abs{M^2} \Pi_{\text{final states j}} \frac{d^3p_j}{(2\pi)^3}\frac{1}{2E_{p_{j}}} (2\pi)^4 \delta(\sum p)$, and $d\sigma$ is related to experiment - roughly, $d\sigma$ is the thing that corresponds to the data that might be obtained from a particle physics experiment; in general, it depends on angles and energies. It tells us the probability that a scattering process will occur. We can write it as:

\[d\sigma =\frac{1}{(2E_1)(2E_2) \abs{\vec{v_1}-\vec{v_2}}} \abs{M^2} \Pi_{\text{final states j}} \frac{d^3p_j}{(2\pi)^3}\frac{1}{2E_{p_{j}}} (2\pi)^4 \delta(\sum p)\]

Most of this is just experimental stuff that is the same for every experiment. The meat of the expression is all in the $\abs{M}^2$ term.
 
 $M$ is called the "scattering amplitude". It is calculated via the S-matrix:
\[\langle f|S|i\rangle_{Heisenberg} = \langle f; t=\infty | i;t=-\infty\rangle_{Schrodinger} = lim_{t_1\to -\infty, t_2\to \infty} \langle f | e^{iH(t_2-t_1)} |i\rangle\]
 
We can calculate S-matrices using the theoretical tools of Feynman diagrams. Actually, the scattering amplitude is really in the T-matrix: 

$S=1+iT$, so we separate $\langle f | i\rangle$, a process that does nothing. The relation for the scattering amplitude is: 

$\langle f |T|i\rangle = (2\pi)^2 \delta(\sum p) M$. The delta function imposes momentum conservation of the process. We calculate this using the LSZ formula. A quick summary (following the approach of \cite{Schw}) with some more details can be found in my slides ``Supersymmetric Juggling" at \url{josephflueg.github.io}.

The theory we will be working with is $\mathcal{N}=4$ supersymmetric Yang-Mills theory.
  
Supersymmetry means that there is a symmetry of the Lagrangian of exchange of fermions (even-dimensional representation, half-integer spin/helicity) and bosons. So there are supersymmetry generators (Lie algebra generators of the Lie group). $\mathcal{N}=4$ is the maximally supersymmetric theory in 4 dimensions, having 4 supersymmetry generator pairs. 
  
In $\mathcal{N}=4$ SYM is a theory, there are 16 particles/states. The maximal supersymmetry implies that every field can be moved to another field via applications of generators. Thus, we can write all 16 states in a single superfield expression $\Omega$. Then, we can write a superamplitude $A[\Omega_1,...,\Omega_n]$ consisting of $n$ external particles, and obtain the correct expression (depending on what particle $q_i$ we want each $\Omega_i$ to correspond to) by apply the correct operator that projects $q_i$ out from $\Omega_i$. 

Our focus in on the case of only planar diagrams, so we restrict to planar $\mathcal{N}=4$ SYM.

We also assume that all particles are massless.

We will not say more about $\mathcal{N}=4$ SYM since it is not central to our discussion, and apologize for the brevity, referring to Chapter 4 of \cite{EH} for more details.

\section{On-Shell Physics}

At a high level, we can think of this textbook quantum field theory approach, encapsulated in Feynman diagrams as calculations on the boundary of spacetime: particles come in at the beginning of spacetime, and then we measure at the end of spacetime, and use this data to calculate the scattering amplitude. In terms of what happens in the middle of spacetime, this is supposed to be represented by Feynman diagrams which depict particles colliding at local spacetime points: it is assumed that there is local evolution in spacetime reflecting the underlying unitary evolution in Hilbert space. 

Much of the motivation for the on-shell approach came from formulas such as the Parke-Taylor formula, where the Feynman diagrammatic approach produces thousands of terms, but the final answer reduces to a single line. This formula suggested that scattering amplitudes have deeper structure than what could be understood in Feynman's approach and opened the question of whether there are easier ways to calculate the scattering amplitude. 

We can think of locality (in spacetime) and unitarity (from quantum mechanics - encoding the conservation of probability) as being the theoretical starting point of the Feynman diagrammatic approach: local, unitary evolution in spacetime are assumed, which becomes encoded in the answer (the scattering amplitude) as certain analytic properties. In contrast, the on-shell scattering amplitudes research, in a sense, seeks to go in the opposite direction: it seeks to derive locality and unitarity from a deeper, more fundamental structure.

One idea that has proved to be very powerful involves positive geometries in kinematic space. The kinematic space is just the data that goes into a scattering amplitude calculation: the labels for the physical on-shell momenta of particles in the scattering process (no virtual particles). It turns out that with no notion of locality or unitarity or even time, we can take the kinematic space and apply mathematical techniques that allow us to do physics, and to eventually see where locality and unitarity come from. These ideas were first discovered in planar N=4 super Yang Mills theory with on-shell diagrams that turned out to be connected to positroids and positroids varieties, and later the amplituhedron. These universal long-distance geometries also capture soft singularities.

Let us explain how a geometry determines a function (that depends on kinematical data) like a scattering amplitude. A positive geometry is roughly a compact region of $\R^n$ for some $n$ that has boundaries of all codimensions. They can be polytopes, or more curved mathematical objects like Grassmannians. This geometry determines a unique differential form called a canonical form, determined by requiring that it have logarithmic singularities on all boundaries (and only the boundaries) of the positive geometry; furthermore, taking resides at a maximal dimension boundary element gives us the canonical form for this boundary (so the residue now has logarithmic singularities only at all boundaries of the boundary element we took the residue of). ``Logarithmic singularity" means that locally there is a coordinate chart where the form looks like the wedge product of a bunch of $dlog$s, e.g. $dlog(x_1)\wedge \cdots \wedge dlog(x_n)=\frac{dx_1}{x_1}\wedge \cdots \wedge \frac{dx_n}{x_n}$. So this is the process:  kinematic space $\to$ positive geometry $\to$ canonical form $=$ scattering amplitude.

We could ask the question of how we calculate the canonical form given a positive geometry. One way to do it would be to triangulate (a notion much more difficult to define for curved spaces): we cut it up into simplices, and since (1) We know what the canonical form is for simplices and (2) we know that summing up the canonical form causes ``spurious poles" that are boundaries between simplices in the triangulation (so not boundaries of the positive geometry) to cancel out, we can just sum up the canonical form for the pieces in the triangulation to get the canonical form for the whole geometry. One particularly nice result is that for a simple polytope (all vertices are locally simplices), the canonical form is $\Omega=\sum_{vertices}sign\abs{X_{a_1}\cdots X_{a_n}}dlogX_{a_1}\cdots dlogX_{a_n}$ (this is related to Feynman diagrams) \cite{AH}.

\section{Leading Singularities}\label{sec:leadsing}

We refer to the book \cite{EH} for more details. You can also refer to notes from my talk on ``What is a Leading Singularity" at the IAS, posted at \url{josephflueg.github.io}. 

What is a leading singularity, and what is its connection to the amplitude in quantum field theory that we gave above?

While there are several related notions of ``leading singularity," the one we use is that a leading singularity is a solution to a maximal cut. To understand this, we should briefly discuss unitary methods.

Roughly, for a certain loop order, one can write the amplitude $A$ in terms of a basis, something like

$A= c_1 B_1 + \cdots c_m B_n$ (this is quite simplified - c.f. Eq 6.26 of \cite{EH}). 

The basis consists of loop integrals that are easier to evaluate, for example what are called ``scalar integrals" that do not have any factors involving momenta in the numerator. Thus, if one can calculate the basis integral for the theory, getting the amplitude becomes a matter of obtaining the coefficients $c_i$. 

It turns out that the coefficients can be obtained using a process called generalized unitarity. One applies maximal cuts. A ``cut" means that you put an internal propagator on-shell ($p^2=m^2$), so you no longer integrate over it; it is as if you have ``cut" an internal line. Removing the propagator from the denominator is equivalent to putting in a delta-function $\delta(p^2-m^2)$ enforcing on-shellness of that momentum in your original expression, which can also be viewed as taking the residue of your original expression at the pole where $p^2=m^2$ (this is sketched in my Leading Singularity notes p3-4). It is called maximal cuts because you are putting as many of the internal propagators on shell as possible. To be precise, one looks at both the left and right side of the equality in $A= c_1 B_1 + \cdots c_m B_n$ and applies these cuts to both sides, and this gives you expressions for the $c$'s in terms of the $B$'s and other $c$'s that you've already solved for. However, in planar $N=4$ supersymmetric Yang-Mills theory, things are much simpler:

(1) You can actually think of each $c_i$ as being a single maximal cut diagram itself (no linear combinations of such diagrams) \\
(2) Most of the basis elements $B_i$ die off - for example, for 1-loop, the bubble and triangle diagrams give 0 contribution to the amplitude, so you only need the box diagrams \\
(3) In fact, all planar loop integrands can be written in a basis of unit leading singularity integrands; this means that all $B_i=$ $0$, $1$, or $-1$. Therefore, the amplitude in planar $N=4$ SYM will literally be just a $\pm$ sum of maximal cut diagrams (and a choice of maximal cut diagram solution is just a leading singularity, as noted in the following paragraph).  

Note that each cut puts a propagator on-shell, which is setting a degree two equation to zero (a quadratic constraint), so it usually has multiple solutions. Thus in general a maximal cut has multiple solutions. Each solution is called a leading singularity. The plabic diagrams (the diagrams with nodes being colored white or black) are simply a representation of a particular solution/leading singularity. This is where positroids enter the story: each positroid thus corresponds to a maximal cut.

The place where positroid varieties enter the story is via the link formula. There is a formula that produces all leading singularities of planar N=4 SYM; the link formula for $n$ external momenta with helicity $k$ is:
\[L_{n,k} =\int_{C} \frac{d^{k\times n} C}{vol(GL(k))}\frac{1}{(1\cdots k)\cdots(n\cdots k-1)} \delta^{k\times 4}(C\cdot \tilde{\eta})\delta^{k\times 2}(C\cdot \tilde{\lambda}) \delta^{2\times (n-k)}(\lambda \cdot C^\perp)
\]

In this formula, $C$ is a matrix the parameterizes the restriction of the (complex) positroid variety to the nonnegative real part. 

We already noted by (3) that in planar N=4 SYM, the amplitude will simply be of the form $A= \pm c_1 + \cdots \pm c_m $, where each $c_i$ is a leading singularity, that is, a maximal cut solution corresponding to a particular on-shell diagram which is labeled by a positroid $f_i$. What we've just said is that then to actually calculate the $c_i$'s, one takes the positroid variety $\Pi_{f_i}$ for $f_i$ and writes a parameterization for $\Pi_{f_i}$ in matrix form that we call $C_i$, which we then plug into the link formula.

When we calculate the link formula, we choose a particular positroid cell $\Gamma_f$ with bounded affine permutation $f$: then the integral is $\oint_{C\subset \Gamma_f}$, and $(i\cdots i+k)$ denotes the cyclic minors of $C$ involving columns $i$ to $i+k$. In the notation of \cite{EH}, writing $M_j$ for these minors instead, switching to angle and square notation ($\lambda=|i\rangle$,$\tilde{\lambda}=|i]$), and writing out the entries of $C$ explicitly:
\[L_{n,k} = \int \frac{d^{k\times n} C}{GL(k) \Pi_{j=1}^n M_j} \bigg[ \Pi_{a=1}^k \delta^2(\sum_i C_{ai} [i|) \delta^{(4)} (\sum_i C_{ai} \eta_{iA}) \bigg] \bigg[ \Pi_{a'=k+1}^n \delta^2(\sum_iC^\perp_{a'i} \langle i |) \bigg] 
\]
Or, taking the first formula of $L_{n,k}$ from \cite{ABCGPT} and parameterizing the cell $\Gamma_f$ with coordinates $\alpha_1,...,\alpha_d$ ($d$ is the dimension of the positroid cell), as in the BCFW bridge decomposition method, the formula becomes:
\[L_{n,k}(f)= \int \frac{d\alpha_1}{\alpha_1} \cdots \frac{d\alpha_d}{\alpha_d} \delta^{k\times 4}(C\cdot \tilde{\eta}) \delta^{k\times 2}(C\cdot \tilde{\lambda}) \delta^{2\times (n-k)}(\lambda \cdot C^\perp)\]

In terms of solving this, it could take some work, but we can do it for simple examples. Our matrix parameterization $C_i$ for the positroid variety $\Pi_{f_i}$ will be written in terms of variables $(\alpha_i)_j$: we have examples below of how to explicitly write this matrix $C_i((\alpha_i)_j)$. We need the data of the external momenta $p_h$ we put into our amplitude calculation.
This should (at least in easier cases), allow us to solve for the $(\alpha_i)_j$ in terms of the $p_h$, using the delta functions in the numerator. We plug these expressions for the $(\alpha_i)_j$ written in terms of the $p_h$ into denominator terms involving determinants. In the nicest cases, the integral is fully localized by the delta functions, so this gives us our leading singularity; however, in general, we may have to perform some integrations. For an explicit example, we refer to Section 9.3 of \cite{EH}; some of the details left to the exercises are done in my partial textbook solutions document at \url{josephflueg.github.io}

\section{Grassmannian Spaces}

Why do we deal with Grassmannian spaces?

The Grassmannian aspect appears right from the start in the notation for spinor helicity. In spinor helicity, momentum conservation is written as: $\sum_j\lambda_\alpha^j\tilde{\lambda}_\alpha^j=0$. It might be easier to visualize if we wrote it out the momenta, rather than as 2-by-2 matrices, as vectors with 4 components: $\lambda^\alpha \lambda^{\dot{\alpha}}=\mqty(\lambda^1 \\ \lambda^2) \mqty(\tilde{\lambda}^{\dot{1}}& \tilde{\lambda}^{\dot{2}})=\mqty( \lambda^1 \tilde{\lambda}^{\dot{1}} & \lambda^1 \tilde{\lambda}^{\dot{2}} \\ \lambda^2 \tilde{\lambda}^{\dot{1}} & \lambda^2 \tilde{\lambda}^{\dot{2}}) \mapsto \mqty(\lambda^1 \\ \lambda^2 \\ \tilde{\lambda}^{\dot{1}} \\ \tilde{\lambda}^{\dot{2}})$. Now momentum conservation says that the sum of $p^0$ for all the particles is equal to zero (same for $p^1$, $p^2$, and $p^3$), and this is where we get $\sum_j\lambda_\alpha^j\tilde{\lambda}_\alpha^j=0$ from (note: we have switched the positions of the $\alpha, \dot{\alpha}$ to subscript to make the matrix below easier to read). And now that we have written it in terms of vectors, for all $n$ particles, we have the matrix displaying the kinematic information:

$\mqty(\lambda_1^1 & \lambda_1^2 & \cdots & \lambda_1^{n-1} & \lambda_1^n \\
\lambda_2^1 & \lambda_2^2 & \cdots & \lambda_2^{n-1} & \lambda_2^n \\
\tilde{\lambda}_{\dot{1}}^1 & \tilde{\lambda}_{\dot{1}}^2 & \cdots & \tilde{\lambda}_{\dot{1}}^{n-1} & \tilde{\lambda}_{\dot{1}}^n \\
\tilde{\lambda}_{\dot{2}}^1 & \tilde{\lambda}_{\dot{2}}^2 & \cdots & \tilde{\lambda}_{\dot{2}}^{n-1} & \tilde{\lambda}_{\dot{2}}^n )$
\\

Thus, the momentum conservation $\sum_j\lambda_\alpha^j\tilde{\lambda}_\alpha^j=0$ is saying that the dot product of any of the first 2 rows with any of the second 2 rows is equal to zero. In other words, the first two rows (the row vectors $\lambda_1^j$ and $\lambda_2^j$) lie in a plane orthogonal to the second two rows (the row vectors $\tilde{\lambda}_{\dot{1}}^j$ and $\tilde{\lambda}_{\dot{2}}^j$). Thus, we can apply $GL(2)$ to either the first 2 rows or the last 2 rows, and this will cause the vectors to stay in the same plane, so if all we wanted to do was to satisfy momentum conservation, we could consider this kinematic data as two elements in $Gr(2,n)$ (the Grassmannian of $2$-planes in $n$ dimensions). In fact, technically it doesn't even have to be GL(2), that is, the transformation does not have to be invertible - for example, in 3 spatial dimensions ($n=3$), one of the planes must be degenerate (at most spans a line); in this case, say $\lambda_1^j$ and $\lambda_2^j$, $j=\{1,2,3\}$, is degenerate, then if we pick the line that $\lambda_1^j$ and $\lambda_2^j$ lie upon, that is, we fix $\lambda_1^j$ or $\lambda_2^j$, then this also uniquely fixes the 2-plane that $\tilde{\lambda}_{\dot{1}}^j$ and $\tilde{\lambda}_{\dot{2}}^j$ span.  
\\

We have shown that as we move the vectors around inside the planes via $GL(2)$-transformations, we continue to satisfy momentum conservation, but one might ask whether the momenta themselves change. Fortunately, they do not. Remember that the Lorentz transformations act on on the $\lambda$ and $\tilde{\lambda}$ as $SL(2)$ transformations. This means that we can move the vectors around in the plane, and we still have the same kinematic information. Thus, the only change that a $GL(2)$ transformation can do is to apply an overall scaling. However, this can be pulled off into the Little group, which can be thought of as a right torus action on the $n$ momenta. (If we have a torus $T=diag(z_1,...,z_n)$ acting on the matrix of $\lambda$ on the right, then the Little group scaling is that it must act as $diag(z_1^{-1},...,z_n^{-1})$ on the matrix of $\tilde{\lambda}$ on the right). 
\\

In fact, the ``C-matrix," which ends up parameterizing positroid varieties (which are subsets of Grassmannians as explained above) in planar N=4 Yang Mills theory, can be seen as arising from linearizing momentum conservation. In other words, rather than dealing with $\sum_{i=1}^n \lambda_i^a \tilde{\lambda}_i^{\dot{a}} =0$, we can introduce an auxiliary $k$-plane in $n$ dimensions:
\[C^\perp \cdot \lambda=0, C\cdot \tilde{\lambda}=0\]

\section{Positroids and Physics}\label{section:PositroidsandPhysics}

(1) The connection between positroids and physics was first discovered by considering on-shell diagrams with 3-particle amplitudes. Perhaps more precisely, by taking diagrams where the only interaction involved 3 particles (in N=4 SYM) and where there were no virtual particles (all particles are on-shell), a couple of important things were noticed: \\

(2) An amazing thing happens \cite{ABCGPT}: if we eliminate so-called ``bubbles" from the diagram (which just involves pulling out dlog variables), we get a \textit{reduced} on-shell diagram. If we have a single physical process in N=4 SYM and we represent the process with an on-shell diagram, nomatter how we do it (since there are multiple ways of gluing together these 3-particle amplitudes that work) or how we do the bubble reduction afterwards - all the reduced diagrams we would get are related by merge/unmerge and square moves that do not affect the physics. It was found that these reduced diagrams are in bijection with the plabic diagrams of \cite{Post}, which equivalently can be labeled by a bounded affine permutation (explained in the sections above; note, this is saying that even if we apply a merge/unmerge or square move we still get the same permutation). Furthermore, for a given number of external particles, there is a finite set of reduced on-shell diagrams (unlike Feynman diagrams, for which we can get an arbitrarily large number). Indeed, for a given reduced diagram, if we apply all the merge/unmerge moves and square moves we can, the set is finite since there are a finite number of vertices at which we can do the merge/unmerge and a finite number of squares for which we can apply the square moves.
\\

(3) \cite{Post} had discovered these diagrams and permutations while dealing with Grassmannians, where he found that these diagrams/permutations form a stratification of the ``positive" (nonnegative) Grassmannian (note: a priori, we are considering the full (complex) Grassmannian, so there could be negative or complex coordinates for our cell - but the positivity is saying that there exists some coordinate transformation such that then these diagrams/permutations are all faces of the Grassmannian restricted to its positive points.) Thus, the physics of on-shell diagrams in N=4 SYM correspond to ways of cutting up Grassmannians. In particular, the BCFW decomposition is given naturally by on-shell diagrams, so BCFW corresponds to the positroid stratification! In particular, the BCFW decomposition is a triangulation of the amplitude (after all, we think of the amplitude as producing a real positive number, and the BCFW decomposition turns it into a sum of smaller pieces); so the BCFW decomposition is given by a collection of positroid cells of the same dimension (we can determine the dimension of a positroid cell by considering the number of inversions, or we can use the formula $dim(C_\sigma)=(\sum_{a=1}^n r[a;\sigma(a)])-k^2$ of \cite{ABCGPT}). \\
Additionally, the boundary structure is important - it gives information on how the cells of the decomposition fit together. The positive Grassmannian has a single ``top cell" and every other cell is of lower dimension: that is, other than one cell, all other cells in the positive Grassmannian are lower dimension - they are all boundaries of this top cell. This makes us think of simplices or of convex polytopes: we have a single top-dimensional region, and a bunch of lower-dimensional boundaries. Indeed, this leads us to the idea of a positive geometry (a geometric space that has boundaries of all codimensions) of which the jewel representative is the amplituhedron is described next. In fact, this is the key idea of a positive geometry, that you can reach each ``facet" by starting from the top with a single object and taking boundaries (or resides). In other words, every facet is a (piece of the) boundary. When we take a sum of a triangulation like the BCFW decomposition, the spurious poles are along common boundaries of multiple pieces, so they cancel out. 
\\

(4) Now, each plabic diagram (or cell of the positive Grassmannian) can be associated with a differential form. We can think of this as a generalization of the canonical form of the simplices. We can quickly recall how we form these, from just the condition that the form has logarithmic singularities (basically, simple poles in the denominator) on all boundaries. So for a standard $n-1$ dimensional simplex given by the convex hull of the $n$ unit basis vectors in $\R^n$, this would be $\Omega_{\Delta^{n-1}}=\frac{dx_1 \cdots dx_n}{(x_1+x_2+\cdots +x_{n-1}-1)(x_1+x_2+\cdots+x_{n-2}+x_n -1)\cdots(x_2+x_3+\cdots +x_n-1)}$. If we go to the facet where $x_n=0$, which we can denote $\Delta^{n-1}_{n=0}$, and we take the residue along this pole we get $\Omega_{\Delta^{n-1}_{n=0}}=lim_{x_n\to 0} \frac{(x_1+x_2+\cdots +x_{n-1}-1)}{dx_n} \frac{dx_1 \cdots dx_n}{(x_1+x_2+\cdots +x_{n-1}-1)(x_1+x_2+\cdots+x_{n-2}+x_n -1)\cdots(x_2+x_3+\cdots +x_n-1)}=$ \\
\noindent $\frac{dx_1 \cdots dx_{n-1}}{(x_1+x_2+\cdots+x_{n-2} -1)\cdots(x_2+x_3+\cdots +x_{n-1}-1)}$, which has singularities exactly on this facet, so is the canonical form for that facet. Note that only for the case of the simplex (or ``simplex-like") are there no zeros (constant numerator), because if there are more facets (simplices have the minimal number of facets for a polytope) then the planes that define them intersect away from the polytope, so there would be additional poles in the denominator that would need to be cancelled. 
\\

(5) Thus, we have to sum up a number of these diagrams to get the amplitude. However, we do not know which combinations work (to get the full amplitude) a priori - rather, to get the right combinations, we need to apply locality and unitarity; this was done by knowing that conditions locality and unitarity place on the singularities of the full integrand and imposing these conditions - locality says that we have the right sort of poles, and unitarity says that the amplitude factorizes on the poles. This was solved for planar N=4 Super Yang-Mills by generalizing the positive Grassmannian with the amplituhedron map, which is a map from a positive Grassmannian to a Grassmannian, where the map is itself an element of a positive Grassmannian given by the external kinematical data. The ``amplituhedron" is the shape that is the image of the positive Grassmannian under this map of the kinematical data. Specifically, the matrix $Z\in Mat^{>0}_{k+m,n}$, the real $k+m \times n$ ($k+m \leq n$) matrices, is determined by the kinematical data and it gives a map 
\[\tilde{Z}:Gr_{k,n}^{\geq 0} \to Gr_{k,k+m}, \tilde{Z}(\expval{v_1,...,v_k}):=\expval{Z(v_1),...,Z(v_k)}\]
 where the angle brackets denote the span (the plane) of the vectors, and the tree amplituhedron $A_{n,k,m}(Z)$ is the image $\tilde{Z}(Gr_{k,n}^{\geq 0})$ inside $Gr_{k,k+m}$. Thus, each scattering experiment gives us a different amplituhedron.
\\

\section{Kinematic Space}

We are used to working with the 4-momentum $p^\mu=(p^0,p^1,p^2,p^3)$. For spinor-helicity, let's recall the basic notions.

We can turn the 4-momentum into a bispinor (still not entirely sure what this word means) using:
$p^{\alpha\dot{\alpha}}=\sigma_\mu^{\alpha\dot{\alpha}}p^\mu=\mqty(p^0-p^3 & -p^1+ip^2 \\ -p^1-ip^2 & p^0+p^3)$. Clearly, we can recover $p^\mu$ if we have the matrix, so we have a one-to-one equivalence between these two notations.

Now, if we take the determinant of this matrix, we get $p_0^2-p_1^2-p_2^2-p_3^2=p_\mu^2=m^2$, so if our particles are massless, this means that the determinant is zero. When a 2-by-2 determinant is zero, we can split the determinant into an outer product of two vectors. We call the two vectors $\lambda^\alpha$ and $\tilde{\lambda}^{\dot{\alpha}}$, and we can write them as:
$\lambda^\alpha=\frac{z}{\sqrt{p^0-p^3}}\mqty(p^0-p^3 \\ -p^1-ip^2)$ and $\tilde{\lambda}^{\dot{\alpha}}=\frac{z^{-1}}{\sqrt{p^0-p^3}}\mqty(p^0-p^3 & -p^1+ip^2)$. It's then easy to check that the outer product $\lambda^\alpha \tilde{\lambda}^{\dot{\alpha}}$ equals $p^{\alpha\dot{\alpha}}$ given above. 

Now notice the presence of the peculiar $z$ in the equations. This is the Little group rescaling: these are exactly the transformations that leave the momentum fixed. [Note: the Little group is a subgroup of the Lorentz transformations - in other words, it acts on the coordinates (on the Lorentz index) of anything in our formulae (amplitudes, S-matrix elements, Lagrangians, etc.) with a 4-vector index; this includes the 4-momentum itself, although this vector cannot change.] In fact, this Little group action is basically equivalent to the masslessness (the determinant being zero). These conditions all cut down on the degrees of freedom by one. More specifically:

We have a 4-momentum $p^\mu=(p^0,p^1,p^2,p^3)$. Masslessness says that $p_\mu^2=0$, so $p_0=\sqrt{p_1^2+p_2^2+p_3^2}$, so $p^\mu=(\sqrt{p_1^2+p_2^2+p_3^2},p^1,p^2,p^3)$, so we have three degrees of freedom. 
Similarly, if we have a matrix, $\mqty(a & b\\ c& d)$ then the zero determinant condition says $ad-bc=0$, so if $d\neq0$ we can write this as $a=\frac{bc}{d}$ or $\mqty(\frac{bc}{d} & b \\ c& d)$ so once again we have 3 degrees of freedom. We said that the zero determinant condition is equivalent to being able to write the matrix as an outer product $\mqty(x \\ y)\mqty(z & w)=\mqty(xz & xw \\ yz& yw)$. Here it's easy to see that if we call $xw=b$, $yz=c$, $yw=d$, then (if $d\neq 0$), $xz=\frac{bc}{d}$ so we get the same relation as before, where only 3 of the entries in the matrix are free. Sure, someone might say, once you start multiplying these numbers together you might get relations, but aren't $x,y,z,w$ independent variables? That depends on what one means by "independent" - the fact is that sure, we could make $x,y,z,w$ whatever we want, but the issue is that it's the matrix (one gets after taking the outer product) that corresponds to the momentum (the actual physical information we have), so two different sets of $x,y,z,w$ that map to the same matrix give the same physical information so are equivalent. And we see [although I don't have a proof] that this happens exactly when two different sets are related by $\mqty(x \\y )\mapsto t\mqty(x\\y), \mqty(z& w)\mapsto t^{-1}\mqty(z & w)$. And if $x\neq0$, we can choose $t=x^{-1}$, so that we get $\mqty(1 \\ x^{-1}y)\mqty( xz & xw)$; relabelling: $\mqty(1  \\\tilde{y})\mqty(\tilde{z} & \tilde{w})$, we see that we indeed have just 3 true degrees of freedom in the $\lambda, \tilde{\lambda}$ that corresponds to a given momentum. 
\\

Spinor helicity is one way to represent the momenta, the kinematical information. There are a couple of other important ways of representing the same information. 

One of them is by momentum twistors. This way is nice because it makes momentum conservation manifest. If we give the momenta of our problem/amplitude calculation a (color) ordering, say $p_1,...,p_n$, then we can imagine a "dual space" where we just take the vector $p_1$ and place its start at the origin and it's arrow/head at the value of $p_1$ so literally $p_1$ is just $p_1$. Let's write the coordinates of $p_1$ as $y_1^\mu$. But now instead of also putting $p_2$ starting at the origin (which would give us $p_2$), let's translate the vector $p_2$ so that it starts at the head of $p_1$ (so it looks like we're adding $p_2$ to $p_1$ in our vector-addition from junior high math). Let's call this result $y_2^\mu$ (so $y_2^\mu$ is the coordinates for $p_1+p_2$). Let's continue this, placing each additional momentum so taht it starts at the head of the previous one. Now, since we have momentum conservation, the total sum of momentum is zero, so this is equivalent to saying that the last momentum $p_n$ comes back to the origin where $p_1$ started. That is, we get a polygon. (This will come in useful later, because we can switch to the kinematical information of consecutive sums of momentum denoted $x_{ij}=p_i+p_{i+1}+\cdots+p_{j-1}+p_j$ and these will be all lines between vertices in the polygon: then we can do things like consider maximally noncrossing sets which triangulate the polygon, producing for us the associahedron, and all that stuff with biadjoint $\phi^3$ theory and even consider crossing chords which gives us the recent research with causal diamonds etc.). (There are also cases where the kinematical information is presented as the Mandelstam invariants) Since these momenta are all null, if we tried to draw it in 2D (say we restrict to two dimensions $x,t$), then all the lines should be at 45 degrees. Our dual space is really the space of the $y^\mu$. Now, let's write the $y^\mu$ in bispinor notation as above (taking the dot product with the Pauli sigma matrices) which gives us $y_{\alpha \dot{\alpha}}$. Now, for consecutive $i$, let's choose $1$ and $2$ for simplicity, we remember that $y_2^\mu-y_1^\mu$ gave us the null momentum $p^\mu$. Therefore, $(y_2^\mu-y_1^\mu)^2=0$ or equivalently $det(y^2_{\alpha\dot{\alpha}}-y^1_{\alpha\dot{\alpha}})=0$. In particular, by the massless Weyl equations, we get $(y^2_{\alpha\dot{\alpha}}-y^1_{\alpha\dot{\alpha}})\lambda^\alpha=0$ (If we had formed the bispinor by using $\ov{\sigma}$ instead, we would get the same equation but with $\tilde{\lambda}$). Pulling one term to the other side, we see that $y^2_{\alpha\dot{\alpha}}\lambda^\alpha=y^1_{\alpha\dot{\alpha}}\lambda^\alpha$. We define $\mu_{\dot{\alpha}}:=y_{\alpha\dot{\alpha}}\lambda^\alpha$. And these are our new momentum twistor variables: $\mqty(\lambda^\alpha \\ \mu_{\dot{\alpha}})$. 
\\

Now, let's note a few things. First, a very important fact. If we scale by $t$: $t\mqty(\lambda^\alpha \\ \mu_{\dot{\alpha}})$, this defines the same exact equation: $(t\mu_{\dot{\alpha}})=y_{\alpha\dot{\alpha}}(t\lambda^\alpha)$. In particular, the Little group scales $\mu_{\dot{\alpha}}$ the exact same way that it scales $\lambda^\alpha$. So, first of all, this means that we can consider our new momentum twistor to be an element of projective space $\PP^3$, so we've dealt with the Little group by projectivizing our 4-component object. Secondly, if we consider all our momenta: $\mqty(\lambda^\alpha_1 & \lambda^\alpha_2 & \cdots & \lambda^\alpha_n \\  \mu_{\dot{\alpha}}^1 & \mu_{\dot{\alpha}}^2 & \cdots & \mu_{\dot{\alpha}}^n)$, then saying that each columns is in projective space is the same as taking a 4-by-n matrix and modding out by the action of a torus acting on the right. This is why such a setup is very close to what has been called the ``configuration space." 
\\

Furthermore, we just saw that if two $y$ were lightlike separated, then they define a null ray (or line). Let's be more precise. Let's take $y_{\alpha\dot{\alpha}}\lambda^\alpha=\mu_{\dot{\alpha}}$. Suppose we are given $y_{\alpha\dot{\alpha}}$; this is a 2-by-2 matrix, so we get two equations (two constraints on $\lambda$ and $\mu$) from this equality. However, the momentum twistor has 3 degrees of freedom, so there's one degree of freedom left that we haven't taken care of. Therefore, points of the dual space map to lines in the momentum twistor space. Alright, so let's pick two points $y_{\alpha\dot{\alpha}}^1$ and $y_{\alpha\dot{\alpha}}^2$ that are lightlike separated. Then the massless Weyl equation $(y_{\alpha\dot{\alpha}}^2-y_{\alpha\dot{\alpha}}^1)\lambda^\alpha=0$ fixes $\lambda^\alpha$ completely, so we can just plug and chug into $y_{\alpha\dot{\alpha}}\lambda^\alpha=\mu_{\dot{\alpha}}$ using either $y^1$ or $y^2$ to get $\mu$. So a null line in the dual space corresponds to a point in momentum twistor space. And we can go the other way. We just saw that a null line in dual space gives us a point in momentum twistor space. If we have a point in momentum twistor space, then we get a null line: look at $y_{\alpha\dot{\alpha}}\lambda^\alpha=\mu_{\dot{\alpha}}$ where $\lambda$ and $\mu$ are fixed - pick a specific value of $y$ that we can calle $\tilde{y}$ satisfying this equation; then any other $y$ that satisfies this equation also satisfies: $(\tilde{y}_{\alpha\dot{\alpha}}-y_{\alpha\dot{\alpha}})\lambda^\alpha=\mu_{\dot{\alpha}}-\mu_{\dot{\alpha}}=0$. Thus $(\tilde{y}_{\alpha\dot{\alpha}}-y_{\alpha\dot{\alpha}})^2=0$ and it is the null momenta producing for us us Weyl's equation for our given $\lambda^\alpha$ so every $y$ for our given momentum twistor $\mqty(\lambda \\ \mu)$ lies along this null line. Similarly, just as a point in the dual space produced a line (not necessairly null) in momentum twistor space, a line in momentum twistor space produces a point in dual space. Specifically, we have 
$y_{\alpha\dot{\alpha}}\lambda_1^\alpha=\mu^1_{\dot{\alpha}}$ and $y_{\alpha\dot{\alpha}}\lambda_2^\alpha=\mu^2_{\dot{\alpha}}$ which gives us 4 equations for the 4 unknowns in $y_{\alpha\dot{\alpha}}$ so we can explicitly solve for the point $y$ and we get: 
\[y_{\alpha\dot{\alpha}}=\frac{\lambda^\alpha_1\mu^2_{\dot{\alpha}}-\lambda^\alpha_2\mu^1_{\dot{\alpha}}}{\langle \lambda_1 \lambda_2\rangle}\]

Of course, our dual space $y_{\alpha\dot{\alpha}}$ acted as sort of a ``messenger" between the spinor-helicity space $\mqty(\lambda\\ \tilde{\lambda})$ and the momentum twistor space $\mqty(\lambda \\ \mu)$. We could eliminate $y$ and pass directly between the two using the relation: 
\[\tilde{\lambda}_i=\frac{\expval{i+,i}\mu_{i-1}+\expval{i,i-1}\mu_{i+1}+\expval{i-1,i+1}\mu_i}{\expval{i-1,i}\expval{i,i+1}}\]
\\

Finally, we note a third way of representing kinematical information. These are said to be twistors in twistor space (as opposed to momentum twistors in momentum twistor space for the previous type). I believe these were the ones used by Witten in his perturbative string theory paper. Twistors are represented by 4-component objects $(\tilde{\mu}_i, \tilde{\lambda}_i)$ where we take the amplitude expressed in $\lambda, \tilde{\lambda}$ and Fourier transform the part with $\lambda$: $\int d^2\lambda_j exp(-i\tilde{\mu}^a_j\lambda_{ja})\cdot$. We note the important fact that under the Little group scaling, $\tilde{\mu}_i\mapsto t_i^{-1}\tilde{\mu}_i$, so just like in the case above with momentum twistors, all the components of the twistor 4-vectors scale the same way (here, as $z_{-1}$), so they are also elements in projective space $\PP^3$. 
This Fourier transform aspect seems a bit weird, but we can actually see that it is dual to the momentum twistor case. Consider solutions of the wave equation $\Box \phi(x)=0$; in Fourier space this is $p^2\tilde{\phi}(p)=0$. So we can write the solution roughly as $\phi(x)=\int d^4p\delta(p^2)dxp(ipx)f(p)=\int \frac{d^3p}{2\abs{\vec{p}}}exp(i(-\abs{p}t-\vec{p}\cdot\vec{x}))f(\vec{p}))$, where we have chosen a frame and used the delta function to eliminate $dp_0$. Note that the condition that $p^2=0$ is automatic if we write the momentum in terms of $\lambda, \tilde{\lambda}$ since this implies the particles are massless. So we write $\phi(x)=\int \frac{d^2\lambda d^2 \tilde{\lambda}}{GL(1)}exp(i\lambda\tilde{\lambda}x)f(\lambda, \tilde{\lambda})$, where the GL(1) we are modding out performs the Little group scaling on $\lambda \mapsto z\lambda, \tilde{\lambda} \mapsto z^{-1} \tilde{\lambda}$. Now we can perform the Fourier transform on $\tilde{\lambda}$. But notice! The thing in front of $\tilde{\lambda}$ in the exponent is $x\lambda$. So if we call the new variable $\mu:=x\lambda$, then we get: $\int \frac{d^2\lambda}{GL(1)}\tilde{f}(\lambda,\mu=x\lambda)=\int \frac{d^2\lambda}{GL(1)}\tilde{f}(Z)|_{\mu=x\lambda}$, then we get nothing other than the momentum twistors (here called "$Z$" following Arkani-Hamed's notation). I believe this is called a Penrose transform. So we now see that analogy with the twistors: we instead perform the half-Fourier transform with $-\tilde{\lambda}x:=\tilde{\mu}$, so the exponent becomes $exp(-i\lambda\tilde{\mu})$ and we integrate with respect ot $\lambda$ to get $\int\frac{d^2\tilde{\lambda}}{GL(1)}f(W)$, $W$ denoting the twistors. 
\\

It's called spinor-helicity because you're representing the information with $p\rangle=\lambda^\alpha=\psi_L$ and $p]=\tilde{\lambda}^{\dot{\alpha}}=\psi_R$ (it could be the other way around for what the angle and square brackets mean, check the conventions of the book/paper). These are clearly Weyl spinors (two-components spinors in the (0,1/2) and (1/2,0) representations of the Lorentz group), and because the particles are massless, the eigenstates of the Weyl equations are the same as the eigenstates of the helicity operator. So positive helicity exactly means that the particle is described by a spinor $\psi_R$ in the (0,1/2) representation. 
Here's where it's not completely clear: $\psi_R$ and $\psi_L$ are describing the particle - the spinors are "quantum states" that we insert into a Lagrangian to get the equations of motion. But now we are writing the \textit{momentum} as spinors in the the left-handed and right-handed helicity/Lorentz representation spaces. Furthermore, the correspondence is not 1-1 because in the case of momentum, it is actually given by a matrix, but we have decomposed it somewhat artificially (using that the determinant is zero in the massless case) into two pieces, both of which are defined only up to a scalar. 
Note that by convention, we usually take all momenta to be incoming in the on-shell stuff (occasionally conventions are that they are all outgoing). This means that for the calculation, since we are flipping the direction of momenta, we also need to flip the helicities (i.e. $\lambda \leftrightarrow \tilde{\lambda}$). Note that this means that things could get quite confusing: \cite{Schw} takes all momenta incoming while \cite{EH} take all momenta outgoing, both use angle brackets for negative helicity, but \cite{Schw} writes the spinors $(\psi_L, \psi_R)$ whereas \cite{EH} write it $(\psi_R, \psi_L)$. The net effect is that both have if you flip the helicities of EH to make momenta incoming to compare the two, you get the negative helicity spinors in both being on the top of the Dirac spinor - although the brackets are still flipped: square brackets on top for \cite{EH} vs. angle brackets on top for \cite{Schw}. 
\\

We do a key calculation show show why the left and right-handed (helicity) spinors are related to a decomposition of the momentum into an outer product.
 \\
 
 Let's start with the helicity operator for general spin. As usual for helicity, we assume the particle is massless or nearly massless. \cite{Schw}writes it as $\frac{\vec{S}\cdot \vec{p}}{s\abs{\vec{p}}}\psi=\pm \psi$. Let's clear denominators: $\vec{S}\cdot \vec{p} \psi=\pm s\abs{\vec{p}}\psi$. In the case of spin one-half, $\vec{S}=\frac{\sigma}{2}$ and $s=\frac{1}{2}$, so the equation is $\vec{\sigma}\cdot \vec{p} \psi=\pm \abs{\vec{p}}\psi$. \cite{Peskin} writes the helicity operator as $h=\hat{p}\cdot \vec{S}$, saying that it is projection of spin along the \textit{direction} of momentum (hence the unit vector), which is the same as \cite{Schw}, just keeping the $\abs{\vec{p}}$ in the denominator. For the case of spin one-half, let's write out the equation explicitly: 
 \begin{equation}
 \vec{\sigma}\cdot \vec{p}\psi=\mqty(p_3 & p_1-ip_2 \\ p_1+ip_2 & -p_3)\psi=\pm p_0 \psi
 \end{equation} 
 In the last step, we uses the masslessness: $p_0^2+\abs{\vec{p}}^2=0$. 
 Now, we know that the solutions $\psi$ are the same as the Weyl equations (the massless Dirac equation; sometimes people say the "massless Weyl equations" as a reminder, but "Weyl equations" already means massless). They are written in various forms:
 \begin{gather}
 i(\del_0-\vec{\sigma}\cdot \vec{\nabla})\psi_L=0, i(\del_0+\vec{\sigma}\cdot \vec{\nabla})\psi_R=0 \\
 \mqty(0 & i\sigma\cdot \del \\ i\ov{\sigma}\cdot \del & 0)\mqty(\psi_L \\ \psi_R)=0 \\
 i\ov{\sigma}\cdot \del \psi_L=0, i\sigma\cdot \del\psi_R=0 \\
 \sigma^\mu p_\mu\psi_R=(E-\vec{\sigma}\cdot \vec{p})\psi_R=0 , \ov{\sigma}^\mu p_\mu \psi_L=(E+\vec{\sigma}\cdot \vec{p})\psi_L=0 \\
 p^{\dot{a}b}p]_b=0, p_{a\dot{b}}p\rangle^{\dot{b}}=0, [p^b p_{b\dot{a}}=0, \langle p_{\dot{b}}p^{\dot{b}a}=0
 \end{gather}
 The first three are from \cite{Peskin}, the last two from \cite{Schw} and \cite{EH}, respectively. 
 The one from \cite{Schw}: $(E-\vec{\sigma}\cdot \vec{p})\psi_R=0$ and $(E+\vec{\sigma}\cdot \vec{p})\psi_L=0$, noting that $E=p_0$ are the closest in form to the helicity equations above, since they clearly show that the spinors are solutions to the helicity operator with eigenvalue $\pm E=\pm p_0$.
 So maybe the best way to think about it, is (for spin one-half), we have the helicity equation $\mqty(p_3 & p_1-ip_2 \\ p_1+ip_2 & -p_3)\psi=\pm p_0 \psi$, and we then move the $E$ over to the other side to get: $\mqty(p_3 \pm p_0 & p_1-ip_2 \\ p_1+ip_2 & -p_3 \pm p_0)\psi=0$, and if we make the substitution $p_\mu=i\del_\mu$ then we get the Weyl equations exactly. 
 You can easily check (or just by inspection), that these equations are all equivalent. 
 \\
 
 Here's the key thing: these helicity eigenstates/solutions to the Weyl equations are exactly the decomposition of the bispinor form of momenta into an outer product.
 \\
 
 The bispinor $(\frac{1}{2},\frac{1}{2})$ form of momenta is written in (27.12) of \cite{Schw}:
\begin{equation}
p^{\alpha\dot{\alpha}}=\sigma^{\alpha\dot{\alpha}}_\mu p^\mu=\mqty(p^0-p^3 & -p^1+ip^2 \\ -p^1-ip^2 & p^0+p^3)
\end{equation}
It's determinant is $(p^0)^2-\abs{\vec{p}}^2=0$ by masslessness, and a 2-by-2 matrix with zero determinant can be decomposed into an outer product of two 2-vectors. \cite{Schw} gives the decomposition:
 $p^{\alpha\dot{\alpha}}=\lambda^\alpha \tilde{\lambda}^{\dot{\alpha}}$, where
\begin{gather*}
 \lambda^\alpha=\frac{z}{\sqrt{p^0-p^3}}\mqty(p^0-p^3 \\ -p^1-ip^2) \\
  \tilde{\lambda}^{\dot{\alpha}}=\frac{z^{-1}}{\sqrt{p^0-p^3}}\mqty(p^0-p^3 & -p^1+ip^2)
 \end{gather*}
 Let's quickly check that this works:
\begin{equation}
\lambda^\alpha \tilde{\lambda}^{\dot{\alpha}}=\frac{1}{p^0-p^3}\mqty((p^0-p^3)^2 & (p^0-p^3)(-p^1+ip^2) \\ (-p^1-ip^2)(p^0-p^3) & (p^1)^2 +(p^2)^2)=\mqty(p^0-p^3 & -p^1+ip^2 \\ -p^1-ip^2 & \frac{(p^1)^2 +(p^2)^2}{p^0-p^3})
\end{equation}
where the last entry is $p^0+p^3$ from masslessness: $(p^1)^2 +(p^2)^2=(p^0)^2-(p^3)^2=(p^0-p^3)(p^0+p^3)$. 
\\

So let's finally show that this decomposition of momenta into these two outer product vectors equals the solutions to the helicity operator/Weyl equations. To match \cite{Schw}, let's start with all raised indices. The helicity operator/Weyl operator is: 
\[ \mqty(p^3 \pm p^0 & p^1-ip^2 \\ p^1+ip^2 & -p^3 \pm p^0)\mqty(x \\ y)=0\] 
Let's calculate for positive $+p^0$ (which actually means the eigenvalue is negative so left-handed):
$(p^1+ip^2)x+(-p^3+p^0)y=0$, so the solution for the eigenvector is $\mqty(\frac{p^3-p^0}{-(p^1+ip^2)} \\ 1)$. Note that if we had used the first row, we would have gotten $\mqty(\frac{-p^1-ip^2}{-(p^3+p^0)} \\ 1)$ and we can check that they are the same thing: 
\[\frac{p^3-p^0}{-(p^1+ip^2)}\cdot 1 \cdot 1=\frac{p^3-p^0}{-(p^1+ip^2)}\frac{p^3+p^0}{p^3+p^0}\frac{-(-p^1-ip^2)}{-(-p^1-ip^2)}=\frac{(p^3)^2-(p^0)^2}{p^3+p^0}\frac{-(-p^1-ip^2)}{-(p^1)^2-(p^2)^2}=\frac{-p^1-ip^2}{-(p^3+p^0)}\]
since again $\frac{(p^3)^2-(p^0)^2}{-(p^1)^2-(p^2)^2}=1$ by masslessness. It's an eigenvalue, so we can scale to get any other eigenvector in the eigenspace, and we will choose: $\mqty(p^3-p^0 \\ -(p^1+ip^2))$. If we let $z=\sqrt{p^0-p^3}$, then this exactly matches \cite{Schw}'s $\lambda^\alpha=\mqty(p^3-p^0 \\ -p^1-ip^2)$.
\\

Similarly, for the negative (right-handed) eigenvector, we get the same thing, except the sign of $p^0$ is flipped: $\mqty(p^3+p^0 \\ p^1+ip^2)$. We need to get this to match \cite{Schw}'s: 
\[ \tilde{\lambda}^{\dot{\alpha}}=\frac{z^{-1}}{\sqrt{p^0-p^3}}\mqty(p^0-p^3 & -p^1+ip^2)\] 
To do this, recall the notation: $\psi=\mqty(\lambda^\alpha \\ \tilde{\lambda}_{\dot{\alpha}})$. That is, $\psi_R=\tilde{\lambda}_{\dot{\alpha}}$. What we need to do is to take $\tilde{\lambda}^{\dot{\alpha}}$, lower the index, and conjugate it. Although it would be better to go the other way, since we have already fixed $z=\sqrt{p^0-p^3}$, so we already have: 
\[\tilde{\lambda}^{\dot{\alpha}}=\frac{1}{p^0-p^3}\mqty(p^0-p^3 & -p^1+ip^2)=\mqty(
1 & \frac{-p^1+ip^2}{p^0-p^3})\] 
So the right thing to do would be to pick the right vector in the eigenspace for $\psi_R$ and to instead raise its index, and conjugate it, and see that it equals this $\tilde{\lambda}^{\dot{\alpha}}$. So let's go back to $\mqty(p^3+p^0 \\ p^1+ip^2)$, which should actually be $\mqty(p_3+p_0 \\ p_1+ip_2)$ since $\psi_R=\tilde{\lambda}_{\dot{\alpha}}$ has lower indices. Let's divide by $p_3+p_0$ to get another eigenvector in the eigenspace of the helicity/Weyl operator, and multiply by $-1$, which hopefully will work out: $\mqty(-1 \\ -\frac{p_1+ip_2}{p_3+p_0})$. Raising indices (recall we are using (+,-,-,-) index), we get: 
\[ \mqty(-1 \\ -\frac{-p^1-ip^2}{-p^3+p^0})=\mqty(-1 \\ \frac{p^1+ip^2}{-p^3+p^0})=\mqty(1 \\ \frac{p^1+ip^2}{p^3-p^0})\] 
Next we conjugate to get: 
\[\mqty(1 & \frac{p^1-ip^2}{p^3-p^0})=\mqty(1 & \frac{-p^1+ip^2}{p^0-p^3})\] 
And this is exactly what we wanted. I think this was confusing for me at first because it feel like there is an abuse of notation: the process of lowering the index on $\tilde{\lambda}^{\dot{\alpha}}$ involved two steps: both the conventional lowering with the Lorentz metric, and the complex conjugation. Even if we tried to lower with the $\ep$ symbol instead, we would still have the wrong type of vector and need to either conjugate or transpose. Let's try this: for the negative one (corresponding to right-handed), since above we simply took the left-handed solutions and switched the sign of $p^0$, let's look at the solutions explicitly, which are $\mqty(p_1-ip_2 \\ -(p_3-p_0))$ or equivalently calculating the eigenvector with the second row: $\mqty(-p_3-p_0 \\ -(p_1+ip_2))$ (knowing that $\psi_R=\tilde{\lambda}_{\dot{\alpha}}$, we've gone ahead and lowered the indices already). Since we will be applying $\ep$ which flips the two entires, and we want to get it into a form like 
\[\tilde{\lambda}^{\dot{\alpha}}=\frac{z^{-1}}{\sqrt{p^0-p^3}}\mqty(p^0-p^3 & -p^1+ip^2)\],
we will use the first one: $\tilde{\lambda}_{\dot{\alpha}}=\mqty(p_1-ip_2 \\ -(p_3-p_0))$. Now using $\ep^{\dot{\alpha}\dot{\beta}}=\mqty(0 & 1 \\ -1 & 0)$, we get: 
\[\tilde{\lambda}^{\dot{\alpha}}=\ep^{\dot{\alpha}\dot{\beta}} \tilde{\lambda}_{\dot{\beta}}=\mqty(0 & 1 \\ -1 & 0) \mqty(p_1-ip_2 \\ -(p_3-p_0)) = \mqty( -(p^3-p^0) \\ -(p^1-ip^2)) =  \mqty( p^0 - p^3 \\ -p^1+ip^2)\] 
So all we have to do is transpose, and we get something of the form: $\tilde{\lambda}^{\dot{\alpha}}=\frac{z^{-1}}{\sqrt{p^0-p^3}}\mqty(p^0-p^3 & -p^1+ip^2)$. This was easier, so maybe this way of doing it is better (of course, there's still that extra transpose or conjugate we have to do). 
\\

This ends the key calculation.

\section{(Differential) Forms on Projective Spaces}

One important change for these integrands (canonical forms) we are defining on these positive-geometry spaces is that they are differential forms on projective or Grassmannian spaces. Thus, they look different than the usual forms in textbooks. This is the reason for the unusual notation: $\Omega=\frac{d^{nk}Y}{VolGL(k)}$. In the case of projective space $k=1$, and this case, say for $n=3$ or $RP^3$, we have coordinates $[x_0:x_1:x_2:x_3]$; but if we write our usual form $dx_0\wedge dx_1\wedge dx_2 \wedge dx_3$, this means that we are integrating 4 variables over a 3-dimensional space, so we are integrating an extra variable, so we are basically carrying along an extra copy of infinity at each point (so the integral will diverge) - what we need to do is to cancel out this extra variable and that is what the VolGL(k) in the denominator is denoting. We could write it more invariantly as $\Omega = \frac{1}{n!}\frac{<x d^m x>}{x_0\cdots x_n}$. The numerator literally means, for $n=3$, formally take the determinant $det\mqty(x_0 & x_1 & x_2 & x_3\\ dx_0 & dx_1 & dx_2 & dx_3 \\ dx_0 & dx_1 & dx_2 & dx_3 \\ dx_0 & dx_1 & dx_2 & dx_3)$, or that we have one row of the $x_i$ variables, and $n$ rows of the $dx_i$ variables. This determinant is not zero because we combine the $dx_i$ variables with the wedge product, so the antisymmetry cancels out. Thus, we end up with a symmetric sum: $\frac{1}{6}(\frac{dx_0dx_1dx_2}{x_0x_1x_2}+\frac{dx_0dx_2dx_3}{x_0x_2x_3}+\frac{dx_0dx_1dx_3}{x_0x_1x_3}+ \frac{dx_1dx_2dx_3}{x_1x_2x_3})$. It's clear how we generalize this to the Grassmannian; this will be the equivalent of $VolGL(k)$ in the denominator. For $Gr(k,n)$, we have terms in the numerator of the form $<X^1X^2\cdots X^k d^{n-k}X^i>$, which looks like $det\mqty(x_1^1 & x_2^1 & \cdots & x_n^1\\ x_1^2 & x_2^2 & \cdots & x_n^2 \\ \vdots & \vdots & \ddots & \vdots \\
x_1^k & x_2^k & \cdots & x_n^k \\ dx_1^i & dx_2^i & \cdots & dx_n^i \\ \vdots & \vdots & \ddots & \vdots \\ dx_1^i & dx_2^i & \cdots & dx_n^i )$. Here, the subscripts give the coordinates in the space, the superscripts denote which row of the matrix representative of $Gr(k,n)$ we have, and we have a copy of our matrix in $Gr(k,n)$ at the top which has $k$ rows, and another $n-k$ rows of $dx^i$ for $i$th row of the matrix representative. Thus, this determinant makes sense because it is of an $n-by-n$ matrix. We multiply together all $k$ of these: $<X^1X^2\cdots X^k d^{n-k}X^1><X^1X^2\cdots X^k d^{n-k}X^2> \cdots <X^1X^2\cdots X^k d^{n-k}X^k>$. Each form in this product has $k(n-k)$ wedges, which is exactly the dimension of the Grassmannian. Now, for the denominator, we saw that for the projective case $\Omega = \frac{1}{n!}\frac{<x d^m x>}{x_0\cdots x_n}$ that the denominator consisted of maximal minors: but in the case of $Gr(1,n)$, all the maximal minors was the same as all cyclic minors starting at each coordinate, so we're not sure a prior how to extend this. It would seem at first that we would need to take all minors. In the projective case, each term in the numerator looks like $x_i dx_0\wedge \cdots \hat{dx_i} \wedge \cdots \wedge dx_n$ (hat denotes omission), while the denominator has all of $x_0\cdots x_n$, so the numerator has a single maximal minor which is cancelled by one of the minors in the denominator, giving us: $\frac{dx_0\wedge \cdots \hat{dx_i} \wedge \cdots \wedge dx_n}{x_0\cdots \hat{x_i} \cdots x_n}$ for this one term (and we symmetrize to get all terms). For the Grassmannian case then, it would seem that we would need to analogously take all determinants because the determinant of the $n-by-n$ matrix above will have terms that involve every single maximal minor of the $k-by-n$ representative. That way, each choice of $n-k$ columns in the lower $n-k$ rows giving us the differential form has its determinant (coefficient), which is given by the determinant of the complement of those $n-k$ columns in the first $k$ columns, cancelling out with precisely one thing in the denominator. We would then need to show the important thing: invariance under GL(k) transformations. We would need to show 2 things: (1) that the terms of the same $i$ (e.g. $dx^i_1 \wedge dx^i_2 \wedge \cdots dx^i_{n-k}$) which now have an extra factor in front due to the transformation, is precisely cancelled by something in the denominator and (2) that the extra "cross-terms" we get like $dx^i_1\wedge dx^{i+1}_1\wedge\cdots$ cancel out in the full product. However, there seems to be two issues with this approach. First, if we include all minors in the denominator, then they are not independent: they satisfy Plucker relations, so we could potentially have double poles. Secondly, the denominator becomes very, very high degree very quickly: the numerator has a lot of terms added together, the denominator has a lot of terms multiplied together. So if we take the most trivial case of a GL(k) transformation, that of scaling by a constant, we would get something in the denominator that is not cancelled out by the numerator. So it seems from this analysis that we would need to choose the second option: making the denominator only consist of the cyclic minors. This definitely has a "cyclic rank" feel (think: positroids) which perhaps was one way that the authors of \cite{ABCGPT} stumbled upon this - because the textbook example of "something wrong with Feynman diagrams" is Parke-Taylor, and Parke-Taylor has a cyclic determinant sort of structure in the denominator, so investigating it would naturally lead to mathematics of subsets of the Grassmannian that have cyclic conditions: that is, positroids. 
Note: some of the minors could be zero; however, this is the same as in the projective case. A generic plane will not have any zero maximal minors. However, if we do, then we are on one of the boundaries, so would should take a residue at that pole. 

This stuff with GL(k) invariance is similar to ideas of the structure sheaf for projective varieties. If you're considering things invariant projective, then you can apply GL(k) transformations if you take the determinant, but if you're considering an affine thing, you can only do (det of) SL(k), or you can divide by something of the same degree in the denominator so that the GL(k) determinants in the numerator and denominator cancel out. 
\\

\chapter{BCFW (Britto Cachazo Feng Witten) Recursion}\label{ch:BCFW}




In this chapter, we discuss the Britto-Cachazo-Feng-Witten (BCFW) relations. We begin with some background, where we derive the BCFW relations. We then turn out attention to the BCFW bridge for on-shell diagrams. We prove a couple of statements regarding the BCFW decomposition method in \cite{ABCGPT}, wring out proofs on how the bridge decomposition algorithm works and that the formula for how the C-matrix changes upon addition of a bridge changes a permutation in the expected way. Finally in this section, we show how to put BCFW bridges on a pipe dream to build the plabic diagram for a permutation $f$; this gives us a way to associate a plabic diagram to a pipe dream.

\section{BCFW for Tree-level Amplitudes}

We can derive the BCFW shift as follows. Suppose we have a collection of momenta $p_1,...,p_n$ satisfying momentum conservation: $\sum_{i=1}^n p_i=0$. Now, let us shift one momentum $p_i$ by some arbitrary complex scalar of some other momentum $q$: $p_i \mapsto \hat{p}_i = p_i+zq$. To be precise, we can include indices so we know which quantities are vectors: $\hat{p}_i^\mu = p_i^\mu+zq^\mu$. The variable $z$ stands for an arbitrary complex number, so $z$ can vary across the complex plane. However, we have freedom to choose $q$, and we will derive the BCFW shift by imposing the two following necessary conditions: (1) momentum conservation (2) null momenta

In other words, we need to ensure that even after the shift (1) the sum of all momenta stays zero, and (2) all momenta are still null. 

To satisfy (1), it is clear that, given that we have $\sum_{i=1}^n p_i=0$, which now after the shift becomes $p_1 + \cdots + \hat{p}_i + \cdots p_n = p_1 + \cdots + p_i + zq + \cdots p_n = \sum_{i=k}^n p_k + zq = 0+zq=zq$, we must have the sum $\sum_{k=1,k\neq i}^n p_k + p_i = -zq$, which is satisfied in the easiest way by letting another external momentum $p_j$ go to $\hat{p}_j = p_j - zq$. 

To satisfy (2), we must have all $p_k \cdot p_k=0$. This was already the case for all $p_k, k\neq i,j$, but now that we have shifted $p_i$ and $p_j$, we need to ensure that they are still null as well. Thus, we first have $\hat{p}_i^2=0 \Rightarrow (p_i + zq)2 = p_i^2 + 2zp_i \cdot q+z^2q^2=0$. Since we already had $p_i^2=0$, this becomes $2zp_i \cdot q+z^2q^2=0$. Because $z$ is free to vary in the complex plane, this equation will only be satisfied for all $z$ if both $p_i \cdot q = 0$, and $q^2=0$. If we do the same thing for $\hat{p_j}$, we get $\hat{p}_j^2=0 \Rightarrow (p_j - zq)2 = p_j^2 - 2zp_j \cdot q+z^2q^2=0$, and again $p_j^2=0$ implies that this becomes $2zp_j \cdot q+z^2q^2=0$. Therefore, we need to satisfy the 3 conditions: (i) $q\cdot q =0$, (ii) $p_i \cdot q = 0$, (iii) $p_j \cdot q = 0$. 

A choice of $q$ that makes this possible will be easy if we go to spinor-helicity notation (this is possible if $q$ is null, which we must impose anyways in order to satisfy (i)). Recall from \cite{EH} (2.24) that $2p\cdot q = (p+q)^2 = \langle p q\rangle [p q]$. Rewriting (i)-(iii), we get (mod factor of $\frac{1}{2}$): (i) $\langle q q\rangle [q q]$, (ii) $\langle p_i q\rangle [p_i q]$, (iii) $\langle p_j q\rangle [p_j q]$. We already noted that the fact that we can write $q$ in spinor-helicity implies that we imposed $q\cdot q =0$, and of course writing $\langle q q\rangle =[q q]=0$ makes this obvious since these are determinants of 2-by-2 matrices where 2 columns are identical. To satisfy (ii) and (iii), we see that letting $q\rangle = p_i \rangle$ and $q] = p_j]$ (or vice-versa switching the $i$ and $j$, or equivalently angle and square brackets) makes both (ii) and (iii) vanish. 

$q= q\rangle [q$ in spinor-helicity; putting in indices: $q_{a\dot{b}}= q\rangle_a [q_{\dot{b}}$. By what we have just done, we let
$q= p_i\rangle [p_j$. Thus, writing the shift in spinor, helicity, we get:

\noindent $\hat{p}_i = p_i+zq = p_i \rangle [p_i + z p_i\rangle [p_j = p_i \rangle \big([p_i + z[p_j \big)$, and \\
$\hat{p}_j = p_j-zq = p_j \rangle [p_j - z p_i\rangle [p_j = \big(p_j \rangle - z p_i\rangle \big) [p_j$.

Therefore, we can describe the BCFW shift in spinor helicity as the shift:

\noindent $[p_i \mapsto [p_i + z [p_j = [\hat{p}_i$ (while $p_i \rangle = \hat{p}_i \rangle$ stays the same), and \\
$p_j \rangle \mapsto p_j \rangle - z p_i \rangle = \hat{p}_j \rangle$ (while $[p_j = [ \hat{p}_j$ stays the same).

(Or, in the language of $\lambda = p]$ and $\tilde{\lambda}=p\rangle$, what we've shown is that the BCFW shift is given by: $\lambda_i \to \lambda_i+z\lambda_j$ and $\tilde{\lambda}_j \to \tilde{\lambda}_j - z \tilde{\lambda}_i$)

Now we think of the amplitude $M$ as an analytic function of $z$, $M(z)$; the amplitude is $M(0)$.

\subsection{BCFW Recursion}

Let us consider a scattering process with $n$ momenta $p_i$. Let us denote the amplitude as $A_n$, which is shorthand for $A_n(p_1,...,p_n)$. Note that we are assuming some quantum field theory (some Lagrangian which specifies the interactions, couplings, etc.), so we think of the amplitude, given the theory, as a function depending on only the momenta $p_i$. 

We choose two momenta $p_i$ and $p_j$ for which we apply the BCFW shift as described in the previous section. Since $q$ is fixed by the choice of $p_i$ and $p_j$ as we already explained, the amplitude now depends on one additional variable $z$: 

$A_n(p_1,...,p_n) \mapsto A_n(p_1,...\hat{p_i},...,\hat{p_j},p_n) = A_n(p_1,...,p_n)(z)$, where the hats denote the shift.
 
 Clearly, the original amplitude $A_n(p_1,...,p_n)$ equals $A_n(p_1,...,p_n)(z=0)$. But we can also write this using an integral as:
 
 \begin{align*}
A_n(p_1,...,p_n) &= A_n(p_1,...,p_n)(0) \\
&= Res_{z=0} \frac{A_n(p_1,...,p_n)(z)}{z} \\
&= \int_C \frac{A_n(p_1,...,p_n)(z)}{z}
 \end{align*}
whre $C$ is a countour that encircles the pole at $z=0$, and no other poles. We can change the orientation of $C$ and add a minus sign, so that:

$\int_C \frac{A_n(p_1,...,p_n)(z)}{z} = \sum_{poles z_j \neq0} Res_{z=z_j} \frac{A_n(p_1,...,p_n)(z)}{z}$

At tree level, this corresponds to looking at all diagrams with an internal propagator put on shell such that $p_i$ and $p_j$ are on opposite sides. 

Perhaps the amazing thing about  BCFW is that it gives a formula for the amplitude where you only sum across channels where momenta $i$ and $j$ (which you can choose) are on opposite sides. There are many other propagators that involve $i$ and $j$ on the same side - but you get these "for free." It seems that the key thing allowing for this magic to happen is that there are no poles as $z\to \infty$; this encapsulates information about certain Feynman diagrams cancelling, and is the first hint of the underlying dual conformal symmetry. 

\subsection{BCFW from on-shell diagrams}

In the case of planar maximally supersymmetric Yang-Mills theory, the leading singularities can be depicted by what are called ``on-shell diagrams" or ``plabic diagrams" (plabic is a shortened form of ``planar bicolored"), consisting of a planar network where all nodes are colored white or black, and in fact all nodes can be made trivalent. The key ideas about the correspondence between these diagrams and positroids was already discussed in Section \ref{section:PositroidsandPhysics} of this thesis, and full details can be found in \cite{ABCGPT}. An example of such a diagram is:

\begin{figure}[htbp]
\includegraphics[scale=0.4,clip=true]{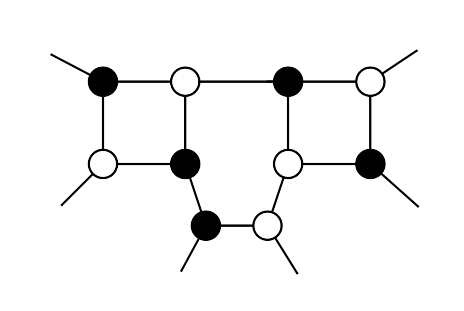}
\end{figure}

In turns out that in this formulation/language, the BCFW shift described above can be captured by the following element:

\begin{tikzpicture}
\node at (0,0) {$\bullet$};
\draw (-0.7,-0.7)--(0,0);
\draw (0,0)--(-0.7,0.7);
\draw (0,0)--(0.9,0);
\draw (1,0) circle (0.1cm);
\draw (1.07,0.07)--(1.7,0.7);
\draw (1.07,-0.07)--(1.7,-0.7);
\end{tikzpicture}

In the following, we refer to this configuration as a ``BCFW bridge."

We first prove an important fact, that any three-edge vertex must be one of two types, which we denote by the white or black circles: if we write the momenta in terms of spinor-helicity so decompose them into $\lambda_i$ and $\tilde{\lambda}_i$ (this of course assumes the momenta are all massless), then either all 3 of the $\lambda$ must be parallel, or all 3 of the $\tilde{\lambda}_i$ must be parallel.
One way of proving this comes from the algebra directly (see \cite{EH} section 2.4, equations (2.58), (2.59)).
Another way comes from thinking about the geometry: at a 3-edge vertex, we can put each $\lambda_i$ as a column vector, and becasue there are 3 momenta, this will be a 2-by-3 matrix:
$\mqty(\lambda_1^x & \lambda_2^x & \lambda_3^x \\ \lambda_1^y & \lambda_2^y & \lambda_3^y)$. 
We have the analogous thing for the $\tilde{\lambda}_i$.

Now, we also have momentum conservation: $p_1+p_2+p_3=0$. In spinor helicity, this says: $\lambda_1 \tilde{\lambda}_1 + \lambda_2 \tilde{\lambda}_2 + \lambda_3 \tilde{\lambda}_3 =0$. In other words, the dot product of $\lambda$ with $\tilde{\lambda}$ is zero; more specifically, the dot product of any of the two rows of the 2-by-3 $\lambda$ matrix above, with any of the two rows of the corresponding  2-by-3 $\tilde{\lambda}$ matrix must also equal zero. This implies that, drawing the $\lambda$ matrix as two vectors in 3-dimensional space, and also drawing the $\tilde{\lambda}$ matrix as another two vectors in same 3-dimensional space, these two pairs of vectors must be orthogonal: the plane containing the two vectors in $\lambda$ must be perpendicualr to the palne containing the two vectors in $\tilde{\lambda}$. Obviously, the space is only 3-dimensional, so one of the two planes must actually be degenerate (the two vectors in this plane must lie along the same line).

\subsection{Permutations from Bridges}

Consider the following figure of a BCFW bridge and what it does to permutations $f$ (this is (5.12) on page 52 of \cite{ABCGPT}).

\begin{figure}[htbp] \centering
\includegraphics[scale=0.8,clip=true]{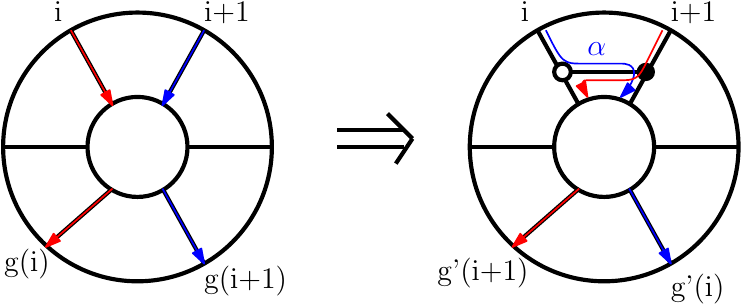}
\caption{}
\end{figure}

These on-shell diagrams correspond to (positroid) varieties in the Grassmannian, which can be represented by sets of matrices. \cite{ABCGPT}) states pictured BCFW bridge application corresponds to the following formula for the columns of the matrix (set):

$v_{a+1} \mapsto \widehat{v_{a+1}} = v_{a+1}+ \alpha v_a$. 

In this section, we write out the proof this formula works as expected; in other words, that this map on a matrix does indeed change the permutation on the matrix (corresponding to the earliest column $f(i)$ after column $i$ such that $v_i\in span\{v_{i+1},v_{i+2},...,v_{f(i)} \}$) in the same way as the change in the permutation on the plabic diagram (given by left-right paths). More precisely, we prove the statement that if I have a matrix whose columns satisfy the bounded affine permutation $\sigma$, and I apply the transformation of $v_{a+1}$ above, that all this does to $\sigma$ is that it switches $\sigma(a)\leftrightarrow \sigma(a+1)$. Let's look at the following picture:

\begin{tikzpicture}
\node[below] at (0,0) {$v_a$}; 
\node[below] at (1,0) {$v_{a+1}$};
\node[below] at (4,0) {$v_{\sigma(a+1)}$}; 
\node[below] at (6,0) {$v_{\sigma(a)}$};
\draw[->] (0,0) .. controls (2,1) and (4,1) .. (6,0);
\draw[->] (1,0) .. controls (2,0.7) and (3,0.7) .. (4,0);
\draw[blue] (-2,-0.5) -- (-2,-0.75) -- (-0.5,-0.75) --(-0.5,-0.5);
\draw[blue] (-0.5,-0.75) -- (3,-0.75) -- (3,-0.5);
\draw[blue] (3,-0.75) -- (5.5,-0.75) -- (5.5,-0.5);
\end{tikzpicture}

Consider the 3 section in blue at the bottom. These are representative. The 1st and 3rd represent vectors $v_b$ such that the entire interval from $b$ to $\sigma(b)$ ($\sigma(b)$ being the place where the ``juggling ball b drops", that is, the first column $v_c$ such that column $v_b$ is in the span of columns $\{v_{b+1},...,v_c\}$.) Thus, these linear dependencies are not affected by the change to $v_{a+1}$ since they do not contain $v_{a+1}$, thus have nothing to do with $v_{a+1}$. 

Now, consider the 2nd section: these are columns $v_b$ such that $b$ is before $a$ but $\sigma(b)$ is after $a+1$. These also are unchanged. This is because we are considering things like span$\{v_{b+1},...,v_c\}$ with $v_a,v_{a+1} \subseteq \{v_{b+1},...,v_c\}$, and this does not change after the transformation because span$\{v_a,v_{a+1}\}=$ span$\{v_a,\hat{v}_{a+1}\}$. Proof: $\alpha_1 v_a + \alpha_2 \hat{v}_{a+1} = \alpha_1 v_a + \alpha_2(v_{a+1} + \alpha v_a) = (\alpha_1+\alpha\alpha_2) v_a+ \alpha_2 v_{a+1}$. 

And clearly any interval $[v,\sigma(b)]$ that does not contain $v_{a+1}$ does not change. Thus, we just have to look at intervals that only contain $v_{a+1}$ but not $v_a$. Theoretically, such intervals could start or end on $v_{a+1}$, but in actuality they cannot end on $v_{a+1}$ because then they would contain $v_a$ as well, thus falling into the 2nd case we just described above.

Thus, we just have to look at $\sigma(a+1)$, and also $\sigma(a)$ (because $\sigma(a)$ corresponds to span$\{v_{a+1},...,v_{\sigma(a)}$, which does not contain both $v_a$ and $v_{a+1}$ so does not fall into the 2nd case above). 

Considering what happens to $\sigma(a)$ first, the span we are analyzing is: span$\{\hat{v}_{a+1},...,v_{\sigma(a)}\}$. Now, $\hat{v}_{a+1}$ equals $v_{a+1}+\alpha v_a$; in other words, we now have a whole copy of $v_a$ in the next column. So all we have to do is to get a span that allows us to eliminate $v_{a+1}$, and we will be left with a scalar multiple of $v_a$. Clearly, we can get that with $\{v_{a+1},...,v_{\sigma(a+1)} \}$. The details are just elementary linear algebra. 

Next, for $\sigma(a+1)$, the important thing to note is that we added $+\alpha v_a$; the important thing is that we have a variable $\alpha$: it has the freedom to be any value, so this is basically saying that we have to find a linear combination that gets us 2 linear combinations: one to kill $v_{a+1}$ and one to kill $v_a$. To kill the $v_{a+1}$, we can use span$\{v_{a+2},...,v_{\sigma(a+1)} \}$. However, to kill the $v_a$, we need span$\{v_{a+1},...,v_{\sigma(a)} \}$. Here are the details:

Letting $v_{a+1} = \beta_{a+2} v_{a+2} + \cdots + \beta_{\sigma(a+1)} v_{\sigma(a+1)}$, we have:

$v_a=\alpha_{a+1} v_{a+1} + \cdots + \alpha_{\sigma(a+1)} v_{\sigma(a+1)} + \cdots + \alpha_{\sigma(a)} v_{\sigma(a)} \\
 = \alpha_{a+1} (\beta_{a+2} v_{a+2} + \cdots + \beta_{\sigma(a+1)} v_{\sigma(a+1)}) + \cdots + \alpha_{\sigma(a+1)} v_{\sigma(a+1)} + \cdots + \alpha_{\sigma(a)} v_{\sigma(a)}$.
 
 Now, we can just add these expressions for $v_{a+1}$ and $\alpha v_a$, and we will get an expression with basis $v_{a+2},...,v_{\sigma(a)}$.

\section{BCFW Bridge Decomposition Method}

Now, we can use these BCFW bridges to create a C-matrix, which is a positive parametrization of the (real nonnegative) positroid cell that can be used in formulas for the on-shell function, which is the amplitude in planar $N=4$ supersymmetric field theory.

The rule is as follows: take a (bounded affine) permutation $g$. Any time you have $g(i)=i$ or $i+n$, this leg gets locked in place, so $g(i)$ cannot participate in any more transpositions. Starting from the left, find the lexicographically first pair $g(a),g(b)$ where $g(a)<g(b)$ and switch them. (Note, g(a) and g(b) are the permutation numbers (where the juggling balls land)). Keep doing this until you get to the identity permutation.

Example: Suppose we have permutation $3142$. This is ${1234 \choose 3142}$, which is actually ${1234 \choose 3546}$ since we are dealing with bounded affine permutations; thus, $f=3546$. Then we look at the lexicographically first (smallest $a,b$) pair $f(a),f(b)$ where $f(a)<f(b)$ (here $f(1)=3<5=f(2)$ is satisfied), so we switch them. This was a transposition of the first two entries, so we record $(12)$. We now have $5346$. Notice that we have $f(1)=5=1+n$, so this leg becomes locked in place. Also, of course, $f(1),(f2)=5,3$ no longer has $5<3$. The new lexicographic minimal pair satisfying the condition $f(a)<f(b)$ is now $f(2)=3<4=f(3)$, so we can switch them: $5436$ (and record $(23)$). Since we have reached $3\to 3$, the $3$ is now locked in place, so it can no longer move. We keep on going until we get to the ``identity permutation" (mod n=4): ${1234 \choose 5634}$. Now, we keep track of all these moves, and going in reverse ordering (starting from the identity), we construct the matrix with the operation $v_{a+1} \mapsto \widehat{v_{a+1}} = v_{a+1}+ \alpha v_a$ that we just proved,

In this example, the sequence of bridges was $(12),(23),(24)$. We start with the identity matrix where the numbers greater than $n$ are in the final permutation of all $0$'s and $n$'s in siteswap language: $5634$ is $4400$ as a siteswap, so $\mqty(1 & 0 & 0 & 0 \\ 0 &1 & 0 & 0)$. Going backwards, so starting with the bridge $(24)$, we noted above that we need to send $v_4 \mapsto v_4+a_1 v_2$: $\mqty(1 & 0 & 0 & 0 \\ 0 &1 & 0 & a_1)$. Next for $(23)$, $\mqty(1 & 0 & 0 & 0 \\ 0 &1 & 0 & a_1)\to \mqty(1 & 0 & 0 & 0 \\ 0 &1 & a_2 & a_1)$. Finally, $(12)$ sends $\mqty(1 & 0 & 0 & 0 \\ 0 &1 & 0 & a_1)\to \mqty(1 & a_3 & 0 & 0 \\ 0 &1 & a_2 & a_1)$.

A nice illustration of this process for a larger example in on page 25 of \cite{ABCGPT}.

We want to prove 2 things about the BCFW bridge decomposition method:\\
(1) that it works: that the algorithm of switching $i$ and $j$ if $\sigma(i)<\sigma(j)$ (or jumping over $0$'s or $n$'s) will take you to the point (and always be a valid bounded affine permutation)
(2) that the BCFW bridge method gives you a sensible matrix (that is, that in the shift $c_i \mapsto c_i + \alpha c_j$, that $c_j \neq 0$ ever. \\

For (1), there are $k$ numbers that let's call B numbers - these will become $n$'s in the siteswap, and $n-k$ numbers that let's call $A$ numbers - these will become $0$'s. The B numbers must move right to get to the right place; they cannot move left because this would imply that they are already right of their correct spot which would imply at siteswap number greater than $n$, and this is impossible. The opposite is true for the A numbers. Every B number is greater in value than an A number. So whenever we have an $i$ and $j$ where we test if $\sigma(i)<\sigma(j)$, there are 3 possibilities: both are in A, both are in B, or $i$ is in A and $j$ in in B (the final possibility that $i$ is in B and $j$ is in A is impossible). We want to show that each entry reaches its desired spot. Suppose $i$ is in A. Then it wants to move left to reach its spot. If $i$ is in the 1st spot, it cannot move left any more, so it must already be an identity column. Thus suppose that it is not in the lefmost spot. This means that, the first time that it moves it must be switching with something form the left, unless every number on the left is larger. 

The proof is obvious from my diagrammatic picture: since we are precisely filling in square with crosses until we have filled the entire diagram in with crosses, which is exactly a siteswap consisting of entirely 0's and n's.  

For (2), this comes from going backwards in the BCFW bridge method. You want to show that $c_j$ at the beginning is an identity column or, has already been hit by an identity column (has already been a $c_i$ before turning into a $c_j$). But this is immediate from my proofs of how to obtain the BCFW bridges by reading off the affine pipe dream diagram. Notice that by the rule depicted in the following, the bridges which are of the form $(c_i c_j)$ always have one of the identity columns as $c_i$ since the identity columns are the ones which ``jump forward" which are precisely on the vertical line segments in the shape, and the first number in $(c_i c_j)$ always comes from the vertical line segment (of the row in which the box/tile lies):

\begin{figure}[htbp]
\includegraphics[scale=0.7,clip=true]{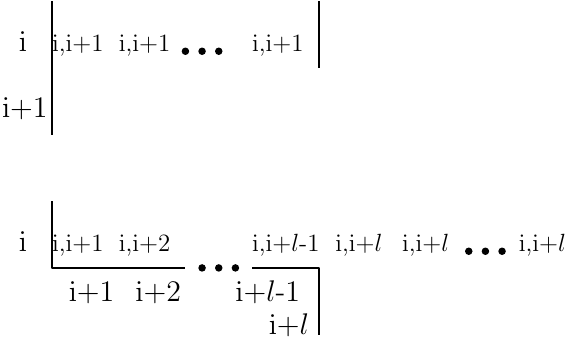}
\end{figure}


Here is a better way to think about this. We have two statements we are trying to show: the first is that the bridge decomposition algorithm does indeed produce the permutation where each $f(i)=i$ or $i+n$ (with a minimal number of transpositions). The second statement is that when doing the change on the matrix given by the formula $v_{a+1} \mapsto \widehat{v_{a+1}} = v_{a+1}+ \alpha v_a$, that we never have $v_a= \vec{0}$, the zero column (otherwise, the matrix would have no record of this transposition).

Split the permutation $f$ into two sets: set $A$, which consists of numbers that are less than or equal to $n$: these are the ones that will end up with $f(i)=i$. And a second set $B$, consisting of numbers that will end up in $f(i)=i+n$. The permutation will initially look like $a_1,...,a_i,b_1,...,b_r, a_{i+1},...,a_{j},b_{r+1},...,b_s,a_{j+1},....,a_k,...,b_{last}$. That is, the permutation will consists of pieces of numbers in $A$ and segments of numbers in $B$. Each number in $A$ will be to the left (including the spot itself) of where it will end up in the end, so it will move to the right in the algorithm. Whereas each number in $B$ will be to the right of where it needs to be, so it will move left. What the algorithm does is it reorganizes each block $a_1,...,a_i$ in decreasing order; it will tend to move the smallest $a$ in this block to the right (or until it gets locked in place). Since the flips must happen to adjacent (unlocked) numbers, this is the most efficient way to do it (another way to see it comes from looking at this as a rotated wiring diagram, as in a pipe dream - this is the inspiration for the technique given at the end of this chapter); in other words, if we move the smallest $a$ in this segment to the right, and it does not get locked in place, then we would have needed to keep moving this $a$ to the right anways, and no other $a$ in this set would have gotten locked in place by moving the $a$'s around in this segment. Each $b_i$ that you get to (e.g., the first is $b_1$) is larger than every other $a$, so it will get moved past the whole segment of $a$'s to the beginning, so it must get locked in place through this process. For example, here if $b_1=1+n$, then we will end up with $b_1,ord(a_1,...,a_i)$, where $ord(a_1,...,a_i)$ stands for $a_1,...,a_i$ being ordered from largest to smallest; otherwise, $b_1$ must get locked in place somewhere within the segment consisting of $a_1,...,a_i$. By the same logic, we know that the general form must be $a_1,...,a_i,b_1,...,b_r, a_{i+1},...,a_{j},b_{r+1},...,b_s,a_{j+1},....,a_k,...,b_{last}$, unless there are entries already locked in place at the start: any $a\in A$ at the end, or $b\in B$ in the beginning must already be locked in place, so cannot move. Since the bridge decomposition method start from the last transposition, it starts with each $b$ where it finally ends up. So we can think of this process we described in reverse: a process where each $b$ starts in the spot where it is $f(i)=i+n$ (corresponding to an identity one-hot column in the matrix) and being moved gradually to where it ends up in the permutation $f$. As such, we will always be adding for $v_{a+1} \mapsto \widehat{v_{a+1}} = v_{a+1}+ \alpha v_a$ a column $v_a$ that is nonzero.

\section{Convention on Vertex Coloring and Grassmannian}

In this section, we discuss the matter of which $k$ in Gr(k,n) we have.

It is a matter of convention whether we have 2 arrows coming in and 1 going out for the black or white vertices. 

We show in this section that this is connected to the choice of left-right paths, so that if we have made a choice of whether we turn right or left at each vertex, this determines whether we have 2 arrows coming in or going out. 

It is most natural to think of $k$ as the number of incoming particles (this relates to the perfect orientation way of doing things): this is because the C-matrix we get from perfect orientations will have a row for every particle coming in. 

If we write the particles on the boundary clockwise, then if we draw the 3-particle black vertex as turning right, then the affine permutation is $123\choose 345$ so the siteswap is $222$. If we then draw the 3-particle white vertex and turning left, then the affine permutation is $111$. Then, the black vertex corresponds to a Grassmannian in Gr(2,3), versus a Grassmannian in Gr(1,3) for the white vertex. Thus, by the previous paragraph, we must have 2 edges incoming (and 1 outgoing) for the black vertex, and vice-versa for the white vertex.

\begin{figure}[htbp]
\includegraphics[scale=0.5,clip=true]{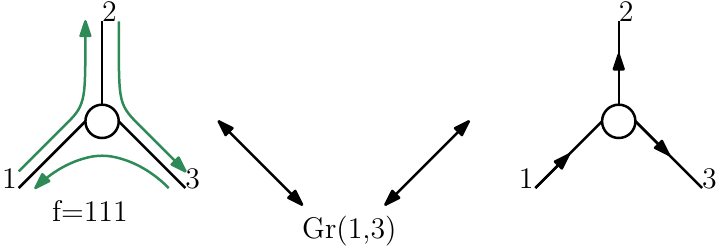}
\end{figure}

And dually for the black vertex (turning right at each verttex, siteswap 222, in Gr(2,3), two incoming arrows to each vertex)

Note that on the right, there is choice/freedom that does not exist in the left figure: namely, the choice of which of the 3 legs is incoming. Thus, the perfect orientation story actually has many, many choices (which should all be equivalent via GL(k) transformations). 

We will see later a conjecture that being able to contract leg $i$ of an on-shell diagram is equivalent ot there existing a perfect orientation where leg $i$ isi ncoming and being able to delete leg $i$ of an on-shell diagram is equivalent ot there existing a perfect orientation where leg $i$ is outgoing.

We also comment that Grassmannian duality occurs if we just exchange all the black and white vertices. On the left side, this just sends all left-right paths in the opposite direction - that is, it produces the inverse permutation, which we can easily check is in $Gr(n-k,n)$. On the right side, we can think of inverting all the arrows (after exchangign all the colors), which clealry preserves the condition that each black vertex has 2 edges incoming, and clearly also sends Gr(k,n) to Gr(n-k,n) since the external edges that were incoming are now outgoing and vice-versa.

\section{Pipe dreams for the BCFW bridge decomposition}\label{section:pipestoplabic}

We can read out the sequence of bridges above by drawing out a pipe dream and following certain rules. It may be easier to follow and start with the example in the next Section \ref{sec:pipesplabicexample}, since this section was written awhile back. 

Intuitively, a pipe dream is just a rotated wiring diagram, and the BCFW bridge method above is just getting to the largest permutation in Bruhat order by filling everything in with crosses. Thus, the transpositions in the BCFW bridge method are simply coming from a choice of filling in squares with elbow tiles in a pipe dream. This is the high-level explanation of where the ideas in this section are coming from. 

We first state the lemma that the patch we must draw corresponds to the first element $I_1$ of the Grassmannian necklace for a bounded affine permutation $g$.

\begin{proof} Note that, by definition, $I_1$ corresponds to the first $k$ linearly independent columns in a matrix representative of $g$; any set of columns that is smaller in Gale order must have determinant that vanishes.

Write out the bounded affine permutation for $1$ to $n$; there will be $k$ excedences (numbers greater than $n$), that we call $J_1$ (thus, $\forall j\in J_1, j>n$). Take $\tilde{J_1} = J_1 (mod n)$. From the juggling way of viewing things, the first $k$ linearly independent columns are the first $k$ spots where balls land, starting from $t=1$, that were not thrown on or after $t=1$; if we consider the periodic interval starting from $1+n$ rather than $1$, these $k$ spots are precisely the balls coming in from before $t=n+1$, which is precisely the excedences in the interval $[1,n]$. Equivalently, in the diagonal diagram, these are the lines that cross the horizontal line at $t=1$ coming down, e.g.:

These numbers are exactly $\tilde{J_1}$. And when we draw the shape of the patch (the boundary), the vertical lines are labeled by the numbers in $\tilde{J_1}$ (if we filled in the shape entirely with crosses, the horizontal pipes would from from $\tilde{J_1} to J_1$), whereas the horizontal parts of the boundary are labeled by the other number (in $[n]\backslash \tilde{J_1}$). 

\end{proof}

We then obtain the pipe dream for the siteswap. Note that we need to create the top pipe dream, as detailed in the section on affine pipe dreams, or Appendix B of \cite{Flu}. The reason is because we want these crosses to be as far away from the bottom pipe dream as possible, and the reason for this is because the BCFW method introduces a BCFW bridge say (12) corresponding to a transposition $s_1$ when $\sigma(1)<\sigma(2)$. The $s_1$ is acting on the \textit{positions}, not the values (in other words, in double-line notation $123 \cdots n \choose \sigma(1) \sigma(2) \sigma(3)...\sigma(n)$ it is the columns corresponding to $s_1$ ($1$ and $2$) in the top row that are switched, the places where the 2nd row equals $1$ and $2$ are irrelevant). Since the \textit{positions} correspond to the order of the pipes coming from the left and bottom (or the lower) boundary, we need, say, $s_1$ to actually switch the pipes that emerge from the 1st and 2nd entries along the left boundary. However, this correspondence will only happen if the tiles corresponding to the bottom pipe dream are ``free"; that is, crosses need to be moved toward the top boundary because if they are lower, then they will interfere with the numbering, causing for example an $s_1$ to switch pipes other than the 1st and 2nd pipes. In summary/short, single the transpositions are acting on the order of the pipes coming from the bottom boundary, the tiles along the bottom boundary need to be kept clear.

We can prove our ordering of reading the word by matching up the BCFW bridge algorithm with what's going on with the pipes. Any time we have 2 adjacent numbers on the left that both have the same number of crosses stretching all the way from the left boundary, they can be switched (the base case is when everything is empty, in which case we satisfy this condition (with the same = 0 number of crosses), and you place a cross in the upper row of these 2 rows with the same number of crosses (see following figure, which adds the red cross via the rule):
\begin{figure}[htbp]
\includegraphics[scale=0.6,clip=true]{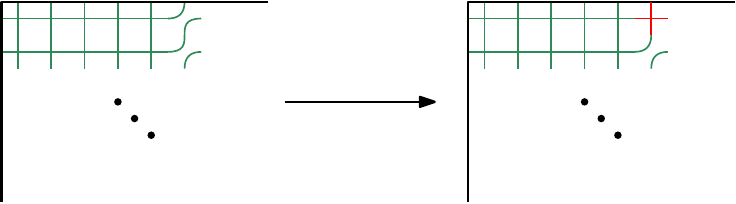}
\end{figure}
In the case where they are empty, this is the northwest corner. Going from top to bottom rows, if there is a mismatch between adjacent rows, then you move one row lower and look at the next pair of rows. Call the rule (*).
Notice that if we have an empty square pipe dream, i.e. the top cell, then this creates the word $w_0$ - this is the permutation that is longest in Bruhat order, e.g. $w_0$ on $S_4$ is $1234 \choose 4321$, and it does this by applying transpositions on the left, so the word is $s_1 s_2 s_1 s_3 s_2 s_1 s_4, s_3 s_2 s_1...$. Getting to $w_0$ consists of filling in everything above the southwest to northeast diagonal. Once this is filled in, one continues filling in SW-NE diagonals going from left to right. For example, in the Gr(4,8) example, we have the following labelling, where we read the word along these SW-NE diagonals:

\begin{tabular}{|c|c|c|c|}
\hline
$s_1$ & $s_1$ & $s_1$ & $s_1$ \\
\hline
$s_2$ & $s_2$ & $s_2$ & $s_2$ \\
\hline
$s_3$ & $s_3$ & $s_3$ & $s_3$ \\
\hline
$s_4$ & $s_{4,6}$ & $s_{4,7}$ & $s_{4,8}$ \\
\hline
\end{tabular}
However, because of the rule for the BCFW bridge decomposition that we skip over entire rows or columns that have been filled in, this process on the empty pipe dreams fills in column by column, and these need to be skipped for: this is the reason for the notation e.g. $s_{46}$ (the column labeled 5 has been filled in entirely with crosses, so can be skipped, so the next $s_4$ after that is actually flipping the 4th and 6th pipes, not the 4th and 5th pipes. See also the next comment. 

One small modification is that once a row or column is filled in completely with crosses, it must be removed (think: full (non-projectionless) deletion or contraction). Since all the pipes which are able to pass through this row or column without change in order, this row or column can essentially be ignored; this is a reason for why the BCFW bridge construction works the way it does with ignoring entries that are loops or coloops. This point is often times a non-issue if, for example, one is dealing with a square pipe dream diagram, where using the BCFW bridge method, the columns are getting filled in starting from the left, because removing the leftmost column does not change the relative ordering of any of the remaining transpositions $s_i$. In contrast, if we are doing a deletion or contraction of a row or column in the middle, removing this row or column of crosses can change the order in which we read the $s_i$. Note that this is just for convenience, if we want to think about reading such that we are going along southwest to northeast diagonals; we don't actually need to make this comment if we just apply the rules (*) above. 

We can think about the grid above as a "base case." Deletion (fillig in column with crosses) does not change any of the numbers, it just eliminates some of the tiles, for example if we contracted the 3rd column, we get:
\begin{tabular}{|c|c|c|c|}
\hline
$s_1$ & $s_1$ & $+$ & $s_1$ \\
\hline
$s_2$ & $s_2$ & $+$ & $s_2$ \\
\hline
$s_3$ & $s_3$ & $+$ & $s_3$ \\
\hline
$s_4$ & $s_{4,6}$ & $+$ & $s_{4,8}$ \\
\hline
\end{tabular}
The reason is because the $s_i$ generally correspond to the rows; the only place where the column shows up is in the $s_{4,col}$ at the end, and the one with $s_{4,7}$ is already eliminated. You do need to change the order, which you can do by continuing to think about SW-NE diagonals along:
\begin{tabular}{|c|c|c|}
\hline
$s_1$ & $s_1$& $s_1$ \\
\hline
$s_2$ & $s_2$ & $s_2$ \\
\hline
$s_3$ & $s_3$ & $s_3$ \\
\hline
$s_4$ & $s_{4,6}$ & $s_{4,8}$ \\
\hline
\end{tabular}

Now, if we contracted (filled in a row with crosses), for example the 2nd row, we get:
\begin{tabular}{|c|c|c|c|}
\hline
$s_1$ & $s_1$ & $s_1$ & $s_1$ \\
\hline
$+$ & $+$ & $+$ & $+$ \\
\hline
$s_3$ & $s_3$ & $s_3$ & $s_3$ \\
\hline
$s_4$ & $s_{4,6}$ & $s_{4,7}$ & $s_{4,8}$ \\
\hline
\end{tabular}

Here, we have to make the modification that, since the 2nd pipe should be considered eliminated, pipes 1 and 3 are now considered adjacent, so $s_1$ switches 1 and 3. To be explicit, we could write:
\begin{tabular}{|c|c|c|c|}
\hline
$s_{13}$ & $s_{13}$ & $s_{13}$ & $s_{13}$ \\
\hline
$+$ & $+$ & $+$ & $+$ \\
\hline
$s_3$ & $s_3$ & $s_3$ & $s_3$ \\
\hline
$s_4$ & $s_{4,6}$ & $s_{4,7}$ & $s_{4,8}$ \\
\hline
\end{tabular}
and we read the word going along SW-NE diagonals of:
\begin{tabular}{|c|c|c|c|}
\hline
$s_{13}$ & $s_{13}$ & $s_{13}$ & $s_{13}$ \\
\hline
$s_3$ & $s_3$ & $s_3$ & $s_3$ \\
\hline
$s_4$ & $s_{4,6}$ & $s_{4,7}$ & $s_{4,8}$ \\
\hline
\end{tabular}

If we put the 2 previous examples together, such that we contract the 2nd row and 3rd column (the 2nd and 7th particles), we get:
\begin{tabular}{|c|c|c|c|}
\hline
$s_{13}$ & $s_{13}$ & $+$ & $s_{13}$ \\
\hline
$+$ & $+$ & $+$ & $+$ \\
\hline
$s_3$ & $s_3$ & $+$ & $s_3$ \\
\hline
$s_4$ & $s_{4,6}$ & $+$ & $s_{4,8}$ \\
\hline
\end{tabular}
and we would be reading:
\begin{tabular}{|c|c|c|}
\hline
$s_{13}$ & $s_{13}$ & $s_{13}$ \\
\hline
$s_3$ & $s_3$  & $s_3$ \\
\hline
$s_4$ & $s_{4,6}$ & $s_{4,8}$ \\
\hline
\end{tabular}

Thus, we can, in some sense, view the bridges as some sort of give-and-take with the pipe dream crosses: the crosses that are filled in when creating the pipe dream are exactly the complement of the ones that are used in creating the BCFW bridge decomposition. Because deletion and contraction put in more crosses, they can be viewed as reducing BCFW bridges (or setting some equal to zero). In some cases we can literally just view deletion and contraction as setting $\alpha$'s equal to zero in a C-matrix parameterization. However, one must be careful because of: (1) the change in ordering of the way the tiles is read if middle columns or rows are deleted or contracted (2) if the pipe dream is not rigid, the moves that can happen in deletion or contraction, which can change the order as well (note that, however, things are a bit easier in this case compared to general, because we are always dealing with the top pipe dream, so we have a unique pipe dream both before and after deletion/contraction). 

We can recap/summarize the material above in the following 2 points:

(1) The order in which the elements are read:

The rule is that we are always reading along southwest to northeast diagonals, starting from the northwest corner. 

\begin{figure}[htbp]
	\includegraphics[scale=0.5,clip=true]{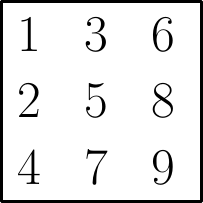}
\end{figure}

In this example the shape was a square, but in practice the shape will be determined by the Grassmann necklace element $I_1$ of our choice of base

(2) What kinds of flips fill the squares:

In terms of what letters or flips to fill in the squares of the diagram with, the rules are the following (1st rule can be seen as a special case of 2nd):

\begin{figure}[htbp]
	\includegraphics[scale=0.6,clip=true]{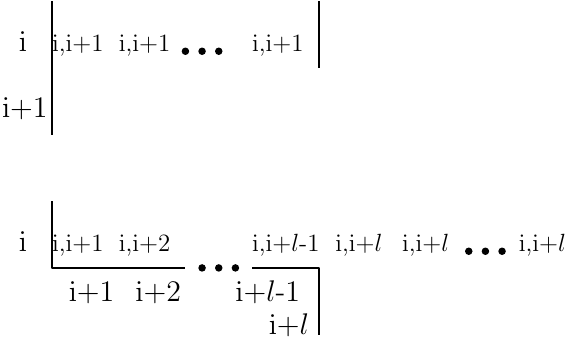}
\end{figure}

Here are some examples: 

\begin{figure}[htbp]
	\includegraphics[scale=0.6,clip=true]{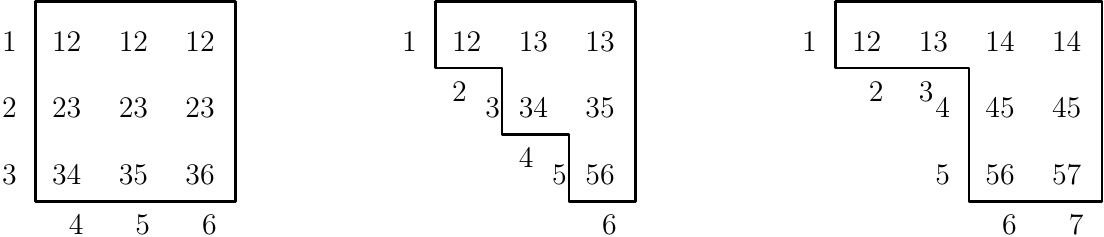}
\end{figure}

E.g. the 2nd example would consist of the sequence: 12,13,34,13,35,56

Things get more complicated if the pipe dream does not consist entirely of elbows - if there are crosses. When this happens, there are 2 additional rules to be taken into account, one primary rule and a secondary rule. The primarily rule is that in the process of filling in the squares (reading the word), whenever a full row is fully filled in with crosses, we must eliminate that row and continue on filling in/reading according to the order given in point (1) above in the \textit{new} grid (with the filled in row eliminated). Additionally, the new diagram must be relabeled, with the labeling of the letters still following point (2) on the new grid; for example, if we had two rows $i$ and $i+1$ after $i$, and $i+1$ is filled in with crosses, and say the next vertical line is $i+2$, then the elements in row $i$ that were originally $(i i+1)$ now becomes $(i i+2)$. Another way to write this rule is that if we have a string of crosses coming all the way from the wall on the right, consider each empty square immediately above this string on crosses: suppose that square is in row $i$. Then this square will be filled with $(i i+j)$, where $j$ is the next box below it that is in a row where not all boxes to its right (including it) are filled in with crosses. 

Usually, if the row becomes entirely filled in with crosses by filling in the rightmost tile in the row last, this step is not important, because this implies that all tiles in rows higher up in the diagram already consist entirely of crosses, so this row will no longer be featured in any transpositions, but if the row already consisted of crosses in the rightmost part and then became a fully crossed out row when a tile in the center was made a cross, then following this step will be essential to getting the right bridge sequences. We already gave an example of this for contraction above, where an entire row is filled in with crosses. We will illustrate a more general example in the next section. Note that we can think of this elimination of a row as the second step of contraction in Section \ref{sec:DelContr}, where we project to an isomorphic positroid variety in $Gr(k-1,n-1)$.

The second rule is that, when reading along a southwest to northeast diagonal line, for entries that are adjacent, we read in the way described above, north southwest or northeast. However, if there are crosses along the diagonal, this causes the transpositions to be disjoint, so the lexicographic-first rule in the BCFW bridge decomposition method will choose the more northeastern piece first. Thus, we can think of crosses along a diagonal as separating this diagonal into ``islands" or blocks, of which, the islands are chosen from northeast to southwest, whereas within an island, the ordering is southwest to northeast. See the following figure as an example, where the numbers along the diagonal give the order in which the tiles are read. 

\begin{figure}[htbp] \centering
	\includegraphics[scale=0.4,clip=true]{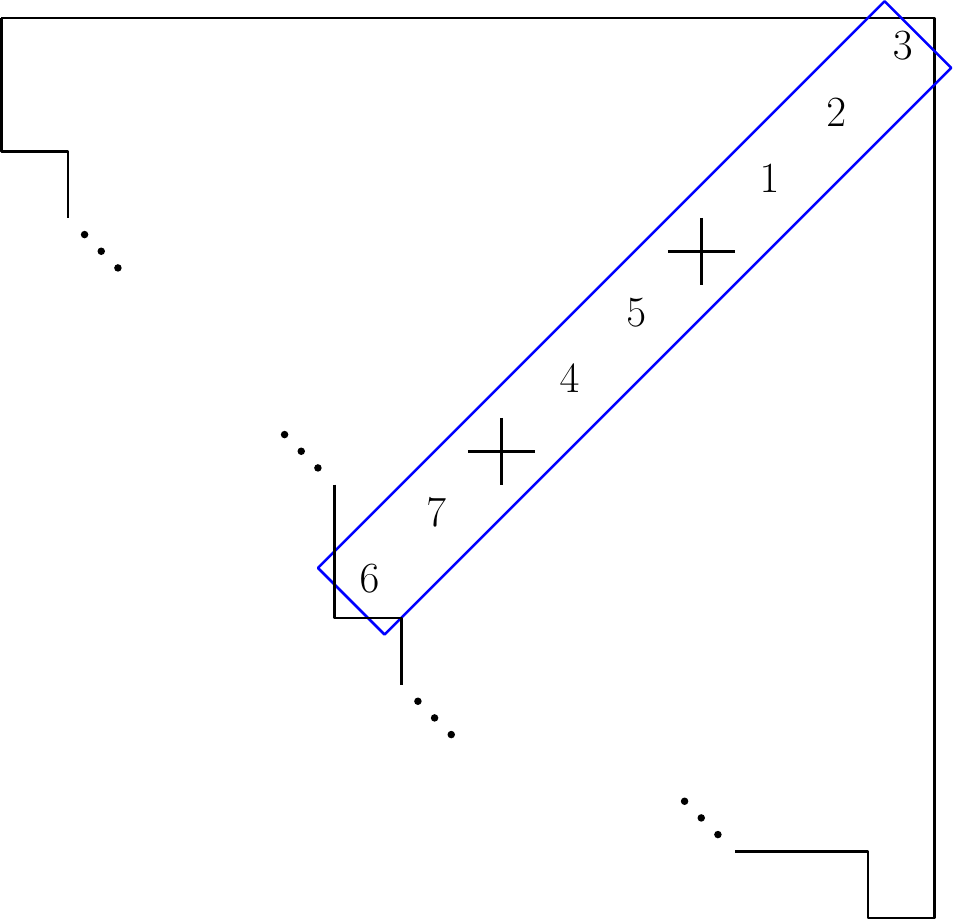}
	\caption{}\label{fig:islands}
\end{figure}

Finally, once we have the sequence of transpositions, with each transposition corresponding to an elbow tile on the pipe dream, we can add BCFW bridges for each transposition (we can draw a BCFW bridge in the given tile), in the order prescribed by the process that we just described. Note that this gives us a map from: pipe dreams to plabic diagrams. \cite{Post} already gave a map from pipe dreams to plabic diagrams for the top pipe dream. Our method also uses the top pipe dream. However, unlike Postnikov, our method gives the map to plabic diagrams using the BCFW bridge decomposition. The example in the following section may help to clarify how this works.

\section{Example: from (Positroid,Point) to Plabic Digram using BCFW Bridges}\label{sec:pipesplabicexample}

If the discussion above was confusing, here we hope to make everything clearer by providing an explicit example of the material above that contains all the relevant features. Specifically, this example with show how to associate a plabic diagram to a bounded affine permutation by placing BCFW bridges on a pipe dream. 
Now, Figure 3 of \cite{ARW}, or Figure 20.1 of \cite{Post}, already provided a way of associating a plabic diagrams to a pipe dream; however, it does not use BCFW bridges as we do here (whether there is a nice way to apply moves to get from one to another would be an interesting question to investigate) so does not produce a matrix representative we can plug into the formula for scattering amplitudes.

As an example, let's take the siteswap $f=2,5,2,7,3,3,6,2,6$. This has 9 numbers, and $2+5+2+7+3+3+6+2+6=36$, so the average is $36/9=4=k$, which means this labels a positroid variety in $Gr(4,9)$. If we write out the first 9 numbers is the bounded affine permutation, it is:

\begin{center}
\begin{tabular}{cccccccccc}
1&2&3&4&5&6&7&8&9& \\
2&5&2&7&3&3&6&2&6&+ \\
\hline
3&7&5&11&8&9&13&10&15&
\end{tabular}
\end{center}

Thus, the bounded affine permutation $g$ corresponding to $f$ has its first $n=9$ entries as $g(1,2...,9)=(3,7,5,11,8,9,13,10,15)$. It's easy to see that we can easily go backwards (from $g$ to $f$) by subtracting $1,2,3,4,5,6,7,8,9,$ instead of adding (this is just reiterating earlier material for those who are only reading particular sections). 

Let's draw out the patch. Recall from Section \ref{ch:affinefinitepipes} that we proved that the finite pipe dream uses the patch given by the elements of the Grassmann necklace $I_1$ of $g$. This is the first $4$ linearly independent columns in a matrix representative of a generic point on the positroid variety. It is easy to get this from the bounded affine permutation $(3,7,5,11,8,9,13,10,15)$. These are the entries that exceed $n$ (mod $n$). Here, the numbers that exceed $n=9$ are $10=1(mod 9), 11=2(mod 9), 13=4(mod 9), 15=6(mod 9)$. Thus, $I_1=\{1,2,4,6\}$. These will be the vertical lines in the boundary of the pipe dream; the rest of the lines will be horizontal lines. This was called the (a) ``distinguished path" in Section \ref{affinepipedreams}. From this, we can create the (b) affine pipe dream and the (b) nonaffine pipe dream, see the following figure: 

  \begin{figure}[htbp]
	\includegraphics[scale=0.5,clip=true]{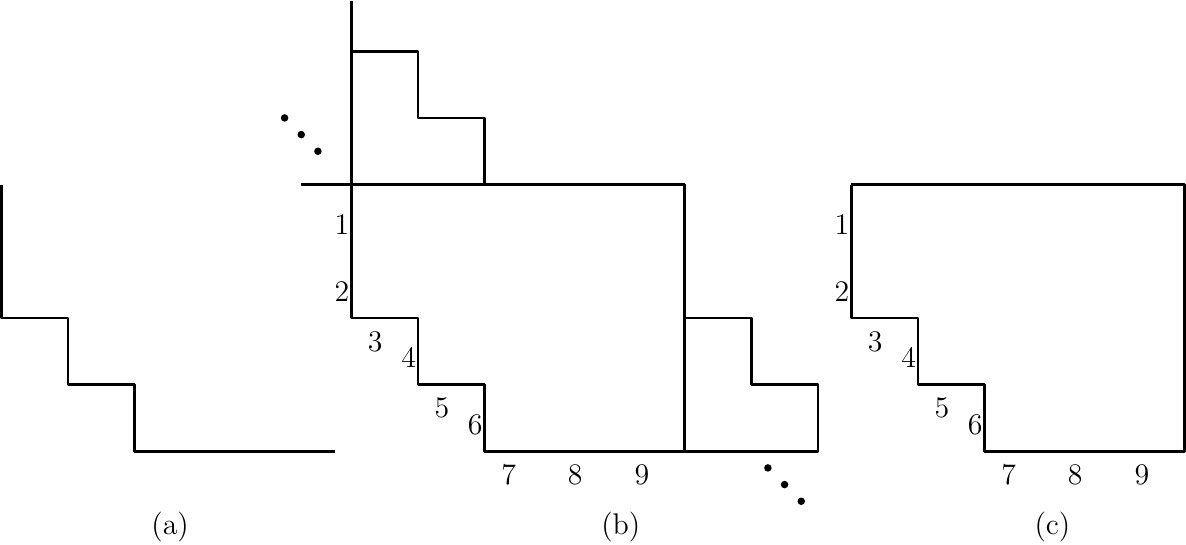}
\end{figure}

The first step is to obtain the top affine pipe dream: since the BCFW bridge method will involve the top pipe dream. Affine pipe dreams are covered in Section \ref{affinepipedreams}. In that section, we noted that the process of obtaining the top pipe dream proceeds by applying transpositions to the inverse permutation (in contrast to the bottom pipe dream, which uses the non-inverse permutation). Thus, to get the inverse permutation, it is perhaps easiest to start with the siteswap $f$ and obtain $f^{-1}(i)$ for each $i$. We can take two windows: $2,5,2,7,3,3,6,2,6 | 2,5,2,7,3,3,6,2,6$, and for the latter window, simply look at which element of the previous window comes in. For example $2,5,2,7,3,3,6,\rd{2},6 | \teal{2},5,2,7,3,3,6,2,6$ shows that the red lands on the teal: $f^{-1}(1)=-2$; similarly, $2,5,2,\rd{7},3,3,6,2,6 | 2,\teal{5},2,7,3,3,6,2,6$, so $f^{-1}(2)=-7$, and $2,5,2,7,3,3,6,2,6 | \rd{2},5,\teal{2},7,3,3,6,2,6$, so $f^{-1}(3)=-2$, etc. We get $f^{-1} =(-2,-7,-2,-6,-2,-6,-5,-3,-3)$. We then add $n=9$ to this sequence, which gives us: $(-2,-7,-2,-6,-2,-6,-5,-3,-3)+9=(7,2,7,3,7,3,4,6,6)$. As a check, the sum $7+2+7+3+7+3+4+6+6=45$ and $45/9=5$, which is exactly what we want since the inverse permutation should live in the dual Grassmannian $Gr(n-k,n)=Gr(9-4,9)=Gr(5,9)$. The inverse bounded affine permutation is then: $(7,2,7,3,7,3,4,6,6)+(1,2,3,4,5,6,7,8,9)=(8,4,10,7,12,9,11,14,15)$. 

For forming the affine pipe dream, we consider just one block - see (a) of the following figure (the full affine pipe dream repeats infinitely to the northwest and southeast). Note that since we using the inverse permutation to obtain the top pipe dream, we need to align the transposition with the upper boundary, not the lower boundary. For example, the transpositions $s_6$ (which swap the 6th and 7th pipes) will be along the southwest to northeast diagonal, aligned to be between where the $6$ and the $7$ are along the upper boundary of the affine pipe dream; this is displayed in the following diagram (a), along with the transpositions $s_7$. In the following, we will simplify the diagrams and just write, for example, ``$6$" rather than $s_6$; the full list of transpositions is in (b) of the figure:

  \begin{figure}[htbp] \centering
	\includegraphics[scale=0.5,clip=true]{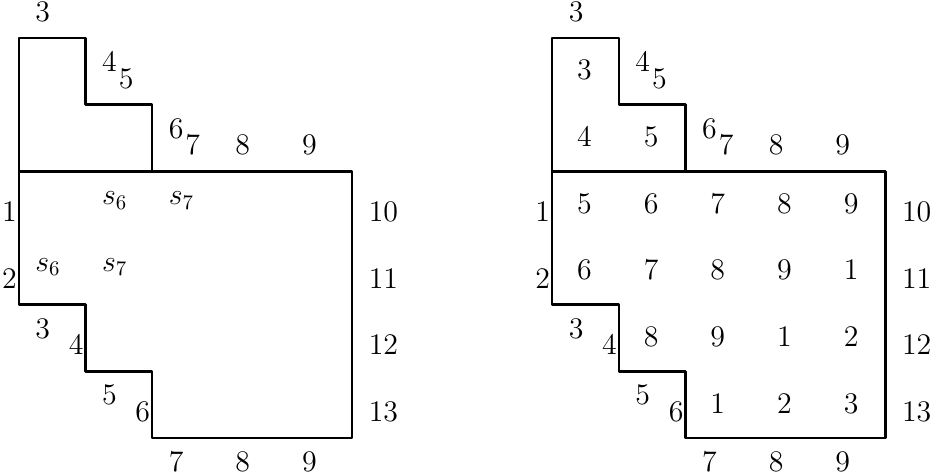}
	\caption{}
\end{figure}

Next, there is the question of what order we read these transpositions in. The order will simply be the reverse of the order for the bottom pipe dream given in Section \ref{affinepipedreams}. It will be the direction of the blue arrows in (a) of the following figure, or, written out explicitly, the ordering in (b) in blue. Thus, we will read off the transpositions (like $s_4$ or $s_7$) in the previous figure, in the order given by the following figure.

  \begin{figure}[htbp] \centering
	\includegraphics[scale=0.5,clip=true]{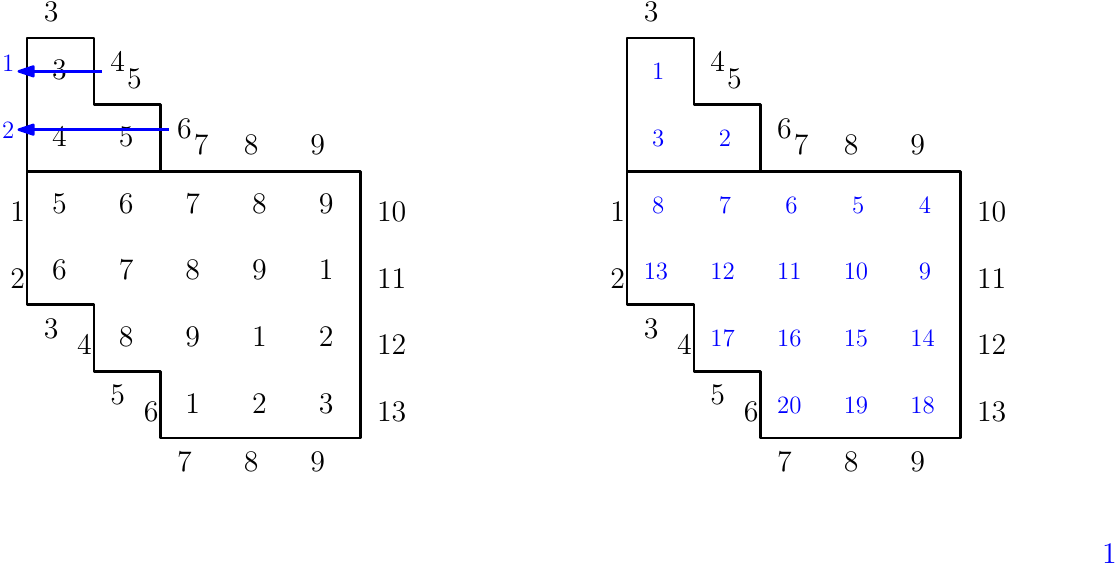}
	\caption{}
\end{figure}

For example, we start with the first entry, which is an $s_3$. The question is: Is $g(3)>g(4)$? If so, then we make the switch, and this first entry gets a cross:

   \begin{figure}[htbp] \centering
	\includegraphics[scale=0.5,clip=true]{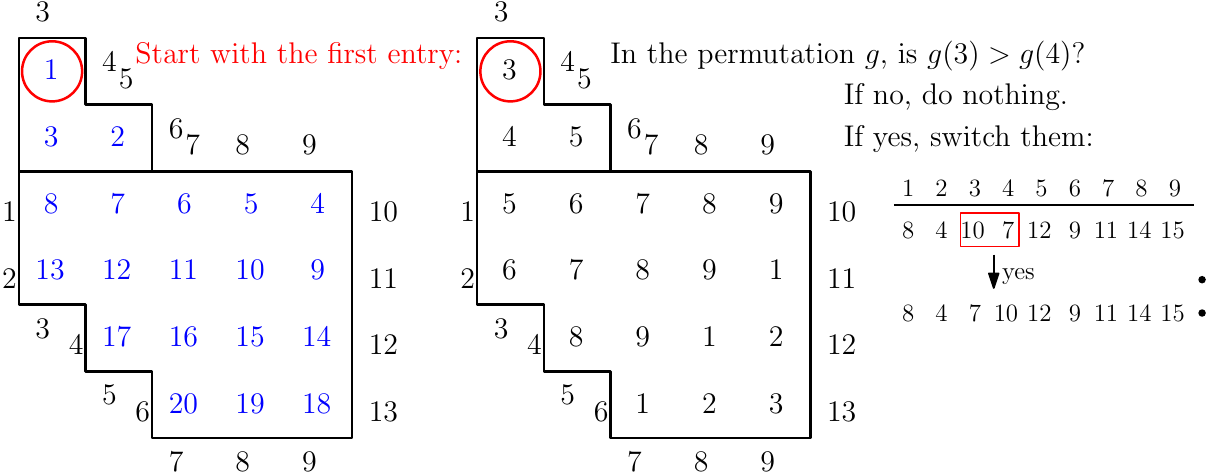}
		\includegraphics[scale=0.5,clip=true]{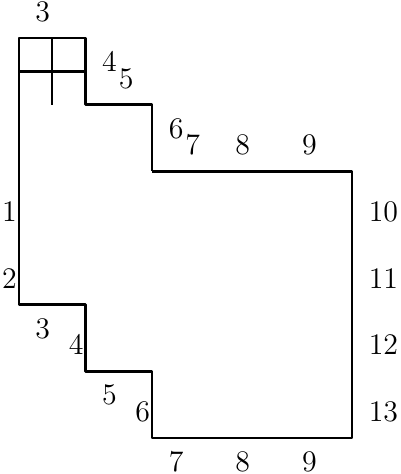}
		\caption{}
\end{figure}

Keeping track of where we are at:

\begin{center}
\begin{tabular}{c|c|ccccccccc|c}
& 0 & 1 & 2 & 3 & 4 & 5 & 6 & 7 & 8 & 9 & 10 \\ \hline
& 6 & 8 & 4 & 10 & 7 & 12 & 9 & 11 & 14 & 15 & 17 \\
\textcolor{red}{$s_3$} & 6 & 8 & 4 & \textcolor{red}{7} & \textcolor{red}{10} & 12 & 9 & 11 & 14 & 15 & 17 
\end{tabular}
\end{center}

We continue on with the 2nd entry, which corresponds to $s_5$, and look to see if $g(5)>g(6)$. This also is true, so we switch these as well. 

\begin{center}
\begin{tabular}{c|c|ccccccccc|c}
& 0 & 1 & 2 & 3 & 4 & 5 & 6 & 7 & 8 & 9 & 10 \\ \hline
& 6 & 8 & 4 & 10 & 7 & 12 & 9 & 11 & 14 & 15 & 17 \\
$s_3$ & 6 & 8 & 4 & \textcolor{red}{7} & \textcolor{red}{10} & 12 & 9 & 11 & 14 & 15 & 17 \\
\textcolor{red}{$s_5$} & 6 & 8 & 4 & 7 & 10 & \textcolor{red}{12} & \textcolor{red}{9} & 11 & 14 & 15 & 17 
\end{tabular}
\end{center}

Recall from Section \ref{ch:affinefinitepipes} that finite pipe dreams are embedded in affine pipe dreams by letting the piece below the distinguished path (in this example, this corresponds to the \begin{tikzpicture}[scale=0.5]
\draw (0,0)--(2,0)--(2,1)--(1,1)--(1,2)--(0,2)--(0,0);
\node at (0.5,1.5) {$3$}; 
\node at (0.5,0.5) {$4$}; 
\node at (1.5,0.5) {$5$}; 
\end{tikzpicture} part) be filled in entirely with crosses, so we expect that $s_3,s_5,s_4$ will all be chosen.

Continuing on, as in Section \ref{affinepipedreams}, we fill in the rest of the table:

\begin{center}
\begin{tabular}{c|c|ccccccccc|c}
& 0 & 1 & 2 & 3 & 4 & 5 & 6 & 7 & 8 & 9 & 10 \\ \hline
& 6 & 8 & 4 & 10 & 7 & 12 & 9 & 11 & 14 & 15 & 17 \\
$s_3$ & 6 & 8 & 4 & \textcolor{red}{7} & \textcolor{red}{10} & 12 & 9 & 11 & 14 & 15 & 17 \\
$s_5$ & 6 & 8 & 4 & 7 & 10 & \textcolor{red}{9} & \textcolor{red}{12} & 11 & 14 & 15 & 17 \\
$s_4$ & 6 & 8 & 4 & 7 & \textcolor{red}{9} & \textcolor{red}{10} & 12 & 11 & 14 & 15 & 17 \\
$s_6$ & 6 & 8 & 4 & 7 & 9 & 10 & \textcolor{red}{11} & \textcolor{red}{12} & 14 & 15 & 17 \\
$s_1$ & 6 & \textcolor{red}{4} & \textcolor{red}{8} & 7 & 9 & 10 & 11 & 12 & 14 & 15 & \textcolor{red}{13} \\
$s_9$ & \textcolor{red}{4} & \textcolor{red}{6} & 8 & 7 & 9 & 10 & 11 & 12 & 14 & \textcolor{red}{13} & \textcolor{red}{15} \\
$s_8$ & \textcolor{red}{5} & 6 & 8 & 7 & 9 & 10 & 11 & 12 & \textcolor{red}{13} & \textcolor{red}{14} & 15 \\
$s_2$ & 5 & 6 & \textcolor{red}{7} & \textcolor{red}{8} & 9 & 10 & 11 & 12 & 13 & 14 & 15 \\
\end{tabular}
\end{center}

Following along on the pipe dream, we get: 

\begin{figure}[htbp] \centering
	\includegraphics[scale=0.5,clip=true]{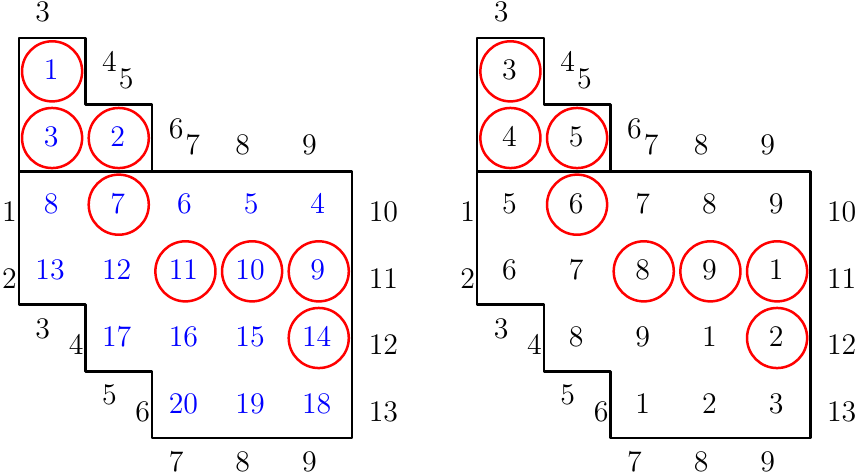}
		\includegraphics[scale=0.5,clip=true]{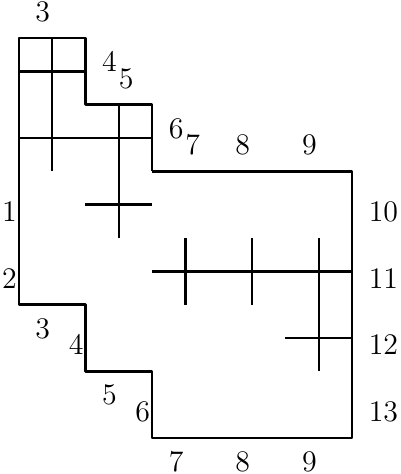}
		\caption{}
\end{figure}

The rightmost part of this figure is the top pipe (affine) pipe dream for this permutation on this patch (the tiles that have no crosses are assumed to be elbows). The corresponding finite pipe dream is as follows. We will work with finite pipe dreams for the rest of this section; however, forming the plabic diagram on affine pipe dreams is essentially no different than the finite pipe dream case, since the only different is a set of tiles consisting entirely of crosses, and these sorts of tiles play no role in the rest of the process.

\begin{figure}[htbp] \centering
	\includegraphics[scale=0.5,clip=true]{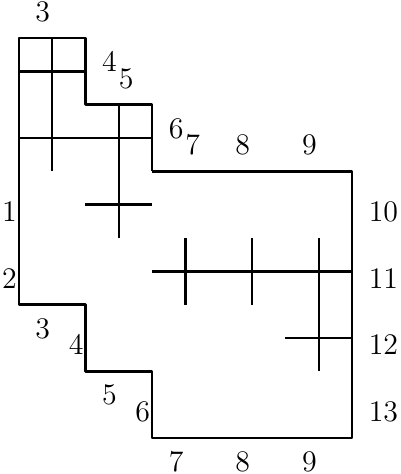}
		\caption{}
\end{figure}

Now, it's time for the 2nd step: obtaining the BCFW bridges. 

What the BCFW bridge method does, is it fills in the rest of the pipe dream with crosses - it provides a particular choice of ordering to replace all the elbows in the previous pipe dream with crosses. Each time one of the elbow tiles becomes filled in with a cross, there is a transposition that must be recorded that will be used to create the plabic diagram at the end. 

As noted in the previous diagram, we write the transpositions on the pipe dream according to the rule that we have transposition $t=(ij)$ where $i$ corresponds to the row that $t$ is in (the vertical boundary to the very left of this tile), and $j$ is the horizontal line segment directly below $t$ if there is a horizontal line segment directly bordering this tile; if there is no horizontal line segment directly bordering this tile, then $j$ is equal to the next vertical line below $i$. In this case, we have:

\begin{figure}[htbp] \centering
	\includegraphics[scale=0.5,clip=true]{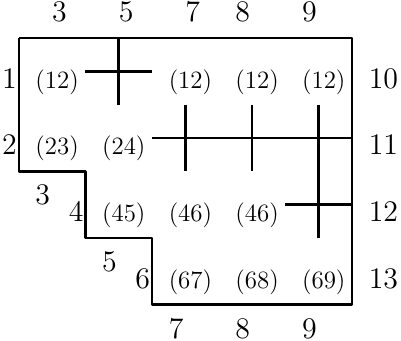}
		\caption{}
\end{figure}

As a first pass, we would read the sequence of transpositions (here, the transpositions are not necessary adjacent: before, we only had elements $s_i=(i,i+1)$, but now we have $(i,j)$ where $j\neq i+1$) as in the following figure, according to the direction given by the arrows in (a), or with the order explicitly written out in the tiles in (b):

\begin{figure}[htbp] \centering
	\includegraphics[scale=0.5,clip=true]{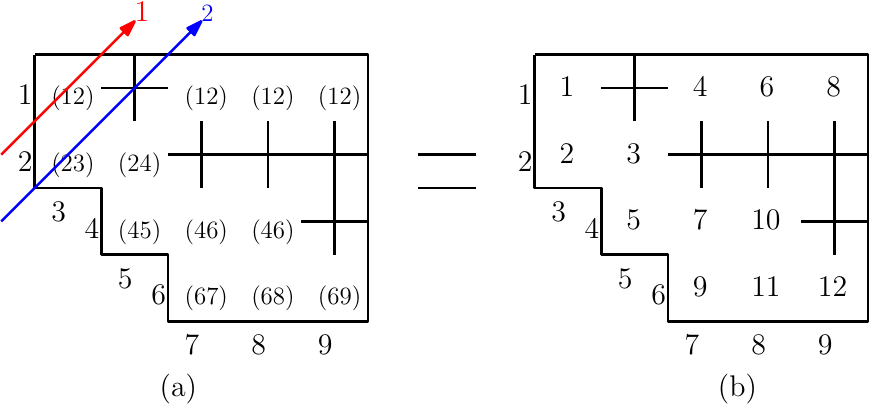}
		\caption{}
\end{figure}

That is, this ordering would produce the sequence: \\
$(12),(23),(24),(12),(45),(12),(46),(12),(67),(46),(68),(69)$. However, this is not quite right, since, as noted in the previous section, there are 2 rules we must follow. The first is the primarily rule, which is that when a row becomes completely filled in with crosses, it must be eliminated from the diagram, and the diagram must be relabeled. 

If we apply these rules, we process is illustrated by the following figure. Note that when we have applied rule 1 to change the shape, we have drawn the new contracted shape in gray (since these gray figures do not advance the list of bridge decompositions):

\begin{figure}[htbp] \centering
	\includegraphics[scale=0.7,clip=true]{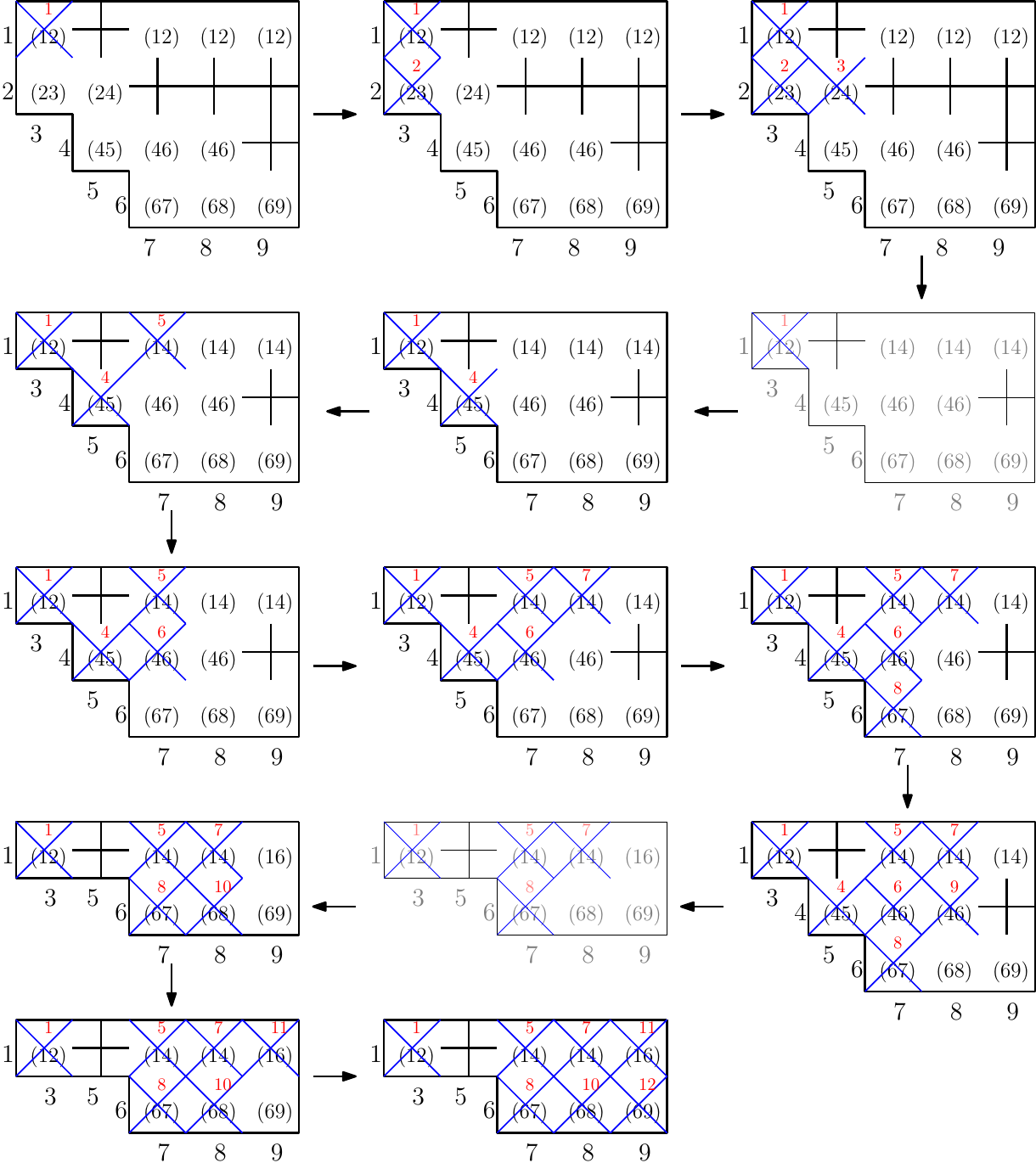}
		\caption{}
\end{figure}

\newpage

Thus, going back to our original pipe dream, we have obtained the following order (and transpositions):

\begin{figure}[htbp] \centering
	\includegraphics[scale=0.7,clip=true]{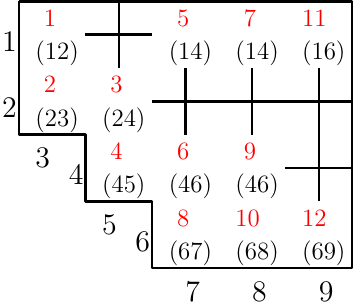}
		\caption{}
\end{figure}

Obtaining this wasn't nicer or easier than simply running the BCFW bridge decomposition method, but it allows us to draw a plabic diagram associated to our pipe dream by placing BCFW bridges on the tiles, which we show in the following figure:

\begin{figure}[htbp] \centering
	\includegraphics[scale=0.7,clip=true]{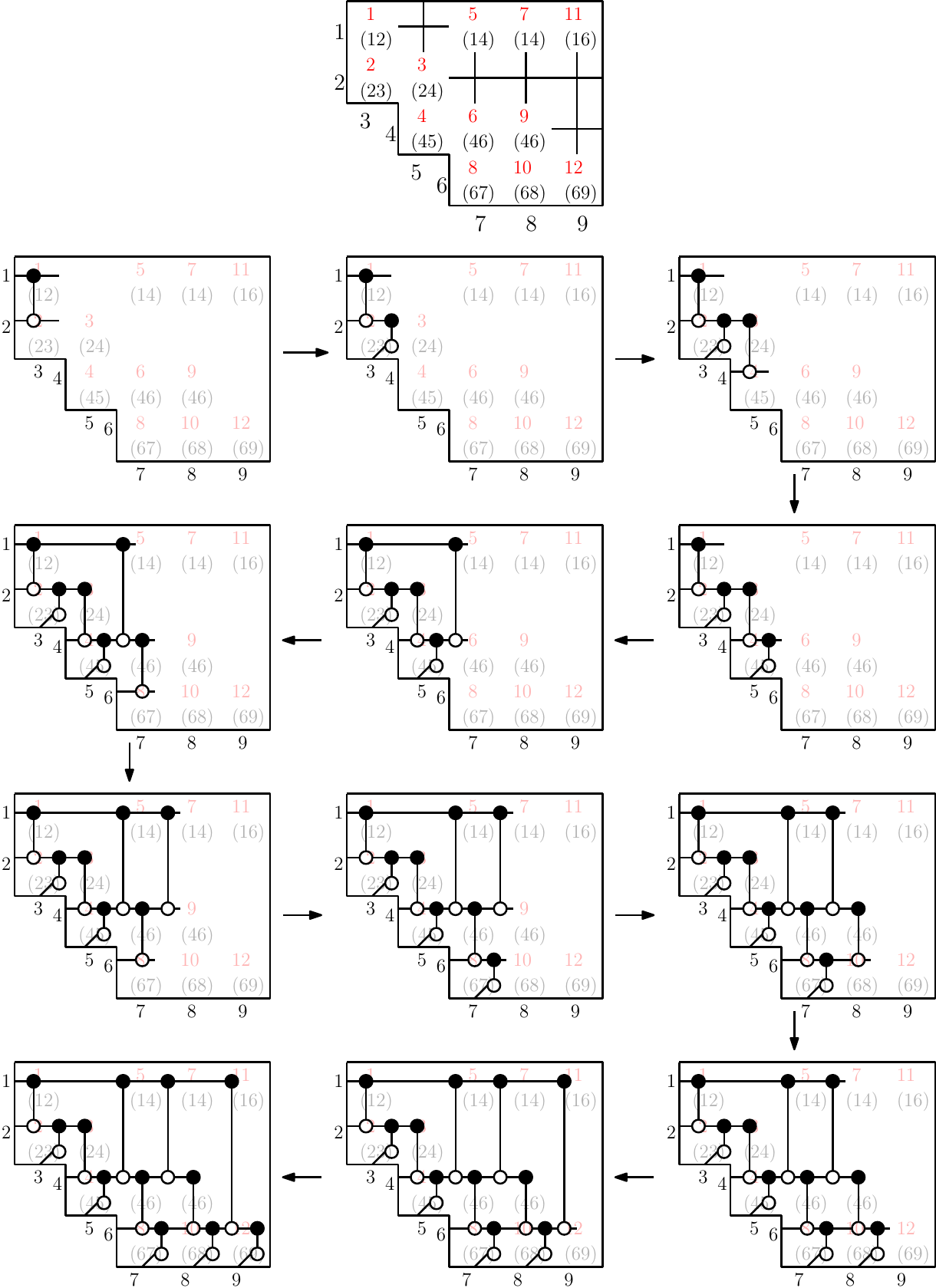}
		\caption{}
\end{figure}

\newpage

This example did not use Rule 2 for the order in which the transpositions are read on the pipe dream, which was illustrated in Figure \ref{fig:islands}. For an example that uses both Rules 1 and 2, we encourage the reader to go through the process in this section for the positroid variety with bounded affine permutation $f=\{5,7,8,10,6,9,11,12\}$ on the point $\lambda=\{1,2,3,4\}$ in Gr(4,8).

Partial answer: going through the process of creating the top pipe dream as in the first part of this section should produce:

\begin{figure}[htbp] \centering
	\includegraphics[scale=0.5,clip=true]{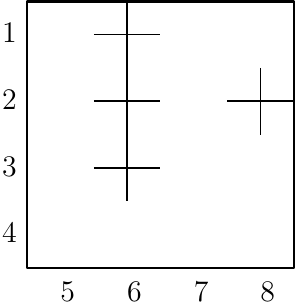}
		\caption{}
\end{figure}

Going through the same steps as Section 5.2 (page 54) of \cite{ABCGPT}, we get the following matrix representative for the positroid variety \cite{BourjMath}:

\[
\mqty(1&\alpha_{12}&0&-\alpha_{6}-\alpha_{8}&(-\alpha_{6}-\alpha_{8})\alpha_{9}&\alpha_{2}-\alpha_{6}\alpha_{7}&\alpha_{2}\alpha_{5}&\alpha_{2}\alpha_{3}&0 \\
0&1&\alpha_{11}&\alpha_{10}&0&0&0&0&0 \\
0&0&0&1&\alpha_{9}&\alpha_{4}+\alpha_{7}&\alpha_{4}\alpha_{5}&0&0 \\
0&0&0&0&0&1&\alpha_{5}&\alpha_{3}&\alpha_{1})
\]

\subsection{Other Schemes and Patches}

Obviously, there was a choice that was made in the BCFW bridge decomposition method of \cite{ABCGPT}, namely, that we began with $g(1)$, so we were dealing with the set of numbers $(3,7,5,11,8,9,13,10,15)$. However, there was nothing stopping us from starting with, say, $g(2)$, such that our permutation is $(6,4,10,7,8,12,9,14,11)$, rather than $(3,7,5,11,8,9,13,10,15)$, and doing the BCFW bridge decomposition for this shifted-by-one permutation. Ultimately, we would end up with a different matrix representative for the same positroid variety. This ends up being useful because, as noted in Section 6.4 (page 65) of \cite{ABCGPT}, all the boundaries of any positroid variety can be found by sending a variable $\alpha_i$ corresponding to a BCFW bridge to zero. 

However, there is something else we can do differently, other than take cyclic rotations of our permutation. We can also go backwards. In other words, we chose a starting or base point: this was our choice of patch/coordinates. In our example above, with $(3,7,5,11,8,9,13,10,15)$, this is like choose the choice of boundaries or endpoints, given by the bars: $...,2,0,6|3,7,5,11,8,9,13,10,15|12,16,14,20,...$. However, our algorithm also had a choice of direction: it started from the $g(1)=3$ and went $\rightarrow$ rightwards. There's nothing stopping us from modifying the BCFW bridge method so that we instead start with $g(9)=15$ and go $\leftarrow$ leftwards. This has a sort of ``opposite Schubert variety" (used in forming Richardson varieties, for instance) feel to it.

This has the slight convenience that, in contrast to the first step in the previous subsection where we had to get the top pipe dream by looking at the inverse permutation, here, if we go backwards, we want the \textit{bottom} pipe dream, so we don't need to deal with the inverse permutation. Overall, the result will be like flipping the original method across the $x=y$ line. To take a simple example, for the case of the positroid variety $\Pi_f, f=(2222)$ in siteswap language, at the point $\lambda=\{1,2\}$, the BCFW bridge decomposition gives the sequence of bridges: $(12),(23),(12),(34)$, whereas the reverse ordering gives $(34),(23),(34),(13)$. The picture is as follows. Rather than labelling according to the row numbers, from top to bottom, the transpositions are labelled according to the column numbers, from right to left. The word is read moving along the southwest to northeast diagonals, but starting from the corner of the diagram in the bottom right:

\begin{figure}[htbp] \centering
	\includegraphics[scale=0.7,clip=true]{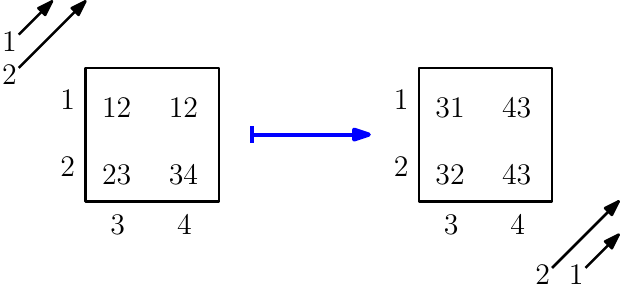}
		\caption{}
\end{figure}

\chapter{Soft Limits}\label{ch:softlimits}

\subsection{Diagonal Diagrams Recap}

In Appendix A of \cite{Flu}, we described diagrams that here we will call \textbf{diagonal diagrams.} These diagrams depict bounded affine permutations or siteswaps where, if we have a bounded affine permutation $g(i)=j$, then row $i$ of the diagram has a single $1$ in the $j$th column. Note that the diagram has $i$ increasing going south, and $j$ increasing going right. Because the bounded condition of bounded affine permutations says that $i\leq j\leq i+n$, as we increase $i$ (moving down), the places where the $1$ can be moves right, creating the northwest to southeast diagonal strip. 

We can also depict what deletion and contraction do to these diagrams. See \cite{Flu} for details - there we showed that deletion of the $i$th coordinates operates on such a diagram by sending the $1$ in the $i$th row to the left, and the new diagram with the post-deletion bounded affine permutation is obtained by a series of transpositions in a sort of L-shaped box as shown below.

 \begin{figure}[htbp]
	\includegraphics[scale=0.6,clip=true]{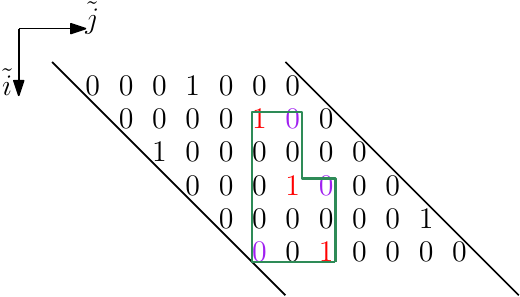}
\end{figure}

Or, step by step:

 \begin{figure}[htbp]
	\includegraphics[scale=0.6,clip=true]{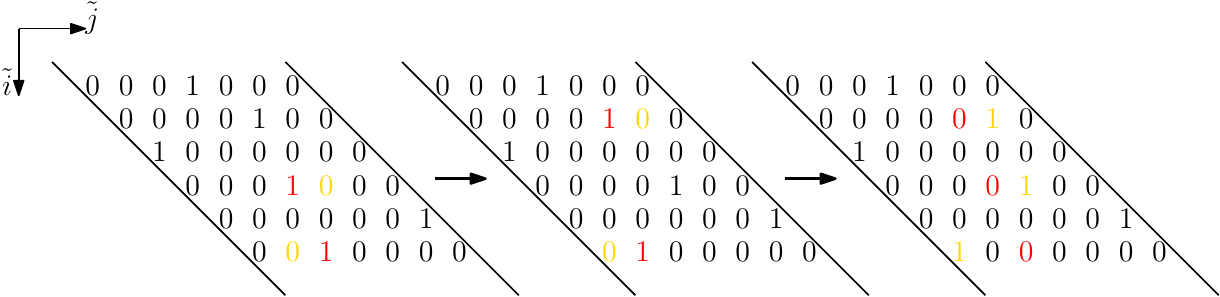}
\end{figure}
 
Unlike the pipe dreams above, these diagonal diagrams are not patch dependent: there is a single such diagram for a single bounded affine permutation.

\subsection{T-duality}\label{subsec:Tduality}

It is especially nice to see what the T-duality map is for these diagonal diagrams. This duality was first pointed out by
\cite{LukWilPar}, although we use the definition from \cite{HuYeats}:

\begin{Def}
The T-duality map from loopless decorated permutations on $[n]$ to co-loopless decorated permutations on $[n]$ is defined as \\
\[\pi\mapsto \hat{\pi}\]
\[(a_1, a_2,..., a_n ) \mapsto (a_n,a_1,...,a_{n-1})\]
where the permutations are written in one-line notation, and any fixed points in $\hat{\pi}$ are declared to be loops. That is, for a given loopless $\pi$, we have $\hat{\pi}(i)=\pi(i-1)$ where all fixed points are loops, and we call $\hat{\pi}$ the T-dual decorated permutation.
\end{Def}

As an example of what the T-duality map looks like on the diagonal diagrams, let's take an example of a bounded affine permutation given $3,8,6,5,10,7$. Writing the BAP and siteswap below it, T-duality sends:

$\mqty(1&2&3&4&5&6\\ 3&8&6&5&10&7\\2&6&3&1&5&1) \mapsto \rd{\mqty(1&2&3&4&5&6\\1&3&8&6&5& 10 \\0&1&5&2&0&4)} $

Diagramatically: 
'
\begin{figure}[htbp]
	\includegraphics[scale=0.8,clip=true]{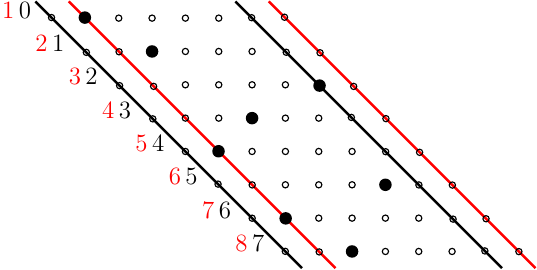}
\end{figure}

In this figure, rather than using $1$'s to depict the entries of the permutation (i.e. the spots where ``the ball lands") as above, we use the larger dots (the smaller dots are $0$'s). 

In other words, the T-duality map shifts the permutation so that it starts $1$ entry earlier (i.e. in the cyclic permutation, we start at $n$ instead of $1$). This is like changing $i$ in $f(i)$ by subtracting $-1$ from every $i$, which would be shown on the diagram by shifting the boundary of the diagram to the right by $1$: that is, the boundary moves from the black diagonal lines to the red diagonal lines. The inverse of the T-duality map would then be given by the opposite process: moving the boundary of the diagonal diagrams one unit to the left (and choosing to ``start" the permutation one row lower when writing it). It is easy to see why the T-duality map must go from loopless decorated permutation to co-loopless deocrated permutations: since any loops (that is, $f(i)=0$'s for the siteswap would get mapped to $f'(i)=-1$ for the T-dual permutation, which would not be a valid positroid; and similarly since the largest an entry can be is $f(i)=n$, this will get mapped to $n-1$ in the T-dual permutation, so it will be co-loopless.

In fact, I realized later that the T-duality map is just the map from momentum space to momentum twistor space give by mapping certain poistroid varieties in Gr(k,n) to positroid vareities in Gr(k-2,n) via a projection to the orthogonal complement of the momentum $\lambda$" (the $n$ 2-component angle spinors in the spinor helicity formalism for $n$ massless momenta, which thus form a 2-plane in $n$ dimensions). See e.g. Chapter 8 in \cite{ABCGPT}. 

It would be interesting to see what happens for the T-duality map at the level of amplitudes; that is, after mapping to the amplituhedron. The diagram below show the map from the momentum space (positive) Grassmannian to the momentum amplituhedron (c.f. \cite{DFLP19}). The T-duality map is labeled $T$ on the arrow on the left; the question is what the map $T'$ is. The original amplituhedron map (c.f. \cite{AT13}) was on momentum twistor space and is labeled by $Z$ in the diagram below. It would be nice to understand what the images of positroid varieties are from $Gr_+(k,n)$ to $Gr(k,k+m)$ as well as what, say, BCFW tiles in $Gr(k,k+m)$ pull back to, under the maps. 

$\begin{tikzcd}
Gr_+(k,n) \ar[r,"\Lambda\times\tilde{\Lambda}"] \ar[d, "T"] & Gr(k-2,k-2+\frac{m}{2})\times Gr(n-k-2,n-k-2+\frac{m}{2})  \ar[d, "T'"] \\
Gr_+(k-2,n) \ar[r, "Z"] & Gr(k,k+m)
\end{tikzcd}$ 

To understand this, we might want to work out the maps in the following diagram:

$\begin{tikzcd}
\tilde{Z} \hookrightarrow Fl(2,k,n-2;n) \ar[r,"J"] \ar[d, "Gr(k-2;n-4)"] & Gr(k;n) \\
 (\lambda, \tilde{\lambda}^\perp) \in Fl(2,n-2;n) &
\end{tikzcd}$ 

In the upper right, we have the 3-step flag variety $Fl(2,k,n-2;n)$, and $\tilde{Z} = \{ V^k: \lambda^2 \subseteq V \subseteq \tilde{\lambda}^{\perp,n-2} \}$ is the set of all $k$ planes $V$ that contain a generic 2-plane $\lambda$ and are contained in the $n-2$-dimensional plane $\tilde{\lambda}^\perp$. For a given momentum, $(\lambda, \tilde{\lambda})$, the fiber over a point in the 2-step flag variety (from forgetting $V$) $(\lambda, \tilde{\lambda}^\perp) \in Fl(2,n-2;n)$ is $Gr(k-2;n-4)$. 
We want to take a positroid variety $\Pi_f$ for a leading singularity (that is, of dimension $2n-4$) in $Gr(k;n)$ and consider its pullback $J^{-1}(\Pi_f)$ to $Fl(2,k,n-2;n)$. The map $J$ simply forgets the 2-plane and $n-2$ plane. An important question is whether the intersection $J^{-1}(\Pi_f) \cap \tilde{Z}$ is a positroid variety, but the fiber $Gr(k-2,n-4)$ does not have coordinates, and there doesn't seem to be a way to do give it coordinates invariantly, so whether the intersection is a positroid variety may be indeterminable. If the positroid variety $\Pi_f$ is $2n-4$ dimensional, then the preimage $J^{-1}(\Pi_f)$ mapped down to $Fl(2,n-2)$ will cover it- the question is what the cardinality of this cover is (how many times does it cover it); equivalently, this is the number of components of $J^{-1}(\Pi_f) \cap \tilde{Z}$.

We do not discuss the amplituhedron in any detail in this thesis; we refer the reader to \cite{AHTr}; a simple introduction to the amplituhedron can also be found in the Amplituhedron notes at \url{josephflueg.github.io}.

\section{Soft factors}

One place where deletion and contraction for positroid varieties has occurred in the literature is in the discussion of inverse soft limits (see Section 16 of \cite{ABCGPT} and \cite{HeZh}). The idea of soft limits comes from the Weinberg soft theorems, which gives formulas for the amplitude as the momentum of one particle (say $p_i$) goes to zero in terms of the amplitude of the rest of the particles, multiplied by a so-called soft-factor (which goes to $1$ as $p_i \to 0$). The idea then, is whether one can go in reverse. 

A ``particle going soft" can be seen as a sort of deletion or contraction kind of phenomenon. One question then, is whether it is possible to rebuild the positroid variety (more specifically, the siteswap $f$) from the data of all, or just some, of its deletions and contractions. Theoretically, it should be possible to reconstruct the pipe dream in most cases from all its deletions and contractions. In fact, if there is a pair $a,b$, where $a$ and $b$ can be either a row or column on a pipe dream that are ``disjoint," meaning that every maximal rectangle (or just rectangle) that contains at least one tile from $a$ is disjoint from every rectangle that contains at least one tile from $b$, then it should be relatively straightforward to rebuild the original pipe dream from the two pipe dreams obtained from deleting/contracting $a$ and $b$, because this means that all the crosses in column $a$ are present in the deletion/contraction of $b$ and vice-versa. Extending this logic, then, the only case that is ambiguous is the case coming from the atomic pairs: a square pipe dream that has a cross along the main diagonal and completely elbow in the diagonal both directly above and below is indeterminate, but every deletion and contraction eliminates that tile, so it cannot be determined whether that cross existed before the deletion/contraction or not (note that these pipe dreams do not have to be atomic since they can involve other moves outside of these 3 main diagonals). However, doing this is practice for larger, complicated pipe dreams may be more difficult, and we do not have a technique for explicitly reconstructing pipe dreams from their deletion and contractions.

If the discussion in the previous paragraph was a bit hard to follow, we give an example here. Let's suppose that we are in $n=4$, on the patch $\lambda=\{1,2\}$. The pipe dream looks as follows (for the affine pipe dream, it will be an infinite number of such squares connected along the northwest and southeast corners; we will drop the numbering in our later figures):
\scalebox{0.6}{
\begin{tikzpicture}
\draw[black, thick] (0,0) -- (2,0) -- (2,2) -- (0,2) -- (0,0);
\draw[dotted] (0,1) -- (2,1);
\draw[dotted] (1,0) -- (1,2);
\node[left] at (0,1.5) {$1$};
\node[left] at (0,0.5) {$2$};
\node[below] at (0.5,0) {$3$};
\node[below] at (1.5,0) {$4$};
\end{tikzpicture}
}

Let's call the original pipe dream, prior to any deletions or contraction, $P_{\lambda,f}$ for $\lambda=\{1,2\}$ and some siteswap $f$). We want to recover $P_{\lambda,f}$, which will tell us what $f$ is. Suppose we know that the deletion of the 3rd entry produces (the empty boxes are elbow tiles) $Del_3(P_{\lambda,f})=$ \scalebox{0.6}{\begin{tikzpicture}
\draw[black, thick] (0,0) -- (2,0) -- (2,2) -- (0,2) -- (0,0);
\draw[dotted] (0,1) -- (2,1);
\draw[dotted] (1,0) -- (1,2);
\draw[black,thick] (0.5,0)--(0.5,2);
\draw[black,thick] (0,0.5)--(2,0.5);
\draw[black,thick] (0,1.5)--(1,1.5);
\draw[red,thick] (1.5,0)--(1.5,1);
\draw[red,thick] (1,0.5)--(2,0.5);
\end{tikzpicture}}
The cross in red is not movable (in particular it can't move to the 3rd column) and it does not participate in the deletion process for the 3rd entry, so it must be there in the original pipe dream prior to $Del_3$. This also implies that the top right corner of $P_{\lambda,f}$ must be an elbow, because the red cross prevents it from moving. Then, if we performed a contraction of leg 4 and got \scalebox{0.6}{\begin{tikzpicture}
\draw[black, thick] (0,0) -- (2,0) -- (2,2) -- (0,2) -- (0,0);
\draw[dotted] (0,1) -- (2,1);
\draw[dotted] (1,0) -- (1,2);
\draw[black,thick] (1.5,0)--(1.5,2);
\draw[black,thick] (1,0.5)--(2,0.5);
\draw[black,thick] (1,1.5)--(2,1.5);
\end{tikzpicture}}, then we would know that $P_{\lambda,f}=$ \scalebox{0.6}{\begin{tikzpicture}
\draw[black, thick] (0,0) -- (2,0) -- (2,2) -- (0,2) -- (0,0);
\draw[dotted] (0,1) -- (2,1);
\draw[dotted] (1,0) -- (1,2);
\draw[black,thick] (1.5,0)--(1.5,1);
\draw[black,thick] (1,0.5)--(2,0.5);
\end{tikzpicture}}, and this is general in this case: that is, since the red cross blocks any cross from moving between the 3rd and 4th column in $P_{\lambda,f}$, whatever pattern we see in the 2 boxes on the left in $Del_4(P_{\lambda,f})$ will be exactly the pattern for $P_{\lambda,f}$. 

Now, let's consider a different case where $Del_3(P_{\lambda,f})=$ \scalebox{0.6}{\begin{tikzpicture}
\draw[black, thick] (0,0) -- (2,0) -- (2,2) -- (0,2) -- (0,0);
\draw[dotted] (0,1) -- (2,1);
\draw[dotted] (1,0) -- (1,2);
\draw[black,thick] (0.5,0)--(0.5,2);
\draw[black,thick] (0,0.5)--(1,0.5);
\draw[black,thick] (0,1.5)--(2,1.5);
\draw[red,thick] (1.5,1)--(1.5,2);
\draw[red,thick] (1,1.5)--(2,1.5);
\end{tikzpicture}}
Here, the difference from the one above is that the red cross is in the upper right box. As before we know that the structure of $Del_3(P_{\lambda,f})$ implies that the red cross is in $P_{\lambda,f}$ as well. However, it cannot be \scalebox{0.6}{\begin{tikzpicture}
\draw[black, thick] (0,0) -- (2,0) -- (2,2) -- (0,2) -- (0,0);
\draw[dotted] (0,1) -- (2,1);
\draw[dotted] (1,0) -- (1,2);
\draw[black,thick] (1.5,1)--(1.5,2);
\draw[black,thick] (1,1.5)--(2,1.5);
\end{tikzpicture}}, otherwise $Del_3$ would move this cross to the bottom left so that $Del_3(P_{\lambda,f})=$ \scalebox{0.6}{\begin{tikzpicture}
\draw[black, thick] (0,0) -- (2,0) -- (2,2) -- (0,2) -- (0,0);
\draw[dotted] (0,1) -- (2,1);
\draw[dotted] (1,0) -- (1,2);
\draw[black,thick] (0.5,0)--(0.5,2);
\draw[black,thick] (0,0.5)--(1,0.5);
\draw[black,thick] (0,1.5)--(1,1.5);
\end{tikzpicture}} (see \ref{section:delcontrpipedream} for proof/details). It also cannot be \scalebox{0.6}{\begin{tikzpicture}
\draw[black, thick] (0,0) -- (2,0) -- (2,2) -- (0,2) -- (0,0);
\draw[dotted] (0,1) -- (2,1);
\draw[dotted] (1,0) -- (1,2);
\draw[black,thick] (1.5,1)--(1.5,2);
\draw[black,thick] (1,1.5)--(2,1.5);
\draw[black,thick] (0,0.5)--(1,0.5);
\draw[black,thick] (0.5,0)--(0.5,1);
\end{tikzpicture}}
since this would introduce a double crossing, so the pipe dream would not be reduced.
Thus, the two possibilities are (A) \scalebox{0.6}{\begin{tikzpicture}
\draw[black, thick] (0,0) -- (2,0) -- (2,2) -- (0,2) -- (0,0);
\draw[dotted] (0,1) -- (2,1);
\draw[dotted] (1,0) -- (1,2);
\draw[black,thick] (0,1.5)--(2,1.5);
\draw[black,thick] (0.5,1)--(0.5,2);
\draw[black,thick] (1.5,1)--(1.5,2);
\end{tikzpicture}} or (B) \scalebox{0.6}{\begin{tikzpicture}
\draw[black, thick] (0,0) -- (2,0) -- (2,2) -- (0,2) -- (0,0);
\draw[dotted] (0,1) -- (2,1);
\draw[dotted] (1,0) -- (1,2);
\draw[black,thick] (0,1.5)--(2,1.5);
\draw[black,thick] (0.5,0)--(0.5,2);
\draw[black,thick] (1.5,1)--(1.5,2);
\draw[black,thick] (0,0.5)--(1,0.5);
\end{tikzpicture}} As before, if we look at $Del_4(P_{\lambda,f})$, this will leave the two boxes on the left untouched, so this will be sufficient to determine $P_{\lambda,f}$. We could have also used $Contr_1(P_{\lambda,f})$ instead, since (A) and (B) are the same along the top row and only differ in the lower left square, which would be seen by $Contr_1(P_{\lambda,f})$. This is in contrast to the first example, where we determined $P_{\lambda,f}$ via $Del_3$ and $Del_4$, but no single contraction would have been sufficient, as all possible combinations of crosses in the left column were possible. 

Obviously this example looked at pairs of positroid varieties (or their permutations) and points $(f,\lambda)$. We could consider the question just for positroid varieties $f$. We know that in the special case of $f$ (siteswaps) with $f(i)=2$ or $f(i)-n-2$ that there is a decomposition of the amplitude into a smaller permutation along with a ``soft factor" that we recap next. However, even here, the decomposition is diagram-by-diagram, rather than at the level of the whole amplitude. It would be interesting to investigate whether there are properties that hold at the level of the whole amplitude, in the consideration of soft limits.

We note that the deletions/contractions are nonempty when deleting $i$ for $f(i)=n$ or contracting $i$ for $f(i)=0$. This was noted in \cite{SOh} or simply by noting that any switch in the diagonal diagram above would result in a siteswap number larger than $n$. One might wonder whether we might, for example trying to delete $f(i)=n-1$; one could for example have $f(i-1)=0$ and thus have to do to a switch with a higher row that would result in a number larger than $n$; however, this is impossible in a valid siteswap because this would imply that in the diagonal diagram the $1$ for the $n-1$ lies above the $1$ for $f(i-1)=0$, etc. Perhaps more interestingly, there are also conditions on where we can place an inverse soft factor. We cannot put a k-preserving inverse soft factor on the $i$th slot if the $i-1$th slot has siteswap $f(i-1)=0$ or $f(i-1)=1$. If it's $0$, then there no way to get $f(i-1)=2$, and if it's $1$, then adding in the extra $i$th particle automatically puts $f(i-1)$ as $2$, so the reverse deletion process just involves one swap (a single BCFW bridge).

The inverse soft factors can be applied to an on-shell (plabic) diagram through one of the two following additions.

For the k-preserving inverse holomorphic soft factor (which is undone by deletion), the addition looks like the following:
 \begin{figure}[htbp]
	\includegraphics[scale=0.6,clip=true]{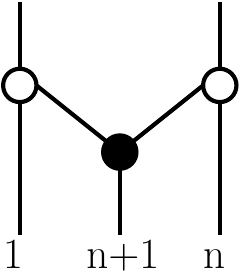}
\end{figure}

For the k-increasing inverse holomorphic soft factor (which is undone by contraction), the addition looks like the following:

 \begin{figure}[htbp]
	\includegraphics[scale=0.6,clip=true]{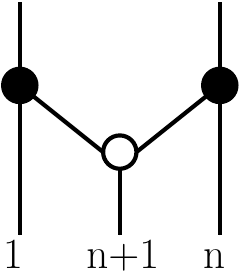}
\end{figure}

We give the rule of inverse soft factors, showing that they are undone by deletion and contraction, and that they correspond to the sequence of two BCFW bridges displayed in \cite{ABCGPT} or \cite{HeZh}.

The rule is that a $k$-preserving inverse-soft factor was applied when $\sigma(i-1)=i+1$. This corresponds to $f(i-1)=2$ in siteswap language. A $k$-increasing inverse-soft factor was applied when $\sigma(i+1)=i-1$. This corresponds to $f(i+1)=n-2$ in siteswap language. Thus, if we undo the $i$th leg via deletion or contraction, this corresponds to being able to do the deletion or contraction in 2 steps or flip on a diagonal diagram, which accords with the notion of adding in 2 degrees of freedom with 2 BCFW bridges:

 \begin{figure}[htbp]
	\includegraphics[scale=0.6,clip=true]{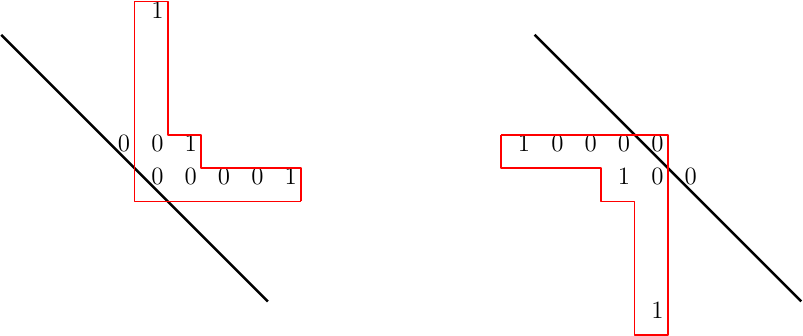}
\end{figure}

In the following discussion, we prove these statements; specifically that deletion of the $i$th leg is a left inverse of a $k$-preserving soft factor (that the deletion undoes the addition of the soft factor, returning to the original permutation), and that contraction of the $i$th leg is a left inverse of a $k$-increasing soft factor. 

For notational purposes, let's call such a deletion or contraction a \textbf{soft-deletion} or \textbf{soft-contraction}.

We can analyze what happens to the bounded affine permutation upon addition of an inverse soft factor:

\begin{figure}[htbp]
	\includegraphics[scale=0.6,clip=true]{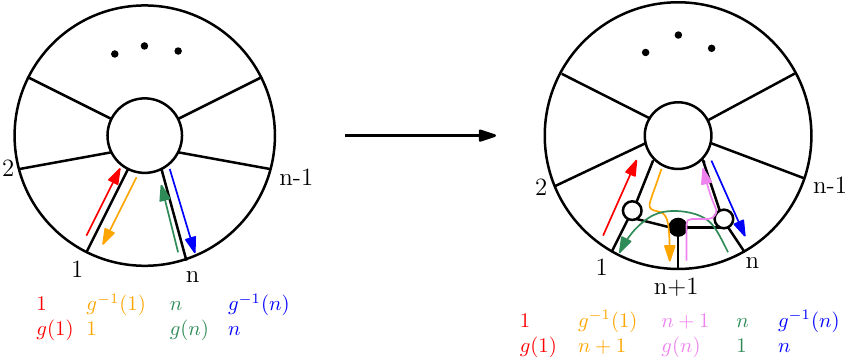}
\end{figure}

However, it might be easier to see what's going on if we think about the removal of a soft factor: 

\begin{figure}[htbp]
	\includegraphics[scale=0.6,clip=true]{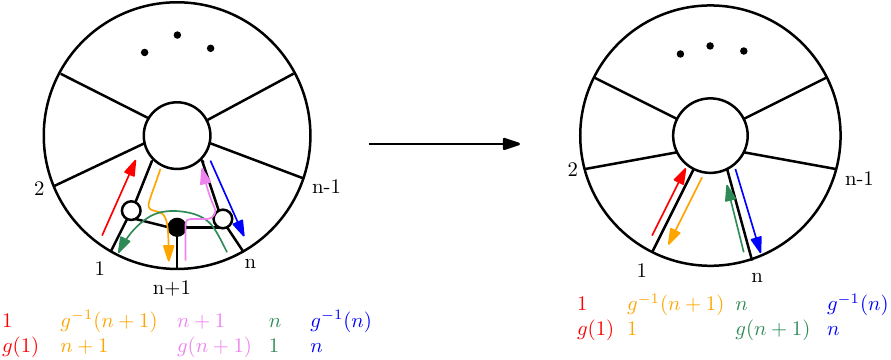}
\end{figure}

Here's an examle of what the addition of a soft factor does on the diagonal diagrams to the siteswap 340362 from earlier: 

\begin{figure}[htbp]
	\includegraphics[scale=0.45,clip=true]{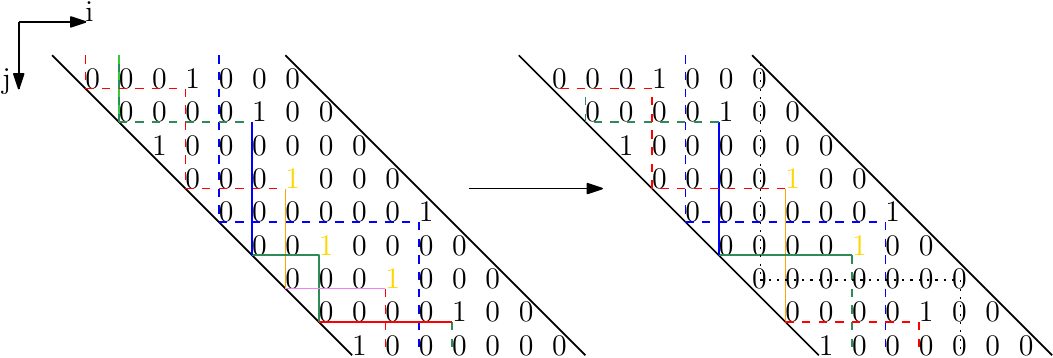}
\end{figure}

We now prove/describe what is happening in the diagrams above in more detail. 

Consider the on-shell diagram where leg $n+1$ has been added between legs $n$ and $1$, by attaching black vertices to $n$ and $1$ that are then attached to a white vertex coming off of leg $n+1$. This ends of being a k-increasing holomorphic soft factor. We have the following: 

$\mqty(1 & \cdots & i & \cdots & j & \cdots & n \\ \sigma(1) & \cdots & 1 & \cdots & n & \cdots &  \sigma(n))$ where we think of $\sigma$ as being the bounded affine permutation. Here, $i$ an $j$ are representative since they could be in reversed order or one of them equal to 1 or n.

For the siteswap, we would have $\mqty(1 & \cdots & i & \cdots & j & \cdots & n \\ f(1) & \cdots & 1-i & \cdots & n-j & \cdots &  f(n))$, obvioulsy with $f(1)=\sigma(1)-1, f(n)=\sigma(n)-n$, and everything considered modulo $n$. 

Consider the left-right paths (c.f. \cite{Post}) for the plabic diagram. After we add in these two BCFW bridges (the k-increasing holomorphic soft factor), the numbers that stay the same, based on tracing out the left-right paths, are $\sigma(i)=1$ and $\sigma(n)$. The others change as follows:

$\mqty(1 & \cdots & i & \cdots & j & \cdots & n & n+1 \\  n & \cdots & 1 & \cdots & n+1 & \cdots &  \sigma(n) & \sigma(1) )$
 
If we shift to siteswap notation, this is:

$\mqty(1 & \cdots & i & \cdots & j & \cdots & n & n+1 \\  n-1 & \cdots & 1-i & \cdots & f(j)+1 & \cdots &  f(n) & f(1)+1 )$

Additionally, every number that was less than $n$ that had $\sigma$ greater than $n$ (every excedence) now needs to have $+1$ added to it, since it now goes past leg $n+1$ as well. 

We see that this is precisely the inverse process of contraction of leg $n+1$: If we contract leg $n+1$, the $j$ with $\sigma(j) = n+1$ goes back by $1$ to $\sigma(j)=n$; this is the lower side of the upside down $L$ shape. Leg $1$ used to have $\sigma(1)=n$ and now becomes $\sigma(1)=\sigma(n+1)$, whereas $n+1$ goes to $n+1$. 

If we want to do a series of inverse soft factors, it is nice to go into siteswap notation, at least temporarily, because it is permutation invariant, so we can add in soft factors between different legs arbitrarily. Here, we go through an example that illustrates how to obtain the sitewap or bounded affine permutation from a siteswap $f$ after inserting an inverse soft factor between legs $i$ and $i+1$. 

Given that we are already familiar with deletion and contraction, let us do an example of removing an inverse soft factor in way that will make it easy to see how to go in reverse. Let us take $f'=234252$ in $Gr(3,6)$. There is a $2$ in the $4th$ slot, so we can think of it as a $k$-preserving inverse soft factor that would be removed by deleting the entry after it (the $5$th slot). To delete, let us first move the $5th$ slot to the end: $223425$.  Now let's get it's bounded affine permutation by adding 123456: $3,4,6,8,7,11$. Let's spot the one that comes into the last (here $n=6$) entry: this would be the 3rd: $3,4,\rd{6},8,7,11$. Because the entry before the last one is equal to $2$, other than the 3rd, this is the only one that changes, and it happens via a sort of rotation of the numbers (as described next); this is the content of Appendix A of \cite{Flu}: $3,4,\rd{6},8,\rd{7},\rd{11} \mapsto 3,4,\rd{7},8,\rd{11},\rd{6}$. We project to the smaller Grassmannian by eliminating the last entry and decreasing all the entries that jump past the last one (these are the numbers greater than $6$ in this example) by $1$: $3,4,6,7,10$. The siteswap is $2,2,3,3,5$. We can now delete or contract any additional columns. 

If we wanted to do contraction, let us summarize what is different than contraction. We need to choose the $n-2$ siteswap number and now change the one before it. In the example above with $f=223425$, the $4$ is in the 4th slot, so we would need to contract the 3rd entry. Rotate now to the beginning: $f$ becomes $342522$. The bounded affine permutation is $465978$. Now, we want they first entry to go to $n+1$, which here is $7$: $4659\rd{7}8$. Contraction now rotates the first two with the red, moving cyclically to the right (instead of to the left, as in the case of deletion): $\rd{46}59\rd{7}8 \mapsto \rd{74}5968$. The siteswap is $622512$, and if we project to the smaller Grassmannian, we will need to subtract $1$ from every entry that reaches or goes past the first entry: $22411$.

Let take the bounded affine permutation $f=22335$. Suppose that we want to insert a k-preserving inverse soft factor (so it would be undone by deletion). Let's insert it at the end. Since it's k-preserving, we add a $0$: $223350$. We go to the bounded affine permutation notation: $3,4,6,7,10,6$. Let's add $1$ to every number that is greater than or equal to $n=6$, since we stuck in an extra particle, so we now need to move them past it: $3,4,7,8,11,6$. Let's highlight the one that is $n+1$: $3,4,\rd{7},8,11,6$, so this one rotates with the last two numbers, but now in the opposite direction (to the right): $3,4,\rd{7},8,\rd{11},{6} \mapsto 3,4,\rd{6},8,\rd{7},{11}$.

Here is how we can only do this for the siteswap directly. Take $f=22335$ again, and again add $0$ to the end to insert a k-preserving inverse soft factor: $223350$. Add $1$ to every number that goes past the $0$ we just added: $223460$. The last $2$ numbers are $6$ and $0$; we switch the two, and change the $6\mapsto 6-1=5$ and $0\mapsto 2$. The final answer is $223425$. (Note, following the process above in detail, we would have need to add $+1$ to every number \textit{including} the one landing on the last entry: $22\rd{4}460$, and then do the two swaps $22\rd{4}4\rd{60} \mapsto 22\rd{4}4\rd{15} \mapsto 22\rd{3}4\rd{25}$, but it may be simpler to think about it in the abbreviated way we just described.) In summary, if we are adding a k-preserving inverse soft factor to the end (let's say this is the $n+1$ entry) [note that this adds the extra particle at position $n+1$, so rotating/unrotating maybe necessary depending on where we want to add the extra particle]:: \\
(A) Tack a $0$ on the end (that is, $f(n+1)=0$) \\
(B) Add $+1$ to every $f(i)$ that goes past $n+1$ \\
(C) Change $f(n+1)$ to be $f(n)-1$ \\
(D) Change $f(n)$ to be $f(n)=2$ \\
Recall that the necessary condition for this algorithm to work (that is, for deletion to be a left inverse of adding a k-preserving inverse soft factor) is that if we are inserting the soft factor on leg $n+1$, then we need $f(n) \geq 2$, so that once we add the extra particle, it becomes at least $3$ and thus 2 switches are necessary to turn $f(n)$ into $2$. Otherwise, if $f(n)=1$, it will become $f(n)=2$ when we add the extra particle so that no switches are possible, or else if $f(n)=0$, it will not participate in the deletion process.


For a k-increasing inverse holomorphic soft factor (contraction) [note that this adds the extra particle at position $1$, so rotating/unrotating may be necessary depending on where we want to add the extra particle]: \\
(A) Tack a $n+1$ at the beginning (that is, $f'(1)=n+1$, and the original $f(1)$ is now $f'(2)$ ) \\
(B) Add $+1$ to every $f(i)$ that reaches or goes past $1$ \\
(C) Find the one that lands on $n$ (the bounded affine permutation lands on $n$) and add $+1$ to this one \\
(D) Change $f'(1)$ to be $f'(2)+1$ \\
(E) Change $f'(2)$ to be $n-1$ (so it lands on $n+1$)\\
Note that this is, equivalent but not exactly the same as the dual to the method above. A completely dual method would go by tacking on $n$ to the beginning, and then adding $+1$ to all numbers (this would be like enlargening the diagonal strip by extending it to the left on the left boundary), and then subtracting $1$ from the entries that \textit{do not} jump over $1$. \\
Recall that the necessary condition for this algorithm to work (that is, for contraction to be a left inverse of adding a k-increasing inverse soft factor) is that if we are inserting the soft factor on leg $1$, then we need $f(2) \leq n-2$; the addition of the extra particle widens the diagonal diagram so that $n-2 = (n+1)-3$, and thus we can apply two switches to get $f'(2)=(n+1)-2$.

By our theorem on the irreducibility of the deletion and contraction, deletions or contractions of positroid variety $\Pi_f$ is just intersections of $\Pi_f$ with various $\Pi_{contr,i}$ and $\Pi_{del,i}$.

Then, we are just intersecting with a bunch of these $\Pi_{del,i}$ and $\Pi_{contr,i}$, and taking intersections is commutative. Based on this notation, we do not project to a smaller Grassmannian after each contraction or deletion (that is, in the case of removing inverse (anti)holomorphic soft factors, we usually project to a smaller Grassmannian, e.g. in the example above, we took the siteswap 223425 and did the contraction on the 5th slot which is the intersection $\cap \Pi_{5,contr}$ which produces: $\Pi_{223425} \cap \Pi_{5,contr}=\Pi_{224460}$, which at this point is still in Gr(3,6); the projection to Gr(3,5) brings it to $\Pi_{22335}$). However, commutativity should still hold even if we project out since the projection is just removing entire rows and columns of zeros which are skipped over when doing the switches in the diagonal diagrams anyways. 

Because we showed that the inverse soft factors are an inverse process to deletion/contraction; that is, if we call deletion $del_i$ an operation on a siteswap $f$, and similarly $contr_i$ for contraction, while we call $pres_i$ the process that adds a k-increasing inverse holomorphic soft factor on leg $i$, and $inc_i$ the process that adds a k-increasing inverse holomorphic soft factor on leg $i$. We have shown that $del_{i}(pres_i(f))=f$, and, if $f(i-1)=2$ (siteswap notation), then $pres_i(del_i(f))=f$. We have proven that $del_{i_1} \cdots contr_{j_1} \cdots del_{i_r} \cdots contr_{j_s}$ can be done in any order (presumably, all are different, but $del_i,del_j$ for $i=j$ does nothing the second time, and $del_i, contr_i$ for the same $i$ would just be empty. Therefore, if we have a set $f' = pres_{i_1} \cdots inc_{j_1} \cdots pres_{i_r} \cdots inc_{j_s} (f)$, then we know $del_{i_1} \cdots contr_{j_1} \cdots del_{i_r} \cdots contr_{j_s} (f') = f$. We can scramble the order of the $del_i$ and $contr_i$, and this does not change. We can then ask the question if this implies that we can similarly scramble the order of the $pres_i$ and $inc_i$. Unfortunately, the answer is no.

The insertion of soft factors is not a commutative operation. Let's take the example of $f=2222$, where we add a k-increasing soft factor at the beginning and a k-preserving soft factor at the end. We do this in both orders:

$inc_1(f): 2222 \mapsto 52222 \mapsto 52233 \mapsto 52333 \mapsto 33333$ \\
$pres_6(inc_1(f)): 33333 \mapsto 333330 \mapsto 333440 \mapsto 333423$

$pres_5(f): 2222 \mapsto 22220 \mapsto 22230 \mapsto 22222$ \\
$inc_1(pres_5(f)): 22222 \mapsto 622222 \mapsto 622233 \mapsto 622333 \mapsto 342333$

These are different positroid varieties (they are isomorphic, but not the same, as we could add in an additional soft factor at the same place in both that makes them no longer isomorphic). The problem is that changing the order of deletion $del_i$ and contraction $contr_i$ does not preserve the entries $f(i-1)$ or $f(i+1)$. And we know that we only have, for example, $pres_i(del_i(f))=f$ if $f(i-1)=2$. This is evident in the example above. Take the second result: $342333$ - we can go backwards and contract the entry $1$ because entry $f(2)$ equals $n-2=6-2=4$; however, if we try to delete the last (6th) entry, because the $f(5)$ is not equal to $2$ it does not corresponding to an inverse soft factor. Now, commutativity of deletion and contraction still holds: we can take $\Pi_f$ for $342333$ and delete the last entry before contracting the first entry, and the result will again be $\Pi_f$ for $f=2222$; the problem is it will not be possible to go in reverse to go back to $342333$ via a k-increasing and then a k-preserving inverse holomorphic soft factor: the calculation above shows that this results in a different positroid variety, and the issue stems from the fact that $342333$ does not have $f(5)=2$ so deletion of the last entry is not inverse to a k-preserving soft factor.

We note that the methods above for obtaining the image positroid vareity after adding in an inverse soft factor are valid as long as the two legs to which the soft facdtor are being attached are not bubbles: in other words, if we are adding in a new particle at $n+1$, so that the soft factor is being attached on legs $n$ and $1$, then neither leg $n$ or $1$ can have a siteswap number of $0$ or $n$. In other words, the soft factor must be attached to a connected piece. Otherwise, if one of the legs is a $0$ or $n$, then adding the soft factor would be like adding just a BCFW bridge.

There is another case where the two legs that the soft factor is attached to are both connected at the same black or white vertex - this special case will be fully fleshed out later...

In terms of getting the form, we can ask whether we are ``building up" or "breaking down." On one hand, the form for inverse soft factors is obtained by taking a bounded affine permutation, and finding spots (where $\sigma(i)=2$ or $n-2$) and doing deletion and contraction until we get to a bounded affine permutation for which we know the differential form; this seems to be a notion of "breaking down." However, we then obtain the form by taking the differential form we know, and wedging with the inverse soft factors until we get to a differential form of the original bounded affine permutation we were interested in. The BCFW bridge construction is similar: we take the bounded affine permutation and breaking down by applying swaps until we get to a point, but then build up the C-matrix by applying the BCFW bridges in reverse to the point. In fact, we can even think about the BCFW bridges as a breaking down, where we eliminate the boxes that are leftover from the top pipe dream one-by-one in the prescribed order. 

\subsection{Soft Factors and Pipe Dreams}

Here, we prove a simple proposition regarding the removal of a soft factor using pipe dreams.

Let us remind ourselves that the number of crossings in the diagonal diagrams vs. pipe dreams are related in the opposite way. That is, for the diagonal diagrams, each crossing represents an increase in the dimension of the associated variety, whereas for pipe dreams, each additional cross represents a decrease in the dimension of the associated variety. The affine pipe dreams always have $k(n-k)$ (the dimension of the Grassmannian) number of tiles, and when every tile is a cross, this is a point, a 0-dimensional variety (for the non-affine pipe dreams, as noted above (in section????) it is embedded in an affine pipe dreams such that the portion of the pipe dream underneath the cut-out boundary consists entirely of crosses). The diagonal diagrams have a minimum of $k^2$ crossings for the 0-dimensional points - where all $k$ lines are of length $n$ so all $k$ lines intersect 

\begin{Prop}
The removal of an anti =holomorphic inverse soft factor of the type mentioned above (that is, performing a soft-deletion) must decrease the dimension of the variety by at least two.
\end{Prop}
\begin{proof}
We noted above that each additional crossing decreases the dimension of the variety by one in a pipe dream. So this proposition is saying that the addition of the soft factor must add at least two crosses to a pipe dream.

We've showed in previous sections of this thesis (c.f. Proposition \ref{prop:delpipedream}) that deletion is implemented on a pipe dream by moving the maximum number of crosses that that column (and note that we are assuming, since we are trying to remove an antiholomorphic soft factor, that if the soft factor is on leg $i$, that the $i$th boundary of the affine pipe dream is horizontal (so that filling in that corresponding column of the pipe dream creates a pipe that takes $i$ vertically upwards to $i$)); after moving the maximum number of crosses to that column, we then fill in the remainder of the column with crosses. 

So our state then becomes, once we have chosen a pipe dream representative that maximizes the number of crosses in teh $i$th column, this representative will have at least 2 elbow tiles in that column. We prove this as follows: 

Look at the $i$th column. It cannot consist of all crosses, otherwise $f(i)=i$ already, so there would be nothing to delete (no soft factor to be removed). So there must be an elbow somewhere along the column, e.g.:

 \begin{figure}[htbp]
	\includegraphics[scale=0.6,clip=true]{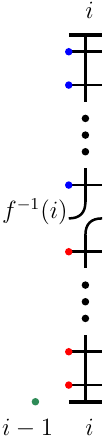}
  \end{figure}

Consider the pipe starting at $i-1$. Its starting point is represented by the green dot in the figure. We know that $f(i-1)=i+1$ (as a bounded affine permutation) since we are performing a soft-deletion. And because $f(i)\neq i$, we must have $f(i)>f(i-1)$. Since two pipes cannot cross and then uncross in a reduced pipe dream, this implies that the $i-1$ and $i$ pipes cannot cross. Therefore, pipe $i-1$ cannot enter any of the red dots. 

Additionally, $f(i-1)\neq i$, so the $i-1$ pipe cannot enter the single elbow tile in the column, otherwise it would end up at $i$. Therefore, we can conclude that pipe $i-1$ must enter through one of the blue dots. However, this produces a contradiction: since $f^{-1}(i)<i-1$, in order to reach one of the blue dots, the pipe $i-1$ must cross the $f^{-1}(i)$ pipe, and then entering via one of the blue dots means that it will cross the $f^{-1}(i)$ pipe again, producing a double crossing. 

Therefore, we must have a 2nd elbow tile along the $i$-th column if we are performing a soft-deletion.

\end{proof}

The proof for the holomorphic (k increasing) soft factor is analogous, with everything reflected across the x=y axis.

\chapter{Smoothness and Physics}\label{ch:smoothphysics}

One can ask the question of whether the singularity structure of positroid varieties is connected to the physics of amplitudes. If so, this would be a big motivation towards doing more work to understand the structure of the singular locus of positroid varieties. There are a couple calculations that suggests the answer is no. 

\section{Equivariant Cohomology Calculation}

There are so-called ``homological identities" for the amplitudes. The easiest textbook example is 6-point NMHV (see \cite{EH} chapter 9), which is, in terms of positroid varieties, `333342"+``423333" + ``334233" = ``342333" + ``333423" + ``233334". Each side involves the positroid variety $\Pi_f$ where $f$ is some cyclic rotation of $333342$ along with the cyclic rotations of $f$ by $2$ slots.
If such an identity holds at the level of the positroid varieies, then it should appear when we calculate equivariant cohomology. If, using the techniques from earlier in the thesis, we calcualte equivariant cohomology restricted to the point $\{1,3,5\}$, the result for the 3 positroid varieties on the left side of the equality are: $z_2-z_3-z_5+z_6$, $z_1-z_2-z_4+z_5$, and $z_1+z_3-z_4-z_6$, respectively. When we restrict to the circle acting by dilation, each of these becoems the polynomial $2h$, so the point $\{1,3,5\}$ is a singular point with multiplicity 2 on all 3 of these positroid vareities. Their total sum is $2(z_1 - z_3 + z_3 - z_4 + z_5 - z_6)$ and 6h after restriction. ON the right side, the calculation of equivariatn cohomology restrict to $\{1,3,5\}$ gives us $(-z_2+z_5)$, $(z_3-z_6)$, and $(z_1-z_4)$ for the 3 positroid varieties. Each term beomes $h$ when we restrict ot the $\C^\times$ cohomology, so $\lambda =\{1,3,5\}$ is a smooth point on all 3 of the positroid varieties on the right of the equality. The total sum is $z_1 - z_3 + z_3 - z_4 + z_5 - z_6$, and $3h$ after restriction. 

This is not surprising, since there are two main problems with this calculation: \\
(1) The homological identity is a sum on the leading singularities, which were calculated using the link formula (for $n$ external momenta with helicity $k$, and $\Gamma_f$ the positroid cell with positroid $f$, see Section \ref{sec:leadsing}):
\[L_{n,k}(f) =\oint_{C\subset \Gamma_f} \frac{d^{k\times n} C}{vol(GL(k))}\frac{1}{(1\cdots k)\cdots(n\cdots k-1)} \delta^{k\times 4}(C\cdot \tilde{\eta})\delta^{k\times 2}(C\cdot \tilde{\lambda}) \delta^{2\times (n-k)}(\lambda \cdot C^\perp)
\]
So the identity holds after doing the calculation where the positroid variety with matrix parameterization $C$ is put into this formula and the integrals evalulated; it is not an identity regarding positroid varieties themselves. \\
(2) The matrix $C$ is a parameterization for the positive real part of the positroid variety, not the positroid variety in the complex Grassmannian. We discussed the connection of the complex Grassmannian and the restriction to the real positive part in Section \ref{sec:positivereal}. There, we noted that the positroid cells are parametrizable and isomorphic to $\R^n$ for some $n$; since $\R^n$ is smooth, that means that all the singularities die off when we go to the nonnegative Grassmannian. This is another reason why we are unlikely to see any connection between the algebraic geometry of positroid varieties (in particular their singularity structure) and the physics in planar $N=4$ coming from these leading singularities.

\section{Spurious Poles}

Another direction one could consider would be whether the singularities of positroid varieties might have a correspondence to the notion of ``spurious poles," which are poles that appear in individual terms in the amplitude but which cancel out in the final answer. For the aforementioned reasons, it would seem that the answer is no; here, we provide an explicit counterexample showing that spurious poles do not correspond to singularities in the Grassmannian. 

Consider the $k=2, m=2, n=5$ amplituhedron. 

We briefly mentioned the amplituhedron in sections \ref{section:PositroidsandPhysics} and \ref{subsec:Tduality}. It is defined as the image of the map $\tilde{Z}$ of the kinematical information in momentum twistor form, from the nonnegative Grassmannian $Gr_{k,n}^{\geq 0}$ to the Grassmannian $Gr_{k,k+m}$. The point of the amplituhedron is that, if one tries to calculate the amplitude (in planar N=4 SYM) using the positive Grassmannian, in general one gets contributions from positroid cells of different dimensions, with complicated homological identities between them - the amplituhedron projects out the extraneous information so that the amplitude only gets contributions from the positroid cells you want. In particular, the BCFW decomposition of the amplitude becomes a triangulation $T_1 + \cdots + T_m$ of the image (amplituhedron) under the amplituhedron map $Z$ where each piece $T_i$ comes from a positroid cell in the preimage with all these positroid cells of the same dimension (this is what is called a ``tiling" of the amplituhedron). The poles in the amplitude correspond to the boundaries of the variety - lower dimensional positroid varieties (which correspond on the physics side to taking limits and factorization of an amplitudes into lower-point amplitudes); and the spurious poles are the boundaries in the image that cancel out in the sum of all terms forming hte amplitude. Of course, BCFW is just one choice of triangulation of the amplituhedron; the nice thing about the amplituhedron is that it (or technically, it's canonical form) gives the amplitude for any choice of triangulation or tiling by positroid cells.   

There is a diagrammatic way to obtain the positroid varieties that form a tiling of this amplituhedron, given by the following (\cite{LukWilPar})

 \begin{figure}[htbp] \centering
	\includegraphics[scale=0.7,clip=true]{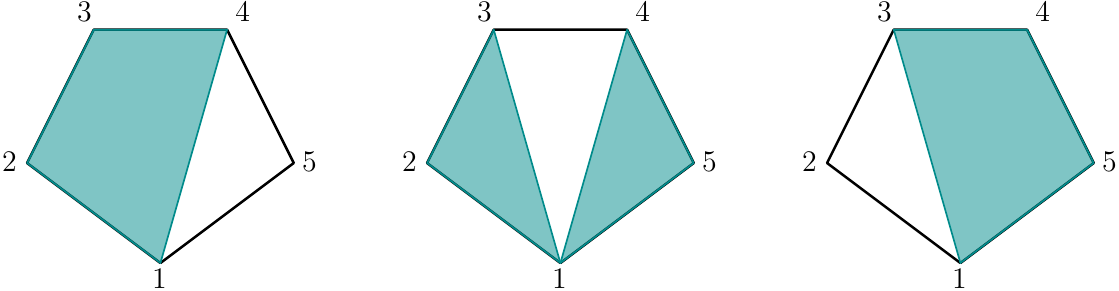}
  \end{figure}
  
The first one is given by the equations $z_{14}\geq 0$, the second by the equations: $z_{12}\geq0,z_{23}\geq0,z_{13}\leq0,z_{14}\leq0,z_{45}\geq0,z_{15}\geq0$, and the last by $z_{13}\geq 0$.

In the mapping to the amplituhedron, the first two share the $z_{14}$ boundary, and the second and third share the $z_{13}$ boundary. 
  
We can form the associated plabic diagrams as follows (we thank Lauren Williams for explaining how this works):

 \begin{figure}[htbp] \centering
	\includegraphics[scale=0.7,clip=true]{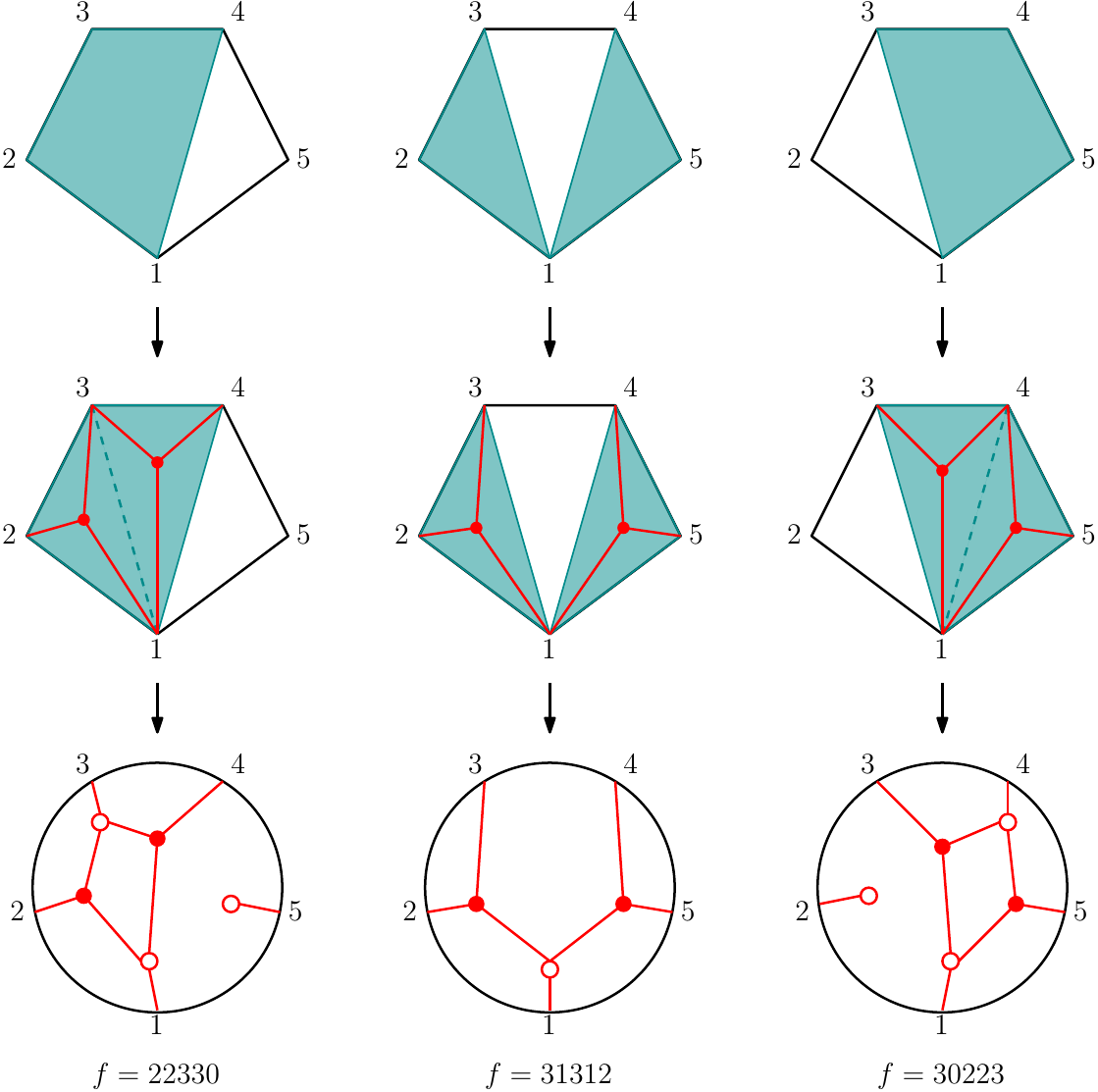}
  \end{figure}
  
That is, we triangulate the green regions (the choices we make here do not matter) and place a dark node at the center of each such triangle, that we connect to the vertices of the triangles. Two such lines that meet become a hollow node going to the corresponding boundary vertex, and we add a hollow lollipop node for each boundary vertex with no connection. 

We can simply trace out the left-right paths to obtain the 3 positroid varieties that these 3 tiles correspond to: 22330, 31312, and 30223 (siteswaps).

The first and third positroid varieties are isomorphic to subgrassmannian (Gr(2,4)), so are smooth. Only the middle one is singular. Yet, in the mapping to the amplituhedron, all three of them have spurious poles (which cancel upon summation). Thus, this example shows that the singularities of the positroids involved in the tiling of the amplituhedron do not correspond to spurious poles.

\chapter{Appendices}

\section{Appendix A: Biadjoint $\phi^3$ Theory and CEGM - Generalized Color Orderings and Reduced Words}

This appendix contains a writeup of a short exercise (part of an email sent to Freddy Cachazo in November 2023) based on the notion of generalized color orderings (GCO) \cite{CEZgco} and \cite{CEZlocplan}. I saw that the figures for an arrangement of 6 lines in the plane can be ``unwound" into the pattern of a wiring diagram. Furthermore, it not hard to see that these are reduced words for the word $w_0$ that is longest in Bruhat order. The first question I investigated was what reduced words fit into which type of GCO. Another question came from looking at the triangles, whose numbers are 6,6,7, and 10, and looking at why there are not triangle numbers of 8 or 9; this begs the question of what sorts of triangle flips take you to other GCO types. I was wondering if the language of reduced words could shed light on these questions. In particular, the mathematical literature on reduced words is vast, so making touch with Coexeter groups could potentially open up new connections. I created an accompanying Mathematica file that inputs a reduced word and outputs the triangles in the corresponding arrangement of lines.

\noindent-For the GCO diagrams, we can let the wires on one side be the labels on the lines on one half of the boundary circle, and the wires coming out the other side as the label on the other half of the circle. \\
-The lines are straight, so no 2 lines cross more than once - this is exactly the statement that the wiring diagram (and it associated word) is \textit{reduced}. \\
-The lines start in one order, and end in the opposite order, so the crossing must form a reduced word for the permutation that is longest in Bruhat order, usually denoted $w_0$ (in particular, the fact that they end up in opposite order means that each line crosses every other line (once, by the previous point)) \\
-Based on my understanding, these classes of GCO's are classes under all permutations of the labels; therefore, they can be said to be orbits under the full group of permutations $S_n$; indeed, action of $S_n$ just changes the labels on lines along the boundary circle, so it does nothing to the arrangement of lines, therefore it should not change the class of GCO. 

In terms of associating to a wiring diagram its word, the point is that the symmetric group $S_n$ is generated by transpositions $s_i, 1\leq i\leq n-1$. These transpositions can be put together to make up a word, which represents a permutation. For example, in one-line notation, the permutation 132 just switches the 2nd and 3rd elements, so its word is just a single letter $s_2$, but some permutations have multiple reduced words representing them. This is certainly the case for $w_0$, which corresponds to the permutation $n(n-1)(n-2)\cdots 321$. These transpositions satisfy the Coxeter relations: 

\noindent(1) $s_i^2=1$ \\
(2) $s_i s_j= s_j s_i$, $\abs{i-j}>1$ \\
(3) $s_i s_{i+1} s_i = s_{i+1} s_i s_{i+1}$.

Relation (1) is unnecessary in this context because it says 2 lines cross and uncross, which is not allowed by reducedness. Relation (2), which I'll call ``commuting moves" does nothing to the arrangement of lines. In fact, if we have, say, lines 1 and 2 crossing, and lines 3 and 4 crossing near each other, there is ambiguity about which pair of lines to write down first if we want to record the crossing of the lines in a linear order. However, for something like relation (3), which we can call ``braid moves," since $s_i$ is the transposition (i, i+1) and $s_{i+1}$ is the transposition (i+1,i+2), both involve the lines at position $i+1$, so the order matters. It turns out that the braid moves correspond precisely to the flips of triangles. 

We should of course note that what's been said so far is for the specific case of (3,n) GCO's, with $n<9$. Clearly, we're in (3,n) because we're dealing with lines (or projectivized planes) - this case is already very interesting to me because even the case of nonplanarity for NMHV in N=4 SYM is not understood, so in terms of connecting the biadjoint scalar theory with N=4, understanding the (3,n) case would be very helpful. In terms of specifying $n\leq 8$, this follows from Theorem 6.3.1 \cite{Bjor99}, stating that every arrangement of up to 8 lines is stretchable (we can homeomorph it to straight lines), or realizable (it turns out after I made the connection to reduced words, that I saw there is a discussion of wiring diagrams in their book, but I did not see a discussion about the triangle flips in terms of reduced words in there). So we can use some of the techniques for n=9, but there may have to be modifications. 

In summary, for the case of $(3,n), n\leq 8$, wiring diagrams up to moving around lines that do not cross is equivalent to considering reduced words up to commuting moves, which is also equivalent to considering arrangements of lines. In other words, we are simply capturing the information of when lines cross, keeping track of the order only when it's necessary (since it would change the connectivity of the diagram).

By considering the arrangement of lines as reduced words, we can immediately deduce a couple of facts. First, we cannot have any 2 triangles sharing a side. This would require a sequence like $s_i s_j s_i s_j$ which involves 2 wires crossing twice; since each triangle involves 3 different lines $a,b,c$, given a side of a triangle $a$, for the adjacent shape to be a triangle as well, that means that lines $b$ and $c$ must form the other 2 sides as well, so $b$ and $c$ must intersect again - this is the second $s_j$ in the sequence above. We also have a result of \cite{Elnit97} that the reduced words for $w_0$ up to commuting moves is connected and bipartite - this means that for any $n$ (although we have to consider pseudoline arrangements as noted above for large $n$), we can obtain all pseudoline arrangements by doing the triangle flips. 

My Mathematica program takes in a reduced word and produces the GCO and the triangles of the associated arrangement of lines.

In terms of producing the GCO, this comes simply from knowing how to map back and forth between the reduced word and the which lines intersect. Let's look at an example. The type 0 GCO had a representative given in Figure 2 of \cite{CEZlocplan}. Its reduced word is:
(2, 3, 5, 1, 2, 4, 3, 4, 5, 4, 2, 1, 2, 3, 2). I am reading from left to right. I have also switched the numbering on the lines so they go in order along the boundary starting from the $1$ at the top of the circle: $1,2,3,4,5,6$. We start with the permutation $(123456)$. The first $2$ switches the $2$ and the $3$, so the permutation becomes $(132456)$, and we record $(23)$. Next, we apply the $3$, which switches the (now) $2$ and $4$, so the permutation becomes $(134256)$. Our record is now $\{(23),(24)\}$. Continuing on, we will get a sequence of 15 pairs $(a,b)$ denoting when lines cross. Here in the type 0 GCO example, $\{(23),(24),...\}$ denotes that lines 2 and 3 cross, and then lines 2 and 4 cross etc. From here, we can read off the GCO very easily; for example, for the second permutation, we simply look at the order of the pairs containing $2$, so $\sigma^{(2)}=34...$. 

For creating a list of the triangles, the goal is to spot when we have something of the form $s_i s_j s_i$ in our reduced word. The Mathematica code does what is pictured in this ``heap" diagram, for the reduced word we just had: (2, 3, 5, 1, 2, 4, 3, 4, 5, 4, 2, 1, 2, 3, 2). We start at the bottom with the $2$, and continually stack on elements of the words. 

\begin{tikzpicture}
\draw node at (2,1*0.5) {2};
\draw node at (3,2*0.5) {3};
\draw node at (5,3*0.5) {5};
\draw node at (1,4*0.5) {1};
\draw node at (2,5*0.5) {2};
\draw node at (4,6*0.5) {4};
\draw node at (3,7*0.5) {3};
\draw node at (4,8*0.5) {4};
\draw node at (5,9*0.5) {5};
\draw node at (4,10*0.5) {4};
\draw node at (2,11*0.5) {2};
\draw node at (1,12*0.5) {1};
\draw node at (2,13*0.5) {2};
\draw node at (3,14*0.5) {3};
\draw node at (2,15*0.5) {2};
\draw[blue, thick] (3,7*0.5) -- (4,6*0.5);
\draw[blue, thick] (3,7*0.5) --(4,8*0.5);
\draw[red, thick] (3,7*0.5) -- (2,5*0.5);
\draw[red, thick] (3,7*0.5) --(2,11*0.5);
\draw node at (3,7*0.5) {3};
\draw[blue, dashed, thick] (2.8,4.1)--(5.2,4.1) --(5.2,2.9)--(2.8,2.9)--(2.8,4.1);
\draw[red, dashed, thick] (3.2,2.4)--(0.8,2.4) --(0.8,5.6)--(3.2,5.6)--(3.2,2.4);
\end{tikzpicture}

In other words, for each number $s_i$, we find the number $s_{i-1}$ that are the closest to it above and below (this is denoted by the two solid red lines). We then draw a box containing these two numbers $s_{i-1}$ of size such that the box extends exactly one unit to the left of $s_{i-1}$ on the left side, and one unit to the right on the other side (this is the depicted by the dotted box). If the box is empty except for the 3 numbers $s_{i-1}, s_i, s_{i-1}$, then they form a triangle. (In our example, both the red and blue form triangles, but if we for example took the $3$ at the very bottom, the $2$ that is above it has a $1$ sandwiched in between, so it does not form a triangle) We do the same on the other side for $s_{i+1}, s_i, s_{i+1}$ (in blue). In this example, both sides form triangles. We can draw the triangles:

\begin{tikzpicture}
\draw node at (2,1*0.5) {2};
\draw node at (3,2*0.5) {3};
\draw node at (5,3*0.5) {5};
\draw node at (1,4*0.5) {1};
\draw node at (2,5*0.5) {2};
\draw node at (4,6*0.5) {4};
\draw node at (3,7*0.5) {3};
\draw node at (4,8*0.5) {4};
\draw node at (5,9*0.5) {5};
\draw node at (4,10*0.5) {4};
\draw node at (2,11*0.5) {2};
\draw node at (1,12*0.5) {1};
\draw node at (2,13*0.5) {2};
\draw node at (3,14*0.5) {3};
\draw node at (2,15*0.5) {2};
\draw[blue] (3,7*0.5) -- (4,6*0.5)--(4,8*0.5)--(3,7*0.5);
\draw[blue] (3,7*0.5) -- (2,5*0.5)--(2,11*0.5)--(3,7*0.5); 
\draw[blue] (2,11*0.5)--(1,12*0.5)--(2,13*0.5)--(2,11*0.5);
\draw[blue] (2,13*0.5)--(3,14*0.5)--(2,15*0.5)--(2,13*0.5);
\draw[blue] (4,8*0.5) -- (5,9*0.5)--(4,10*0.5) --(4,8*0.5);
\end{tikzpicture}

However, there should be 6 triangles. The problem is that when you projectivize, there's a sort of Mobius thing that happens to the reduced words: it turns out that the triangles that occur by going through the boundary circle can also be described in terms of the occurrence of potential braid moves $s_i s_j s_i$, but you have to reverse the orientation when you come back around (I guess this is expected since you do the antipodal identification on the boundary of projective space).

Thus, if we look at (\rd{2}, 3, 5, 1, 2, 4, 3, 4, 5, 4, 2, 1, 2, 3, 2), this is equal to (3, 5, 1, 2, 4, 3, 4, 5, 4, 2, 1, 2, 3, 2,\rd{4}). So we actually have the triangle: (3, 5, 1, 2, 4, 3, 4, 5, \rd{4}, 2, 1, 2, \rd{3}, 2,\rd{4}), which we can visualize as: 

\begin{tikzpicture}
\draw node at (2,1*0.5) {2};
\draw node at (3,2*0.5) {3};
\draw node at (5,3*0.5) {5};
\draw node at (1,4*0.5) {1};
\draw node at (2,5*0.5) {2};
\draw node at (4,6*0.5) {4};
\draw node at (3,7*0.5) {3};
\draw node at (4,8*0.5) {4};
\draw node at (5,9*0.5) {5};
\draw node at (4,10*0.5) {4};
\draw node at (2,11*0.5) {2};
\draw node at (1,12*0.5) {1};
\draw node at (2,13*0.5) {2};
\draw node at (3,14*0.5) {3};
\draw node at (2,15*0.5) {2};
\draw[red] node at (4,16*0.5) {4};
\draw[blue] (3,7*0.5) -- (4,6*0.5)--(4,8*0.5)--(3,7*0.5);
\draw[blue] (3,7*0.5) -- (2,5*0.5)--(2,11*0.5)--(3,7*0.5); 
\draw[blue] (2,11*0.5)--(1,12*0.5)--(2,13*0.5)--(2,11*0.5);
\draw[blue] (2,13*0.5)--(3,14*0.5)--(2,15*0.5)--(2,13*0.5);
\draw[blue] (4,8*0.5) -- (5,9*0.5)--(4,10*0.5) --(4,8*0.5);
\draw[red]  (4,16*0.5) --(3,14*0.5)--(4,10*0.5)--(4,16*0.5);
\end{tikzpicture}

Thus, we can think about the wiring diagram as being along a (Mobius) strip that twists before closing back on itself. 

There's one more issue, which is that we can have a potential braid move involving lines 1 and 6. If we drew out the wiring diagram, this would be a triangle connecting the bottom of the wiring diagram and the top of the wiring diagram. This is due to a cyclic symmetry that is evident in the circle diagram, but that is lost when we go to a wiring diagram - we have artificially chosen the gap between lines 1 and 6 and stretched that out infinitely in forming a strip. Thus, if we look at our reduced word and consider all cases of $s_1$ and $s_{n-1}$ and we find that there are only a total of 3 occurrences, once with $s_1$ and twice with $s_{n-1}$ (or vice-versa), then we have one additional triangle. I have put this in by-hand in the code. Perhaps the best way to deal with this is to use the cyclic symmetry to shift everything by 1; then, due to the fact noted above that no 2 triangles are adjacent, the gap between the bottom and top of the wiring diagram will no longer be a triangle, and we can use the ordinary method. 

For example, one reduced word for the type III (3,6) GCO is (5,2,4,3,2,4,1,2,3,4,5,4,2,3,2). If we look at all $s_1$ and $s_5$, we just have $5,1,5$, so we have such a triangle. If we shift everything by $+1(mod5)$, we get (1,3,5,4,3,5,2,3,4,5,1,5,3,4,3), and this is the same arrangement of lines, but has no such triangle if we consider the wiring diagram with wire 1 entering at the bottom and wire 6 entering at the top. 

Therefore, perhaps the right way to think about the topology of the wiring diagram coming from a projectivized arrangement of planes in $\mathbb{R}^3$ is as a cylinder with a twist before matching the ends (this should just be a Klein bottle). 

One goal I had was to see if there was anything about the formulation in terms of reduced words that made the structure of the GCOs easier to see. Putting reduced words into the computer to get the answer was pretty easy, but seeing which class of GCO (in the (3,6) case) the reduced words falls into directly from the reduced word itself turns out to be difficult, so I'm not sure that much is gained. 

Even for the case of type III, which is probably the simplest case, this is not easy to see. 

For type III, every letter in the word participates in exactly 2 different triangles. Type III has this beautiful pattern with 6 pentagons and 10 triangles, where every face of the pentagon is part of a triangle, and all triangles touch 3 other triangles at its 3 vertices. In particular, every crossing of 2 lines is a vertex for 2 different triangles, so every letter is a part of 2 possible braid moves. We can see this for an example of type III, given by the reduced word:

354321354325345. Taking the first 4 letters (3543) as an example, we have the obvious potential braid moves: \rd{3}5\rd{43}21354325345, 
3\rd{54}3213\rd{5}4325345, 354\rd{32}1\rd{3}54325345.

Thus, we have the 2 triangles that the 43 (3rd and 4th letters) participate in, but we only have 1 triangle for the first two letter. And this not particularly easy to see from the word; one of them comes from wrapping around: \rd{3}54321354325\rd{34}5. 

The second triangle involves the 5 is most difficult to see: 5 becomes a 1 on the other side:
3\rd{5}4321354325345$\to$43213543253453\rd{1}, and then forms a 121 triangle: 432\rd{1}3543\rd{2}53453\rd{1}, so this one is especially hard to see. 

In terms of the question of which GCO types go to which other ones under the triangle/braid moves, it's clear from the picture for type III that all possible moves take it to type II with 7 triangles. In the ones I tested on the computer, type 0 only goes to type 0 or type I (the number of triangles stay the same at 6 triangles). However, beyond that I would need to do more investigation. (Or there might also be some information in some of the portions of your paper that I haven't had a chance to read as carefully yet).

\section{Appendix B: Nonplanar On-Shell Diagrams}

\subsection{Removability and Dimension}\label{subsec:removable}

It is an important question whether an edge is removable or not. Most of the time we put labels on the edges via the BCFW bridge method (we create the representative on-shell diagram using the BCFW method and accordingly label the edges). However, we cannot access all boundaries in any one chart - as noted in \cite{ABCGPT}, any boundary can be accessed within a chart for (the BCFW bridge decomposition) for one of the $n$ rotations. In any one chart, setting some of the variables/bridges equal to zero drops the dimension of the variety by more than 1 - these are non-removable edges.

It is important to check whether an edge is removable or not. The way we do this is we check the dimension of the tangent space after setting this variable equal to zero and checking that it decreased by 1. Note that this is a concept that applies to both planar and nonplanar diagrams.

We do an example here. Let's take the top-cell for the (planar) positroid variety in Gr(3,6). A representative on-shell diagram (with a perfect orientation) is provided by the following from \cite{BFGW}

\begin{figure}[htbp]
\includegraphics[scale=0.1,clip=true]{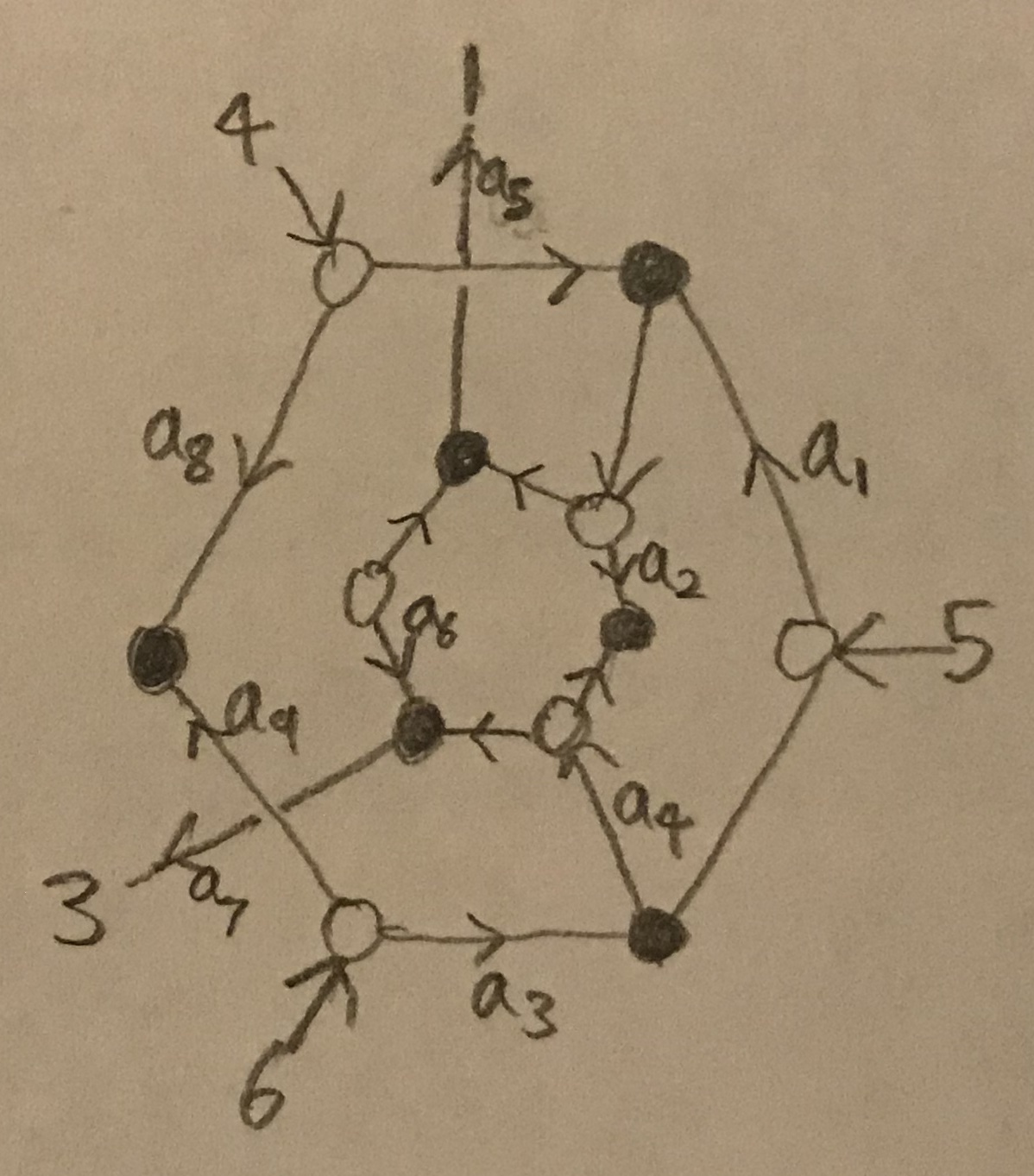}
\end{figure}

The matrix is: $C-matrix = \mqty( -a_5 (1+a_8)& -a_2 & a_6 a_7 a_8 & 1& 0 & 0 \\
a_1 a_5 & a_1 a_2+a_4 & a_4 a_7 & 0 & 1 & 0 \\
-a_5 a_9 & a_3 a_4 & a_7 (a_3 a_4+a_6 a_9) & 0 & 0 & 1 \\
)$

We can consider the space as a variety determined by the entires in the first 3-by-3 grid. This variety (and we can check this using the technique we are currently in the process of describing) is 9-dimensional, and has 9 variables. We can write the tangent space by taking this $3\times 3$ piece and stretching it into a 9-by-1 vector, and computing successive partial derivatives for all 9 variables and putting this into a 9-by-9 matrix, and then computing its rank to get the dimension of the tangent space. In other words, we can do:

\[M_{tang} = \frac{\del}{\del x_i} \mqty( -a_5 (1+a_8)\\ -a_2 \\ a_6 a_7 a_8 \\ 
a_1 a_5 \\ a_1 a_2+a_4 \\ a_4 a_7 \\ 
-a_5 a_9 \\ a_3 a_4 \\ a_7 (a_3 a_4+a_6 a_9) 
)\]

The first $\frac{\del}{\del x_1}$ is: $\mqty( 0\\ 0 \\ 0 \\ 
a_5 \\ a_2 \\ 0 \\ 
0 \\ 0 \\ 0
)$

The full matrix is: 

\begin{figure}[htbp]
\includegraphics[scale=0.5,clip=true]{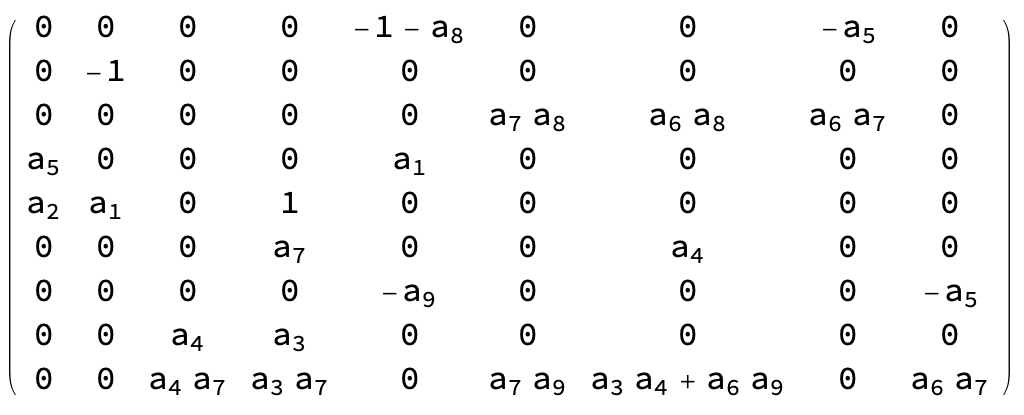}
\end{figure}

Now, the definition of removable is that removing that edge only lowers the dimension of the variety by 1. So the test of this is now automatic: we take the matrix $Cmatrix$ above and set a variable $a_i=0$ to see if it is removable. Take this $Cmatrix|_{a_i=0}$ matrix, and now do the process to get $M_{tang}$ from the modified matrix with $a_i=0$ and finally calculate the rank to see if it just dropped by 1.

\subsection{MHV nonplanar counting}

From (4.35) and (4.36) of \cite{ABCGPT}, we have the equations:

$k = 2n_b+n_w-n_I$ and $n=3n_V-2n_I$.

where $n_V = n_b+n_w$, the total number of vertices is the number of black plus the number of white vertices; $n_I$ is the number of internal lines.

Note that $n=3n_V-2n_I$ requires that all edges be trivalent, e.g. it does not hold for one of the diagrams in \cite{ABCPT}, which has 4-valent vertices: 

\begin{figure}[htbp]
\includegraphics[scale=0.5,clip=true]{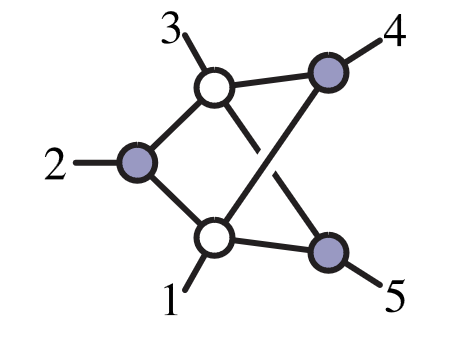}
\end{figure}

(we must use the equivalence moves to expand the white 4-vertices into a pair of trivalent vertices, and then the formula will hold)

They furthermore note in \cite{ABCPT} that the number of delta functions imposed on the external kinematics beyond momentum conservation is $n_\delta = 4 n_V - 3n_I - 4$. For an on-shell function/leading singularity, we require $n_\delta =0$.

Substituting $n=3n_V-2n_I$ into $n_\delta = 4 n_V - 3n_I - 4$ and setting it equal to zero, we get: $n + n_V - n_I - 4 =0$, or $n + n_b + n_w - n_I - 4 =0$.

If we subtract $2n_b+n_w-n_I -k =0$: $n + n_b + n_w - n_I - 4 -(2n_b+n_w-n_I -k) =0$, we get

$n -n_b -4 +k =0$.

Thus, if $k=2$ (MHV), we get $n_b=n-2$, and if $k=3$ (NMHV), we get $n_b=n-1$. 

However, the interpretation in terms of the form is not as clear in the NMHV case. In the MHV case, we can place every triplet for every black vertex into a matrix as a set of 2-brackets, and then there is an interpretation in terms of eliminating any 2 columns, computing the determinant of the remaining matrix and dividing by the bracket of the 2 eliminated columns to obtain a Jacobian factor of the differential form. For NMHV, the brackets should now be 3-dimensional. The process for writing the differential form does not carry over immediately.

\subsection{Deletion and Contraction for Nonplanar Diagrams}

We are considering nonplanar bicolored diagrams, such as described in \cite{{BFGW}}. 


Note that these are the top-dimensional varieties obtained using the nonplanar diagrams. However, the leading singularities/on-shell functions are obtained by using the integral formula, which has $2n-4$ delta functions. Here, $n=6$ particles, and $k=3$ since the NMHV next to maximal helicity violating case is being considered. Thus, to obtain a function, we need a $2n-4=12-4=8$ dimensional variety. It is conjectured that these (let's call them "on shell") varieties are the codimension-1 boundaries of the top-dimensional varieties. 


Of course, to be precise, the diagram represents (loosely used: more precisely there is a way to assign variables to edges of the diagram and produce a ``C-matrix") a variety in a Grassmannian; the physics comes from applying the integral formula using this variety (putting the C-matrix into the delta functions and evaluating, etc.)

However, we have the additional fact from \cite{ABCPT} that the MHV case (Gr(2,n)) has a nice description of the differential form of nonplanar diagrams as sums of Parke-Taylor factors. In other words, MHV amplitudes can be decomposed in terms of ``permuted positroid varieties." They give an easy criterion for seeing if a nonplanar diagram is MHV, which is that there are $n-2$ black vertices (n legs). 

Now, we can think about applying deletion or contraction to the Grassmannian. For example, if we delete the 3rd leg, we are considering the subset of the Grassmannian that does not use the 3rd coordinate (in the orthogonal complement to the 3rd coordinate), and if we contract the 3rd leg, we are considering the subset of planes that contains the 3rd coordinate. In \cite{Flu}, this was described as a 2-step process: first, intersecting with the positroid variety $\Pi_{del,i}$ or  $\Pi_{contr,i}$, and then projecting to the isomorphic Grassmannin in $n-1$. 

For deletion, the first step would still be a variety in Gr(3,6), while the second step would take it to Gr(3,5). For contraction, the first step would also be a variety in Gr(3,6), but the second step would drop it to Gr(2,5), which would be the MHV case. The deletion case of Gr(3,5) would also be isomorphic to a variety in Gr(2,5) by Grassmann duality. It would thus be interesting to see what deletion and contraction do to the top-dimensional Gr(3,6) varieties listed above. It seems possible for these operations (particularly contraction) to take a diagram representing an on-shell variety in $Gr(3,6)$ to an MHV or anti-MHV diagram.

Thus, the questions are: what does deletion/contraction do to the nonplanar diagram? 
And, what does del/contr do to the differential form? If it ends up as a sum of Parke-Taylor forms, which sum does it land on? Do multiple contractions or deletions (of different columns) of the same top-dim variety land on the same form? What about contractions of different top-dimensional varieties landing on the same MHV form? 

Below, we took the simplest nonplanar diagram and chose a perfect orientation (we are not aware of any proof that the choice of perfect orientation doesn't matter for these nonplanar diagrams (the case for planar ones was proven in \cite{Post}) 

\begin{figure}[htbp]
\includegraphics[scale=0.06,clip=true]{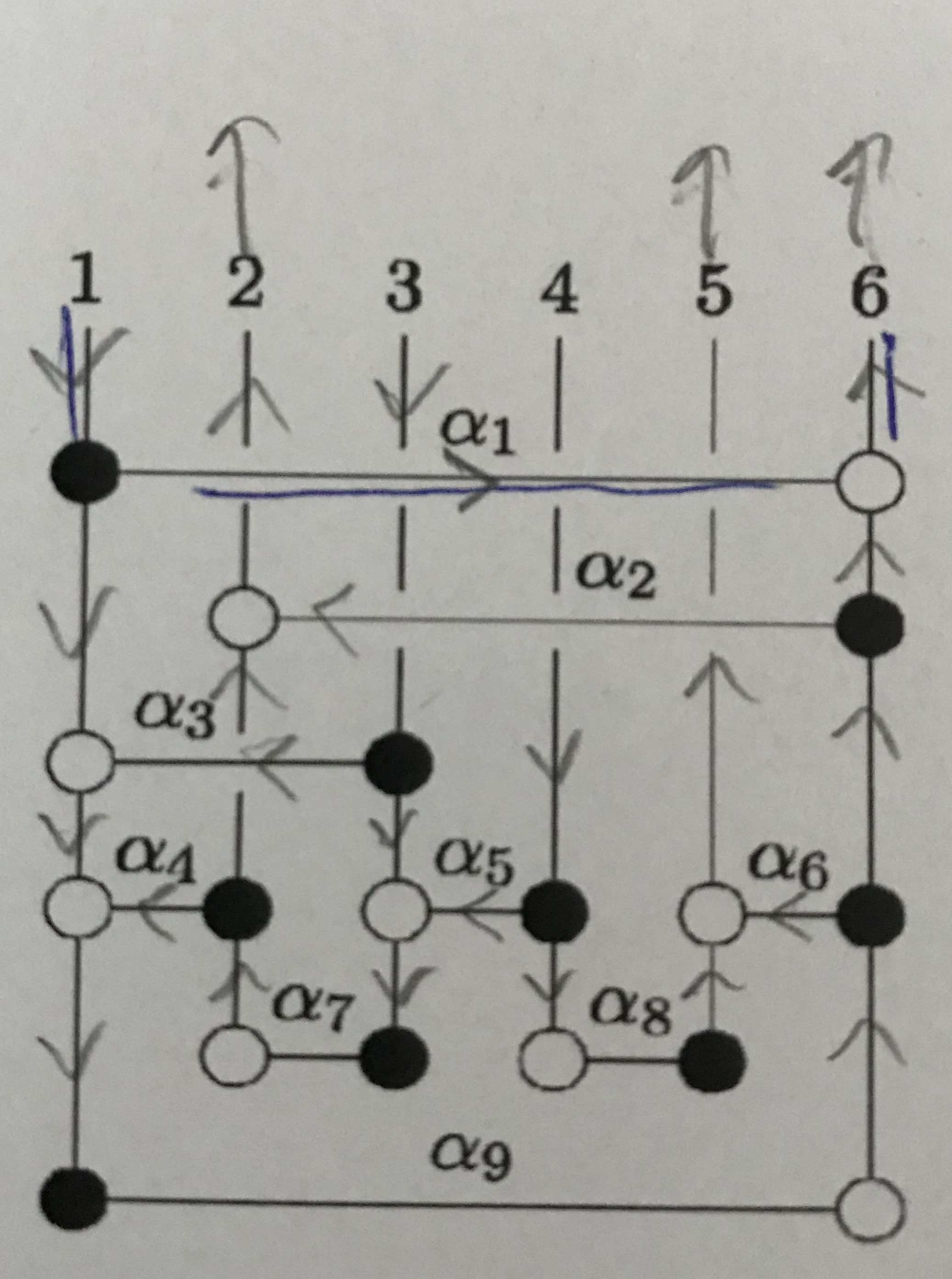}
\end{figure}

The important thing was to choose orientation such that there is no internal loop (that makes the matrix have no denominators, thus be easier to deal with). Note that the conventions for arrows coming in and out are opposite those of \cite{ABCGPT}. Here, we initially used those conventions with each white vertex having two edges pointed out, but this resulted in a C-matrix in Gr(2,6), so it needed to be changed such that each white vertex has 2 edges pointing in.

$\mqty(
a_2 a_9 & 1 & a_2 a_9 (a_3+a_4 a_7) & a_2 a_9 (a_4 a_5 a_7) & 0 & 0 \\
a_6 a_9 & 0 & a_7 + a_3 a_6 a_9 + a_4 a_6 a_7 a_9 & a_8 + a_5 a_7 + a_4 a_5 a_6 a_7 a_9 & 1 & 0 \\
a_1 + a_9 & 0 & a_3 a_9 + a_4 a_7 a_9 & a_4 a_5 a_7 a_9 & 0 & 1 
)$

There are 9 variables, corresponding to the fact that we have a 9-dimensional variety. Now, suppose that we set $a_2=0$, which would correspond to removing that bridge/edge in the graph. Then we get:

$\mqty(
0 & 1 & 0 & 0 & 0 & 0 \\
a_6 a_9 & 0 & a_7 + a_3 a_6 a_9 + a_4 a_6 a_7 a_9 & a_8 + a_5 a_7 + a_4 a_5 a_6 a_7 a_9 & 1 & 0 \\
a_1 + a_9 & 0 & a_3 a_9 + a_4 a_7 a_9 & a_4 a_5 a_7 a_9 & 0 & 1 
)$

This is in the contracted form. We can project to Gr(2,5):

$\mqty(
a_6 a_9 & a_7 + a_3 a_6 a_9 + a_4 a_6 a_7 a_9 & a_8 + a_5 a_7 + a_4 a_5 a_6 a_7 a_9 & 1 & 0 \\
a_1 + a_9 & a_3 a_9 + a_4 a_7 a_9 & a_4 a_5 a_7 a_9 & 0 & 1 
)$

This gets rid of the 2nd leg. The elimination of $a_2$ already got rid of one black vertex. So that left us with 5 black vertices. So the removal of the 2nd leg leaves us with 4 black vertices for 5 particles. So we have 1 too many black vertices. We should be able to do a small change to the diagram to get an equivalent diagram with 4 vertices. One thing we should do is that eliminating leg 2 sets $a_4$ equal to $a_7$. It seems like the contraction should be a 6-dimensional variety, so it seems like we should eliminate one more variable/bridge. 

It is also interesting to try to obtain the form using the trick on page 7 of \cite{ABCPT} We first have to obtain $C^\perp$, which is:

$\mqty(
1 & 0 & 0 &  \\
0 & 1 & 0 &  \\ 
0 & 0 & 1 \\
-a_6 a_9 & -(a_7 + a_3 a_6 a_9 + a_4 a_6 a_7 a_9) & -(a_8 + a_5 a_7 + a_4 a_5 a_6 a_7 a_9) \\
-(a_1 + a_9) & -(a_3 a_9 + a_4 a_7 a_9) & -(a_4 a_5 a_7 a_9) 
)$

We should take the transpose of this before eliminating columns as on page 7 of that paper. If we can get the correct on-shell diagram with 3 black vertices, we can write down the form for that, and match up with the form in terms of the $a$ variables that we get using this method.  

(In the planar/positroid case, the C-matrix that we get by doing the boundary measurement map should be the same as the one from using the BCFW decomposition, up to SL(k) transformation. However, here we are sort of mixing the two: we are doing the boundary measurement map, but we are only using variables associated with the bridges.) 


If we set a2=0, which is removing that bridge, we still have 5 black vertices, rather than 6-2=4, which is what we should have. (But this will be taken care of by the leg removal later, hopefully)

This contracted thing should already be MHV, in particular, we can take the 2 by 5 matrix, and  now it has a1, a3,a4,a5,a6,a7,a8,a9 all still in there, so it should be 8-dimensional, although this doesn't make sense since Gr(2,5) is 6-dimensional, so some of these degrees of freedom are redundant. If we removed 2 edges, we lose 2 more black vertices, so we're down to 3 black edges, and we have to remove that leg that we're contracting. Removing the 2nd leg here forces a4 and a7 to be equal (could this be related to the ``amalgamation" procedure in the \cite{ABCGPT})

Ultimately, the purpose of considering deletion and contraction is because of the hope that we can get an "inverse soft" construction that allows us to build up amplitudes from smaller pieces. 
\\

In summary, the idea is that if we take an on-shell diagram for Gr(3,6) and do a deletion or contraction, this should take us to Gr(3,5) or Gr(2,5), which are either MHV or anti-MHV (equivalent to MHV under Grassmannian duality). 

In terms of getting the $n-2$ black vertices (as noted in \cite{ABCPT}), this does not work, because for the Gr(3,6) nonplanar top-cells, the equivalence moves are still valid; however, since an equivalence move on a black vertex (switching between a 4-valent black vertex and two 3-valent black vertices) adds one to both the edges and vertices (similarly for the white vertex), the equations above: $k = 2n_b+n_w-n_I$ cannot stay constant under a black equivalence move. This could imply that we cannot have any two black vertices that are connected to each other. 
For white vertices, this is ok, if we take $n_2-n_I$ and change both by the same amount; this is precisely what happens in one of the equivalent moves, as long as none of the edges connecting to the white vertex are external edges. (i.e. the white vertex must have all edges connecting to other vertices (black or white) and none connecting to external edges)

A main issue with examples from top cells is that where deletion/contraction seems to not send them to one of the $n_b=n-2$ type diagrams.

So this seems to suggest one of two things: (1) there's an issue with the delta functions/leading singularities, so we actually needed to use one of the Gr(3,6) leading singularities, which have dimension 8; or (2) there is an issue with the process itself.

Evidence in favor of (2) would come from the fact that deletion and contraction are difficult to do on plabic diagrams (there's no clean rule) - they are nice on pipe dreams, as explained in this thesis, but not for plabic diagrams.

There's a tension that deleting a nonremovable edge lowers the number of variables by less (but lowers the dimension of the variety by more), so perhaps it is necessary to only remove removable edges: this means that we need to move 3 (set 3 variables equal to zero) for Gr(3,n). This may be the best way to do deletion and contraction for these diagrams. However, the pipe dream story also gives you good information - it tells you how many extra crosses you have to put in (so this may be of dimension less than $k$ (or less than $n-k$ for contraction)). 

The pipe dream is patch dependent. The question is how many crosses we put in. The cleanest story is when there are $k$ for deletion (no crosses), since then we can do it by successively doing removable edges. If not, then we need to do fewer than $k$, which means that we will be able to set the column equal to zero with fewer than $k$ variables, which implies that we will be removing at least one non-removable edge. We should try this for the positroid case: test if removing all removable edges does indeed get you the right positroid variety when we do left-right paths after deleting the edges/bridges at the diagram level. 


We should answer the question of how many delta functions does a deletion or contraction change.

The main ideas/questions in the previous discussion were summarized in my 2024 Research Statement (available at \url{josephflueg.github.io}, but pasted here for convenience):

\subsection{Research Statement 2024}

Here we discuss how to apply deletion and contraction to gain a better understanding of \textit{nonplanar on-shell diagrams}, which we call \textit{nonplabic diagrams}. These are like plabic diagrams, in that they are bicolored and can be drawn with entirely trivalent vertices, but they cannot be drawn without edges crossing each other if all ``external edges" are on the boundary of a disc. 

One of the cases that is very interesting for the physics is the 6-point NMHV case, corresponding to $Gr(3,6,\R)$. By calculating a very large number of examples, \cite{BFGW} determined that the nonplabic diagrams produce 24 inequivalent top-dimensional ``varieties" in $Gr(3,6,\R)$. However, these are open sets. We conjecture that they are all of $Gr(3,6,\R)$ complemented by a divisor, which is the vanishing locus of a section of the anticanonical bundle. 

These nonplabic diagrams are difficult to deal with (they have a number of moves that preserve the equivalence class, including so-called ``sphere moves" (\cite{CEGM}) and possibly others). It may be best to start working with the varieties themselves, for which we have a $C$-matrix parameterization obtained by using the boundary measurement matrix after picking a perfect orientation \cite{Post}. Deletion and contraction of plabic diagrams is not a nice story (we know how it works on pipe dreams, but the map from pipe dreams to plabic diagrams does not interact nicely with moves). 

We know that deletion takes $Gr(3,6)$ to $Gr(3,5)$. Via Grassmannian duality (which can be obtained on the plabic or nonplabic diagrams by switching the white and black vertices), $Gr(3,5)$ is isomorphic to $Gr(2,5)$, so we can map to a variety within $Gr(2,5)$. Contraction takes $Gr(3,6)$ to $Gr(2,5)$. Thus, in both cases, we get a map to (something isomorphic to) $Gr(2,5)$. Now, the $Gr(2,n)$ case is the MHV (maximally helicity violating case). We understand the nonplanar MHV amplitudes \cite{ABCPT}: the differential forms consists of sums of positroid varieties with rotated coordinates. Thus, the question is whether the deletions and contractions give us MHV amplitudes. If so, what is the map? And can it shed light on the geometry of nonplanar NMHV amplitudes? One of the difficulties will be determining when a $C$-matrix corresponds to a nonvanishing rational function for generic momenta, since this is the case for which the nice results of \cite{ABCPT} hold.   

We do a quick example here. Let's take the top-cell for the (planar) positroid variety in $Gr(3,6)$ stated in the Removability and Dimension subsection \ref{subsec:removable} above (\cite{BFGW}). We restate the C-matrix here:

The matrix is: $C-matrix = \mqty( -a_5 (1+a_8)& -a_2 & a_6 a_7 a_8 & 1& 0 & 0 \\
a_1 a_5 & a_1 a_2+a_4 & a_4 a_7 & 0 & 1 & 0 \\
-a_5 a_9 & a_3 a_4 & a_7 (a_3 a_4+a_6 a_9) & 0 & 0 & 1 \\
)$

The contraction of column 5 gives us the vanishing set of the ideal $I=\langle a_1 a_5 , a_1 a_2+a_4 , a_4 a_7 \rangle$, which has top-dimensional components $a_1=a_4=0$. We conjecture that a point $\lambda$ exists on the open set parameterized by a C-matrix if and only if there is a perfect orientation such that all external edges corresponding to $\lambda$ point outwards. 

Once we go to nonplabic diagrams, we start dealing with minors that are negative (e.g. see\cite{PTZ}), so it makes sense to use oriented matroids to describe some of the combinatorial structure of these varieties. In this case, using an analogue of deletion and contraction for oriented matroids could uncover interesting results. The end goal would be to try to obtain some sort of combinatorial story for the varieties in $Gr(3,6)$ given by these nonplabic diagrams. If we are able to delete/contract these varieties, we could ask the additional question of whether we still have closure for smoothness under going down. One problem here is that the non-irreducibility of the contraction, as we saw in the example above, suggests that these varieties are not torus-invariant.


\bibliography{FluegemannBibliography}

@misc{KLS,
  author       = {A. Knutson and T. Lam and D. Speyer},
  title        = {Positroid Varieties: Juggling and Geometry},
  eprint       = {1111.3660},
  archivePrefix= {arXiv},
  primaryClass = {math.AG}
}

@misc{KLS2,
  author       = {A. Knutson and T. Lam and D. Speyer},
  title        = {Projections of Richardson Varieties},
  eprint       = {1008.3939},
  archivePrefix= {arXiv},
  primaryClass = {math.AG}
}

@phdthesis{Sni,
  author = {M. Snider},
  title  = {Affine Patches on Positroid Varieties and Affine Pipe Dreams},
  note   = {arXiv:1011.3705 [math.CO]},
  year   = {2011},
  school={Cornell University}
}

@misc{KM03,
  author       = {A. Knutson and E. Miller},
  title        = {Subword Complexes in Coxeter Groups},
  eprint       = {math/0309259},
  archivePrefix= {arXiv},
  primaryClass = {math.CO}
}

@article{AJSBilley,
  author    = {Sara C. Billey},
  title     = {Kostant polynomials and the cohomology ring for $G/B$},
  journal   = {Duke Mathematical Journal},
  volume    = {96},
  number    = {1},
  pages     = {205--224},
  year      = {1999},
  publisher = {Duke University Press},
  doi       = {10.1215/S0012-7094-99-09606-0}
}

@article{BB93,
  author  = {N. Bergeron and S. Billey},
  title   = {RC-graphs and Schubert Polynomials},
  journal = {Experimental Mathematics},
  volume  = {2},
  number  = {4},
  pages   = {257--269},
  year    = {1993}
}

@misc{BilBra03,
  author       = {S. Billey and Tom Braden},
  title        = {Lower Bounds for Kazhdan-Lusztig Polynomials from Patterns},
  eprint       = {math/0202252},
  archivePrefix= {arXiv},
  primaryClass = {math.RT}
}

@misc{SOh,
  author       = {S. Oh},
  title        = {Contraction and Restriction of Positroids in Terms of Decorated Permutations},
  eprint       = {0804.0882},
  archivePrefix= {arXiv},
  primaryClass = {math.CO}
}

@misc{Oh08,
  author       = {Suho Oh},
  title        = {Positroids and Schubert Matroids},
  eprint       = {0803.1018},
  archivePrefix= {arXiv},
  primaryClass = {math.CO}
}

@misc{Post,
  author       = {A. Postnikov},
  title        = {Total Positivity, Grassmannians, and Networks},
  eprint       = {math/0609764},
  archivePrefix= {arXiv},
  primaryClass = {math.CO}
}

@misc{Rossmann,
  author = {W. Rossmann},
  title  = {Equivariant Multiplicities on Complex Varieties},
  note   = {\url{https://mysite.science.uottawa.ca/rossmann/Papers_files/Equivariant_multiplicities.pdf}}
}

@book{Hartshorne,
  author    = {R. Hartshorne},
  title     = {Algebraic Geometry},
  publisher = {Springer},
  year={1977}
}

@misc{KnTao01,
  author       = {A. Knutson and T. Tao},
  title        = {Puzzles and (Equivariant) Cohomology of Grassmannians},
  eprint       = {math/0112150},
  archivePrefix= {arXiv},
  primaryClass = {math.CO}
}

@book{Kumar,
  author    = {S. Kumar},
  title     = {Kac-Moody Groups, their Flag Varieties and Representation Theory},
  series    = {Progress in Mathematics},
  volume    = {204},
  publisher = {Springer Science},
  year      = {2002}
}

@article{Kleiman,
  author  = {S. Kleiman},
  title   = {The Transversality of a General Translate},
  journal = {Compositio Mathematica},
  volume  = {28},
  pages   = {287--297},
  year    = {1974},
  note    = {MR 0360616}
}

@misc{AH,
  author = {N. Arkani-Hamed},
  title  = {Quantum Mechanics and Spacetime Lectures},
  note   = {Harvard, Fall 2019}
}

@book{EH,
  author    = {H. Elvang and Y. Huang},
  title     = {Scattering Amplitudes in Gauge Theory and Gravity},
  publisher = {Cambridge University Press},
  edition   = {1st},
  year      = {2015}
}

@book{Peskin,
  author    = {M. Peskin and D. Schroeder},
  title     = {An Introduction to Quantum Field Theory},
  publisher = {CRC Press},
  edition   = {1st},
  year      = {1995}
}

@book{Schw,
  author    = {M. Schwartz},
  title     = {Quantum Field Theory and the Standard Model},
  publisher = {Cambridge University Press},
  edition   = {1st},
  year      = {2013}
}

@misc{AHTr,
  author       = {N. Arkani-Hamed and J. Trnka},
  title        = {The Amplituhedron},
  eprint       = {1312.2007},
  archivePrefix= {arXiv},
  primaryClass = {hep-th}
}

@misc{ARW,
  author       = {F. Ardila and F. Rinc{\'o}n and L. Williams},
  title        = {Positroids and Non-Crossing Partitions},
  eprint       = {1308.2698},
  archivePrefix= {arXiv},
  primaryClass = {math.CO}
}

@misc{Flu,
  author       = {J. Fluegemann},
  title        = {Smooth Points on Positroid Varieties},
  eprint       = {2407.21116},
  archivePrefix= {arXiv},
  primaryClass = {math.CO}
}

@misc{ABCGPT,
  author       = {N. Arkani-Hamed and J. Bourjaily and F. Cachazo and A. Goncharov and A. Postnikov and J. Trnka},
  title        = {Scattering Amplitudes and the Positive Grassmannian},
  eprint       = {1212.5605},
  archivePrefix= {arXiv},
  primaryClass = {math.CO}
}

@misc{HeZh,
  author       = {S. He and C. Zhang},
  title        = {Notes on Scattering Amplitudes as Differential Forms},
  eprint       = {1807.11051},
  archivePrefix= {arXiv},
  primaryClass = {hep-th}
}

@misc{LukWilPar,
  author       = {T. Lukowski and M. Parisi and L. K. Williams},
  title        = {The Positive Tropical Grassmannian, the Hypersimplex, and the m=2 Amplituhedron},
  eprint       = {2002.06164},
  archivePrefix= {arXiv},
  primaryClass = {math.CO}
}

@misc{HuYeats,
  author       = {S. Hu and K. Yeats},
  title        = {Reconstructing a Bijection on the Level of Le Diagrams},
  eprint       = {2308.06307},
  archivePrefix= {arXiv},
  primaryClass = {hep-th}
}

@article{DFLP19,
  author  = {D. Damgaard and L. Ferro and T. Lukowski and M. Parisi},
  title   = {{The Momentum Amplituhedron}},
  journal = {JHEP},
  volume  = {08},
  pages   = {042},
  year    = {2019},
  eprint  = {1905.04216},
  archivePrefix = {arXiv}
}

@article{AT13,
  author  = {N. Arkani-Hamed and J. Trnka},
  title   = {{The Amplituhedron}},
  journal = {JHEP},
  volume  = {10},
  pages   = {030},
  year    = {2014},
  doi     = {10.1007/JHEP10(2014)030},
  eprint  = {1312.2007},
  archivePrefix = {arXiv},
  primaryClass   = {hep-th}
}

@article{AtiBott84,
  author  = {M. F. Atiyah and R. Bott},
  title   = {The Moment Map and Equivariant Cohomology},
  journal = {Topology},
  volume  = {23},
  number  = {1},
  pages   = {1--28},
  year    = {1984}
}

@misc{KnutQuiver,
  author = {A. Knutson},
  title  = {Schubert Calculus and Quiver Varieties},
  note   = {Unpublished notes, 2019. \url{https://pi.math.cornell.edu/~allenk/courses/19fall/7310/quiv.pdf}}
}

@article{CEZgco,
  author  = {F. Cachazo and N. Early and Y. Zhang},
  title   = {Generalized Color Orderings: CEGM Integrands and Decoupling Identities},
  journal = {Nucl. Phys. B},
  volume  = {1002},
  pages   = {116552},
  year    = {2024},
  doi     = {10.1016/j.nuclphysb.2024.116552},
  eprint  = {2304.07351},
  archivePrefix = {arXiv},
  primaryClass   = {hep-th}
}

@article{CEZlocplan,
  author  = {F. Cachazo and N. Early and Y. Zhang},
  title   = {Color-Dressed Generalized Biadjoint Scalar Amplitudes: Local Planarity},
  journal = {SIGMA},
  volume  = {20},
  pages   = {016},
  year    = {2024},
  doi     = {10.3842/SIGMA.2024.016},
  eprint  = {2212.11243},
  archivePrefix = {arXiv},
  primaryClass   = {hep-th}
}

@book{Bjor99,
  author    = {A. Bj{\"o}rner and M. Las Vergnas and B. Sturmfels and N. White and G. M. Ziegler},
  title     = {Oriented Matroids},
  series    = {Encyclopedia of Mathematics and its Applications},
  volume    = {46},
  publisher = {Cambridge University Press},
  address   = {Cambridge},
  edition   = {2nd},
  year      = {1999}
}

@article{Elnit97,
  author  = {S. Elnitsky},
  title   = {Rhombic Tilings of Polygons and Classes of Reduced Words in Coxeter Groups},
  journal = {Journal of Combinatorial Theory, Series A},
  volume  = {77},
  number  = {2},
  pages   = {193--221},
  year    = {1997}
}

@incollection{Mnev,
  author    = {N. E. Mn{\"e}v},
  title     = {The Universality Theorems on the Classification Problem of Configuration Varieties and Convex Polytopes Varieties},
  booktitle = {Topology and Geometry --- Rokhlin Seminar},
  series    = {Lecture Notes in Mathematics},
  volume    = {1346},
  pages     = {527--544},
  publisher = {Springer},
  address   = {Berlin, Heidelberg},
  year      = {1988}
}

@article{BFGW,
  author  = {J. L. Bourjaily and S. Franco and D. Galloni and C. Wen},
  title   = {Stratifying On-Shell Cluster Varieties: the Geometry of Non-Planar On-Shell Diagrams},
  journal = {JHEP},
  volume  = {10},
  pages   = {003},
  year    = {2016},
  doi     = {10.1007/JHEP10(2016)003},
  eprint  = {1607.01781},
  archivePrefix = {arXiv},
  primaryClass   = {hep-th}
}

@misc{BourjMath,
  author = {J. L. Bourjaily},
  title  = {Positroids, Plabic Graphs, and Scattering Amplitudes in Mathematica},
  eprint = {1212.6974},
  archivePrefix = {arXiv},
  primaryClass   = {hep-th}
}

@article{ABCPT,
  author  = {N. Arkani-Hamed and J. L. Bourjaily and F. Cachazo and A. Postnikov and J. Trnka},
  title   = {On-Shell Structures of MHV Amplitudes Beyond the Planar Limit},
  journal = {JHEP},
  volume  = {06},
  pages   = {179},
  year    = {2015},
  doi     = {10.1007/JHEP06(2015)179},
  eprint  = {1412.8475},
  archivePrefix = {arXiv},
  primaryClass   = {hep-th}
}

@article{PTZ,
    author = "Paranjape, Shruti and Trnka, Jaroslav and Zheng, Minshan",
    title = "{Non-planar BCFW Grassmannian geometries}",
    eprint = "2208.02262",
    archivePrefix = "arXiv",
    primaryClass = "hep-th",
    doi = "10.1007/JHEP12(2022)084",
    journal = "JHEP",
    volume = "12",
    pages = "084",
    year = "2022"
}

@article{CEGM,
   title={$\Delta$-algebra and scattering amplitudes},
   volume={2019},
   ISSN={1029-8479},
   url={http://dx.doi.org/10.1007/JHEP02(2019)005},
   DOI={10.1007/jhep02(2019)005},
   number={2},
   journal={Journal of High Energy Physics},
   publisher={Springer Science and Business Media LLC},
   author={Cachazo, Freddy and Early, Nick and Guevara, Alfredo and Mizera, Sebastian},
   year={2019},
   month=Feb }

\bibliographystyle{alpha}

\end{document}